\documentclass[11pt]{amsdtx}
\usepackage{amscd,amsmath,amssymb,epic}
\usepackage[dvips]{graphicx}
\usepackage[mathscr]{eucal}

\makeatletter \@addtoreset{equation}{section}\makeatother
\renewcommand{\theequation}{\thesection.\arabic{equation}}

\newtheorem{theorem}{Theorem}[section]
\newtheorem{lemma}[theorem]{Lemma}
\newtheorem{proposition}[theorem]{Proposition}
\newtheorem{corollary}[theorem]{Corollary}

\title{\LARGE \bf LECTURES ON q-ANALOGUES OF CARTAN DOMAINS AND ASSOCIATED
HARISH-CHANDRA MODULES}

\author{\\ \bf \Large L. Vaksman (Ed.)}
\date{}

\begin{document}
\begin{titlepage}
\bigskip \bigskip \bigskip
\maketitle
\vfill

\unitlength=1mm \linethickness{.5mm}
\begin{center}
\begin{picture}(111,160)(0,0)

\multiput(3,25)(15,0){6}{\circle*{2}} \put(78,25){\circle{4}}
\put(33,10){\circle*{2}} \put(3,25){\line(1,0){75}}
\put(33,25){\line(0,-1){15}}

\multiput(3,50)(15,0){5}{\circle*{2}} \put(3,50){\circle{4}}
\put(33,35){\circle*{2}} \put(3,50){\line(1,0){60}}
\put(33,50){\line(0,-1){15}}

\multiput(3,70)(15,0){7}{\circle*{2}}
\multiput(108,60)(0,20){2}{\circle*{2}} \put(3,70){\circle{4}}
\dashline[70]{3}(63,70)(78,70) \put(3,70){\line(1,0){60}}
\put(78,70){\line(1,0){15}} \put(93,70){\thicklines \line(3,2){15}}
\put(93,70){\thicklines \line(3,-2){15}}

\multiput(3,90)(15,0){5}{\circle*{2}}
\multiput(93,90)(15,0){2}{\circle*{2}} \put(3,90){\circle{4}}
\put(93.5,89.4){\framebox(14,1.3){}} \dashline[70]{3}(63,90)(93,90)
\put(3,90){\line(1,0){60}} \put(100,89.1){$\big >$}

\multiput(3,110)(15,0){7}{\circle*{2}}
\multiput(108,100)(0,20){2}{\circle*{2}} \put(108,120){\circle{4}}
\dashline[70]{3}(63,110)(78,110) \put(3,110){\line(1,0){60}}
\put(78,110){\line(1,0){15}} \put(93,110){\thicklines \line(3,2){15}}
\put(93,110){\thicklines \line(3,-2){15}}

\multiput(3,135)(15,0){5}{\circle*{2}}
\multiput(93,135)(15,0){2}{\circle*{2}} \put(108,135){\circle{4}}
\put(93.5,134.4){\framebox(14,1.3){}} \dashline[70]{3}(63,135)(93,135)
\put(3,135){\line(1,0){60}} \put(100,134.1){$\big <$}

\multiput(3,150)(15,0){8}{\circle*{2}} \put(63,150){\circle{4}}
\put(3,150){\line(1,0){30}} \put(48,150){\line(1,0){30}}
\put(93,150){\line(1,0){15}} \dashline[70]{3}(33,150)(48,150)
\dashline[70]{3}(78,150)(93,150)

\end{picture}
\end{center}

\vfill

\centerline{\LARGE \bf Kharkov, Ukraine --- 2001}
\thispagestyle{empty}
\end{titlepage}
\newpage
\thispagestyle{empty}
\bigskip
\hfill {\bf \LARGE To Vladimir Drinfeld}
\newpage
\setcounter{tocdepth}{0}

\makeatletter
\renewcommand{\@oddhead}{CONTENTS \hfill \thepage}
\renewcommand{\@evenhead}{\thepage \hfill CONTENTS  }
\makeatother

{\bf\boldmath\raggedright

\centerline{Contents}\medskip

\noindent PREFACE\dotfill 5\medskip

\noindent INTRODUCTION

{\sl \ S. Sinel'shchikov and L. Vaksman}\dotfill 6\bigskip

\noindent Part I \ THE SIMPLEST EXAMPLE\hfill 12\medskip

\noindent A NON-COMMUTATIVE ANALOGUE OF THE FUNCTION THEORY IN THE UNIT DISC

{\sl D. Shklyarov, S. Sinel'shchikov, and L. Vaksman}\dotfill\ 13\medskip

\noindent QUANTUM DISC: THE BASIC STRUCTURES

{\sl D. Shklyarov, S. Sinel'shchikov, and L. Vaksman}\dotfill 30\medskip

\noindent QUANTUM DISC: THE CLIFFORD ALGEBRA AND THE DIRAC OPERATOR

{\sl K. Schm\"udgen, S. Sinel'shchikov, and L. Vaksman}\dotfill 40\medskip

\noindent ON UNIQUENESS OF COVARIANT \ DEFORMATION WITH SEPARATION OF
VARIABLES \ OF \ THE QUANTUM DISC

{\sl D. Shklyarov}\dotfill 51\bigskip

\noindent Part II \ UNITARY REPRESENTATIONS AND NON-COMMUTATIVE HARMONIC
ANALYSIS\hfill 62\medskip

\noindent HARISH-CHANDRA EMBEDDING AND q-ANALOGUES OF BOUNDED SYMMETRIC
DOMAINS

{\sl S. Sinel'shchikov and L. Vaksman}\dotfill 63\medskip

\noindent q-ANALOGUES OF SOME BOUNDED SYMMETRIC DOMAINS

{\sl D. Shklyarov, S. Sinel'shchikov, and L. Vaksman}\dotfill 67\medskip

\noindent ON A q-ANALOGUE OF THE FOCK INNER PRODUCT

{\sl D. Shklyarov}\dotfill 73\medskip

\noindent GEOMETRIC REALIZATIONS FOR SOME SERIES \ OF \ REPRESENTATIONS OF
THE QUANTUM GROUP $SU_{2,2}$

{\sl D. Shklyarov, S. Sinel'shchikov, and L. Vaksman}\dotfill 94\bigskip

\noindent Part III \ QUANTUM HARISH-CHANDRA MODULES ASSOCIATED TO q-CARTAN
DOMAINS\hfill 111\medskip

\noindent NON-COMPACT QUANTUM GROUPS AND HARISH-CHANDRA \ MODULES

{\sl D. Shklyarov, S. Sinel'shchikov, A. Stolin, and L. Vaksman}\dotfill
112\medskip

\noindent ON A q-ANALOGUES OF THE PENROSE TRANSFORM

{\sl D. Shklyarov, S. Sinel'shchikov, A. Stolin, and L. Vaksman}\dotfill
117\medskip

\noindent SPHERICAL PRINCIPAL NON-DEGENERATE SERIES OF REPRESENTATIONS FOR
THE QUANTUM GROUP $SU_{2,2}$

{\sl S. Sinel'shchikov, A. Stolin, and L. Vaksman}\dotfill 125\bigskip

\noindent Part IV \ ADDITIONAL RESULTS ON \ SOME \ QUANTUM \ VECTOR \
SPACES\hfill 135\medskip

\noindent HIDDEN SYMMETRY OF THE DIFFERENTIAL CALCULUS ON THE QUANTUM MATRIX
SPACE

{\sl S. Sinel'shchikov and L. Vaksman}\dotfill 136\medskip

\noindent q-ANALOGS OF CERTAIN PREHOMOGENEOUS VECTOR SPACES: COMPARISON OF
SEVERAL APPROACHES

{\sl D. Shklyarov}\dotfill 141\medskip

\noindent HIDDEN SYMMETRY OF SOME ALGEBRAS OF q-DIFFERENTIAL OPERATORS

{\sl D. Shklyarov, S. Sinel'shchikov, and L. Vaksman}\dotfill 155\medskip

}



\title{\bf PREFACE}\author{}\date{}

\newpage
\setcounter{section}{0}
\large

\makeatletter
\renewcommand{\@oddhead}{PREFACE \hfill \thepage}
\addcontentsline{toc}{chapter}{\@title \dotfill}
\makeatother

\maketitle

This volume contains a mildly expanded version of lectures and talks at
seminars and conferences, as well as review papers on subjects listed in the
title of the volume. A great deal of these texts have already been published
or sent to press. However, the only way to provide a good exposition of the
field we are interested in is to collect all of those papers together.

Bounded symmetric domains form a favorite subject of research in function
theory, non-commutative harmonic analysis, and representation theory of real
reductive Lie groups. The authors introduce the notions of q-analogues of
bounded symmetric domains and q-Harish-Chandra modules. For that, they
follow the traditions of quantum group theory in replacing Lie groups with
the quantum universal enveloping Drinfeld-Jimbo algebras and representations
of groups with associated Harish-Chandra modules.

The purpose of this volume is to convince the reader in exceptional
attraction of the deduced class of quantum homogeneous spaces and the
related class of modules over quantum universal enveloping algebras.

The first part of the volume contains 4 lectures on the quantum unit disc.
The first of them does not assume the knowledge of quantum group theory. It
contains a list of problems which we find specific for this part of
'non-commutative function theory', together with explanations about what
kind of results would be reasonable to expect when studying q-analogues of
bounded symmetric domains.

The second part of the volume contains 3 lectures devoted to 'function
theory' in q-analogues of bounded symmetric domains and the related problems
of non-commutative harmonic analysis. A feature of these fields is that the
algebras they deal with are involutive, and the modules over universal
enveloping algebras involved are normally unitarisable.

The third part of the volume includes 3 lectures devoted to q-analogues of
Harish-Chandra modules. These modules arise naturally as soon as one
disregards involution and unitarisability. It is worthwhile to emphasize an
important open problem in the quantum group theory, the problem of
construction and classification of simple quantum Harish-Chandra modules
discussed in the first lecture of this part.

The fourth part of the volume collects some additional and auxiliary
results. In particular, the last lecture contains a discussion on
equivalence of several approaches to construction of q-analogues for vector
spaces which admit natural embeddings of the bounded symmetric domains in
question.

These lectures were delivered within the period 1996 -- 2001. In some of
them it would be reasonable to provide references to papers which contain
complete proofs of the results announced in the lectures. Also, in some
other parts of the volume, we would like to mention some applications of our
results. This information is attached as notes of the Editor to particular
lectures.

\hfill Editor




\title{\bf INTRODUCTION}
\author{S. Sinel'shchikov \and L. Vaksman}
\date{\tt Institute for Low Temperature Physics \& Engineering\\
47 Lenin Avenue, 61103 Kharkov, Ukraine}

\newpage
\setcounter{section}{0}
\large

\makeatletter
\renewcommand{\@oddhead}{INTRODUCTION \hfill \thepage}
\renewcommand{\@evenhead}{\thepage \hfill S. Sinel'shchikov and L. Vaksman}
\let\@thefnmark\relax \@footnotetext{This lecture has been delivered at the
conference dedicated to the 40-th anniversary of the Institute for Low
Temperature Physics and Engineering, August 2000, Kharkov}
\addcontentsline{toc}{chapter}{\@title \\ {\sl S. Sinel'shchikov and L.
Vaksman}\dotfill}
\makeatother

\maketitle

A great deal of a background of the quantum group theory was worked out by
V.~Drinfeld and expounded in his report at International Congress of
Mathematicians \cite{uni0cdqa.texD}.

Just after its appearance in 80-th the theory of quantum groups attracts an
attention of a large number of specialists due to substantial new problems
that arise here as well as due to perfectly surprising applications and
links to various fields of mathematics and theoretical physics. By now,
dozens of monographs and hundreds of research papers have been dedicated to
the theory of quantum groups and related topics. It is too cumbersome to
describe all the areas in quantum group theory. We are going to dwell here
on the theory of 'functions' in q-Cartan domains.

Consider a bounded domain $\mathscr{D}$ in a finite dimensional complex
vector space. It is said to be symmetric if every point $p \in \mathscr{D}$
is an isolated fixed point for a biholomorphic involutive automorphism
$$
\varphi_p:\mathscr{D}\to \mathscr{D},\qquad \varphi_p \circ
\varphi_p=\mathrm{id}.
$$

$\mathscr{D}$ is called reducible if it is biholomorphically isomorphic to a
Cartesian product of two nontrivial domains. A classification of irreducible
bounded symmetric domains up to isomorphism was obtained by Eli Cartan.
There exist four series of such domains (A, C, D, BD) and the two
exceptional domains. It was demonstrated by Harish-Chandra that every
irreducible bounded symmetric domain admits a so called {\sl standard}
realization as the unit ball in a finite dimensional {\sl normed} vector
space. The irreducible bounded symmetric domains of this sort are called
Cartan domains.

\medskip

{\sc Series A.} Let $m,n \in \mathbb{N}$ and $m \le n$. Consider the vector
space $\mathrm{Mat}_{mn}$ of complex matrices with $m$ rows and $n$ columns.
Every such matrix $\mathbf{z}$ determines a linear map of Hermitian spaces
$\mathbb{C}^n \to \mathbb{C}^m$. Equip $\mathrm{Mat}_{mn}$ with the operator
norm $\|\mathbf{z}\|$. The unit ball $\mathscr{D}=\{\mathbf{z}\in
\mathrm{Mat}_{mn}|\,\|\mathbf{z}\|<1 \}$ is a Cartan domain of series A. It
is easy to produce an automorphism with a fixed point $\mathbf{z}_0$ in the
special case $\mathbf{z}_0=0$ (just $\mathbf{z}\mapsto-\mathbf{z}$). The
general case reduces to this special case by considering an action of the
group $\mathrm{Aut}(\mathscr{D})$. This group admits an easy description if
one passes from matrices $\mathbf{z}$ to graphs of the associated linear
maps, thus using the embedding $\mathscr{D}\hookrightarrow
\mathrm{Gr}_m(\mathbb{C}^{m+n})$ (similar embeddings exist for all Cartan
domains; these are called Borel embeddings).

\medskip

{\sc Example (series C, D).} The series C and D are defined in a similar
way, except that $m=n$, and in the case of series C one should consider the
symmetric matrices ($\mathbf{z}=\mathbf{z}^t$) while in the case of series D
antisymmetric matrices ($\mathbf{z}=-\mathbf{z}^t$). For all those series
the automorphisms of domains are just fractionally linear maps
$$
\mathbf{z}\mapsto(\mathbf{a}\mathbf{z}+\mathbf{b})
(\mathbf{c}\mathbf{z}+\mathbf{d})^{-1}.
$$

The following areas attract a permanent attention of specialists in function
theory and theory of representations of real reductive Lie groups (see a
list of references in \cite{uni0cdqa.texAr}):
\begin{itemize}
\item construction of unitary representations of $\mathrm{Aut}(\mathscr{D})$
in Hilbert spaces of functions in the domain $\mathscr{D}$ or on its
boundary,

\item the theory of Toeplitz and Hankel operators in those Hilbert spaces,

\item non-commutative harmonic analysis,

\item description of algebras of $\mathrm{Aut}(\mathscr{D})$-invariant
differential operators and their common eigenfunctions,

\item studying of generalized hypergeometric functions related to Cartan
domains.
\end{itemize}


Our principal observation is that the methods of quantum group theory allows
one to embed a great deal of the results on Cartan domains into a
one-parameter family of 'q-analogues' (everywhere in the sequel $q
\in(0,1)$).

The term 'q-analogue' can be easily illustrated by an example of a
q-analogue of the binomial formula:
$$(a+b)^k=\sum_{j=0}^k\binom{k}{j}b^ja^{k-j}.$$
Consider the free algebra $\mathbb{C}\langle a,b \rangle$ over the field of
complex numbers with two generators $a,b$, together with the two-sided ideal
generated by $ab-qba$. The associated factor algebra $\mathbb{C}[a,b]_q$ is
called the algebra of functions on the quantum plane. The monomials
$\{b^ja^l \}_{j,l \in\mathbb{Z}_+}$ form a basis in the vector space
$\mathbb{C}[a,b]_q$. Thus we get well defined numbers
$\genfrac{[}{]}{0pt}{0}{k;q}{j}$ such that
$$(a+b)^k=\sum_{j=0}^k \genfrac{[}{]}{0pt}{0}{k;q}{j}b^ja^{k-j}.$$
Of course, $\lim \limits_{q \to 1}\genfrac{[}{]}{0pt}{0}{k;q}{j}
=\dbinom{k}{j}$. The numbers $\genfrac{[}{]}{0pt}{0}{k;q}{j}$ are called
q-binomial coefficients or Gaussian polynomials as they are really
polynomials of the indeterminate $q$. They are well studied and are among
the simplest q-special functions \cite{uni0cdqa.texGR}.

After introducing 'Planck's constant' $h$ by $q=e^{-h/2}$, we observe a
similarity between the passage from the relation $ab-ba=0$ to
$ab-e^{-h/2}ba=0$ and the quantization procedure.

The quantum group theory \cite{uni0cdqa.texCP} allows one to obtain non-commutative
algebras which are q-analogues for algebras of functions in Cartan domains
and q-analogues of some results of function theory in those domains. To
rephrase, the quantum group theory leads to a variety of results of
non-commutative function theory in q-Cartan domains.

Consider a Cartan domain $\mathscr{D}\subset V$. The algebra
$\mathrm{Pol}(V)_q$ of 'polynomial functions' on the quantum vector space
$V$ for $0<q<1$ is introduced as a $*$-algebra determined by its generators
$z_1,z_2,\ldots,z_M$ and the relations
\begin{eqnarray*}
z_iz_j&=&q^{k_{ij}}z_jz_i+\sum_{M \ge i>i'>j'>j \ge
1}a_{ii'jj'}(q)z_{j'}z_{i'},\qquad i>j
\\ z_i^*z_j&=&\sum_{i',j'=1}^Mb_{ii'jj'}(q)z_{j'}z_{i'}^*+
(1-q^2)\delta_{ij}.
\end{eqnarray*}
The choice of $k_{ij}$, $a_{ii'jj'}(q)$, $b_{ii'jj'}(q)$ is determined by
quantum symmetry arguments; in particular, $\lim \limits_{q \to
1}a_{ii'jj'}(q)=0$, $\lim \limits_{q \to
1}b_{ii'jj'}(q)=\delta_{ii'}\delta_{jj'}$.

\medskip

{\sc Example.} The simplest Cartan domain is just the unit disc in
$\mathbb{C}$. In this special case the $*$-algebra
$\mathrm{Pol}(\mathbb{C})_q$ is determined by a single generator $z$ and a
single relation:
$$z^*z=q^2zz^*+1-q^2.$$

In the classical case (i.e. for $q=1$) the $*$-algebra $\mathrm{Pol}(V)$ can
be equipped with the norm
$$
\|f \|=\max_{\boldsymbol{z}\in \overline{\mathscr{D}}}|f(\boldsymbol{z})|,
$$
where $\overline{\mathscr{D}}$ is the closure of the domain $\mathscr{D}$.
Each element $\boldsymbol{z}_0 \in V$ determines a $*$-homomorphism
$$
\mathrm{Pol}(V)\to \mathbb{C},\qquad f(\boldsymbol{z})\mapsto
f(\boldsymbol{z}_0).
$$
In this setting the points of $\overline{\mathscr{D}}$ are exactly the
points which determine {\sl bounded} linear functionals.

That is, one can use $\|f \|$ in $\mathrm{Pol}(V)$ instead of
$\|\boldsymbol{z}\|$ in $V$ to distinguish the domain $\mathscr{D}\subset
V$. This approach can be transferred onto the quantum case. By definition,
$$\|f \|=\sup_{\widetilde{T}}\|\widetilde{T}(f)\|,$$
with $\widetilde{T}$ is any irreducible $*$-representation of
$\mathrm{Pol}(V)_q$ in a Hilbert space. (An experience shows that there
exists a unique\footnote{Up to unitary equivalence} faithful irreducible
$*$-representation $T$. This representation provides a supremum in the above
equality. In the example $\mathrm{Pol}(\mathbb{C})_q$ such
$*$-representation can be defined in $l^2(\mathbb{Z}_+)$ by
$T(z)e_k=(1-q^{(2k+1)})^{1/2}e_{k+1}$, $T(z^*)e_k=(1-q^{2k})^{1/2}e_{k-1} $,
$k \ne 0$, $T(z^*)e_0=0$.)

To conclude, we restrict ourselves to the case of unit disc in $\mathbb{C}$
in presenting a list of those formulae of classical function theory and
classical harmonic analysis where we were lucky to 'insert the index q'.
First, recall that
\begin{gather*}
\begin{pmatrix}a & b \\ c & d \end{pmatrix},\qquad z \mapsto
\frac{az+b}{cz+d},
\\ \mathscr{D}=\ [z \in \mathbb{C}|\,|z|<1 \},
\\ \mathrm{Aut}(\mathscr{D})=\left \{\left.\begin{pmatrix}a & b \\ c & d
\end{pmatrix}\right|\,\overline{a}=d,\;\overline{b}=c,\;ad-bc=1 \right
\}\left/\left \{\pm \begin{pmatrix}1 & 0 \\ 0 & 1 \end{pmatrix} \right
\}\right.,
\\ \begin{pmatrix}a & b \\ c & d \end{pmatrix},\qquad z \mapsto
\frac{az+b}{cz+d}.
\end{gather*}

The ring of $\mathrm{Aut}(\mathscr{D})$-invariant differential operators is
generated by
$$
\square=-(1-|z|^2)^2\frac{\partial}{\partial z}\cdot
\frac{\partial}{\partial \overline{z}}
$$
(the invariant Laplace operator). An invariant integral is of the form
$$d \nu=(1-|z|^2)^{-2}\cdot \frac{d \mathrm{Re}z\,d\mathrm{Im}z}{\pi}$$
(the normalizing multiple is not essential here).

The spectrum of $\square$ in $L^2(d \nu)$ is purely continuous and fills the
semiaxis $(-\infty,-\frac14]$. Recall an explicit form of an eigenfunction
expansion for $\square$.

The function $P(z,\zeta)=\dfrac{1-|z|^2}{|1-z \overline{\zeta}|^2}$ is
called the Poisson kernel. If $f$ is a 'good' function on the circle
$\partial \mathscr{D}=\{z \in \mathbb{C}|\,|z|=1 \}$, then for all $\rho \in
\mathbb{R}$ the function
$$
u(z)=\int \limits_{\partial \mathscr{D}}P(z,e^{i\theta})^{i
\rho+\frac12}f(e^{i\theta})\frac{d \theta}{2 \pi}
$$
is a solution of the differential equation:
$$\square u=\left(-\frac14-\rho^2 \right)u.$$
An eigenfunction expansion of a 'good' function in the disc $\mathscr{D}$
has the form
$$
u(z)=\int \limits_0^\infty \left \{\int \limits_0^{2 \pi}P(z,e^{i\theta})^{i
\rho+\frac12}\widetilde{u}(\rho,e^{i \theta})\frac{d \theta}{2 \pi}\right
\}s(\rho)d \rho,
$$
with $\widetilde{u}(\rho,e^{i \theta})$ being an analogue of the Fourier
transform,
$$
\widetilde{u}(\rho,e^{i \theta})=\int
\limits_{\mathscr{D}}P(z,e^{i\theta})^{-i \rho+\frac12}u(z)d \nu(z),
$$
and $s(\rho)d \rho=2 \rho \mathrm{th}(\pi \rho)d \rho$ the Plancherel
measure.

One more result is related to the invertibility of (the closure of) the
operator $\square$ in $L^2(d \nu)$:
$$
(\square^{-1}f)(z)=\int \limits_{\mathscr{D}}G(z,\zeta)f(\zeta)d \nu(\zeta).
$$
Here $G(z,\zeta)=\ln \left(\left|\dfrac{z-\zeta}{1-z
\overline{\zeta}}\right|^2 \right)$ is the Green function of $\square$.

\medskip

For a 'good' function $u(z)$ in the disc $\mathscr{D}$ one has
$$
u(z)=\frac1{2 \pi i}\int
\limits_{\mathscr{D}}\frac1{(1-\overline{\zeta}z)^2}u(\zeta)d
\overline{\zeta}\wedge d \zeta+\frac1{2 \pi i}\int
\limits_{\mathscr{D}}\frac{1-|z|^2}{(z-\zeta)(1-\overline{\zeta}z)}\cdot
\frac{\partial u}{\partial \overline{\zeta}}d \overline{\zeta}\wedge d
\zeta.
$$

All the above results admit q-analogues. Here is the idea of some of the
proofs. For $q=1$ the radial part $\square^{(0)}$ of the operator $\square$
is a differential operator in a suitable Hilbert space of functions on the
interval $(0,1)$. For $0<q<1$ the spectrum of the operator $T(1-zz^*)$ is a
union of $\{0 \}$ and the geometric progression $\{1,q^2,q^4,q^6,\ldots \}$.
This progression substitutes the interval $(0,1)$, and the operator
$\square^{(0)}$ is replaced by a difference operator (a difference
approximation of the differential operator $\square^{(0)}$) in the space of
functions on the above progression. What remains is to obtain formulae for
eigenfunction expansion for this difference operator. A passage from the
radial part $\square^{(0)}$ to the entire operator $\square$ is accessible
via an application of quantum symmetry arguments.

Describe briefly a passage from the disc to generic Cartan domains and their
q-analogues. Let $K \subset \mathrm{Aut}(\mathscr{D})$ be a subgroup of all
linear automorphisms of $\mathscr{D}$. The subalgebra of $K$-invariant
elements of the algebra of continuous functions in $\overline{\mathscr{D}}$
admits a q-analogue. One can prove that this subalgebra is commutative, and
its space of maximal ideals works as the set $\mathrm{spec}\,T=\{0 \}\cup
\{1,q^2,q^4,q^6,\ldots \}$. In this setting a single difference operator is
replaced by a family of pairwise commuting self-adjoint difference operators
in the associated Hilbert space. The problems in question are just the
expansion problems in the common eigenfunctions.

\section*{Notes of the Editor}

The recent works \cite{uni0cdqa.texSSV1,uni0cdqa.texSSV2} deals with the $*$-representations of
$\mathrm{Pol}(V)_q$ in the special case of series A, that is
$V=\mathrm{Mat}_{m,n}$. The results of those works are applied in some of
the lectures of this volume as well as in \cite{uni0cdqa.texT,uni0cdqa.texPT}. The work \cite{uni0cdqa.texT}
contains a description of all irreducible $*$-representations of
$\mathrm{Pol}(\mathrm{Mat}_{2,2})_q$ up to a unitary equivalence. In
\cite{uni0cdqa.texPT} one can find a proof of the fact that for $f \in
\mathrm{Pol}(\mathrm{Mat}_{2,2})_q$ the norms $\|f \|$, $\|Tf \|$ are the
same: $\|f \|=\|Tf \|$. These norms and the identity have been discussed in
the introduction.




\title{\bf A NON-COMMUTATIVE ANALOGUE OF THE FUNCTION THEORY IN THE UNIT
DISC}

\author{D. Shklyarov \and  S. Sinel'shchikov \and L. Vaksman}

\date{\tt Institute for Low Temperature Physics \& Engineering \\
National Academy of Sciences of Ukraine\\ e-mail:
sinelshchikov@ilt.kharkov.ua, vaksman@ilt.kharkov.ua}

\makeatletter \@addtoreset{equation}{section}\makeatother
\renewcommand{\theequation}{\thesection.\arabic{equation}}

\newpage
\setcounter{section}{0}
\large

\thispagestyle{empty}

\ \vfill \begin{center}\LARGE \bf PART I \\ THE SIMPLEST EXAMPLE
\end{center}
\vfill \addcontentsline{toc}{chapter}{Part I \ \ THE SIMPLEST EXAMPLE}
\newpage

\makeatletter
\renewcommand{\@oddhead}{A NON-COMMUTATIVE ANALOGUE OF THE UNIT DISC \hfill
\thepage}
\renewcommand{\@evenhead}{\thepage \hfill D. Shklyarov, S. Sinel'shchikov,
and L. Vaksman}
\let\@thefnmark\relax
\@footnotetext{This lecture has been delivered at the monthly seminar
'Quantum Groups', Paris, February 1996.}
\addcontentsline{toc}{chapter}{\@title \\ {\sl D. Shklyarov, S.
Sinel'shchikov, and L. Vaksman}\dotfill} \makeatother

\maketitle

\begin{quotation}\small {\bf Abstract.} The present work considers one of the
simplest homogeneous spaces of the quantum group $SU(1,1)$, the q-analogue
of the unit disc in ${\mathbb C}$. We present q-analogues of Cauchy-Green
formulae, integral representations of eigenfunctions of the Laplace-Beltrami
operator, Green functions for Poisson equation and an inversion formula for
Fourier transform. It is also demonstrated that the two-parameter
quantization of the disc introduced before by S. Klimek and A. Lesniewski,
can be derived via the Berezin method.
\end{quotation}

\bigskip

\section{Introduction}

The theory of von Neumann algebras is a non-commutative analogue of the
function theory of real variable and its far-reaching generalization. This
approach is also applicable for developing a non-commutative analogue of the
function theory of complex variable, as one can see in the classical work of
W. Arveson \cite{uni1FTQDIR.TEXAr}. However, a progress in non-commutative complex
analysis was inhibited by the absence of substantial examples. The break was
made possible by the fundamental works of L. Faddeev and his collaborators,
V.~Drinfeld, M. Jimbo, S. Woronovicz in quantum group theory. The subject of
this paper is the simplest example of a quantum homogeneous complex
manifold, the quantum disc. Meanwhile, a part of the exposed results are
extensible onto the case of quantum bounded symmetric domains \cite{uni1FTQDIR.TEXSV}.

We obtain non-commutative analogues for some integral representations of
functions in the disc and find a non-commutative analogue of the Plancherel
measure. The Berezin method is applied to produce a formal deformation for
our quantum disc. Our results provide an essential supplement to those of
the well known work \cite{uni1FTQDIR.TEXKL} in this field. The initial sections contain
the necessary material on differential and integral calculi in the quantum
disc.

In the cases when the statement in question to be proved via the quantum
group theory is available in the mathematical literature (in particular, in
our electronic preprints \cite{uni1FTQDIR.TEXSSV2},\cite{uni1FTQDIR.TEXSSV3},\cite{uni1FTQDIR.TEXSSV4},\cite{uni1FTQDIR.TEXSSV5}),
we restrict ourselves to making exact references. This approach allows us to
expound all the main results without applying the quantum group theory, and
thus making this text intelligible for a broader class of readers.

The authors are grateful to V. Drinfeld for helpful discussions of a draft
version of this work.

\bigskip

\section{Functions in the quantum disc}

We assume in the sequel all the vector spaces to be complex, and let $q$
stand for a real number, $0<q<1$.

The works \cite{uni1FTQDIR.TEXNN,uni1FTQDIR.TEXKL,uni1FTQDIR.TEXF,uni1FTQDIR.TEXCK} consider the involutive algebra
$\mathrm{Pol}({\mathbb C})_q$ given by its generator $z$ and the commutation
relation
$$z^*z=q^2zz^*+1-q^2.$$
It arises both in studying the algebras of functions in the quantum disc
\cite{uni1FTQDIR.TEXNN,uni1FTQDIR.TEXKL,uni1FTQDIR.TEXF} and in studying the q-analogues of the Weil algebra
(oscillation algebra) \cite{uni1FTQDIR.TEXCK}.

Any element $f \in \mathrm{Pol}({\mathbb C})_q$ admits a unique
decomposition
$$f=\sum_{jk}a_{jk}(f)z^jz^{*k},\quad a_{jk}\in{\mathbb C}.$$
Moreover, at the limit $q \to 1$ we have $z^*z=zz^*$. This allows one to
treat $\mathrm{Pol}({\mathbb C)}_q$ as a polynomial algebra on the quantum
plane.

It is worthwhile to note that the passage to the generators
$a^+=(1-q^2)^{-1/2}z^*$, $a=(1-q^2)^{-1/2}z$ realizes an isomorphism between
the $*$-algebra $\mathrm{Pol}({\mathbb C})_q$ and the q-analogue of the Weil
algebra considered in \cite{uni1FTQDIR.TEXCK}: $a^+a-q^2aa^+=1$.

Impose the notation $y=1-zz^*$. It is straightforward that
\begin{equation}z^*y=q^2yz^*,\quad zy=q^{-2}yz.\end{equation}
It is also easy to show that any element $f \in \mathrm{Pol}({\mathbb C)}_q$
admits a unique decomposition
\begin{equation}\label{uni1FTQDIR.TEXzycr} f=\sum_{m>0}z^m
\psi_m(y)+\psi_0(y)+\sum_{m>0}\psi_{-m}(y)z^{*m}.
\end{equation}
The passage to the decomposition (\ref{uni1FTQDIR.TEXzycr}) is similar to the passage from
Cartesian coordinates on the plane ${\mathbb C}\simeq{\mathbb R}^2$ to polar
coordinates.

We follow \cite{uni1FTQDIR.TEXKL,uni1FTQDIR.TEXNN} in completing the vector space
$\mathrm{Pol}({\mathbb C)}_q$ to obtain the function space in the quantum
disc.

Consider the $\mathrm{Pol}({\mathbb C)}_q$-module $H$ determined by its
generator $v_0$ and the relation $z^*v_0=0$. Let $T:\mathrm{Pol}({\mathbb
C)}_q \to \mathrm{End}_{\mathbb C}(H)$ be the representation of
$\mathrm{Pol}({\mathbb C)}_q$ in $H$, and $l_{m,n}(f)$, $m,n \in{\mathbb
Z}_+$, stand for the matrix elements of the operator $T(f)$ in the basis
$\{z^mv_0 \}_{m \in {\mathbb Z}_+}$.

Impose three topologies in the vector space $\mathrm{Pol}({\mathbb C)}_q$
and prove their equivalence.

Let ${\cal T}_1$ be the weakest among the topologies in which all the linear
functionals $l_{j,k}':f \mapsto a_{jk}(f)$, $j,k \in{\mathbb Z}_+$, are
continuous,\\ ${\cal T}_2$ the weakest among the topologies in which all the
linear functionals $l_{m,n}'':f \mapsto \psi_m(q^{2n})$, $m \in{\mathbb Z}$,
$n \in{\mathbb Z}_+$, are continuous,\\ ${\cal T}_3$ the weakest among the
topologies in which all the linear functionals $l_{m,n}$, $m,n \in{\mathbb
Z}_+$, are continuous.

\medskip

\begin{proposition}
The topologies ${\cal T}_1, {\cal T}_2, {\cal T}_3$ are equivalent.
\end{proposition}

\smallskip

{\bf Proof.} Remind the standard notation:
$$
(t;q)_m=\prod_{j=0}^{m-1}(1-tq^j);\quad
(t;q)_\infty=\prod_{j=0}^\infty(1-tq^j).
$$
It follows from the definitions that with $m \ge k$
$$z^{*k}z^mv_0=(q^{2m};q^{-2})_k \cdot z^{m-k}v_0.$$
Hence
$$
l_{m,n}(f)=\sum_{j=0}^{\min(m,n)}(q^{2m};q^{-2})_{m-j} \cdot
l_{n-j,m-j}';\quad m,n \in{\mathbb Z}_+.
$$
Thus the topology ${\cal T}_1$ is stronger than ${\cal T}_3$. In fact, these
topologies are equivalent since the linear span of the functionals
$\{l_{m-j,n-j}\}_{j=0}^{\min(m,n)}$ coincides with that of
$\{l_{m-j,n-j}'\}_{j=0}^{\min(m,n)}$. The equivalence of ${\cal T}_2$ and
${\cal T}_3$ follows from the relations
$l_{j,k}=(q^{2j};q^{-2})_{\max(0,j-k)}l_{k-j,\min(j,k)}''$. \hfill $\square$

\medskip

We denote the completion of the Hausdorff topological vector space
$\mathrm{Pol}({\mathbb C)}_q$ by $D(U)_q'$ and call it the space of
distributions in the quantum disc. $D(U)_q'$ may be identified with the
space of formal series of the form (\ref{uni1FTQDIR.TEXzycr}) whose coefficients
$\psi_j(y)$ are defined in $q^{2{\mathbb Z}_+}$. The linear functionals
$l_{m,n}$ are extendable by continuity onto the topological vector space
$D(U)_q'$. Associate to each distribution $f \in D(U)_q'$ the infinite
matrix $T(f)=(l_{m,n}(f))_{m,n \in{\mathbb Z}_+}$.

A distribution $f \in D(U)_q'$ is said to be finite if
$\#\{(j,k)|\,\psi_j(q^{2k})\ne 0\}<\infty$. Evidently, a distribution $f$ is
finite iff the matrix $T(f)$ has only finitely many non-zero entries. The
vector space of finite functions in the quantum disc is denoted by $D(U)_q$.
There exists a non-degenerate pairing
$$
D(U)_q' \times D(U)_q \to {\mathbb C};\quad f_1 \times f_2 \mapsto
\mathrm{tr}\,T(f_1)T(f_2).
$$

The extension by continuity procedure allows one to equip $D(U)_q$ with a
structure of $^*$-algebra, and $D(U)_q'$ with a structure of
$D(U)_q$-bimodule.

Consider the algebra $\mathrm{Pol}({\mathbb C})_q^{op}$ derived from
$\mathrm{Pol}({\mathbb C})_q$ via a replacement of the multiplication law by
the opposite one. The elements $z \otimes 1,\,z^*\otimes 1,\,1 \otimes z,\,1
\otimes z^*$ of $\mathrm{Pol}({\mathbb C})_q^{op}\otimes
\mathrm{Pol}({\mathbb C})_q$ are denoted respectively by
$z,z^*,\zeta,\zeta^*$. To avoid confusion in the notation, we use braces to
denote the multiplication in $\mathrm{Pol}({\mathbb
C})_q^{op}\otimes\mathrm{Pol}({\mathbb C})_q$, e.g.
$\{zz^*\}=q^2\{z^*z\}+1-q^2$. The module $H^{op}$ over
$\mathrm{Pol}({\mathbb C})_q^{op}$ is defined by its generator $v_0^{op}$
and the relation $zv_0^{op}=0$. Apply the above argument to
$\mathrm{Pol}({\mathbb C})_q^{op}\otimes \mathrm{Pol}({\mathbb C})_q$-module
$H^{op}\otimes H$ in order to introduce the algebra $D(U \times U)_q$ of
finite functions in the Cartesian product of quantum discs, together with
$D(U \times U)_q$-bimodule $D(U \times U)_q'$.

The reason for replacement the multiplication law in $\mathrm{Pol}({\mathbb
C})_q$ by the opposite one is cleared in \cite{uni1FTQDIR.TEXSSV2,uni1FTQDIR.TEXV1,uni1FTQDIR.TEXV2}.

The linear functional
$$
\mu(f)=(1-q^2)\sum_{m \in{\mathbb Z}_+}\psi_0(q^{2m})q^{2m},\quad f \in
D(U)_q,
$$
is called the (normalized) Lebesgue integral in the quantum disc, since
under the formal passage to the limit as $q \to 1$ one has $\mu(f)\to
\frac{1}{\pi}\mathop{\int \int}\limits_Ufd \mathrm{Im}\,z \cdot d
\mathrm{Re}\,z$.

Let $K \in D(U \times U)_q'$; the integral operator $f \mapsto
\mathrm{id}\otimes \mu(K(1 \otimes f))$ with kernel $K$ maps $D(U)_q$ into
$D(U)_q'$. We are interested in solving an inverse problem which is in
finding out the explicit formulae for kernels $K \in D(U \times U)_q'$ of
well known linear operators. In this field, an analogue of the Bergman
kernel for the quantum disc was obtained in the work of S. Klimek and A.
Lesniewski \cite{uni1FTQDIR.TEXKL}:
$$K_q(z,\zeta)=(1-z \zeta^*)^{-1}(1-q^2z \zeta^*)^{-1}.$$

Finally, equip $H$ with the structure of a pre-Hilbert space by setting
$$(z^jv_0,z^mv_0)=\delta_{jm}(q^2;q^2)_m;\quad j,m \in{\mathbb Z}_+.$$
It is easy to show that $T(z^*)=T(z)^*$, $I-T(z)T(z^*)\ge 0$. Thus we get a
$*$-representation of $\mathrm{Pol}({\mathbb C})_q$ in the completion
$\overline{H}$ of $H$ (see \cite{uni1FTQDIR.TEXKL}).

\bigskip

\section{Differential forms and $\overline{\partial}$-problem}

Let $\Omega_q({\mathbb C})$ stand for the involutive algebra determined by
its generators $z,dz$ and the relations
\begin{gather*}
1-z^* \cdot z=q^2(1-z \cdot z^*),\quad dz \cdot z^*=q^{-2}z^*\cdot dz,\quad
dz \cdot z=q^2z \cdot dz,
\\ dz \cdot dz^*=-q^{-2}dz^*\cdot dz,\quad dz \cdot dz=0.
\end{gather*}
(Also, an application of the involution $^*$ to the above yields
$$
dz^* \cdot dz^*=0,\quad dz^*\cdot z=q^2z \cdot dz^*,\quad dz^* \cdot
z^*=q^{-2}z^*\cdot dz^*.)
$$
Equip $\Omega_q({\mathbb C})$ with the grading as follows:
$$
\mathrm{deg}\,z=\mathrm{deg}\,z^*=0,\quad \mathrm{deg}\,dz=
\mathrm{deg}\,dz^*=1.
$$
There exists a unique linear map $d:\Omega_q({\mathbb C})\to
\Omega_q({\mathbb C})$ such that
$$
d:z \mapsto dz,\quad d:z^*\mapsto dz^*,\quad d:dz \mapsto 0,\quad
d:dz^*\mapsto 0,
$$
and
$$
d(\omega' \cdot \omega'')=d \omega'\cdot
\omega''+(-1)^{\mathrm{deg}\,\omega'}\cdot \omega'd \omega'';\quad
\omega',\omega'\in \Omega_q({\mathbb C}).
$$
Evidently, $d^2=0$, and $(d \omega)^*=d \omega^*$ for all $\omega \in
\Omega_q({\mathbb C})$.

Turn to a construction of operators $\partial,\overline{\partial}$. For
that, we need a bigrading in $\Omega_q({\mathbb C})$:
$$
\mathrm{deg}\,z=\mathrm{deg}\,z^*=(0,0);\quad \mathrm{deg}\,(dz)=(1,0);
\quad \mathrm{deg}\,(dz^*)=(0,1).
$$
Now $d$ has a degree 1 and admits a unique decomposition into a sum
$d=\partial+\overline{\partial}$ of operators $\partial$,
$\overline{\partial}$ with bidegrees respectively (1,0) and (0,1). A
standard argument allows one to deduce from $d^2=0$ that
$\partial^2=\overline{\partial}^2=\partial
\overline{\partial}+\overline{\partial}\partial=0$. It is also easy to show
that $(\partial \omega)^*=\overline{\partial}\omega^*$ for all $\omega \in
\Omega_q({\mathbb C})$.

Each element $\omega \in \Omega_q({\mathbb C})$ is uniquely decomposable
into a sum
$$
\omega=f_{00}+dz \,f_{10}+f_{01}dz^*+dz \,f_{11}dz^*,\quad f_{ij}\in
\mathrm{Pol}({\mathbb C})_q,\quad i,j=0,1.
$$
Equip $\Omega_q({\mathbb C})$ with a topology corresponding to this
decomposition:
$$
\Omega_q({\mathbb C})\simeq \mathrm{Pol}({\mathbb C})_q \oplus
\mathrm{Pol}({\mathbb C})_q \oplus \mathrm{Pol}({\mathbb C})_q \oplus
\mathrm{Pol}({\mathbb C})_q.
$$
Pass as above via a completion procedure from $\mathrm{Pol}({\mathbb C})_q$
to the space of distributions $D(U)_q'$ and then to the space of finite
functions to obtain the bigraded algebra $\Omega(U)_q$. The operators
$d,\partial,\overline{\partial}$ are transferred by continuity from
$\Omega_q({\mathbb C})$ onto the algebra $\Omega(U)_q$ of differential forms
with finite coefficients in the quantum disc.

The subsequent constructions involve essentially q-analogues of type (0,*)
differential forms with coefficients in sections of holomorphic bundles. The
latter carry a structure of bimodules over algebras of type (0,*)
differential forms as above. Remind the notion of a differentiation for such
a bimodule.

Let $\Omega$ be a ${\mathbb Z}_+$-graded algebra and $M$ a ${\mathbb
Z}_+$-graded $\Omega$-bimodule. A degree 1 operator is said to be a
differentiation if for all $m \in M,\,\omega \in \Omega$ one has
$\overline{\partial}(m
\omega)=(\overline{\partial}m)\omega+(-1)^{\mathrm{deg}\,m}m \cdot
\overline{\partial}\omega$, $\overline{\partial}(\omega
m)=(\overline{\partial}\omega)\cdot m+(-1)^{\mathrm{deg}\,\omega}\omega
\cdot \overline{\partial}m$.

Let $\lambda \in {\mathbb C}$. Consider the graded bimodule over
$\Omega({\mathbb C})_q^{(0,*)}=\Omega({\mathbb C})_q^{(0,0)}+\Omega({\mathbb
C})_q^{(0,1)}$ determined by its generator $v_ \lambda$ with
$\mathrm{deg}(v_ \lambda)=0$ and the relations
$$
z \cdot v_ \lambda=q^{-\lambda}v_ \lambda \cdot z,\quad z^*\cdot v_
\lambda=q^{\lambda} v_ \lambda \cdot z^*,\quad dz^*\cdot v_
\lambda=q^{\lambda} v_ \lambda \cdot dz^*.
$$
Denote this graded bimodule by $\Omega({\mathbb C})_{\lambda,q}^{(0,*)}$. It
possesses a unique degree 1 differentiation $\overline{\partial}$ such that
$\overline{\partial}v_ \lambda=0$. Pass (via an extension by continuity)
from polynomial coefficients to finite ones to obtain the graded bimodule
$\Omega(U)_{\lambda,q}^{(0,*)}$ over $\Omega(U)_q^{(0,*)}$, together with
its differentiation $\overline{\partial}$.

We restrict ourselves to the case $\lambda \in{\mathbb R}$ and equip the
spaces $\Omega(U)_{\lambda,q}^{(0,0)}$, $\Omega(U)_{\lambda,q}^{(0,1)}$ with
the scalar products
\begin{eqnarray}
(f_1 \cdot v_ \lambda,f_2 \cdot v_ \lambda)&=&\int
\limits_{U_q}f_2^*f_1(1-zz^*)^{\lambda-2}d \mu,
\\ (f_1v_ \lambda dz^*,f_2v_ \lambda dz^*)&=&\int
\limits_{U_q}f_2^*f_1(1-zz^*)^{\lambda}d \mu.
\end{eqnarray}

The completions of the pre-Hilbert spaces $\Omega(U)_{\lambda,q}^{(0,0)}$
and $\Omega(U)_{\lambda,q}^{(0,1)}$ can be used in the formulation of
$\overline{\partial}$-problem. Specifically, we mean finding a solution of
the equation $\overline{\partial}u=f$ in the orthogonal complement to the
kernel of $\overline{\partial}$. In the classical case ($q=1$) such a
formulation is standard \cite{uni1FTQDIR.TEXBFG}, and the solution is very well known. If
$\lambda$ stand for a real number and $\lambda>1$ than one has
\begin{equation}\label{uni1FTQDIR.TEXsol}
u(z)=\frac{1}{2 \pi i}\int \limits_U \frac{1}{z-\zeta}
\left(\frac{1-|\zeta|^2}{1-\overline{\zeta}z}\right)^{\lambda-1}f(\zeta)d
\overline{\zeta}\wedge d \zeta.
\end{equation}

(\ref{uni1FTQDIR.TEXsol}) implies the "Cauchy-Green formula":
\begin{multline}\label{uni1FTQDIR.TEXCGf}
u(z)=\frac{\lambda-1}{2 \pi i}\int \limits_U
\frac{(1-|\zeta|^2)^{\lambda-2}}{(1-\overline{\zeta}z)^{\lambda}}u(\zeta)d
\overline{\zeta}\wedge d \zeta+
\\ +\frac{1}{2 \pi i}\int \limits_U
\frac{(1-|\zeta|^2)^{\lambda-1}}{(z-\zeta)(1-\overline{\zeta}z)^{\lambda-1}}
\frac{\partial u}{\partial \overline{\zeta}}d \overline{\zeta}\wedge d
\zeta.
\end{multline}
Our purpose is to obtain the q-analogues of (\ref{uni1FTQDIR.TEXsol}), (\ref{uni1FTQDIR.TEXCGf}) for
$\lambda=2$.

The standard way of solving the $\overline{\partial}$-problem is to solve
first the Poisson equation $\square \omega=f$ with
$\square=-\overline{\partial}^*\overline{\partial}$. In the case $\lambda=2$
the kernel in (\ref{uni1FTQDIR.TEXsol}) is derived by a differentiation in $z$ of the
Green function
\begin{equation}\label{uni1FTQDIR.TEXGff}
G(z,\zeta)=\frac{1}{\pi}\ln(|z-\zeta|^2)-\frac{1}{\pi}\ln(|1-z
\overline{\zeta}|^2).
\end{equation}
In its turn, (\ref{uni1FTQDIR.TEXGff}) can be obtained by the d'Alembert method: the first
term is contributed by a real source, and the second one is coming from an
imaginary source.

Note that the differential calculus for the quantum disc we use here is well
known (see, for example, \cite{uni1FTQDIR.TEXZu}). Its generalization onto the case of an
arbitrary bounded symmetric domain was obtained in \cite{uni1FTQDIR.TEXSV} via an
application of a quantum analogue of the Harish-Chandra embedding \cite{uni1FTQDIR.TEXH1}.

\bigskip

\section{Green function for Poisson equation}

With $q=1$ the measure $d \nu=(1-|z|^2)^{-2}d \mu$ is invariant with respect
to the M\"obius transformations. In the case $q \in(0,1)$ impose an
"invariant integral" $\nu:D(U)_q \to{\mathbb C}$, $f \mapsto \int
\limits_{U_q}fd \nu$ by setting
$$
\int \limits_{U_q}fd \nu \stackrel{def}{=}\int \limits_{U_q}f
\cdot(1-zz^*)^{-2}d \mu.
$$
The Hilbert spaces $L^2(d \nu)_q$, $L^2(d \mu)_q$ are defined as completions
of the vector spaces $D(U)_q=\Omega(U)_q^{(0,0)}$, $\Omega(U)_q^{(0,1)}$
with respect to the norms
$$
\| f \|=\left(\int \limits_{U_q}f^*fd \nu \right)^{1/2}, \quad \|fdz^*
\|=\left(\int \limits_{U_q}f^*fd \mu \right)^{1/2}.
$$

Proofs of the following statements are to be found in \cite[Proposition 5.7,Corollary 5.8]{uni1FTQDIR.TEXSSV2}, \cite[Corollary 4.2]{uni1FTQDIR.TEXSSV4}.

\medskip

\begin{lemma}
There exist $0<c_1 \le c_2$ such that
\begin{equation}
\label{uni1FTQDIR.TEXc1c2} c_1 \le \overline{\partial}^*\overline{\partial}\le c_2.
\end{equation}
\end{lemma}

\medskip

\begin{proposition}
The exact estimates for $\overline{\partial}^*\overline{\partial}$ are of
the form
\begin{equation}\label{uni1FTQDIR.TEXeef}
\frac{1}{(1+q)^2}\le \overline{\partial}^*\overline{\partial}\le
\frac{1}{(1-q)^2}.
\end{equation}
\end{proposition}

\medskip

The inequalities (\ref{uni1FTQDIR.TEXc1c2}) allow one to extend by continuity the
operators $\overline{\partial}$,
$\square=-\overline{\partial}^*\overline{\partial}$ from the dense subspace
of finite functions $D(U)_q$ onto the entire $L^2(d \nu)_q$. They also imply
that for any $f \in L^2(d \nu)_q$ there exists a unique solution $u$ of
Poisson equation $\square u=f$. Now it follows from (\ref{uni1FTQDIR.TEXeef}) that $\| u
\| \le (1+q)^2 \| f \|$.

One obtains easily from the definitions

\medskip

\begin{lemma}
The series
$$
\sum_{i,j \in{\mathbb Z}_+}z^{*i}\zeta^{i}\psi_{ij}(y,\eta)z^{}\zeta^{*j},
\quad y=1-z^*z,\quad \eta=1-\zeta \zeta^*
$$
converges in $D(U \times U)_q'$ for any family $\{\psi_{ij}(y,\eta)\}_{i,j
\in{\mathbb Z}_+}$ of functions defined on $q^{2{\mathbb Z}_+}\times
q^{2{\mathbb Z}_+}$.
\end{lemma}

\medskip

\begin{corollary}
For all $m \ge 0$ there exists a well defined generalized kernel
\begin{equation}\label{uni1FTQDIR.TEXGm}
G_m=\left \{\left((1-\zeta \zeta^*)(1-z^*\zeta)^{-1}\right)^m
\left((1-z^*z)(1-z \zeta^*)^{-1}\right)^m \right \}.
\end{equation}
\end{corollary}

\medskip

To state the principal result of the section, we need an expansion of the
Green function (\ref{uni1FTQDIR.TEXGff}):
\begin{multline*}
\ln \frac{|z-\zeta|^2}{|1-z \overline{\zeta}|^2}=\ln
\left(1-\frac{(1-|z|^2)(1-|\zeta|^2)}{|1-z \overline{\zeta}|^2}\right)=
\\ =-\sum_{m=1}^\infty \frac{1}{m}\left(\frac{(1-|z|^2)(1-|\zeta|^2)}{|1-z
\overline{\zeta}|^2}\right)^m.
\end{multline*}

Evidently, a formal passage to a limit yields
$$
\lim_{q \to 1}\,G_m=\left(\frac{(1-|z|^2)(1-|\zeta|^2)}{|1-z
\overline{\zeta}|^2}\right)^m.
$$

A proof of the following result one can find in \cite[Theorem 1.2]{uni1FTQDIR.TEXSSV4}.

\medskip

\begin{theorem}
The continuous operator $\square^{-1}$ in $L^2(d \nu)_q$ coincides on the
dense linear subspace $D(U)_q \subset L^2(d \nu)_q$ with the integral
operator whose kernel is $G=- \sum \limits_{m=1}^\infty
\frac{q^{-2}-1}{q^{-2m}-1}G_m$:
$$\square^{-1}f=\int \limits_{U_q}G(z,\zeta)f(\zeta)d \nu.$$
Here $G_m \in D(U \times U)_q'$ is given by (\ref{uni1FTQDIR.TEXGm}).
\end{theorem}

\medskip

Note in conclusion that the operators $\partial$, $\overline{\partial}$,
$\square$ admit an extension by continuity onto the space $D(U)_q'$ of
distributions in the quantum disc.

\bigskip

\section{The Cauchy-Green formula}

One can use the differentials $\partial:\Omega_q^{(0,0)}\to
\Omega_q^{(1,0)}$, $\overline{\partial}:\Omega_q^{(0,0)}\to
\Omega_q^{(0,1)}$ to define the partial derivatives
$\frac{\partial^{(l)}}{\partial z}$, $\frac{\partial^{(r)}}{\partial z}$,
$\frac{\partial^{(l)}}{\partial z^*}$, $\frac{\partial^{(r)}}{\partial
z^*}$. Specifically, we set up $\partial f=dz \cdot
\frac{\partial^{(l)}f}{\partial z}=\frac{\partial^{(r)}f}{\partial z}dz$,
$\overline{\partial}f=dz^* \cdot \frac{\partial^{(l)}f}{\partial
z^*}=\frac{\partial^{(r)}f}{\partial z^*}dz^*$. It is easy to show that
these operators admit extensions by continuity from $D(U)_q$ onto $D(U)_q'$.

Let $f \in D(U)_q$. Define the integral of the (1,1)-form $dz \cdot f \cdot
dz^*$ over the quantum disc by $\int \limits_{U_q}dz \cdot f \cdot dz^*=-2i
\pi \int \limits_{U_q}fd \mu$.

A proof of the following proposition can be found in \cite[Theorem2.1]{uni1FTQDIR.TEXSSV4}.

\medskip

\begin{proposition}
Let $f \in D(U)_q$. Then
\begin{enumerate}

\item There exists a unique solution $u \in L^2(d \mu)_q$ of the
$\overline{\partial}$-problem $\overline{\partial}u=f$, which is orthogonal
to the kernel of $\overline{\partial}$.

\item $u=\frac{1}{2 \pi i}\int \limits_{U_q} d \zeta
\frac{\partial^{(l)}}{\partial z}G(z,\zeta)fd \zeta^*$, with $G \in D(U
\times U)_q'$ being the Green function of the Poisson equation.

\item $f=-\frac{1}{2 \pi i}\int \limits_{U_q}(1-z \zeta^*)^{-1}(1-q^{-2}z
\zeta^*)^{-1}d \zeta f(\zeta)d \zeta^*-\frac{1}{2 \pi i}\int \limits_{U_q} d
\zeta \frac{\partial^{(l)}}{\partial z}G(z,\zeta) \cdot
\frac{\partial^{(r)}f}{\partial \zeta^*}d \zeta^*$.

\end{enumerate}
\end{proposition}

\bigskip

\section{Eigenfunctions of the operator $\square$}

Let ${\mathbb C}[\partial U]_q$ stand for the algebra of finite sums of the
form
\begin{equation}\label{uni1FTQDIR.TEXfs}
\sum_{m \in{\mathbb Z}}a_me^{im \theta},\quad \theta \in {\mathbb R}/2
\pi{\mathbb Z}
\end{equation}
with complex coefficients. The ${\mathbb C}[\partial U]_q$-module of formal
series of the form (\ref{uni1FTQDIR.TEXfs}) is denoted by ${\mathbb C}[[\partial U]]_q$.
We also denote the algebra of finite sums like (\ref{uni1FTQDIR.TEXfs}) with coefficients
from $D(U)_q$ by $D(U \times \partial U)_q$, and the module of formal series
(\ref{uni1FTQDIR.TEXfs}) with coefficients from $D(U)_q'$ by $D(U \times \partial U)_q'$.
This vector space will be equipped by the topology of coefficientwise
convergence.

The use of the index $q$ in the notation for the above vector spaces is
justified by the fact that, as one can show, they are in fact modules over
the quantum universal enveloping algebra.

Recall the notations \cite{uni1FTQDIR.TEXGR}:
$$
(a;q^2)_\infty=\prod_{j \in{\mathbb Z}_+}(1-aq^{2j}),\quad(a;q^2)_\gamma=
\frac{(a;q^2)_\infty}{(aq^{2 \gamma};q^2)_\infty},\quad \gamma \in{\mathbb
C}.
$$

With $q=1$, the integral
$$
u(z)=\int \limits_{\partial U}\left(\frac{1-|z|^2}{(1-z \overline
\zeta)(1-\overline z \zeta)}\right)^{l+1}f(\zeta)d \nu,\quad d \nu=\frac{d
\theta}{2 \pi},
$$
represents an eigenfunction of $\square$ (see \cite{uni1FTQDIR.TEXH2}):
$$
\square u=\lambda(l)u,\quad
\lambda(l)=\left(l+\frac{1}{2}\right)^2-\frac{1}{4}.
$$

With $q \in(0,1)$, the power $P^\gamma$ of the Poisson kernel
$P=\frac{1-|z|^2}{|1-z \overline \zeta|^2}$ is replaced by the element $P_
\gamma \in D(U \times \partial U)_q$:
\begin{equation}
P_ \gamma=(1-zz^*)^\gamma(z
\zeta^*;q^2)_{-\gamma}\cdot(q^2z^*\zeta;q^2)_{-\gamma}.
\end{equation}

Here $(z \zeta^*;q^2)_{-\gamma}$, $(q^2z^*\zeta;q^2)_{-\gamma}$ are the
$q$-analogues of the powers $(1-z \overline \zeta)^{-\gamma}$, $(1-\overline
z \zeta)^{-\gamma}$, and the $q$-binomial theorem (see \cite{uni1FTQDIR.TEXGR}) implies
\begin{eqnarray*}
(z \zeta^*;q^2)_{-\gamma}&=&\sum_{n \in{\mathbb Z}_+}\frac{(q^{2
\gamma};q^2)_n}{(q^2;q^2)_n}(q^{-2 \gamma}z \zeta^*)^n,
\\ (q^2z^*\zeta;q^2)_{-\gamma}&=&\sum_{n \in{\mathbb Z}_+}\frac{(q^{2
\gamma};q^2)_n}{(q^2;q^2)_n}(q^{2-2 \gamma}z^*\zeta)^n.
\end{eqnarray*}

The following proposition is proved in \cite[Theorem 3.1]{uni1FTQDIR.TEXSSV4}.

\medskip

\begin{proposition}\label{uni1FTQDIR.TEXevir}
For all $f \in{\mathbb C}[\partial U]_q$ the element
\begin{equation}\label{uni1FTQDIR.TEXevf}
u=\int \limits_{\partial U}P_{l+1}(z,e^{i \theta})f(e^{i \theta})\frac{d
\theta}{2 \pi}
\end{equation}
of $D(U)_q'$ is an eigenvector of $\square$:
$$
\square u=\lambda(l)u,\quad
\lambda(l)=-\frac{(1-q^{-2l})(1-q^{2l+2})}{(1-q^2)^2}.
$$
\end{proposition}

\medskip

We need the following standard notation (\cite{uni1FTQDIR.TEXGR}):
\begin{multline*}
_r \Phi_s \genfrac{[}{]}{0pt}{0}{a_1,a_2,\ldots,a_r;q;z}{b_1,\ldots,b_s}=
\\ =\sum_{n \in{\mathbb Z}_+} \frac{(a_1;q)_n \cdot(a_2;q)_n \cdot \ldots
\cdot(a_r;q)_n} {(b_1;q)_n \cdot(b_2;q)_n \cdot \ldots
\cdot(b_s;q)_n(q;q)_n}\left((-1)^n \cdot
q^\frac{n(n-1)}{2}\right)^{1+s-r}\cdot z^n.
\end{multline*}

\medskip

\begin{corollary} (cf. \cite{uni1FTQDIR.TEXV1}).
The series
$$
\varphi_l(y)=\,_3 \Phi_2
\genfrac{[}{]}{0pt}{0}{y^{-1},q^{-2l},q^{2(l+1)};q^2;q^2}{q^2;0}
$$
converges in $D(U)_q'$, and its sum is an eigenfunction of $\square$:\
$\square \varphi_l=\lambda(l)\varphi_l$.
\end{corollary}

\smallskip

{\bf Proof.} The convergence of the series is due to the fact that it breaks
for each $y \in q^{2{\mathbb Z}_+}$. So it suffices to establish the
relation
$$\varphi_l(y)=\int \limits_{\partial U}P_{l+1}(z,\zeta)d \nu.$$

It follows from the definitions that the above integral equals to
\begin{multline*}
\sum_{n \in{\mathbb Z}_+}
\frac{(q^{2l+2};q^2)_n^2}{(q^2;q^2)_n^2}q^{-2(2l+1)n}y^{l+1}z^nz^{*n}=
\\ =y^{l+1}\, _3 \Phi_1
\genfrac{[}{]}{0pt}{0}{q^{2+2l},q^{2+2l},y^{-1};q^2;q^{-2(2l+1)}y}{q^2}.
\end{multline*}

Now it remains to apply the identity (see \cite{uni1FTQDIR.TEXGR}):
$$
b^n \,_3 \Phi_1 \genfrac{[}{]}{0pt}{0}{q^{-n},b,\frac{q}{z};q,
\frac{z}{c}}{\frac{bq^{1-n}}{c}}= \,_3 \Phi_2
\genfrac{[}{]}{0pt}{0}{q^{-n},b,
\frac{bzq^{-n}}{c};q,q}{\frac{bq^{1-n}}{c},0},
$$
with $q$ being replaced by $q^2$, $y$ by $q^{2n}$, $b$ by $q^{2l+2}$, $z$ by
$q^{-2l}$, and $c$ by $q^{2+2l-2n}$. \hfill $\square$

\medskip

Note that $\varphi_l(y)$ is a $q$-analogue of a spherical function on a
hyperbolic plane (see \cite{uni1FTQDIR.TEXH2}).

For each $l \in{\mathbb C}$ a linear operator has been constructed from
${\mathbb C}[\partial U]_q$ into the eigenspace of $\square$, associated to
the eigenvalue $\lambda(l)$. Now we try to invert this linear operator.

For that, we need a $q$-analogue of the operator $b_r:f(z)\mapsto f(re^{i
\theta})$ which restricts the function in the disc onto the circle $|z|=r$
of radius $r \in(0,1)$. Let $r>0$, $1-r^2 \in q^{2{\mathbb Z}_+}$. Define a
linear operator $b_r:D(U)_q'\to{\mathbb C}[[\partial U]]_q$ by
\begin{multline*}
b_r:\sum_{j>0}z^j \cdot \psi_j(y)+\psi_0(y)+\sum_{j>0}\psi_{-j}(y) \cdot
z^{*j}\mapsto
\\ \sum_{j>0}(re^{i \theta})^j \cdot
\psi_j(1-r^2)+\psi_0(1-r^2)+\sum_{j>0}\psi_{-j}(q^{-2j}\cdot (1-r^2))(re^{-i
\theta})^j.
\end{multline*}
(It is implicit that the functions $\psi_j(y)$, $j \in{\mathbb Z}$, vanish
at $y \notin q^{2{\mathbb Z}_+})$.

Recall the definition of the $q$-gamma-function (\cite{uni1FTQDIR.TEXGR}):
$$\Gamma_q(x)=\frac{(q;q)_\infty}{(q^x;q)_\infty}(1-q)^{1-x}.$$

One may assume without loss of generality that
$$
0 \le \mathrm{Im}\,l<\frac{\pi}{2 \ln(q^{-1})},\quad \mathrm{Re}\,l
\ge-\frac{1}{2}.
$$

\medskip

\begin{proposition}\label{uni1FTQDIR.TEXef}
Let $\mathrm{Re}\,l>-\frac{1}{2}$, and $u \in D(U)_q'$ is an eigenfunction
of $\square$ given by (\ref{uni1FTQDIR.TEXevf}). Then in ${\mathbb C}[\partial U]_q$ one
has
$$
f=\frac{\Gamma_{q^2}^2(l+1)}{\Gamma_{q^2}(2l+1)}\lim_{1-r^2 \in q^{2{\mathbb
Z}_+},\,r \to 1}(1-r^2)^lb_ru.
$$
\end{proposition}

\medskip

The proof of this proposition is based on the following result which was
communicated to the authors by L. I. Korogodsky:

\medskip

\begin{lemma}\label{uni1FTQDIR.TEXphas}
$$
\lim_{x \in q^{-2{\mathbb Z}_+},\,x \to \infty}\varphi_l
\left(\frac{1}{x}\right)\left/
\left(\frac{\Gamma_{q^2}(2l+1)}{\Gamma_{q^2}^2(l+1)}x^l \right)=1
\right.\leqno 1).
$$
if $\mathrm{Re}\,l>-\frac{1}{2}$.

$$
\lim_{x \in q^{-2{\mathbb Z}_+},\,x \to \infty}\varphi_l
\left(\frac{1}{x}\right)\left/
\left(\frac{\Gamma_{q^2}(-2l-1)}{\Gamma_{q^2}^2(-l)}x^{-l-1}\right)=1
\right.\leqno 2).
$$
if $\mathrm{Re}\,l<-\frac{1}{2}$.
\end{lemma}

{\bf Proof.} It follows from the relation $\varphi_l(y)=\varphi_{-1-l}(y)$
that one may restrict oneself to the case $\mathrm{Re}\,l>-\frac{1}{2}$. An
application of the identity (\cite{uni1FTQDIR.TEXGR})
$$
_2 \Phi_1 \genfrac{[}{]}{0pt}{0}{q^{-n},b;q;z}{c}=\frac{\left(\frac{c}{b};q
\right)_n}{(c;q)_n}\;_3 \Phi_2 \genfrac{[}{]}{0pt}{0}{q^{-n},b,
\frac{bzq^{-n}}{c};q;q}{\frac{bq^{1-n}}{c},0},
$$
with $q$, $b$, $c$, $z$ being replaced respectively by $q^2$, $q^{-2l}$,
$q^{-2l-2n}$, $q^{2l+2}$, yields
\begin{multline*}
\varphi_l(q^{2n})=\frac{(q^{-2l-2n};q^2)_n}{(q^{-2n};q^2)_n}\cdot \,_2
\Phi_1 \genfrac{[}{]}{0pt}{0}{q^{-2n};q^{-2l};q^2;q^{2l+2}}{q^{-2l-2n}}\sim
\\ \sim q^{-2nl}\frac{(q^{2(l+1)};q^2)_\infty}{(q^2;q^2)_\infty}\cdot \,_1
\Phi_0[q^{-2l};q^2;q^{2(l+1)}]=
\\ =q^{-2nl}\frac{(q^{2(l+1)};q^2)_\infty}{(q^2;q^2)_\infty}\cdot
\frac{(q^{2(l+1)};q^2)_\infty}{(q^{2(2l+1)};q^2)_\infty}.
\end{multline*}
Now it remains to refer to the definition of the $q$-gamma-function. \hfill
$\square$

\medskip

In the special case $f=1$ proposition \ref{uni1FTQDIR.TEXef} follows from lemma
\ref{uni1FTQDIR.TEXphas}. The general case reduces to the above special case via an
application of a quantum symmetry argument, which is described in
\cite[Theorem 3.7]{uni1FTQDIR.TEXSSV4}.

\bigskip

\section{Fourier transform}

It follows from the definitions that the integral operator with kernel
$K=\sum \limits_ik_i''\otimes k_i'$ is conjugate to the integral operator
with the kernel $K^t=\sum \limits_ik_i^{\prime*}\otimes k_i^{\prime
\prime*}$. Note that the conjugate to the unitary is an inverse operator.

Recall \cite{uni1FTQDIR.TEXH2} the heuristic argument that leads to the Fourier transform.
Proposition \ref{uni1FTQDIR.TEXevir} allows one to obtain eigenfunctions of $\square$. It
is natural to expect that "any" function $u$ admits a decomposition in
eigenfunctions of $\square$, and that the associated Fourier operator is
unitary.

Impose the notations: $h=-2 \ln q$,
\begin{eqnarray*}
P_{\gamma}^t&=&(q^2z^*\zeta;q^2)_{-\gamma}(z \zeta^*;q^2)_{-\gamma}(1-\zeta
\zeta^*)^{\gamma},
\\ c(l)&=&\Gamma_{q^2}(2l+1)/(\Gamma_{q^2}(l+1))^2.
\end{eqnarray*}

It is shown in \cite[section 5]{uni1FTQDIR.TEXSSV4} that, just as in the standard
representation theory (see \cite{uni1FTQDIR.TEXH2}), one has

\medskip

\begin{proposition}
Consider the Borel measure $d \sigma$ on $[0,\frac{\pi}{h}]$, given by
$$
d \sigma(\rho)=\frac{1}{2 \pi}\cdot \frac{h \cdot
e^h}{e^h-1}c(-\frac{1}{2}+i \rho)^{-1}\cdot c(-\frac{1}{2}-i \rho)^{-1}d
\rho.
$$
The integral operators
\begin{eqnarray*}
u(z)&\mapsto &\int \limits_{U_q}P_{\frac{1}{2}-i \rho}^t(z,\zeta)u(\zeta)d
\nu,
\\ f(e^{i \theta},\rho)&\mapsto &\int \limits_0^{ \pi/h}\int
\limits_0^{2 \pi}P_{\frac{1}{2}+i \rho}(z,e^{i \theta})f(e^{i
\theta},\rho)\frac{d \theta}{2 \pi}d \sigma(\rho)
\end{eqnarray*}
are extendable by continuity from the dense linear subspaces
$$
D(U)_q \subset L^2(d \nu)_q,\quad C^\infty[0,\frac{\pi}{h}]\otimes{\mathbb
C}[\partial U]_q \subset L^2(d \sigma)\otimes L^2(\frac{d \theta}{2 \pi})
$$
up to mutually inverse unitaries $F$, $F^{-1}$.
\end{proposition}

\medskip

{\sc Remark.} The function $c(l)$, the measure $d \sigma(\rho)$ and the
operator $F$ are the quantum analogues for c-function of Harish-Chandra,
Plancherel measure and Fourier transform respectively (see \cite{uni1FTQDIR.TEXH2}).

\bigskip

\section{The Berezin deformation of the quantum disc}

We are going to use in the sequel bilinear operators $L:D(U)_q \times D(U)_q
\to D(U)_q$ of the form
\begin{equation}\label{uni1FTQDIR.TEXL}
L:\,f_1 \times f_2 \to
\sum_{ijkm=0}^{N(L)}a_{ijkm}\left(\left(\frac{\partial^{(r)}}{\partial
z^*}\right)^if_1 \right)z^{*j}z^k \left(\left(\frac{\partial^{(l)}}{\partial
z}\right)^mf_2 \right),
\end{equation}
with $a_{ijkm}\in{\mathbb C}$. Such operators will be called
$q$-bidifferential.

Our principal purpose is to construct the formal deformation of the
multiplication law in $D(U)_q$. The new multiplication is to be a bilinear
map
\begin{align*}
*:&\,D(U)_q \times D(U)_q \to D(U)_q[[t]],
\\ *:&\,f_1 \times f_2  \mapsto  f_1 \cdot f_2+\sum_{i=1}^\infty
t^iC_i(f_1,f_2),
\end{align*}
which satisfies the formal associativity condition
$$\sum_{i+k=m}C_i(f_1,C_k(f_2,f_3))=\sum_{i+k=m}C_i(C_k(f_1,f_2),f_3)$$
(cf. \cite{uni1FTQDIR.TEXL}). When producing the new multiplication $*$, we follow F.
Berezin \cite{uni1FTQDIR.TEXBe}. The bilinear operators $C_j:D(U)_q \times D(U)_q \to
D(U)_q$, $j \in{\mathbb N}$, will turn out to be $q$-bidifferential, and we
shall give explicit formulae for them.

To begin with, choose a positive $\alpha$ and consider a linear functional
$\nu_ \alpha:\mathrm{Pol}({\mathbb C})_q \to {\mathbb C}$;
$$
\int \limits_{U_q}fd \nu_ \alpha \stackrel{def}{=}\frac{1-q^{4
\alpha}}{1-q^2}\cdot \int \limits_{U_q}f \cdot (1-zz^*)^{2 \alpha+1}d
\nu=(1-q^{4 \alpha})\,\mathrm{tr}\,T(f \cdot(1-zz^*)^{2 \alpha}).
$$

Impose a norm $\| f \|_ \alpha=\left(\int \limits_{U_q}f^*fd \nu_ \alpha
\right)^{1/2}$ on $\mathrm{Pol}({\mathbb C})_q$. Let $L_{q,\alpha}^2$ stand
for the completion of $\mathrm{Pol}({\mathbb C})_q$ with respect to the
above norm, and $H_{q,\alpha}^2$ for the linear span of monomials $z^j \in
L_{q,\alpha}^2$, $j \in{\mathbb Z}_+$.

\medskip

\begin{lemma}\label{uni1FTQDIR.TEXorth}
The monomials $\{z^m\}_{m \in{\mathbb Z}_+}$ are pairwise orthogonal in
$H_{q,\alpha}^2$, and $\| z^m \|_ \alpha=((q^2;q^2)_m/(q^{4
\alpha+2};q^2)_m)^{1/2}$.
\end{lemma}

\smallskip

{\bf Proof.} The pairwise orthogonality of the monomials $z^m$ is obvious;
\begin{multline*}
\| z^m \|_ \alpha^2=(1-q^{4 \alpha})\cdot
\mathrm{tr}\,T(z^{*m}z^m(1-zz^*)^{2 \alpha})=
\\ =\frac{1-q^{4 \alpha}}{1-q^2}\int \limits_0^1(q^2y;q^2)_m \cdot y^{2
\alpha-1}d_{q^2}y=
\\ =\frac{1-q^{4 \alpha}}{1-q^2}\cdot \frac{\Gamma_{q^2}(2 \alpha)\cdot
\Gamma_{q^2}(m+1)}{\Gamma_{q^2}(m+2 \alpha+1)}=\frac{(q^2;q^2)_m}{(q^{4
\alpha+2};q^2)_m}.
\end{multline*}

We have used the well known \cite[\S 1.11]{uni1FTQDIR.TEXGR} identity
$$
\int \limits_0^1t^{\beta-1}\cdot(tq^2;q^2)_{\alpha-1}d_{q^2}t=
\frac{\Gamma_{q^2}(\beta)\Gamma_{q^2}(\alpha)}{\Gamma_{q^2}(\alpha+\beta)}.
\eqno \square
$$

\medskip

\begin{corollary}\label{uni1FTQDIR.TEXzh}
Let $\widehat{z}$ be the operator of multiplication by $z$ in
$H_{q,\alpha}^2$, and $\widehat{z}^*$ the conjugate operator. Then
$\widehat{z}$, $\widehat{z}^*$ are bounded, and
\begin{equation}\label{uni1FTQDIR.TEXzhcr}
\widehat{z}^*\widehat{z}=q^2 \widehat{z}\widehat{z}^*+1-q^2+q^{4
\alpha}\cdot \frac{1-q^2}{1-q^{4
\alpha}}(1-\widehat{z}\widehat{z}^*)(1-\widehat{z}^*\widehat{z}).
\end{equation}
\end{corollary}

\smallskip

{\bf Proof} \ follows from
\begin{gather*}
\widehat{z}:z^m \mapsto z^{m+1},\qquad m \in {\mathbb Z}_+;
\\ \widehat{z}^*:1 \mapsto 0,\qquad \widehat{z}^*:z^m \mapsto
\frac{1-q^{2m}}{1-q^{4 \alpha+2m}}z^{m-1},\qquad m \in{\mathbb N}.
\end{gather*}

In fact,
\begin{align*}
(1-\widehat{z}\widehat{z}^*)^{-1}:&\,z^m \mapsto((q^{-2m}-q^{4
\alpha})/(1-q^{4 \alpha})z^m,
\\ (1-\widehat{z}^*\widehat{z})^{-1}:&\,z^m \mapsto((q^{-2m-2}-q^{4
\alpha})/(1-q^{4 \alpha})z^m.
\end{align*}
Hence
$(1-\widehat{z}\widehat{z}^*)^{-1}=q^2(1-\widehat{z}^*\widehat{z})^{-1}-q^{4
\alpha}\frac{1-q^2}{1-q^{4 \alpha}}$. \hfill $\square$

\medskip

Lemma \ref{uni1FTQDIR.TEXorth} and corollary \ref{uni1FTQDIR.TEXzh} were proved in the work by S. Klimek
and A. Lesniewski \cite{uni1FTQDIR.TEXKL} on two-parameter quantization of the disc beyond
the framework of perturbation theory.

To every element $f=\sum a_{ij}z^iz^{*j}\in \mathrm{Pol}({\mathbb C})$ we
associate the linear operator $\widehat{f}=\sum a_{ij}\widehat{z}^i
\widehat{z}^{*j}$ in $H_{q,\alpha}^2$. The formal deformation of the
multiplication law in the algebra of functions in the quantum disc will be
derived via an application of "Berezin quantization procedure" $f \mapsto
\widehat{f}$ to the ordinary multiplication in the algebra of linear
operators.

More exactly, (\ref{uni1FTQDIR.TEXzhcr}) allows one to get a formal asymptotic expansion
$$
\widehat{f_1}\cdot \widehat{f_2}=\widehat{f_1 \cdot f_2}+\sum_{k=1}^\infty
q^{4 \alpha k}{C_k \widehat{(f_1,f_2)}}. \quad f_1,f_2 \in
\mathrm{Pol}({\mathbb C})_q,
$$
with $C_k:\mathrm{Pol}({\mathbb C})_q \times \mathrm{Pol}({\mathbb C})_q \to
\mathrm{Pol}({\mathbb C})_q$, $k \in{\mathbb N}$, bilinear maps. In this
way, we get a formal deformation
\begin{gather*}
*:\,\mathrm{Pol}({\mathbb C})_q \times \mathrm{Pol}({\mathbb C})_q \to
\mathrm{Pol}({\mathbb C})_q[[t]];
\\ f_1*f_2=f_1 \cdot
f_2+\sum_{k=1}^\infty t^k \cdot C_k(f_1,f_2);\quad f_1,f_2 \in
\mathrm{Pol}({\mathbb C})_q.
\end{gather*}

We present an explicit formula for the multiplication $*$, and thus also for
bilinear maps $C_k$, $k \in{\mathbb N}$. Let $\stackrel{\sim}{\square}$ be a
linear operator in $\mathrm{Pol}({\mathbb C})_q^{op}\otimes
\mathrm{Pol}({\mathbb C})_q$ given by
$$
\stackrel{\sim}{\square}=q^{-2}(1-(1+q^{-2})z^*\otimes z+q^{-2}z^{*2}\otimes
z^2)\frac{\partial^{(r)}}{\partial z^*}\otimes
\frac{\partial^{(l)}}{\partial z},
$$
and $m:\,\mathrm{Pol}({\mathbb C})_q \times \mathrm{Pol}({\mathbb C})_q \to
\mathrm{Pol}({\mathbb C})_q$, $m:\,\psi_1 \otimes \psi_2 \to \psi_1 \psi_2$
the multiplication in $\mathrm{Pol}({\mathbb C})_q$.

\medskip

\begin{theorem}\label{uni1FTQDIR.TEXcvef}
For all $f_1,f_2 \in \mathrm{Pol}({\mathbb C})_q$
$$
f_1*f_2=(1-t)\sum_{j \in{\mathbb Z}_+}t^j \cdot
m(p_j(\stackrel{\sim}{\square})f_1 \otimes f_2),
$$
with
$$
p_j(\stackrel{\sim}{\square})=\sum_{k=0}^j
\frac{(q^{-2j};q^2)_k}{(q^2;q^2)^2_k}q^{2k}\cdot
\prod_{i=0}^{k-1}(1-q^{2i}((1-q^2)^2 \cdot
\stackrel{\sim}{\square}+1+q^2)+q^{4i+2}).
$$
\end{theorem}

\smallskip

The proof can be found in \cite[Theorem 8.4]{uni1FTQDIR.TEXSSV5}.

\medskip

{\sc Example.} For all $f_1,f_2 \in \mathrm{Pol}({\mathbb C})_q$
$$
f_1*f_2=f_2 \cdot f_2+t \cdot(q^{-2}-1)\frac{\partial^{(r)}f_1}{\partial
z^*}(1-z^*z)^2 \frac{\partial^{(l)}f_2}{\partial z}+O(t^2).
$$

\medskip

It is worthwhile to note that the formal associativity of the multiplication
$*$ follows from the associativity of multiplication in the algebra of
linear operators.

\medskip

\begin{corollary}
The bilinear operators $C_k$ are of the form (\ref{uni1FTQDIR.TEXL}) and are extendable by
continuity up to $q$-bidifferential operators $C_k:\,D(U)_q \times D(U)_q
\to D(U)_q$.
\end{corollary}

\medskip

The above $q$-bidifferential operators determine a formal deformation of the
multiplication in $D(U)_q$. The formal associativity of the newly formed
multiplication $*:\,D(U)_q \times D(U)_q \to D(U)_q[[t]]$ follows from the
formal associativity of the previous multiplication
$*:\,\mathrm{Pol}({\mathbb C})_q \times \mathrm{Pol}({\mathbb C})_q \to
\mathrm{Pol}({\mathbb C})_q[[t]]$.

Finally, note that our proof of theorem \ref{uni1FTQDIR.TEXcvef} is based on the
properties of some $q$-analogue for Berezin transform \cite{uni1FTQDIR.TEXUU}.

\bigskip

\section{Appendix. On q-analogue of the Green formula}

Consider the two-sided ideal $J \in \mathrm{Pol}({\mathbb C})_q$ generated
by the element $1-zz^*\in \mathrm{Pol}({\mathbb C})_q$, and the commutative
quotient algebra ${\mathbb C}[\partial U]_q
\stackrel{def}{=}\mathrm{Pol}({\mathbb C})_q/J$. Its elements will be
identified with the corresponding polynomials on the circle $\partial U$.
The image $f|_{\partial U}$ of $f \in \mathrm{Pol}({\mathbb C})_q$ under the
canonical homomorphism $\mathrm{Pol}({\mathbb C})_q \to{\mathbb C}[\partial
U]_q$ will be called a restriction of $f$ onto the boundary of the quantum
disc.

Define the integral $\Omega({\mathbb C})_q^{(1,0)}\to{\mathbb C}$ by
$$
\int \limits_{\partial U}dz \cdot f \stackrel{def}{=}2 \pi i \int
\limits_{\partial U}(z \cdot f)|_{\partial U}d \nu,\quad f \in
\mathrm{Pol}({\mathbb C})_q,
$$
with
$$
\int \limits_{\partial U}\psi d \nu \stackrel{def}{=}\int \limits_0^{2
\pi}\psi(e^{i \theta})\frac{d \theta}{2 \pi}.
$$

\medskip

\begin{proposition}\label{uni1FTQDIR.TEXGf}
For all $\psi \in \Omega({\mathbb C})_q^{(0,1)}$ one has
\begin{equation}\label{uni1FTQDIR.TEXGfe}
\int \limits_{U_q}\overline{\partial}\psi=\int \limits_{\partial U}\psi.
\end{equation}
\end{proposition}

\smallskip

{\sc Remark.} The integral $\int \limits_{U_q}dz \cdot f \cdot dz^*=-2i \pi
\int \limits_{U_q}fd \mu$ introduced in section 4 for $f \in D(U)_q$, is
extendable by continuity onto all (1,1)-forms $dz \cdot f \cdot dz^*$ with
$$
f=\sum_{m>0}z^m \psi_m(y)+\psi_0(y)+\sum_{m>0}\psi_{-m}(y)z^{*m}\in D(U)_q',
$$
such that $\sum \limits_{m \in{\mathbb Z}_+}|\psi_0(q^{2m})|q^{2m}<\infty$.
Under these assumptions one also has
\begin{equation}
\int \limits_{U_q}dzdz^*f=\int \limits_{U_q}dz \cdot f \cdot dz^*=\int
\limits_{U_q}fdzdz^*.
\end{equation}

\smallskip

{\bf Proof.} Recall that (see (\ref{uni1FTQDIR.TEXzycr}))
\begin{equation}\label{uni1FTQDIR.TEXpd}
\psi=dz \left(\sum_{m>0}z^m
\psi_m(y)+\psi_0(y)+\sum_{m>0}\psi_{-m}(y)z^{*m}\right).
\end{equation}
We can restrict ourselves to the case $\psi=dz \psi_{-1}(y)z^*$, since this
is the only term in (\ref{uni1FTQDIR.TEXpd}) which could make a non-zero contribution to
(\ref{uni1FTQDIR.TEXGfe}).

It follows from the definitions that $\overline{\partial}\psi=dz \cdot
f(y)dz^*$, with
$$
f(y)=\psi_{-1}(y)-q^{-2}
\frac{\psi_{-1}(q^{-2}y)-\psi_{-1}(y)}{q^{-2}y-y}(1-y).
$$
In fact, $\overline{\partial}y=\overline{\partial}(1-zz^*)=-zdz^*$. Hence
$\overline{\partial}y^m=\sum
\limits_{j=0}^{m-1}y^j(-zdz^*)y^{m-1-j}=-\frac{1-q^{2m}}{1-q^2}zy^{m-1}dz^*$.

That is, for any polynomial $p(y)$ one has
\begin{equation}\label{uni1FTQDIR.TEXdpy}
\overline{\partial}p(y)=-z \frac{p(y)-p(q^2y)}{y-q^2y}\cdot dz^*.
\end{equation}
(Note that the validity of (\ref{uni1FTQDIR.TEXdpy}) for polynomials already implies its
validity for all distributions). Finally,
$$
\overline{\partial}(dz \psi_{-1}(y)z^*)=dz \left(-z
\frac{\psi_{-1}(y)-\psi_{-1}(q^2y)}{y-q^2y}z^*+\psi_{-1}(y)\right)dz^*.
$$
On the other hand, $-z \frac{\psi_{-1}(y)-\psi_{-1}(q^2y)}{y-q^2y}z^*+
\psi_{-1}(y)=f(y)$, since $zy=q^{-2}yz$, $zz^*=1-y$.

If one assumes $\psi_{-1}(0)=0$, it is easy to show that $\sum \limits_{n
\in{\mathbb Z}_+}f(q^{2n})q^{2n}=0$. Hence, in this case $\int
\limits_{U_q}\overline{\partial}\psi=\int \limits_{\partial U}\psi=0$. Thus,
Proposition \ref{uni1FTQDIR.TEXGf} is proved for all (1,0)-forms from some linear subspace
of codimensionality 1. Now it remains to prove (\ref{uni1FTQDIR.TEXGfe}) in the special
case $\psi=dz \cdot z^*$. \hfill $\square$

\medskip

\begin{corollary}
If $\psi \in \Omega(U)_q^{(1,0)}$, then $\int
\limits_{U_q}\overline{\partial}\psi=0$.
\end{corollary}




\makeatletter \@addtoreset{equation}{section} \makeatother

\renewcommand{\theequation}{\thesection.\arabic{equation}}

\title{\bf QUANTUM DISC: THE BASIC STRUCTURES}
\author{D. Shklyarov \and S. Sinel'shchikov \and L. Vaksman}
\date{}

\newpage
\setcounter{section}{0}
\large

\makeatletter
\renewcommand{\@oddhead}{QUANTUM DISC: THE BASIC STRUCTURES \hfill \thepage}
\renewcommand{\@evenhead}{\thepage \hfill D. Shklyarov, S. Sinel'shchikov,
and L. Vaksman}
\let\@thefnmark\relax
\@footnotetext{This is an expanded version of a lecture delivered at a
seminar for undergraduate students in Kharkov, April 2000.}
\addcontentsline{toc}{chapter}{\@title \\ {\sl D. Shklyarov, S.
Sinel'shchikov, and L. Vaksman}\dotfill}
\makeatother

\maketitle

\section{Introduction}

The theory of quantum groups and their homogeneous spaces is in an extensive
progress since mid-80's \cite{uni1k_qdbs.texDr}. An important class of such homogeneous
spaces is formed by q-analogues of Cartan domains. The simplest Cartan
domain is the unit disc $\mathbb{U}$ in $\mathbb{C}$; the present work is
devoted to its q-analogue, the quantum disc. We intend to introduce the
basic notions of the theory of q-Cartan domains while restricting ourselves
to this simplest case. We hope this text could facilitate reading works on
general q-Cartan domains and will allow better understanding the results of
non-commutative function theory in quantum disc \cite{uni1k_qdbs.texSSV1}. The concluding
part of our work contains references to the papers of other authors devoted
to the related topics. Among those one should emphasize the work by K.
Schm\"udgen and A. Sch\"uler \cite{uni1k_qdbs.texSchSch} which introduces the differential
calculus in quantum disc which we use below.

The next section recalls the basic notions of the quantum group theory after
the lecture notes \cite{uni1k_qdbs.texJantz} and the monograph \cite{uni1k_qdbs.texCh-P}.

\bigskip

\section{The $*$-Hopf algebra $\mathbf{U_q \mathfrak{su}_{1,1}}$}

Everywhere in the sequel we assume $0<q<1$, and $\mathbb{C}$ is implicit as
a ground field. The algebras under consideration are supposed to be unital,
unless the contrary is stated explicitly.

Recall the definition of the Hopf algebra $U_q \mathfrak{sl}_2$. It is
determined by its generators $E$, $F$, $K$, $K^{-1}$ and the following
relations:
\begin{align*}
&KK^{-1}=K^{-1}K=1,&\\ &K^{\pm 1}E=q^{\pm 2}EK^{\pm 1},&K^{\pm 1}F=q^{\mp
2}FK^{\pm 1},\\ &EF-FE=\frac{K-K^{-1}}{q-q^{-1}}.&
\end{align*}

A $U_q \mathfrak{sl}_2$-module $V$ is said to be $\mathbb{Z}$-weight if
$$
V=\bigoplus_{\mu \in \mathbb{Z}}V_\mu,\qquad V_\mu=\{v \in V|\:K^{\pm
1}v=q^{\pm \mu}v \}.
$$

We restrict our considerations to $U_q \mathfrak{sl}_2$-modules of the above
form and define a linear operator $H$ by $H|_{V_\mu}=\mu$. Obviously,
$K^{\pm 1}v=q^{\pm H}v$, $v \in V$, and the following relations are valid in
the algebra of linear operators:
\begin{gather*}
HE-EH=2E,\qquad HF-FH=-2F,
\\ EF-FE=\frac{q^H-q^{-H}}{q-q^{-1}}.
\end{gather*}

A formal passage to a limit as $q \to 1$ leads to the determining relations
of the ordinary universal enveloping algebra $U \mathfrak{sl}_2$.

Within the category of modules over an algebra $A$, the operation of tensor
product is not defined; it is also not clear which $A$-module is to be
treated as trivial and how to define an $A$-module $V^*$ dual to an
$A$-module $V$. These three problems are solvable via equipping $A$ with the
structure of Hopf algebra. Specifically, one has to distinguish
\begin{itemize}
 \item a homomorphism $\triangle:A \to A \otimes A$, called a
comultiplication and used for producing a tensor product of $A$-modules
$V'$, $V''$:
$$
A \underset{\triangle}{\to}A \otimes A \to \mathrm{End(V')}\otimes
\mathrm{End(V'')}\simeq \mathrm{End}(V' \otimes V'');
$$
 \item a homomorphism $\varepsilon:A \to \mathbb{C}$, called a counit, to be
used for producing the trivial $A$-module $\mathbb{C}$;
 \item an antihomomorphism $S:A \to A$, called an antipode, to be used for
producing a dual $A$-module:
$$
(\xi f)(v)\stackrel{\mathrm{def}}{=}f(S(\xi)v),\qquad \xi \in A,\quad v \in
V,\quad f \in V^*.
$$
\end{itemize}

Of course, the definition of a Hopf algebra includes several assumptions on
$(A,\triangle,\varepsilon,S)$ which provide the habitual properties of the
operations of tensor product and passage to a dual module.

Introduce a comultiplication $\triangle:U_q \mathfrak{sl}_2 \to U_q
\mathfrak{sl}_2 \otimes U_q \mathfrak{sl}_2$, a counit $\varepsilon:U_q
\mathfrak{sl}_2 \to \mathbb{C}$, and an antipode $S:U_q \mathfrak{sl}_2 \to
U_q \mathfrak{sl}_2$ as follows:
\begin{align*}
\triangle(E)=E \otimes 1+K \otimes E,&& \triangle(F)=F \otimes K^{-1}+1
\otimes E,&& \triangle(K^{\pm 1})=K^{\pm 1}\otimes K^{\pm 1};
\\ \varepsilon(E)=\varepsilon(F)=0,&& \varepsilon(K^{\pm 1})=1;
\\ S(E)=-K^{-1}E,&& S(F)=-FK,&& S(K^{\pm 1})=K^{\mp 1}.
\end{align*}

Equip $U_q \mathfrak{sl}_2$ with an antilinear involution $*:U_q
\mathfrak{sl}_2 \to U_q \mathfrak{sl}_2$ given by
$$E^*=-KF,\qquad F^*=-EK^{-1},\qquad K^*=K.$$
It follows from the definitions that
$$
\triangle(\xi^*)=\triangle(\xi)^{*\otimes *},\qquad \xi \in U_q
\mathfrak{sl}_2.
$$
The $*$-Hopf-algebra $(U_q \mathfrak{sl}_2,*)$ we thus obtain will be
denoted by $U_q \mathfrak{su}_{1,1}$. It is a q-analogue of the
$*$-Hopf-algebra $U \mathfrak{su}_{1,1}\otimes_\mathbb{R}\mathbb{C}$. A $U_q
\mathfrak{su}_{1,1}$-module $V$ is said to be unitarizable if for some
scalar product (i.e. positive sesquilinear form) in $V$
$$
(\xi v_1,v_2)=(v_1,\xi^*v_2),\qquad \xi \in U_q \mathfrak{su}_{1,1},\quad
v_1,v_2 \in V.
$$
Obviously, a tensor product of unitarizable $U_q
\mathfrak{su}_{1,1}$-modules is a unitarizable $U_q
\mathfrak{su}_{1,1}$-module.

\bigskip

\section{Examples of $U_q \mathfrak{sl}_2$-module algebras: \boldmath
$\mathbb{C}[z]_q$, $\Lambda_q$, $A_{1,q}$}

Consider a $U_q \mathfrak{sl}_2$-module $V$ and its element $v$. This
element is called $U_q \mathfrak{sl}_2$-invariant if the linear map
$\mathbb{C}\to V$, $z \mapsto zv$, is a morphism of $U_q
\mathfrak{sl}_2$-modules (equivalently, $Ev=Fv=Hv=0$).

Consider an algebra $\mathcal{F}$ and a linear map $m:\mathcal{F}\otimes
\mathcal{F}\to \mathcal{F}$, $m:f_1 \otimes f_2 \mapsto f_1f_2$ determined
by the multiplication law in $\mathcal{F}$: $f_1 \times f_2 \mapsto f_1f_2$.
$\mathcal{F}$ is called a $U_q \mathfrak{sl}_2$-module algebra if it is
equipped with a structure of $U_q \mathfrak{sl}_2$-module and the
multiplication $m:\mathcal{F}\otimes \mathcal{F}\to \mathcal{F}$ is a
morphism of $U_q \mathfrak{sl}_2$-modules (equivalently, the following
q-analogue of the Leibnitz rule is valid:
\begin{eqnarray*}
E(f_1f_2)&=&(Ef_1)f_2+(q^Hf_1)(Ef_2),
\\ F(f_1f_2)&=&(Ff_1)(q^{-H}f_2)+f_1(Ff_2),
\\ H(f_1f_2)&=&(Hf_1)f_2+f_1(Hf_2)).
\end{eqnarray*}
Under the presence of the unit element $1 \in \mathcal{F}$ there is an
additional requirement of its $U_q \mathfrak{sl}_2$-invariance:
$E1=F1=H1=0$.

In a similar way, the notion of $U_q \mathfrak{sl}_2^{\mathrm{op}}$-module
coalgebra is introduced, with $U_q \mathfrak{sl}_2^{\mathrm{op}}$ being the
Hopf algebra deduced from $U_q \mathfrak{sl}_2$ via replacing its
comultiplication $\triangle$ with the opposite one $\triangle^{\mathrm{op}}$
(if $\triangle(\xi)=\sum \limits_j \xi_j'\otimes \xi_j''$ then
$\triangle^{\mathrm{op}}(\xi)=\sum \limits_j \xi_j''\otimes \xi_j'$).

An important example of a $U_q \mathfrak{sl}_2^{\mathrm{op}}$-module
coalgebra is built from the Verma module of zero weight $M(0)$. This module
is determined by its generator $v(0)\in M(0)$ and the relations $Ev(0)=0$,
$K^{\pm 1}v(0)=v(0)$. Obviously, there exists a unique morphism of $U_q
\mathfrak{sl}_2^{\mathrm{op}}$-modules $\delta:M(0)\to M(0)\otimes M(0)$
such that $\delta v(0)=v(0)\otimes v(0)$. Coassociativity of this
comultiplication follows from
$$((v(0)\otimes v(0))\otimes v(0)=v(0)\otimes(v(0)\otimes v(0)).$$
Note that the vector spaces $M(0)$, $M(0)\otimes M(0)$ are graded:
\begin{gather*}
M(0)=\bigoplus_jM(0)_j,\qquad M(0)_j=\{v \in M(0)|\,Hv=2jv \},\\
(M(0)\otimes M(0))_j=\bigoplus_{i+k=j}M(0)_i \otimes M(0)_k,
\end{gather*}
and that the comultiplication preserves the degree of homogeneity:
$$\triangle:M(0)_j \to (M(0)\otimes M(0))_j,$$
with $M(0)_0=\mathbb{C}v(0)$. That is why the dual graded vector space
$M(0)^*\stackrel{\mathrm{def}}{=}\bigoplus \limits_j(M(0)_j)^*$ is a unital
$\mathbb{Z}$-graded algebra. By a construction, $M(0)$ is a $U_q
\mathfrak{sl}_2^{\mathrm{op}}$-module coalgebra and $M(0)^*$ is a $U_q
\mathfrak{sl}_2$-module algebra. We are going to describe it by explicit
formulae.

The elements $v(0)$, $S(F)v(0)$, $S(F^2)v(0)$, \ldots, form a basis of the
vector space $M(0)$. Hence there exists a unique element $z \in M(0)^*$ such
that
$$
z(S(F^j)v(0))=
\begin{cases}
 q^{1/2}, & j=1, \\
 0, & j \ne 1.
\end{cases}
$$
It is easy to prove that the algebra $M(0)^*$ is isomorphic to
$\mathbb{C}[z]$ and
\begin{equation}\label{uni1k_qdbs.texac}
Fz=q^{1/2},\qquad Ez=-q^{1/2}z^2,\qquad K^{\pm 1}z=q^{\pm 2}z.
\end{equation}
Thus we have equipped the polynomial algebra $\mathbb{C}[z]$ with a
structure of $U_q \mathfrak{sl}_2$-module algebra, to be denoted
$\mathbb{C}[z]_q$ in the sequel.

Turn to producing a differential calculus on $\mathbb{C}[z]_q$. A duality
argument like that we used to derive $\mathbb{C}[z]_q$ allows one to produce
a differential $U_q \mathfrak{sl}_2$-module algebra $(\Lambda_q,d)$. In this
setting, the embeddings $\mathbb{C}[z]_q \hookrightarrow \Lambda_q$ and the
differential $d:\Lambda_q \to \Lambda_q$ are morphisms of $U_q
\mathfrak{sl}_2$-modules. One can also prove that
$$dz \cdot dz=0,\qquad dz \cdot z=q^2z \cdot dz.$$

See \cite{uni1k_qdbs.texSV} for details. Note that this differential calculus is well
known and can be derived from the Wess-Zumino differential calculus on the
quantum $\mathbb{C}^2$ ($t_1t_2=qt_2t_1$) via the localization procedure
$z=t_2^{-1}t_1$ and further restriction onto the subalgebra of zero
homogeneity degree elements.

To conclude, describe one more $U_q \mathfrak{sl}_2$-module algebra, which
is a q-analogue of the Weyl algebra $A_1$.

Consider the linear operator $\dfrac{d}{dz}$ in $\mathbb{C}[z]_q$ given by
$df=dz \cdot \left(\dfrac{d}{dz}f \right)$, $f \in \mathbb{C}[z]_q$. It is
easy to demonstrate that
\begin{equation}\label{}
\frac{d}{dz}\cdot \widehat{z}-q^{-2}\widehat{z}\cdot \frac{d}{dz}=1,
\end{equation}
with $\widehat{z}$ being the left multiplication operator by $z$
($\widehat{z}:f \mapsto zf$) in $\mathbb{C}[z]_q$. The algebra determined by
the two generators $\dfrac{d}{dz}$ and $\widehat{z}$ and the above relations
is called the quantum Weyl algebra (q-oscillatory algebra) $A_{1,q}$. Our
purpose is to equip it with a structure of a $U_q \mathfrak{sl}_2$-module
algebra. We use the fact that the algebra
$\mathrm{End}_{\mathbb{C}}(\mathbb{C}[z]_q)$ of all linear operators in
$\mathbb{C}[z]_q$ is $U_q \mathfrak{sl}_2$-module:
$$
\xi(T)=\sum_j\xi_j'\cdot T \cdot S(\xi_j'')\qquad \text{for}\quad
\triangle(\xi)=\sum_j\xi_j'\otimes \xi_j''.
$$
We are to prove that the image of $A_{1,q}$ under the canonical embedding
into $\mathrm{End}_{\mathbb{C}}(\mathbb{C}[z]_q)$ appears to be a submodule
of the $U_q \mathfrak{sl}_2$-module
$\mathrm{End}_{\mathbb{C}}(\mathbb{C}[z]_q)$. This is a consequence of the
following relations which describe the action of the generators $E$, $F$,
$K^{\pm 1}$ on the generators $\widehat{z}$, $\dfrac{d}{dz}$:
\begin{multline}\label{uni1k_qdbs.texgg}
E(\widehat{z})=-q^{1/2}\widehat{z}^{\,2},\qquad
F(\widehat{z})=q^{1/2},\qquad K^{\pm 1}(\widehat{z})=q^{\pm 2}\widehat{z},
\\
E \left(\frac{d}{dz}\right)=q^{-3/2}(q^{-2}+1)\widehat{z}\frac{d}{dz},\qquad
F \left(\frac{d}{dz}\right)=0,\qquad K^{\pm
1}\left(\frac{d}{dz}\right)=q^{\mp 2}\frac{d}{dz}.
\end{multline}
We restrict ourselves to proving the most intricate of these relations. It
follows from the definitions that
\begin{gather*}
\frac{df(z)}{dz}=\frac{f(q^{-2}z)-f(z)}{q^{-2}z-z},\qquad Ef(z)=-q^{1/2}z^2
\frac{f(z)-f(q^2z)}{z-q^2z},
\\ Ff(z)=q^{1/2}\frac{f(q^{-2}z)-f(z)}{q^{-2}z-z},\qquad K^{\pm
1}f(z)=f(q^{\pm 2}z).
\end{gather*}
Therefore,
\begin{gather*}
E \frac{d}{dz}f(z)=
-\frac{q^{1/2}}{(1-q^2)^2}(q^2f(q^{-2}z)-(1+q^2)f(z)+f(q^2z)),
\\ \frac{d}{dz}Ef(z)=
-\frac{q^{1/2}}{(1-q^2)^2}(f(q^{-2}z)-(1+q^2)f(z)+q^2f(q^2z)).
\end{gather*}
It follows from the definition of $E \left(\dfrac{d}{dz}\right)$ that
$$
E \left(\frac{d}{dz}\right)=E \frac{d}{dz}-K \frac{d}{dz}K^{-1}E=E
\frac{d}{dz}-q^{-2}\frac{d}{dz}E,
$$
hence
$$
E \left(\frac{d}{dz}\right)f(z)=
-\frac{q^{-3/2}}{(q^{-2}-1)^2}(1-q^{-4})(f(q^{-2}z)-f(z))=
q^{-3/2}(1+q^{-2})\widehat{z}\frac{df(z)}{dz}.
$$
It follows from the above observations that $A_{1,q}$ is a $U_q
\mathfrak{sl}_2$-module algebra and the action of $U_q \mathfrak{sl}_2$ in
$A_{1,q}$ is given by (\ref{uni1k_qdbs.texgg}).

\bigskip

\section{The \boldmath $U_q \mathfrak{su}_{1,1}$-module algebras
$\mathrm{Pol}(\mathbb{C})_q$, $\mathcal{D}(\mathbb{U})_q$}

Consider a $U_q \mathfrak{sl}_2$-module algebra $F$ with an involution $*$.
$F$ is called a $U_q \mathfrak{su}_{1,1}$-module algebra if the involutions
agree:
$$(\xi f)^*=(S(\xi))^*f^*,\qquad \xi \in U_q \mathfrak{su}_{1,1},\;f \in F.$$

We present below examples of such algebras.

Note first that
$$(S(E))^*=q^{-2}F,\qquad (S(F))^*=q^2E,\qquad (S(K^{\pm 1}))^*=K^{\mp 1}.$$
Equip the vector space $\mathrm{Pol}(\mathbb{C})_q=\mathbb{C}[z]_q \otimes
\mathbb{C}[z^*]_q$ with an involution: $(f_2(z)\otimes
f_2(z^*))^*=\overline{f}_2(z)\otimes \overline{f}_1(z^*)$, where bar denotes
complex conjugation for the coefficients of polynomials. The involutions in
$\mathrm{Pol}(\mathbb{C})_q$ and $U_q \mathfrak{su}_{1,1}$ agree in the
above sense. What remains is to equip $\mathrm{Pol}(\mathbb{C})_q$ with a
multiplication $m:\mathrm{Pol}(\mathbb{C})_q \otimes
\mathrm{Pol}(\mathbb{C})_q \to \mathrm{Pol}(\mathbb{C})_q$ which agree with
the action of $U_q \mathfrak{sl}_2$ and the involution $*$ in
$\mathrm{Pol}(\mathbb{C})_q$. It was demonstrated in \cite{uni1k_qdbs.texSV} that such
multiplication can be derived from the morphism of $U_q
\mathfrak{sl}_2$-modules
$$
\check{R}:\mathbb{C}[z^*]_q \otimes \mathbb{C}[z]_q \to \mathbb{C}[z]_q
\otimes \mathbb{C}[z^*]_q
$$
determined by the action of the universal R-matrix in $\mathbb{C}[z^*]_q
\otimes \mathbb{C}[z]_q$ and a subsequent ordinary flip of tensor multiples.
More precisely, $m=m_+\otimes m_-(\mathrm{id}\otimes \check{R}\otimes
\mathrm{id})$, with $m_{\pm}$ being the multiplications in $\mathbb{C}[z]_q$
and $\mathbb{C}[z^*]_q$, respectively.

To describe the action of the universal R-matrix one has well known
Drinfeld's formulae \cite{uni1k_qdbs.texDr}:
$$R=\exp_{q^2}((q^{-1}-q)E \otimes F)q^{-H \otimes H/2},$$
with $\exp_t(x)=\sum \limits_{n=1}^\infty x^n \left(\prod
\limits_{j=1}^n((1-t^j)/(1-t))\right)^{-1}$.

An application of this relation allows one to prove (see \cite{uni1k_qdbs.texSV}) that
\begin{equation}\label{uni1k_qdbs.texqd}
z^*z=q^2zz^*+1-q^2.
\end{equation}
Furthermore, $\mathrm{Pol}(\mathbb{C})_q$ can be determined by the
generators $z$, $z^*$ and the relation (\ref{uni1k_qdbs.texqd}), and the action of $E$,
$F$, $K^{\pm 1}$ on $z$, $z^*$ can be given by (\ref{uni1k_qdbs.texac}).

We need a positive invariant integral in the quantum disc. The problem is
that in the classical case ($q=1$) the positive $SU_{1,1}$-invariant measure
in the unit disc such integral is not defined on the polynomial algebra
since
$$
\int \limits_{\mathbb{U}}\left(1-|z|^2 \right)^{-2}d\mathrm{Re}z \wedge d
\mathrm{Im}z=\infty.
$$
To produce an invariant integral, one uses an extension
$\mathrm{Fun}(\mathbb{U})_q$ of the $U_q \mathfrak{su}_{1,1}$-module
$*$-algebra $\mathrm{Pol}(\mathbb{C})_q$, derived by adding an element $f_0$
with the following properties:
\begin{equation}\label{uni1k_qdbs.texcrf0}
z^*f_0=0,\qquad f_0z=0,\qquad f_0^2=f_0.
\end{equation}
The involution is extended from $\mathrm{Pol}(\mathbb{C})_q$ onto
$\mathrm{Fun}(\mathbb{U})_q$ in such a way that $f_0$ becomes a selfadjoint
element: $f_0^*=f_0$.

To motivate this definition, let us consider a faithful irreducible
$*$-representation $T$ of the $*$-algebra $\mathrm{Pol}(\mathbb{C})_q$. Such
$*$-representation is unique up to a unitary equivalence. In the standard
basis $\{e_k\}_{k=0}^\infty$ of the Hilbert space $l^2(\mathbb{Z}_+)$ it can
be given by
$$
T(z)e_k=\left(1-q^{2(k+1)}\right)^{1/2}e_{k+1},\qquad T(z^*)e_k=
\begin{cases}
\left(1-q^{2k}\right)^{1/2}e_{k-1}, & k>0,
\\ \qquad 0,& k=0.
\end{cases}
$$
It follows from the definition of the $*$-algebra
$\mathrm{Fun}(\mathbb{U})_q$ that $T$ is extendable up to a
$*$-representation $\widetilde{T}$ of $\mathrm{Fun}(\mathbb{U})_q$ and
$\widetilde{T}(f_0)$ appears to be an orthogonal projection onto the 'vacuum
subspace' $\mathbb{C}e_0$.

Note that the system of equation $z^*f_0=0$, $f_0z=0$ has a solution in the
space of formal series $\mathcal{D}(\mathbb{U})_q'=\left\{\left.\sum
\limits_{j,k=0}^\infty c_{jk}z^jz^{*k}\right|\,c_{jk}\in
\mathbb{C}\right\}$. Specifically,
$$
f_0=\sum_{k=0} ^\infty
\frac{(-1)^kq^{k(k-1)}}{(1-q^2)(1-q^4)\ldots(1-q^{2k})}z^kz^{*k}.
$$
We thus obtain an embedding of vector spaces $\mathrm{Fun}(\mathbb{U})_q
\subset \mathcal{D}(\mathbb{U})_q'$. The vector space
$\mathcal{D}(\mathbb{U})_q'$ with the coefficientwise convergence topology
will be called the distribution space in the quantum disc. One can
demonstrate that the involution and the $U_q \mathfrak{sl}_2$-action are
extendable by a continuity from $\mathrm{Pol}(\mathbb{C})_q$ onto
$\mathcal{D}(\mathbb{U})_q'$ and thus equip $\mathrm{Fun}(\mathbb{U})_q$
with a structure of $U_q \mathfrak{su}_{1,1}$-module algebra. The action of
$E$, $F$, $K^{\pm 1}$ onto $f_0$ is described by the explicit formulae:
$$
Ef_0=-\frac{q^{1/2}}{1-q^2}zf_0,\qquad
Ff_0=-\frac{q^{1/2}}{q^{-2}-1}f_0z^*,\qquad K^{\pm 1}f_0=f_0.
$$
The two-sided ideal $\mathcal{D}(\mathbb{U})_q=\mathrm{Fun}(\mathbb{U})_qf_0
\mathrm{Fun}(\mathbb{U})_q$ works as an algebra of finite functions in the
quantum disc. Note that there exists a non-degenerate pairing of
$\mathcal{D}(\mathbb{U})_q$ and $\mathcal{D}(\mathbb{U})_q'$.

Of course, the subalgebra $\mathcal{D}(\mathbb{U})_q$ is a $U_q
\mathfrak{su}_{1,1}$-module algebra. It is just the domain where a positive
$U_q \mathfrak{sl}_2$-invariant integral is to be defined.

Consider the linear span $\mathcal{L}\subset l^2(\mathbb{Z}_+)$ of the
standard basis $\{e_k\}_{k=0}^\infty$ together with a linear operator
$\Gamma$ in $\mathcal{L}$,
$$\Gamma e_k=q^{-2k}e_k,\qquad k \in \mathbb{Z}_+.$$
The linear functional
$$
\int \limits_{\mathbb{U}_q}fd \nu
\stackrel{\mathrm{def}}{=}(1-q^2)\mathrm{Tr}(\widetilde{T}(f)\Gamma)
$$
is well defined on $\mathcal{D}(\mathbb{U})_q$; it is positive and $U_q
\mathfrak{sl}_2$-invariant (more precisely,
\begin{gather*}
\int \limits_{\mathbb{U}_q}f^*fd \nu>0 \quad \text{for}\quad f \ne 0, \quad
\text{and}
\\ \int \limits_{\mathbb{U}_q}(\xi f)d \nu=\varepsilon(\xi) \int
\limits_{\mathbb{U}_q}fd \nu,\qquad f \in \mathcal{D}(\mathbb{U})_q,\quad
\xi \in U_q \mathfrak{sl}_2).
\end{gather*}

Note finally an outward similarity of the above invariant integral in
$\mathcal{D}(\mathbb{U})_q$ and the well known in the quantum group theory
invariant integral $\mathrm{tr}_q:\mathrm{End}_{\mathbb{C}}(V)\to
\mathbb{C}$, with $V$ being a $U_q \mathfrak{sl}_2$-module. This similarity
is due to the fact that $\mathcal{D}(\mathbb{U})_q$ is embedable into a
one-parameter family of $U_q \mathfrak{su}_{1,1}$-module algebras
$\mathrm{End}_{\mathbb{C}}V_t$ as a limit object. In turn, the origin of the
new deformation parameter $t$ is in the Berezin quantization \cite{uni1k_qdbs.texSSV2}.

\bigskip

\section{Differential calculi on \boldmath $\mathrm{Pol}(\mathbb{C})_q$,
$\mathrm{Fun}(\mathbb{U})_q$}

A universal R-matrix was used in the previous section to pass from the $U_q
\mathfrak{sl}_2$-module algebra $\mathbb{C}[z]_q$ to the $U_q
\mathfrak{su}_{1,1}$-module algebra $\mathrm{Pol}(\mathbb{C})_q$. A similar
construction described in \cite{uni1k_qdbs.texSV} allows one to pass from the $U_q
\mathfrak{sl}_2$-module algebra $(\Lambda_q,d)$ as in section 3 to the $U_q
\mathfrak{su}_{1,1}$-module differential algebra described below.

Let $\Omega_q$ be the algebra determined by the generators $z$, $z^*$, $dz$,
$dz^*$ and the relations
\begin{align*}
z^*z&=q^2zz^*+1-q^2,
\\  dz \cdot z&=q^2z \cdot dz, & dz^*\cdot z^*&=q^{-2}z^*\cdot dz^*,
\\  dz \cdot z^*&=q^{-2}z^*\cdot dz, & dz^*\cdot z&=q^2z \cdot dz^*,
\\  dz \cdot dz&=dz^*\cdot dz^*=0, & dz^*\cdot dz&=-q^2dz \cdot dz^*.
\end{align*}

Equip this algebra with an involution $*:z \mapsto z^*$; $*:dz \mapsto dz^*$
and a $\mathbb{Z}_2$-grading $\deg(z)=\deg(z^*)=0$, $\deg(dz)=\deg(dz^*)=1$.
There exists a unique differentiation $d$ of the superalgebra $\Omega_q$
such that $d^2=0$ and $d:z \mapsto dz$, $d:z^*\mapsto dz^*$. The $*$-algebra
$\mathrm{Pol}(\mathbb{C})_q$ is uniquely embedable into $\Omega_q$. There
exists a unique extension of the structure of $U_q
\mathfrak{su}_{1,1}$-module algebra from $\mathrm{Pol}(\mathbb{C})_q$ onto
$\Omega_q$ in such a way that the differential $d$ appears to be a morphism
of $U_q \mathfrak{su}_{1,1}$-modules.

The $*$-algebra $\mathrm{Fun}(\mathbb{U})_q$ was derived from
$\mathrm{Pol}(\mathbb{C})_q$ via adding the element $f_0$ which satisfies
(\ref{uni1k_qdbs.texcrf0}). Besides that, one should also include to the definition of the
algebra of differential forms the relations $f_0 \cdot dz=dz \cdot f_0$,\ \
$f_0 \cdot dz^*=dz^*\cdot f_0$; the definition of the differential $d$
should also include the equality
$$df_0=-\frac1{1-q^2}(dzf_0z^*+zf_0dz^*).$$
In \cite[section 5]{uni1k_qdbs.texSSV3} these relations and the latter equality were
obtained via the embedding $\mathrm{Fun}(\mathbb{U})_q \subset
\mathcal{D}(\mathbb{U})_q'$. Also, it was proved there that the structure of
$U_q \mathfrak{su}_{1,1}$-module algebra is uniquely extendable from
$\mathrm{Fun}(\mathbb{U})_q$ onto the above differential $*$-algebra.

Note that, similarly to the case $q=1$, one has a bigrading of the algebra
of differential forms:
$$
\deg(z)=\deg(z^*)=\deg(f_0)=(0,0),\qquad \deg(dz)=(1,0),\qquad
\deg(dz^*)=(0,1).
$$
The differential $d$ admits a decomposition
$$d=\partial+\overline{\partial}$$
as a sum of holomorphic and antiholomorphic differentials, to be defined in
the standard way. For example,
$$
\partial f_0=-\frac1{1-q^2}dz \cdot f_0 \cdot z^*,\qquad \overline{\partial}
f_0=-\frac1{1-q^2}z \cdot f_0 \cdot dz^*.
$$
Of course, $\partial^2=\partial
\overline{\partial}+\overline{\partial}\partial=\overline{\partial}^2=0$.

Define the operators $\frac{\partial}{\partial z}$,
$\frac{\partial}{\partial z^*}$ in $\mathrm{Pol}(\mathbb{C})_q$ by
$$
\partial f=dz \cdot \frac{\partial f}{\partial z},\qquad
\overline{\partial}f=dz^*\cdot \frac{\partial f}{\partial z^*},
$$
and the operators $\widehat{z}$, $\widehat{z}^*$ by $\widehat{z}f=zf$,
$\widehat{z}^*f=z^*f$. These linear operators are extendable by a continuity
onto the entire space of distributions $\mathcal{D}(\mathbb{U})_q'$ and
satisfies the following commutation relations:
\begin{align*}
&\widehat{z}^*\widehat{z}=q^2\widehat{z}\widehat{z}^*+1-q^2,&&
\\ &\frac{\partial}{\partial
z}\widehat{z}-q^{-2}\widehat{z}\frac{\partial}{\partial z}=1,&&
\frac{\partial}{\partial z}\widehat{z}^*=q^2
\widehat{z}^*\frac{\partial}{\partial z},
\\ &\frac{\partial}{\partial
z^*}\widehat{z}=q^{-2}\widehat{z}\frac{\partial}{\partial z^*},&&
\frac{\partial}{\partial z^*}\widehat{z}^*-q^2
\widehat{z}^*\frac{\partial}{\partial z^*}=1,
\\ &\frac{\partial}{\partial z^*}\frac{\partial}{\partial z}=q^2
\frac{\partial}{\partial z}\frac{\partial}{\partial z^*}.&&
\end{align*}

These commutation relations follow from
\begin{align*}
&z^*z=q^2zz^*+1-q^2,&
\\ &\partial(zf)=dz \cdot f+zdz \frac{\partial f}{\partial z}=dz
\left(f+q^{-2}z \frac{\partial f}{\partial z}\right),&
\\ &\partial(z^*f)=z^*dz \frac{\partial f}{\partial z}=q^2dz\cdot
z^*\frac{\partial f}{\partial z},&
\\ &\overline{\partial}(zf)=zdz^*\frac{\partial f}{\partial
z^*}=q^{-2}dz^*\cdot z \cdot \frac{\partial f}{\partial z^*},&
\\ &\overline{\partial}(z^*f)=dz^*\cdot f+z^*\cdot dz^*\cdot \frac{\partial
f}{\partial z^*}=dz^*\left(f+q^2z^*\frac{\partial f}{\partial z^*}\right),&
\\ &\partial \overline{\partial}f=\partial \left(dz^*\frac{\partial
f}{\partial z^*}\right)=-dz^*dz \frac{\partial}{\partial
z}\frac{\partial}{\partial z^*}f=q^2dzdz^*\frac{\partial}{\partial
z}\frac{\partial}{\partial z^*}f,&
\\ &\overline{\partial}\partial f=\overline{\partial}\left(dz \frac{\partial
f}{\partial z}\right)=-dzdz^*\frac{\partial}{\partial
z^*}\frac{\partial}{\partial z}f,&\partial
\overline{\partial}+\overline{\partial}\partial=0.
\end{align*}

We thus obtain the algebra of differential operators in the quantum disc. An
example of such operator is a q-analogue of the Laplace operator:
$(1-zz^*)^2 \dfrac{\partial}{\partial z}\dfrac{\partial}{\partial z^*}$.

\bigskip

\section{Concluding notes}

The results described above constitute a background for producing
non-commutative function theory in the quantum disc \cite{uni1k_qdbs.texSSV1}. On the
other hand, these results admit an extension onto the case of generic
q-Cartan domains introduced in \cite{uni1k_qdbs.texSV}.

Applications to non-commutative function theory use essentially the fact
that $f_0$ generate the $U_q \mathfrak{sl}_2$-module
$\mathcal{D}(\mathbb{U})_q$. This allows one to reduce proofs of some
relations in algebras of intertwining operators to their verification on
$f_0$.

A great deal of the algebras we consider in this work were mentioned in the
literature (sometimes in different contexts). Nevertheless, most of works on
those algebras either do not mention presence of a $U_q
\mathfrak{sl}_2$-module structure or do not use involution and positive
invariant integral. That is why we treat the authors of \cite{uni1k_qdbs.texKL,uni1k_qdbs.texNN} as our
direct predecessors.

Our approach to constructing the differential calculus is inspired by the
classical result which establishes a duality between the covariant
differential operators and morphisms of generalized Verma modules (see
\cite[\S 11.1]{uni1k_qdbs.texBaE}, \cite{uni1k_qdbs.texHarr}).

Among other works on quantum algebras, related to this text, one should
mention \cite{uni1k_qdbs.texCK,uni1k_qdbs.texFl,uni1k_qdbs.texKMT,uni1k_qdbs.texJak}.




\makeatletter \@addtoreset{equation}{section} \makeatother

\renewcommand{\theequation}{\thesection.\arabic{equation}}

\title{\bf QUANTUM DISC: THE CLIFFORD ALGEBRA AND THE DIRAC OPERATOR}
\author{K. Schm\"udgen \and S. Sinel'shchikov \and L. Vaksman}
\date{}

\newpage
\setcounter{section}{0}
\large

\makeatletter
\renewcommand{\@oddhead}{QUANTUM DISC: THE CLIFFORD ALGEBRA AND THE DIRAC
OPERATOR \hfill \thepage}
\renewcommand{\@evenhead}{\thepage \hfill K. Schm\"udgen, S. Sinel'shchikov,
and L. Vaksman}
\let\@thefnmark\relax
\@footnotetext{This research was supported in part by of the US Civilian
Research \& Development Foundation (Award No UM1-2091) and by Swedish
Academy of Sciences (project No 11293562).}

\addcontentsline{toc}{chapter}{\@title \\ {\sl K. Schm\"udgen, S.
Sinel'shchikov, and L. Vaksman}\dotfill} \makeatother

\maketitle

\section{Introduction}

In his non-commutative geometry A. Connes \cite{uni1l_qdkad.TEXCo} has developed methods
and notions in order to define and to study "non-commutative Riemannian
manifolds". A crucial role in this theory is played by the Dirac operators
and its non-commutative analogs.

One of the simplest (non-compact) Riemannian manifolds of a constant
curvature is the unit disc $\mathbb{U}\subset \mathbb{C}$ whose metric is
given by
$$\left(1-x^2-y^2 \right)^{-2}dxdy$$
(the Poincar\'e model of Lobachevskian geometry). A non-commutative analogue
of this Riemannian manifold was extensively studied recently within the
framework of the quantum group theory \cite{uni1l_qdkad.TEXKlSch}. The aim of the present
paper is to develop the Clifford algebra and the Dirac operator for this
quantum disc.

\bigskip

\section{A q-Clifford bundle}

In what follows $q \in(0,1)$ and $\mathbb{C}$ is the ground field. We assume
some acquaintance with the basics of quantum group theory \cite{uni1l_qdkad.TEXKlSch,uni1l_qdkad.TEXJant}. All algebras in this section are supposed to be unital.

We use the standard notation for the quantum universal enveloping algebra
$U_q \mathfrak{sl}_2$, its generators $E$, $F$, $K$, $K^{-1}$, the
comultiplication $\triangle$, the counit $\varepsilon$, and the antipode $S$
\cite{uni1l_qdkad.TEXJant,uni1l_qdkad.TEXSSV1}. Throughout we shall use the Sweedler notation
$\triangle(f)=f_{(1)}\otimes f_{(2)}$ instead of $\triangle(f)=\sum
\limits_i f_i' \otimes f_i''$. The $*$-Hopf algebra $U_q
\mathfrak{su}_{1,1}=(U_q \mathfrak{sl}_2,*)$ is defined by
$$E^*=-KF,\qquad F^*=-EK^{-1},\qquad \left(K^{\pm 1}\right)^*=K^{\pm 1}.$$

Primarily we are interested in $U_q \mathfrak{sl}_2$-module algebras and
$U_q \mathfrak{su}_{1,1}$-module $*$-algebras\footnote{A $U_q
\mathfrak{sl}_2$-module $*$-algebra $F$ is called $U_q
\mathfrak{su}_{1,1}$-module algebra if the following compatibility condition
on ivolutions is valid:
$$
(\xi f)^*=(S(\xi))^*f^*,\qquad \xi \in U_q \mathfrak{su}_{1,1},\quad f \in
F.
$$}.
We follow \cite{uni1l_qdkad.TEXSSV1} in starting with the $*$-algebra
$\mathrm{Pol}(\mathbb{C})_q$ with single generator $z$ and defining relation
\begin{equation}\label{uni1l_qdkad.TEXqd}
z^*z=q^2zz^*+1-q^2.
\end{equation}
The structure of the $U_q \mathfrak{su}_{1,1}$-module algebra used by the
authors in \cite{uni1l_qdkad.TEXSSV1} attracted an attention about 10 years ago
\cite{uni1l_qdkad.TEXKlLe}. It can be determined by
\begin{equation}\label{}
Fz=q^{1/2},\qquad Ez=-q^{1/2}z^2,\qquad K^{\pm 1}z=q^{\pm 2}z.
\end{equation}
We follow \cite{uni1l_qdkad.TEXSSV1} in recalling the description of a covariant
differential calculus over the $U_q \mathfrak{su}_{1,1}$-module algebra
$\mathrm{Pol}(\mathbb{C})_q$ which initially appeared in \cite{uni1l_qdkad.TEXSch-Sch}.
This is a $U_q \mathfrak{su}_{1,1}$-module algebra $\Omega_q \supset
\mathrm{Pol}(\mathbb{C})_q$ equipped with a (super)differentiation
$d:\Omega_q \to \Omega_q$ which commutes with an involution $*$ and is a
morphism of $U_q \mathfrak{su}_{1,1}$-modules.

The list of relations which determine the $U_q \mathfrak{su}_{1,1}$-module
algebra $\Omega_q$ consists of the above relations determining the $U_q
\mathfrak{su}_{1,1}$-module algebra $\mathrm{Pol}(\mathbb{C})_q$, and the
additional relations
\begin{align}\label{uni1l_qdkad.TEXdzzcr}
&dz \cdot z=q^2 z \cdot dz,&& dz \cdot z^*=q^{-2}z^*\cdot dz,
\\ \label{uni1l_qdkad.TEXdzdzcr} &dz \cdot dz=0,&& dz^*\cdot dz^*=0,
\\ &dz^*\cdot dz+q^2dz \cdot dz^*=0.
\end{align}

In the classical case $q=1$ the passage from the algebra of exterior
differential forms with polynomial coefficients to the algebra of polynomial
sections of the Clifford bundle (associated to an invariant metric) reduces
to a replacement of the relations
$$
d \overline{z}dz+dzd \overline{z}=0 \qquad \left(\Leftrightarrow \frac{d
\overline{z}dz+dzd \overline{z}}{(1-|z|^2)^2}=0\right)
$$
and
$$
d \overline{z}dz+dzd \overline{z}=(1-|z|^2)^2 \qquad \left(\Leftrightarrow
\frac{d \overline{z}dz+dzd \overline{z}}{(1-|z|^2)^2}=1 \right).
$$
A similar argument for $q \in(0,1)$ leads to the following definition of the
'algebra of polynomial sections' of a q-Clifford bundle.

\medskip

{\sc Definition.} Let $\mathrm{Cl}_q$ denote the $*$-algebra with generators
$z$, $dz$ and determining relations (\ref{uni1l_qdkad.TEXqd}), (\ref{uni1l_qdkad.TEXdzzcr}),
(\ref{uni1l_qdkad.TEXdzdzcr}), and
\begin{equation}\label{}
dz^*\cdot dz+q^2dz \cdot dz^*=(1-zz^*)^2.
\end{equation}

\medskip

We equip $\mathrm{Cl}_q$ with the standard filtration given by
$\deg(z)=\deg(z^*)=0$, $\deg(dz)=\deg(dz^*)=1$.

It follows from the definition of $\mathrm{Cl}_q$ that, just as in the
classical case $q=1$, one has

\medskip

\begin{proposition}\label{uni1l_qdkad.TEXcl*}
The graded algebra of $\mathrm{Cl}_q$ associated to its standard filtration
is isomorphic to $\Omega_q$.
\end{proposition}

\medskip

Consider the $U_q \mathfrak{sl}_2$-module subalgebras
\begin{gather*}
\Omega_q^{(*,0)}=\{f+dz \cdot g|\,f,g \in \mathrm{Pol}(\mathbb{C})_q
\}\subset \Omega_q,
\\ \Omega_q^{(0,*)}=\{f+g \cdot dz^*|\,f,g \in \mathrm{Pol}(\mathbb{C})_q
\}\subset \Omega_q.
\end{gather*}

As a consequence of the definition of $\mathrm{Cl}_q$ one has
\begin{equation}\label{uni1l_qdkad.TEXhah}
\Omega_q^{(*,0)}\hookrightarrow \mathrm{Cl}_q,\qquad
\Omega_q^{(0,*)}\hookrightarrow \mathrm{Cl}_q.
\end{equation}

\medskip

\begin{proposition}\label{uni1l_qdkad.TEXuqsuma}
There exists a unique structure of $U_q \mathfrak{su}_{1,1}$-module algebra
in $\mathrm{Cl}_q$ such that the embeddings (\ref{uni1l_qdkad.TEXhah}) appear to be
morphisms of $U_q \mathfrak{sl}_2$-modules.
\end{proposition}

\smallskip

{\bf Proof.} The only non-trivial statement here is about the existence of a
structure of $U_q \mathfrak{sl}_2$-module algebra in $\mathrm{Cl}_q$ such
that the embeddings (\ref{uni1l_qdkad.TEXhah}) are morphisms of $U_q
\mathfrak{sl}_2$-modules.

Consider the $*$-algebra $\mathcal{F}$ with generators $z$, $dz$ and
determining relations (\ref{uni1l_qdkad.TEXqd}), (\ref{uni1l_qdkad.TEXdzzcr}), (\ref{uni1l_qdkad.TEXdzdzcr}).
$\mathcal{F}$ is a $U_q \mathfrak{su}_{1,1}$-module algebra, and the
canonical homomorphism $\mathcal{F}\to \Omega_q$ is a morphism of $U_q
\mathfrak{su}_{1,1}$-module algebras. One has to prove that the two-sided
ideal $J$ of $\mathcal{F}$ generated by $dz^*\cdot dz+q^2dz \cdot
dz^*-(1-zz^*)^2$, is a submodule of the $U_q \mathfrak{su}_{1,1}$-module
$\mathcal{F}$. Obviously, $K^{\pm 1}J \subset J$, one needs only to prove
that $EJ \subset J$, $FJ \subset J$. The proofs of these inclusions are
similar. We restrict ourselves to demonstrating the first of them. It
follows from the following statement.

\medskip

\begin{lemma}\label{uni1l_qdkad.TEXej}\hfill

i) $E(dz dz^*)=-q^{1/2}\left(1+q^2 \right)zdz dz^*$,

ii) $E(dz^*dz)=-q^{1/2}\left(1+q^2 \right)zdz^*dz$,

iii) $E\left((1-zz^*)^2 \right)=-q^{1/2}\left(1+q^2 \right)z(1-zz^*)^2$.
\end{lemma}

\smallskip

{\bf Proof} of lemma \ref{uni1l_qdkad.TEXej}.
\begin{multline*}
i)\;E(dzdz^*)=(Edz)dz^*=d(Ez)dz^*=d(-q^{1/2}z^2)dz^*=-q^{1/2}(dz \cdot z+z
\cdot dz)dz^*=
\\ =-q^{1/2}\left(1+q^2 \right)zdzdz^*.
\end{multline*}
$$
ii)\; E(dz^*dz)=(Kdz^*)(Edz)=q^{-2}dz^*\left(-q^{1/2}\left(1+q^2\right)zdz
\right)=-q^{1/2}\left(1+q^2 \right)zdz^*dz.
$$
\ \ \ \ {\it iii)} Note first that $Ez^*=q^{-3/2}$. (In fact,
$(Fz)^*=(S(F))^*z^*$. Hence
$q^{1/2}=(-FK)^*z^*=-K^*\left(-EK^{-1}\right)z^*=q^2Ez^*$). Furthermore,
$$E(1-zz^*)=q^{1/2}zz^*-q^2zq^{-3/2}=-q^{1/2}z(1-zz^*).$$
Finally,
$$
E \left((1-zz^*)^2 \right)=-q^{1/2}z(1-zz^*)^2-q^{1/2}(1-zz^*)z(1-zz^*)=
-q^{1/2}\left(1+q^2 \right)z(1-zz^*)^2. \eqno \square
$$

\medskip

{\sc Remark.} Propositions \ref{uni1l_qdkad.TEXcl*} and \ref{uni1l_qdkad.TEXuqsuma} suggest that an
invariant Riemannian metric in the quantum disc should be defined by an
expession of the form
\begin{equation}\label{uni1l_qdkad.TEXirm}
(1-zz^*)^{-2}(dz^*\otimes dz+q^2dz \otimes dz^*).
\end{equation}

To begin with, consider the multiplicative closed subset
$(1-zz^*)^{\mathbb{N}}\subset \mathrm{Pol}(\mathbb{C})_q \subset \Omega_q$.
It follows from the relations
\begin{align*}
&(1-zz^*)z=q^2z(1-zz^*),&&(1-zz^*)z^*=q^{-2}z^*(1-zz^*),
\\ &(1-zz^*)dz=dz(1-zz^*),&&(1-zz^*)dz^*=dz^*(1-zz^*),
\end{align*}
that this multiplicative closed system is an Ore set. Consider the
associated localizations $\widetilde{\mathrm{Pol}}(\mathbb{C})_q$ and
$\widetilde{\Omega}_q$ of $\mathrm{Pol}(\mathbb{C})_q$ and $\Omega_q$,
respectively. It is well-known that there embeddings
$\mathrm{Pol}(\mathbb{C})_q \hookrightarrow
\widetilde{\mathrm{Pol}}(\mathbb{C})_q$ and $\Omega_q \hookrightarrow
\widetilde{\Omega}_q$. One can also prove just as in \cite{uni1l_qdkad.TEXSStV} that the
structure of $U_q \mathfrak{su}_{1,1}$-module algebra is uniquely extended
from $\Omega_q$ onto $\widetilde{\Omega}_q$.

Certainly, the grading mentioned above is canonically extendable onto
$\widetilde{\Omega}_q$ and
$$
\widetilde{\Omega}_q^1=\left \{f'dz+f''dz^*\left|\,f',f''\in
\widetilde{\mathrm{Pol}}(\mathbb{C})_q \right.\right \}
$$
is a $U_q \mathfrak{sl}_2$-module
$\widetilde{\mathrm{Pol}}(\mathbb{C})_q$-bimodule, equipped with the
involution $*$ inherited from $\widetilde{\Omega}_q$.

Now we are in a position to introduce the $U_q \mathfrak{sl}_2$-module
$\widetilde{\Omega}_q
\otimes_{\widetilde{\mathrm{Pol}}(\mathbb{C})_q}\widetilde{\Omega}_q$, with
a 'metric tensor' (\ref{uni1l_qdkad.TEXirm}) chosen among the $U_q
\mathfrak{sl}_2$-invariant elements of this module.

\bigskip

\section{A q-Dirac bundle}

Consider the $U_q \mathfrak{sl}_2$-module algebra of 'differential forms'
$\Omega_q^{(0,*)}$, together with its $U_q \mathfrak{sl}_2$-module
subalgebra $\mathrm{Pol}(\mathbb{C})_q$. Obviously, $\Omega_q^{(0,*)}$ is a
free (left) $\mathrm{Pol}(\mathbb{C})_q$-module with generators $1$, $dz^*$.
Of course, $\Omega_q^{(0,*)}$ is a $U_q \mathfrak{sl}_2$-module
$\mathrm{Pol}(\mathbb{C})_q$-module (a q-analogue of the space of polynomial
sections of the vector bundle $T^{(0,*)}\mathbb{C}$).

In the classical case ($q=1$) the Clifford algebra acts in the fibers of
$T^{(0,*)}\mathbb{C}$, which allows us to speak about the Dirac bundle. We
are going to describe a quantum analogue of this action. It follows from the
definition of $\mathrm{Cl}_q$ that, in the category of $U_q
\mathfrak{sl}_2$-module $\mathrm{Pol}(\mathbb{C})_q$-modules there exists a
canonical isomorphism
$$
\Omega_q^{(0,*)}\simeq \mathrm{Cl}_q/\mathrm{Cl}_q \cdot dz.
$$
In particular, $\Omega_q^{(0,*)}$ is a $U_q \mathfrak{sl}_2$-module
$\mathrm{Cl}_q$-module.

Consider a sesquilinear map $\mathbf{h}:\Omega_q^{(0,*)}\times
\Omega_q^{(0,*)}\to \mathrm{Pol}(\mathbb{C})_q$,
$$
\mathbf{h}(\varphi_0+dz^*\varphi_1,\psi_0+dz^*\psi_1)=
\psi_0^*\varphi_0+q^{-2}\psi_1^*(1-zz^*)^2\varphi_1,
$$
for $\varphi_0,\varphi_1,\psi_0,\psi_1 \in \mathrm{Pol}(\mathbb{C})_q$.

We are going to demonstrate that the 'Hermitian metric' $\mathbf{h}$ is
positive definite, invariant and 'respects' the $\mathrm{Cl}_q$-action. That
is, it can be used to produce a q-analogue of the Dirac operator.

\medskip

\begin{proposition}\label{uni1l_qdkad.TEXbh}\hfill

i) $\mathbf{h}(\omega,\omega)\ge 0$ for all $\omega \in \Omega_q^{(0,*)}$
and \ $\mathbf{h}(\omega,\omega)=0 \Rightarrow \omega=0$,

ii) $\xi \mathbf{h}(\omega_1,\omega_2)=
\mathbf{h}(\xi_{(2)}(\omega_1),(S(\xi_{(1)}))^*\omega_2)$ for all $\xi \in
U_q \mathfrak{su}_{1,1}$, $\omega_1,\omega_2 \in \Omega_q^{(0,*)}$.
\footnote{Recall that we use the Sweedler notation. The condition ii) is
just a rephrasing of a $U_q \mathfrak{su}_{1,1}$-invariance of the Hermitian
metric $\mathbf{h}$. The flip of tensor multiples $\xi_{(1)}$, $\xi_{(2)}$
in ii) is due to the fact that a similar flip is normally used implicitly
while constructing scalar products in the spaces of functions and
differential forms. For example, for functions
$$
f_1 \otimes f_2 \mapsto f_2 \otimes f_1 \mapsto f_2^* \otimes f_1 \mapsto
f_2^*f_1 \mapsto \int f_2^*f_1 d \nu.
$$}
\begin{equation}\label{uni1l_qdkad.TEXcinv}
iii)\; \mathbf{h}(c
\omega_1,\omega_2)=\mathbf{h}(\omega_1,c^*\omega_2)\qquad \text{for
all}\quad c \in \mathrm{Cl}_q,\quad \omega_1,\omega_2 \in
\Omega_q^{(0,*)}.\qquad \qquad \qquad
\end{equation}
\end{proposition}

\smallskip

{\bf Proof.} The first statement is evident. To prove the second statement,
it suffices to establish $U_q \mathfrak{su}_{1,1}$-invariance of the
Hermitian forms
\begin{align*}
&\mathbf{h}_1:\mathrm{Pol}(\mathbb{C})_q \times \mathrm{Pol}(\mathbb{C})_q
\to \mathrm{Pol}(\mathbb{C})_q,&&\mathbf{h}_1:f_1 \times f_2: \mapsto
f_2^*f_1,
\\ &\mathbf{h}_2:\Omega_q^{(0,1)}\times \Omega_q^{(0,1)}\to
\mathrm{Pol}(\mathbb{C})_q,&&\mathbf{h}_2:dz^*f_1 \times dz^*f_2 \mapsto
f_2^*(1-zz^*)^2f_1.
\end{align*}

On the other hand, since $\Omega_q$ is a $U_q \mathfrak{su}_{1,1}$-module
algebra, the forms $\mathbf{h}_1$ and
$$
\mathbf{h}_3:\Omega_q^{(0,1)}\times \Omega_q^{(0,1)}\to
\Omega_q^{(1,1)},\qquad \mathbf{h}_3:\omega \times \omega \mapsto
\omega^*\omega
$$
are $U_q \mathfrak{su}_{1,1}$-invariant. (Here
$\Omega_q^{(1,1)}=\{fdzdz^*|\,f \in \mathrm{Pol}(\mathbb{C})_q \}\subset
\Omega_q$.) The $U_q \mathfrak{su}_{1,1}$-invariance of $\mathbf{h}_2$
follows from the invariance of $\mathbf{h}_3$ and the invariance of the
element
\begin{equation}\label{uni1l_qdkad.TEXkf}
(1-zz^*)^{-2}dzdz^*=dzdz^*(1-zz^*)^{-2}\in \widetilde{\Omega}_q.
\end{equation}

One can use a similar argument to prove (\ref{uni1l_qdkad.TEXcinv}) in the special case $c
\in \mathrm{Pol}(\mathbb{C})_q$. (The passage from $\mathbf{h}_3$ to
$\mathbf{h}_2$ involves the fact that the 'K\"ahler form' (\ref{uni1l_qdkad.TEXkf}) is in
the center of $\widetilde{\Omega}_q$.)

We prove (\ref{uni1l_qdkad.TEXcinv}) in the general case. In view of the relations
$$
\mathbf{h}(\omega_1f,\omega_2)=\mathbf{h}(\omega_1,\omega_2)f,\qquad
\mathbf{h}(\omega_1,\omega_2f)=f^*\mathbf{h}(\omega_1,\omega_2),\qquad
\omega_1,\omega_2 \in \Omega_q^{(0,*)},\;f \in \mathrm{Pol}(\mathbb{C})_q,
$$
it remains only to consider the case $\omega_1,\omega_2 \in \{1,dz^* \}$,
$c=dz^*$. In this case (\ref{uni1l_qdkad.TEXcinv}) follows from
$$\mathbf{h}(dzdz^*,1)=\mathbf{h}(1,dzdz^*)=q^{-2}(1-zz^*)^2,$$
$$\mathbf{h}(dz^*,dz^*)=q^{-2}(1-zz^*)^2.$$
The first relation here follows from the fact that in the
$\mathrm{Cl}_q$-module $\Omega_q^{(0,*)}\simeq \mathrm{Cl}_q/\mathrm{Cl}_q
\cdot dz$
$$dz \cdot dz^*=q^{-2}(1-zz^*)^2,$$
while the second one is obvious. \hfill $\square$

\bigskip

\section{An important $*$-representation of \boldmath $\mathrm{Cl}_q$}

Every positive linear functional $\mu:\mathrm{Pol}(\mathbb{C})_q \to
\mathbb{C}$ determines a $*$-representation of $\mathrm{Cl}_q$ in the
pre-Hilbert space $\Omega_q^{(0,*)}$ with scalar product
\begin{equation}\label{uni1l_qdkad.TEXsp}
(\omega_1,\omega_2)=\int \mathbf{h}(\omega_1,\omega_2)d \mu.
\end{equation}

However, to produce a Dirac operator we need a $U_q
\mathfrak{su}_{1,1}$-invariant scalar product. We first recall the
corresponding definition.

Consider a $U_q \mathfrak{su}_{1,1}$-module $V$ and a sesquilinear form $V
\times V \to \mathbb{C}$, $v_1 \times v_2 \mapsto(v_1,v_2)$. This
sesquilinear form is called $U_q \mathfrak{su}_{1,1}$-invariant if for all
$\xi \in U_q \mathfrak{su}_{1,1}$, $v_1,v_2 \in V$
$$
\varepsilon(\xi)\cdot(v_1,v_2)=\left(\xi_{(2)}v_1,(S(\xi_{(1)}))^*v_2
\right).
$$

The invariance of a sesquilinear map $\mathbf{h}$ indicates a naive way of
producing the required scalar product. That is, formula (\ref{uni1l_qdkad.TEXsp}) suggests
to define a $U_q \mathfrak{su}_{1,1}$-invariant integral $\mu$ by
$$
\int(\xi f)d \mu=\varepsilon(\xi)\int fd \mu,\qquad \xi \in U_q
\mathfrak{su}_{1,1},\quad f \in \mathrm{Pol}(\mathbb{C})_q.
$$
In fact this fails even in the classical case ($q=1$), since the positive
$SU_{1,1}$-invariant measure in the unit disc has the form
$$d \nu=\mathrm{const}\left(1-|z|^2 \right)^{-2}dzd \overline{z},$$
and the associated integral is not defined on the polynomial algebra (unless
$\mathrm{const}=0$). That is why in producing an integral calculus, the
algebra $\mathcal{D}(\mathbb{U})$ of smooth functions with support inside
the disc $\mathbb{U}$ is used. In the quantum case one can also replace the
$U_q \mathfrak{su}_{1,1}$-module algebra $\mathrm{Pol}(\mathbb{C})_q$ by the
$U_q \mathfrak{su}_{1,1}$-module algebra $\mathcal{D}(\mathbb{U})_q$ to
produce a $U_q \mathfrak{su}_{1,1}$-invariant integral
$\nu:\mathcal{D}(\mathbb{U})_q \to \mathbb{C}$ (see \cite{uni1l_qdkad.TEXSSV1}).

Recall the definition of the $U_q \mathfrak{su}_{1,1}$-module algebra
$\mathcal{D}(\mathbb{U})_q$. Consider the $*$-algebra
$\mathrm{Fun}(\mathbb{U})_q$ with generators $z$, $f_0$ and determining
relations (\ref{uni1l_qdkad.TEXqd}) and
\begin{equation}\label{uni1l_qdkad.TEXzf0}
z^*f_0=f_0z=0,\qquad f_0^2=f_0,\qquad f_0^*=f_0.
\end{equation}
Obviously, we have an embedding $\mathrm{Pol}(\mathbb{C})_q \hookrightarrow
\mathrm{Fun}(\mathbb{U})_q$.

An argument used in \cite{uni1l_qdkad.TEXSSV1} leads to (\ref{uni1l_qdkad.TEXzf0}) and to the relations
\begin{equation}\label{uni1l_qdkad.TEXuqf0}
Ef_0=-\frac{q^{1/2}}{1-q^2}zf_0,\qquad
Ff_0=-\frac{q^{1/2}}{q^{-2}-1}f_0z^*,\qquad K^{\pm 1}f_0=f_0.
\end{equation}
This argument also allows us to prove that the structure of $U_q
\mathfrak{su}_{1,1}$-module algebra admits an extension from
$\mathrm{Pol}(\mathbb{C})_q$ to $\mathrm{Fun}(\mathbb{U})_q$ which satisfies
(\ref{uni1l_qdkad.TEXuqf0}). The uniqueness of this extension is obvious. The $U_q
\mathfrak{su}_{1,1}$-module algebra $\mathcal{D}(\mathbb{U})_q$ is defined
as a two-sided ideal of $\mathrm{Fun}(\mathbb{U})_q$ generated by $f_0$:
$$
\mathcal{D}(\mathbb{U})_q=\mathrm{Fun}(\mathbb{U})_q \cdot f_0 \cdot
\mathrm{Fun}(\mathbb{U})_q.
$$

Recall the explicit form of the invariant integral $\nu$.

Let $Q$ be the linear operator in $l^2(\mathbb{Z}_+)$ given by
$Qe_k=\left(1-q^2 \right)q^{-2k}e_k$, with $\{e_k \}$ being the standard
basis in $l^2(\mathbb{Z}_+)$. Consider the $*$-representation $T$ of
$\mathrm{Fun}(\mathbb{U})_q$ in the Hilbert space $l^2(\mathbb{Z}_+)$ given
by
\begin{gather*}
T(z)e_k=\left(1-q^{2(k+1)}\right)^{1/2}e_{k+1},\qquad T(z^*)e_k=
\begin{cases}
\left(1-q^{2k}\right)^{1/2}e_{k-1}, & k>0,
\\ \hfill 0 \hfill,& k=0.
\end{cases},
\\ T(f_0)e_k=
\begin{cases}
e_k, & k=0,
\\ 0, & k \ne 0.
\end{cases}
\end{gather*}

Note that the only non-zero matrix element of the operator $T(z^jf_0z^{*k})$
is in line $j$ column $k$. Hence the operators $T(z^jf_0z^{*k})$ are linear
independent. From this one deduces that the elements $z^jf_0z^{*k}$, $j,k
\in \mathbb{Z}_+$, form a basis in $\mathcal{D}(\mathbb{U})_q$ and that the
representation $T$ of $\mathcal{D}(\mathbb{U})_q$ is faithful.

It is easy to prove (see \cite{uni1l_qdkad.TEXSSV2}):

\medskip

\begin{proposition}\label{uni1l_qdkad.TEXii}
The linear functional
$$
\int \limits_{\mathbb{U}_q}fd \nu
\stackrel{\mathrm{def}}{=}\mathrm{tr}(T(f)Q),\qquad f \in
\mathcal{D}(\mathbb{U})_q,
$$
on the $U_q \mathfrak{su}_{1,1}$-module $*$-algebra
$\mathcal{D}(\mathbb{U})_q$ is well defined, positive definite, and $U_q
\mathfrak{su}_{1,1}$-invariant.
\end{proposition}

\medskip

We have thus derived a $U_q \mathfrak{su}_{1,1}$-module algebra
$\mathrm{Fun}(\mathbb{U})_q$ from $\mathrm{Pol}(\mathbb{C})_q$ by adding a
generator and additional relations (\ref{uni1l_qdkad.TEXzf0}), (\ref{uni1l_qdkad.TEXuqf0}).

In a similar way one can extend the $U_q \mathfrak{su}_{1,1}$-module
$*$-algebras $\Omega_q$ and $\mathrm{Cl}_q$ (see \cite{uni1l_qdkad.TEXSSV1}). For that, one
has to add $f_0$ to the list of generators and to complete the list of
relations with (\ref{uni1l_qdkad.TEXzf0}), (\ref{uni1l_qdkad.TEXuqf0}) and
\begin{equation}\label{}
f_0dz=dz \cdot f_0,\qquad f_0 \cdot dz^*=dz^*\cdot f_0.
\end{equation}
Let $\Omega(\mathbb{U})_q$, $\mathrm{Cl}(\mathbb{U})_q$ be the two-sided
ideals of those algebras generated by $f_0$. We call $\Omega(\mathbb{U})_q$
the algebra of differential forms with finite coefficients in the quantum
disc, and $\mathrm{Cl}(\mathbb{U})_q$ the algebra of finite sections of the
q-Clifford bundle.

Obviously, $\Omega(\mathbb{U})_q$ is a $U_q \mathfrak{su}_{1,1}$-module
$\Omega_q$-bimodule and $\mathrm{Cl}(\mathbb{U})_q$ is a $U_q
\mathfrak{su}_{1,1}$-module $\mathrm{Cl}_q$-bimodule. As usual,
\begin{align}\label{uni1l_qdkad.TEXouq}
\Omega(\mathbb{U})_q^{(0,1)}&=\{f+dz^*g|\,f,g \in \mathcal{D}(\mathbb{U})_q
\},\nonumber
\\ \Omega(\mathbb{U})_q^{(0,1)}&=
\mathrm{Cl}(\mathbb{U})_q/\mathrm{Cl}(\mathbb{U})_q \cdot dz.
\end{align}

It follows from (\ref{uni1l_qdkad.TEXouq}) that $\Omega(\mathbb{U})_q^{(0,1)}$ is a
$\mathrm{Cl}_q$-module. We intend to equip $\Omega(\mathbb{U})_q^{(0,1)}$
with an invariant Hermitian metric and an invariant scalar product.

Consider a sesquilinear map
$\mathbf{h}_{\mathbb{U}}:\Omega(\mathbb{U})_q^{(0,*)}\times
\Omega(\mathbb{U})_q^{(0,*)}\to \mathcal{D}(\mathbb{U})_q$,
\begin{equation}
\mathbf{h}_{\mathbb{U}}(\varphi_0+dz^*\cdot \varphi_1,\psi_0+dz^*\cdot
\psi_1)=\psi_0^*\phi_0+q^{-2}\psi_1^*(1-zz^*)^2 \varphi_1,
\end{equation}
with $\varphi_0,\varphi_1,\psi_0,\psi_1 \in \mathcal{D}(\mathbb{U})_q$.


\medskip

\begin{proposition}\label{uni1l_qdkad.TEXbh_}\hfill

i) $\mathbf{h}_{\mathbb{U}}(\omega,\omega)\ge 0$ for all $\omega \in
\Omega(\mathbb{U})_q^{(0,*)}$ and $\mathbf{h}_{\mathbb{U}}(\omega,\omega)=0
\Rightarrow \omega=0$,

ii) $\xi \mathbf{h}_{\mathbb{U}}(\omega_1,\omega_2)=
\mathbf{h}_{\mathbb{U}}(\xi_{(2)}\omega_1,(S(\xi_{(1)}))^*\omega_2)$ for all
$\xi \in U_q \mathfrak{su}_{1,1}$, $\omega_1,\omega_2 \in
\Omega(\mathbb{U})_q^{(0,*)}$.

iii) $\mathbf{h}_{\mathbb{U}}(c
\omega_1,\omega_2)=\mathbf{h}_{\mathbb{U}}(\omega_1,c^*\omega_2)$ for all $c
\in \mathrm{Cl}_q$, $\omega_1,\omega_2 \in \Omega(\mathbb{U})_q^{(0,*)}$.
\end{proposition}

\medskip

This proposition can be proved in a similar manner as proposition \ref{uni1l_qdkad.TEXbh}.

Propositions \ref{uni1l_qdkad.TEXii} and \ref{uni1l_qdkad.TEXbh_} allow us to equip
$\Omega(\mathbb{U})_q^{(0,*)}$ with a $U_q \mathfrak{su}_{1,1}$-invariant
scalar product: $(\omega_1,\omega_2)\stackrel{\mathrm{def}}{=}\int
\limits_{\mathbb{U}_q}\mathbf{h}_{\mathbb{U}}(\omega_1,\omega_2)d \nu$.

Consider the Hilbert space $\mathcal{H}_q$ which is the completion of the
pre-Hilbert space $\Omega(\mathbb{U})_q^{(0,*)}$.

The boundedness of the operators $T(f)$, $f \in \mathrm{Pol}(\mathbb{C})_q$,
and the explicit form of the invariant integral
$\nu:\mathcal{D}(\mathbb{U})_q \to \mathbb{C}$ imply the following

\medskip

\begin{proposition}
The linear operators $\Omega(\mathbb{U})_q^{(0,*)}\to
\Omega(\mathbb{U})_q^{(0,*)}$, $\omega \mapsto c \omega$, are bounded for
all $c \in \mathrm{Cl}_q$.
\end{proposition}

\medskip

\begin{corollary}
There exists a unique $*$-representation $\pi$ of $\mathrm{Cl}_q$ (by
bounded operators) in $\mathcal{H}_q$ such that for all $c \in
\mathrm{Cl}_q$, $\omega \in \Omega(\mathbb{U})_q^{(0,*)}$,
$$\pi(c)\omega=c \omega.$$
\end{corollary}

\medskip

We close this section by recalling that the structure of $U_q
\mathfrak{su}_{1,1}$-module (super)differential algebra can be canonically
transferred from $\Omega_q$ onto $\Omega(\mathbb{U})_q^{(0,*)}$ (see
\cite{uni1l_qdkad.TEXSSV1}):
$$
df_0=-\frac1{1-q^2}(dz \cdot f_0 \cdot z^*+z \cdot f_0 \cdot dz^*).
$$

\bigskip

\section{The Dirac operator in the quantum disc}

The Dirac operator is normally defined in terms of special affine
connections. Nevertheless, in the case of K\"ahler manifold, this operator
admits a definition by a simple explicit formula (see, for example,
\cite{uni1l_qdkad.TEXMich}). In the special case of the unit disc $\mathbb{U}\subset
\mathbb{C}$ and the Clifford bundle as it is considered in this work, the
Dirac operator in $\Omega(\mathbb{U})^{(0,*)}=
\Omega(\mathbb{U})^{(0,0)}+\Omega(\mathbb{U})^{(0,1)}$ is defined by
\begin{equation}\label{uni1l_qdkad.TEXD}
{\not}\mathcal{D}=\begin{pmatrix}0 & (\overline{\partial})^*\\
\overline{\partial} & 0\end{pmatrix},
\end{equation}
where $(\overline{\partial})^*$ being a formal adjoint differential operator
(with respect to the standard $SU_{1,1}$-invariant scalar products in
$\Omega(\mathbb{U})^{(0,0)}$, $\Omega(\mathbb{U})^{(0,1)}$). The operator
${\not}\mathcal{D}$ is essentially selfadjoint in the Hilbert space
completion $\mathcal{H}$ of $\Omega(\mathbb{U})^{(0,*)}$.

It is worth to note that
$$
{\not}\mathcal{D}^2=\begin{pmatrix}\square_0 & 0 \\ 0 & \square_1
\end{pmatrix},
$$
with $\square_0=(\overline{\partial})^*\overline{\partial}$,
$\square_1=\overline{\partial}(\overline{\partial})^*$ being essentially
selfadjoint $SU_{1,1}$-invariant differential operators (invariant Laplace
operators in the spaces of functions and differential type $(0,1)$ forms,
respectively).

Now we are going to pass from the classical disc to the quantum disc. Just
as in the classical case, one has the decompositions
$$
\Omega(\mathbb{U})_q=\bigoplus_{i,j=0}^1 \Omega(\mathbb{U})_q^{(i,j)},\qquad
d=\partial+\overline{\partial},
$$
\begin{flalign*}
\text{with} && &\partial:\Omega(\mathbb{U})_q^{(i,j)}\to
\Omega(\mathbb{U})_q^{(i+1,j)},&&
\\ && &\overline{\partial}:\Omega(\mathbb{U})_q^{(i,j)}\to
\Omega(\mathbb{U})_q^{(i,j+1)}&&
\end{flalign*}
(these decompositions are related to the standard bigrading of the algebra
of differential forms: $\deg(z)=\deg(z^*)=\deg(f_0)=(0,0)$,
$\deg(dz)=(1,0)$, $\deg(dz^*)=(0,1)$).

In some sense the quantum case turns out to be easier than the classical
one, because we have

\medskip

\begin{proposition}\label{uni1l_qdkad.TEXbns}
The linear operator $\overline{\partial}:\Omega(\mathbb{U})_q^{(0,0)}\to
\Omega(\mathbb{U})_q^{(0,1)}$ is bounded.
\end{proposition}
See a proof in \cite{uni1l_qdkad.TEXSSV3}.

\medskip

Let $\mathcal{H}_q^{(0)}$ and $\mathcal{H}_q^{(1)}$ be the closures of the
linear manifolds $\Omega(\mathbb{U})_q^{(0,0)}=\mathcal{D}(\mathbb{U})_q$
and $\Omega(\mathbb{U})_q^{(0,1)}$, respectively, in $\mathcal{H}_q$. From
proposition \ref{uni1l_qdkad.TEXbns} it follows that the linear map
$\overline{\partial}:\Omega(\mathbb{U})_q^{(0,0)}\to
\Omega(\mathbb{U})_q^{(0,1)}$ has an extension by a continuity up to a
bounded operator $\overline{\partial}:\mathcal{H}_q^{(0)}\to
\mathcal{H}_q^{(1)}$. Now the decomposition
$\mathcal{H}_q=\mathcal{H}_q^{(0)}\oplus \mathcal{H}_q^{(1)}$ allows us to
introduce the Dirac operator in the quantum disc by (\ref{uni1l_qdkad.TEXD}). As we have
shown, in the quantum case the operator ${\not}\mathcal{D}$ is bounded in
$\mathcal{H}_q$.

Note that ${\not}\mathcal{D}={\not}\mathcal{D}^*$, and that the Dirac
operator is an endomorphism of the $U_q \mathfrak{su}_{1,1}$-module
$\Omega(\mathbb{U})_q^{(0,*)}$, by this property of the operator
$\overline{\partial}$ and the $U_q \mathfrak{su}_{1,1}$-invariance of the
scalar products.

Introduce the standard notation $\dfrac{\partial}{\partial z}$,
$\dfrac{\partial}{\partial z^*}$ for the linear operators defined by
$$
\partial\psi=dz \frac{\partial \psi}{\partial z},\qquad
\overline{\partial}\psi=dz^*\frac{\partial \psi}{\partial z^*}.
$$

\medskip

\begin{proposition}
For all $\psi \in \mathcal{D}(\mathbb{U})_q$
$$
(\overline{\partial})^*:dz^*\psi \mapsto q^2(1-zz^*)^2 \frac{\partial
\psi}{\partial z}.
$$
\end{proposition}

\smallskip

{\bf Proof.} It suffices to prove that the linear operators $q^2
\dfrac{\partial}{\partial z}$, $\dfrac{\partial}{\partial z^*}$ are formal
adjoints with respect to the scalar product in $\mathcal{D}(\mathbb{U})_q$
$$
\varphi \times \psi \mapsto \int
\limits_{\mathbb{U}_q}\psi^*\varphi(1-zz^*)^2d \nu,
$$
determined by the 'q-Lebesgue measure' $(1-zz^*)^2d \nu$. The latter
assertion can be proved in the same way as the correspoding statement in
\cite{uni1l_qdkad.TEXSSV4}.

\begin{corollary}
For all $\varphi,\psi \in \mathcal{D}(\mathbb{U})_q$
\begin{gather*}
{\not}\mathcal{D}(\varphi+dz^*\psi)=q^2(1-zz^*)^2 \frac{\partial
\psi}{\partial z}+dz^*\frac{\partial \varphi}{\partial z^*},
\\ (\overline{\partial})^*\overline{\partial}\varphi=q^2(1-zz^*)^2
\frac{\partial}{\partial z}\frac{\partial}{\partial z^*},
\\ \overline{\partial}(\overline{\partial})^*(dz^*\psi)=
q^2dz^*\frac{\partial}{\partial z^*}(1-zz^*)^2 \frac{\partial}{\partial
z^*}.
\end{gather*}
\end{corollary}

\bigskip

\section{Concluding remarks}

In A. Connes' non-commutative Riemannian geometry \cite[chapter VI]{uni1l_qdkad.TEXCo} the
differential $d$ is closely related to the Dirac operator
${\not}\mathcal{D}$. Specifically, for any 'function' $f$, the elements $f$
and $df$ are in the Clifford algebra. They act in the same space as
${\not}\mathcal{D}$, and
\begin{equation}\label{uni1l_qdkad.TEXCC}
df \cdot \psi=[{\not}\mathcal{D},f]\cdot \psi.
\end{equation}
A similar relation is also valid in the case of quantum disc as one can see
from the following obvious

\medskip

\begin{proposition}
For all $f \in \mathrm{Pol}(\mathbb{C})_q$, $\psi \in
\mathcal{D}(\mathbb{U})_q$
\begin{equation}\label{uni1l_qdkad.TEXpar}
(\overline{\partial}f)\cdot \psi=[{\not}\mathcal{D},f]\cdot \psi.
\end{equation}
\end{proposition}

\medskip

Of course, one can replace the $\mathrm{Cl}_q$-module $\Omega_q^{(0,*)}$ by
$\mathrm{Cl}_q$ treated as a $\mathrm{Cl}_q$-module, in order to pass from
$\overline{\partial}$ to $d$, i.e. for getting the relation (\ref{uni1l_qdkad.TEXCC})
instead of (\ref{uni1l_qdkad.TEXpar}). This is just what stands for the Dirac bundle in
\cite[chapter 5A]{uni1l_qdkad.TEXMich}.

There exists a unique Hermitian metric $\mathbf{h}$ which has all the
properties stated in proposition \ref{uni1l_qdkad.TEXbh} and satisfies the relation
$\mathbf{h}(dz,dz)=\mathbf{h}(dz^*,dz^*)$:
\begin{multline*}
\mathbf{h}(\varphi_{00}+dz
\varphi_{10}+dz^*\varphi_{01}+dzdz^*\varphi_{11},\psi_{00}+dz
\psi_{10}+dz^*\psi_{01}+dzdz^*\psi_{11})=
\\ =\psi_{00}^*\varphi_{00}+\frac1{1+q^2}\left(\psi_{10}^*(1-zz^*)^2
\varphi_{10}+\psi_{01}^*(1-zz^*)^2
\varphi_{01}\right)+\frac1{q^2(1+q^2)}\psi_{11}^*(1-zz^*)^4 \varphi_{11},
\end{multline*}
with $\varphi_{ij},\psi_{ij}\in \mathrm{Pol}(\mathbb{C})_q$.

{\sc Remark.} The $U_q \mathfrak{su}_{1,1}$-invariance of the Hermitian
metric $\mathbf{h}$ follows from the $U_q \mathfrak{su}_{1,1}$-invariance of
the following elements:
\begin{gather*}
(1-zz^*)^{-2}dzdz^* \in \widetilde{\Omega}_q,\footnotemark
\\ (1-zz^*)^{-2}dz \otimes dz^* \in
\widetilde{\Omega}_q^{(1,0)}\mathop{\otimes}
\nolimits_{\mathrm{Pol}(\mathbb{C})_q}\widetilde{\Omega}_q^{(0,1)},
\\ (1-zz^*)^{-2}dz^*\otimes dz \in
\widetilde{\Omega}_q^{(0,1)}\mathop{\otimes}
\nolimits_{\mathrm{Pol}(\mathbb{C})_q}\widetilde{\Omega}_q^{(1,0)}.
\end{gather*}
\footnotetext{Note that the tilde symbol stands for the localization with
respect to the multiplicative closed system \hbox{$(1-zz^*)^{2
\mathbb{N}}$.}}

In fact, one can use the two latter elements in order to define the
following morphisms of $U_q \mathfrak{su}_{1,1}$-modules:
\begin{align*}
&\widetilde{\Omega}_q^{(0,1)}\mathop{\otimes}
\nolimits_{\widetilde{\mathrm{Pol}}(\mathbb{C})_q}
\widetilde{\Omega}_q^{(1,0)}\to
\widetilde{\mathrm{Pol}}(\mathbb{C})_q,\qquad \psi dz^* \otimes dz \varphi
\mapsto \psi(1-zz^*)^2 \varphi;
\\ &\widetilde{\Omega}_q^{(1,0)}\mathop{\otimes}
\nolimits_{\widetilde{\mathrm{Pol}}(\mathbb{C})_q}
\widetilde{\Omega}_q^{(0,1)}\to
\widetilde{\mathrm{Pol}}(\mathbb{C})_q,\qquad \psi dz \otimes dz^* \varphi
\mapsto \psi(1-zz^*)^2 \varphi;
\\ &\widetilde{\Omega}_q^{(1,0)}\mathop{\otimes}
\nolimits_{\widetilde{\mathrm{Pol}}(\mathbb{C})_q}
\widetilde{\Omega}_q^{(0,1)}
\mathop{\otimes}\nolimits_{\widetilde{\mathrm{Pol}}(\mathbb{C})_q}
\widetilde{\Omega}_q^{(1,0)}\mathop{\otimes}
\nolimits_{\widetilde{\mathrm{Pol}}(\mathbb{C})_q}
\widetilde{\Omega}_q^{(0,1)}\to \widetilde{\mathrm{Pol}}(\mathbb{C})_q,
\\ &\psi dz \otimes dz^* \otimes dz \otimes dz^*\varphi \mapsto
\psi(1-zz^*)^4 \varphi.
\end{align*}
What remains to elaborate is the fact that the corresponding tensor algebra
is a $U_q \mathfrak{su}_{1,1}$-module algebra.

In sections 2 -- 5 we have considered a more simple case ($\Omega_q^{(0,*)}$
instead of $\Omega_q$) in order to reduce the relevant computations and to
make the exposition more plausible.

It should be noted that A. Connes assumed in his theory that the spectrum of
the Dirac operator is discrete. This fails both for classical and quantum
discs because of the non-compactness of $\mathbb{U}$. Even more, the
spectrum of $(\overline{\partial})^*\overline{\partial}$ in the space of
square summable functions is absolutely continuous in the classical case
\cite{uni1l_qdkad.TEXHe} and in the quantum case \cite{uni1l_qdkad.TEXSSV3} as well. In particular, the
spectrum of our Dirac operator is not discrete.




\makeatletter \@addtoreset{equation}{section} \makeatother

\renewcommand{\theequation}{\thesection.\arabic{equation}}

\title{\bf ON UNIQUENESS OF COVARIANT DEFORMATION WITH SEPARATION OF
VARIABLES OF THE QUANTUM DISC}
\author{D. Shklyarov}
\date{}

\newpage
\setcounter{section}{0}
\large

\makeatletter
\renewcommand{\@oddhead}{ON UNIQUENESS OF DEFORMATION OF THE QUANTUM DISC
\hfill \thepage}
\renewcommand{\@evenhead}{\thepage \hfill D. Shklyarov}
\let\@thefnmark\relax
\addcontentsline{toc}{chapter}{\@title \\ {\sl D. Shklyarov}\dotfill}
\@footnotetext{This research was supported in part by of the US Civilian
Research \& Development Foundation (Award No UM1-2091) and by Swedish
Academy of Sciences (project No 11293562).} \makeatother

\maketitle

\section{Introduction}

The concept of deformation quantization was introduced around 1977 by Bayen,
Flato, Fronsdal, Lichnerowicz, and Sternheimer \cite{uni1M_UNIQ.TEXBFFLS}. Since then the
concept has become very popular. In this approach quantization means a
deformation of the usual product of smooth function on a 'phase space' into
a noncommutative associative product $\star_t$ (star-product) with
additional properties.

The problem of constructing star-products explicitly is still of importance.
For some special symplectic manifolds star-products can be constructed by
using the famous Berezin quantization method \cite{uni1M_UNIQ.TEXB}. Specifically, the
method can be explored for a wide class of Kahler manifolds. The
star-product on a Kahler manifold constructed by means of the Berezin scheme
possesses some remarkable properties. First, it respects the complex
structure of the Kahler manifold. Second, it respects the action of the
group of biholomorphic automorphisms of the Kahler manifold.

In general, it seems to be important to look for those star-products on a
symplectic manifold which keep some additional geometric structures on this
manifold. In this connection, we want to mention the result on complete
classification of star-products with separation of variables on Kahler
manifolds \cite{uni1M_UNIQ.TEXKar} and results concerning invariant
star-products on homogeneous symplectic manifolds \cite{uni1M_UNIQ.TEXCG1},
\cite{uni1M_UNIQ.TEXCG2}.

A remarkable class of Kahler manifolds to which the Berezin method is
applicable is the class of bounded symmetric domains. The simplest example
is the unit disc within complex plane $\mathbb{C}$. The star-products on
bounded symmetric domains arising from the Berezin method were studied in
\cite{uni1M_UNIQ.TEXMo} in the simplest case of the disc and in \cite{uni1M_UNIQ.TEXCGR} in general
setting.

Recently, in the framework of the quantum group theory $q$-analogues of
bounded symmetric domains have been constructed \cite{uni1M_UNIQ.TEXSV}. In \cite{uni1M_UNIQ.TEXSSV} we
used a $q$-analog of the Berezin method to produce a deformation of product
in a noncommutative algebra of 'functions on the quantum unit disc' (see
section 2 for definitions). This deformation respects both the complex
structure and the quantum group symmetry action (precise definitions are to
be found in section 3). The main result of this paper is Theorem \ref{uni1M_UNIQ.TEXmain}
which states that the deformation constructed in \cite{uni1M_UNIQ.TEXSSV} is essentially
the unique deformation possessing the above properties (specifically, any
deformation of that kind can be obtained from the Berezin one via change of
parameter).

{\bf Acknowledgement.} I thank L. Vaksman for posing the problem and
fruitful discussions, and E. Karolinsky for many valuable remarks on
improving the text. Unfortunately, I didn't succeed to take into account all
his remarks, and I hope to do it in forthcoming versions of this paper.

\section{The algebra $\mathrm {Pol}({\mathbb C})_q$}

Everywhere in the sequel we suppose that $q\in (0,1)$ and the ground field
is the field ${\mathbb C}$ of complex numbers.

Denote by $\mathrm {Pol}({\mathbb C})_q$ the involutive unital algebra given
by its generator $z$ and the unique commutation relation
\begin{equation}\label{uni1M_UNIQ.TEXpol}
z^*z=q^2zz^*+1-q^2.
\end{equation}

This algebra was studied in \cite{uni1M_UNIQ.TEXNaNi}. It serves as the initial object in
constructing function theory in the quantum unit disc (see \cite{uni1M_UNIQ.TEXD1}).

The algebra $\mathrm {Pol}({\mathbb C})_q$ may be endowed with an important
extra structure, namely, the structure of a $U_q\mathfrak{su}_{1,1}$-module
algebra. Let us recall one what it means.

To start with, remind the definition of the quantum universal enveloping
algebra $U_q\mathfrak{sl}_{2}$ and its "real form" $U_q\mathfrak{su}_{1,1}$.
$U_q \mathfrak{sl}_2$ is a Hopf algebra over ${\mathbb C}$ determined by the
generators $K,K^{-1}, E, F$ and the relations $$KK^{-1}=K^{-1}K=1,\quad
K^{\pm 1}E=q^{\pm 2}EK^{\pm 1},\quad K^{\pm 1}F=q^{\mp 2}FK^{\pm 1},$$
$$EF-FE=\frac{K-K^{-1}}{q-q^{-1}},$$ $$\Delta(K^{\pm 1})=K^{\pm
1}\otimes K^{\pm 1},\quad \Delta(E)=E \otimes 1+K \otimes E,\quad
\Delta(F)=F \otimes K^{-1}+1 \otimes F.$$

 Note that
$$\varepsilon(E)=\varepsilon(F)=\varepsilon(K^{\pm 1}-1)=0,$$
$$S(K^{\pm 1})=K^{\mp 1},\quad S(E)=-K^{-1}E,\quad S(F)=-FK,$$
with $\varepsilon:U_q \mathfrak{sl}_2 \rightarrow{\mathbb C}$ and $S:U_q
\mathfrak{sl}_2 \rightarrow U_q \mathfrak{sl}_2$ being the counit and the
antipode of $U_q \mathfrak{sl}_2$, respectively. Equip $U_q \mathfrak{sl}_2$
with the involution given on the generators by
\begin{equation}\label{uni1M_UNIQ.TEXinvol}
 E^*=-KF,\quad F^*=-EK^{-1},\quad (K^{\pm 1})^*=K^{\pm 1}.
\end{equation}
The pair $(U_q \mathfrak{sl}_2,*)$ is the $*$-Hopf algebra. It is denoted by
$U_q \mathfrak{su}_{1,1}$.

Let ${\cal F}$ stands for a unital algebra equipped also with a structure of
$U_q \mathfrak{sl}_2$-module. ${\cal F}$ is called a $U_q
\mathfrak{sl}_2$-module algebra if the multiplication $m:{\cal F}\otimes
{\cal F}\to {\cal F}$ is a morphism of $U_q\mathfrak{sl}_2$-modules and for
any $\xi\in U_q \mathfrak{sl}_2$
\begin{equation}\label{uni1M_UNIQ.TEXunit}
\xi(1)=\varepsilon(\xi)\cdot1
\end{equation}
(in other words, the unit $1$ of the algebra ${\cal F}$ is $U_q
\mathfrak{sl}_2$-invariant).

An involutive algebra ${\cal F}$ is said to be a $U_q
\mathfrak{su}_{1,1}$-module algebra if it is a $U_q
\mathfrak{sl}_{2}$-module one and
\begin{equation}\label{uni1M_UNIQ.TEXsogl}
(\xi f)^*=(S(\xi))^* f^*
\end{equation}
for any $\xi\in U_q \mathfrak{su}_{1,1}$ and $f\in {\cal F}$ (the first
asterisk in the right hand side of (\ref{uni1M_UNIQ.TEXsogl}) means the involution
(\ref{uni1M_UNIQ.TEXinvol})).

We have explained the meaning of the term '$U_q \mathfrak{su}_{1,1}$-module
algebra'. The following statement is proved in \cite{uni1M_UNIQ.TEXSV}.

\medskip

\begin{proposition}
There exists a unique structure of $U_q \mathfrak{su}_{1,1}$-module algebra
in $\mathrm {Pol}({\mathbb C})_q$ such that
\begin{equation}\label{uni1M_UNIQ.TEXaction1}
 K^{\pm1}z=q^{\pm2}z,\quad Ez=-q^{1/2}z^2,\quad Fz=q^{1/2}.
\end{equation}
\end{proposition}
\medskip
Note that uniqueness of $U_q \mathfrak{su}_{1,1}$-module algebra structure
in $\mathrm {Pol}({\mathbb C})_q$ satisfying (\ref{uni1M_UNIQ.TEXaction1}) is a simple
consequence of definitions. Indeed, after application of the involution to
the both hand sides of equalities (\ref{uni1M_UNIQ.TEXaction1}) we find that the following
relations hold (see (\ref{uni1M_UNIQ.TEXsogl}))
\begin{equation}\label{uni1M_UNIQ.TEXaction2}
 K^{\pm1}z^*=q^{\mp2}z^*,\quad Ez^*=q^{-3/2},\quad Fz=-q^{-5/2}z^{*2}.
\end{equation}
Relations (\ref{uni1M_UNIQ.TEXunit}), (\ref{uni1M_UNIQ.TEXaction1}), and (\ref{uni1M_UNIQ.TEXaction2}) allow one to
apply $E$, $F$, $K^{\pm1}$ to any element of $\mathrm {Pol}({\mathbb C})_q$.
This implies the uniqueness. Existence of $U_q \mathfrak{su}_{1,1}$-module
algebra structure in $\mathrm {Pol}({\mathbb C})_q$ given by (\ref{uni1M_UNIQ.TEXaction1})
is a more difficult fact and we don't adduce its proof (an analogous result
is proved in \cite{uni1M_UNIQ.TEXSV} in a much more general setting).

Let $\mathbb{C}[z]_q$ and $\mathbb{C}[z^*]_q$ stand for the unital
subalgebras in $\mathrm {Pol}({\mathbb C})_q$ generated by $z$ and $z^*$,
respectively. In view of (\ref{uni1M_UNIQ.TEXaction1}) and (\ref{uni1M_UNIQ.TEXaction2}), these
subalgebras in the $U_q \mathfrak{sl}_{2}$-module algebra $\mathrm
{Pol}({\mathbb C})_q$ are $U_q \mathfrak{sl}_{2}$-module algebras themselves
(however, these subalgebras are not $U_q \mathfrak{su}_{1,1}$-module
subalgebras).

The defining relation (\ref{uni1M_UNIQ.TEXpol}) allows one to rewrite any element $f$ of
the algebra $\mathrm {Pol}({\mathbb C})_q$ in the form of a linear
combination of "normally ordered" monomials: $$
f=\sum_{i,j}a_{ij}z^iz^{*j},\qquad a_{ij}\in \mathbb{C}. $$ Thus, there is a
natural isomorphism of vector spaces
\begin{equation}\label{uni1M_UNIQ.TEXisom}
\mathrm {Pol}({\mathbb C})_q\simeq\mathbb{C}[z]_q\otimes\mathbb{C}[z^*]_q.
\end{equation}
Evidently, this isomorphism respects the action of $U_q \mathfrak{sl}_{2}$.

In conclusion of this section, let us comment formulas (\ref{uni1M_UNIQ.TEXaction1}) and
(\ref{uni1M_UNIQ.TEXaction2}). As we mentioned in Introduction, the unit disc is a
homogeneous space of the group $SU(1,1)$. Formulas (\ref{uni1M_UNIQ.TEXaction1}) and
(\ref{uni1M_UNIQ.TEXaction2}) are just $q$-analogs of the corresponding 'infinitesimal'
action of $U \mathfrak{su}_{1,1}$ in the space of polynomials.

\bigskip
\section{Deformations of the algebra $\mathrm {Pol}({\mathbb C})_q$.\\
Real, and covariant deformations}

First of all, let us explain what the term 'deformation of the algebra
$\mathrm {Pol}({\mathbb C})_q$' means. Let $\mathbb{C}[[t]]$ and $\mathrm
{Pol}({\mathbb C})_q[[t]]$ stand for the algebra of formal power series in
$t$ over $\mathbb{C}$ and the $\mathbb{C}[[t]]$-module of formal series in
$t$ over $\mathrm {Pol}({\mathbb C})_q$, respectively. According to the
traditional approach \cite{uni1M_UNIQ.TEXGer}, a deformation of the algebra $\mathrm
{Pol}({\mathbb C})_q$ is an associative $\mathbb{C}[[t]]$-bilinear product
(called star-product) $$ \star_t: \mathrm {Pol}({\mathbb C})_q[[t]]\times
\mathrm {Pol}({\mathbb C})_q[[t]]\rightarrow\mathrm {Pol}({\mathbb
C})_q[[t]] $$ given for $f,g\in\mathrm {Pol}({\mathbb C})_q\subset\mathrm
{Pol}({\mathbb C})_q[[t]]$ by
\begin{equation}\label{uni1M_UNIQ.TEXstar}
f\star_tg=f\cdot g+t\cdot m_1(f,g)+t^2\cdot m_2(f,g)+\ldots
\end{equation}
with $m_i:\mathrm {Pol}({\mathbb C})_q\times\mathrm {Pol}({\mathbb
C})_q\rightarrow \mathrm {Pol}({\mathbb C})_q$ being some ${\mathbb
C}$-bilinear maps. The product of arbitrary elements $f,g\in\mathrm
{Pol}({\mathbb C})_q[[t]]$ is defined via the conditions of ${\mathbb
C}$-bilinearity and $t$-adic continuity.

In subsequent sections we consider deformations of certain class only. In
this section we present all necessary definitions.

A deformation $\star_t$ of the algebra $\mathrm {Pol}({\mathbb C})_q$ is
said to be deformation with separation of variables (see
\cite{uni1M_UNIQ.TEXKar, uni1M_UNIQ.TEXBor}) if
$$f\star_tg=f\cdot g$$
with $f\in{\mathbb C}[z]_q$ or $g\in{\mathbb C}[z^*]_q$.

\medskip

{\bf Remark 1.} Deformations with separation of variables possess the
following property: any such deformation can be reconstructed from the
series $z^*\star_tz$.

\medskip

A deformation $\star_t$ of the algebra $\mathrm {Pol}({\mathbb C})_q$ is
called real if for any $f,g\in\mathrm {Pol}({\mathbb C})_q$ $$
(f\star_tg)^*=g^*\star_tf^* $$ (we suppose the formal parameter $t$ is real:
$t^*=t$). In this case the algebra $\mathrm {Pol}({\mathbb C})_q[[t]]$ is an
involutive algebra.

Finally, let us define the notion of a covariant deformation of the algebra
$\mathrm {Pol}({\mathbb C})_q$. Its classical counterpart is the notion of a
$SU(1,1)$-invariant deformation of the unit disc. Endow $\mathrm
{Pol}({\mathbb C})_q[[t]]$ with a structure of $U_q\mathfrak{sl}_2$-module
via ${\mathbb C}[[t]]$-linearity (and $t$-adic continuity). A deformation
$\star_t$ of the algebra $\mathrm {Pol}({\mathbb C})_q$ is said to be
covariant if the algebra $\mathrm {Pol}({\mathbb C})_q[[t]]$ is a
$U_q\mathfrak{sl}_2$-module algebra with respect to the product $\star_t$
(see section 2 for definitions). \footnote{Note that if a deformation
$\star_t$ is real and covariant then the algebra $\mathrm {Pol}({\mathbb
C})_q[[t]]$ with the product $\star_t$ {\it automatically} become a
$U_q\mathfrak{su}_{1,1}$-module algebra.}

{\it Further we restrict ourselves of studying real covariant deformations
with separation of variables (abbreviated to RCW-deformations) of the
algebra $\mathrm {Pol}({\mathbb C})_q$ only}.

\bigskip
\section{Examples of RCW-deformations of the algebra \boldmath
$\mathrm{Pol}({\mathbb C})_q$}

The simplest example of a RCW-deformation is the trivial deformation (the
trivial deformation is the deformation $\star_t$ such that $f\star_tg=f\cdot
g$ for any $f,g\in\mathrm {Pol}({\mathbb C})_q$). The remaining part of this
section is devoted to a more substantial example, namely, the Berezin
deformation. This deformation is constructed in \cite{uni1M_UNIQ.TEXSSV} by a q-analogue
of the Berezin quantization method. This approach is applicable to all
q-Cartan domains \cite{uni1M_UNIQ.TEXSV}. In the simplest case of quantum disc this method
could be replaced by a plausible procedure, to be sketched below.

Let us begin with some non-rigorous arguments. Suppose $t\in[0,1)$. Consider
the unital involutive algebra given by its generator and the following
relation \cite{uni1M_UNIQ.TEXKL}:
\begin{equation}\label{uni1M_UNIQ.TEXdefpol}
z^*z=q^2zz^*+1-q^2+\frac{(1-q^2)t}{1-q^2t}(1-z^*z)(1-zz^*).
\end{equation}
We denote it by $\mathrm {Pol}({\mathbb C})_{q,t}$. Note that application of
the involution to the both hand sides of (\ref{uni1M_UNIQ.TEXdefpol}) leads to the
relation
\begin{equation}\label{uni1M_UNIQ.TEXdefpol1}
z^*z=q^2zz^*+1-q^2+\frac{(1-q^2)t}{1-q^2t}(1-zz^*)(1-z^*z).
\end{equation}
In particular, the elements $1-zz^*$ and $1-z^*z$ commute.



We describe the procedure of producing the Berezin deformation \cite{uni1M_UNIQ.TEXSSV}
starting from the family of algebras $\mathrm {Pol}({\mathbb C})_{q,t}$.
Consider a polynomials $f$ in (non-commuting) variables $z$, $z^*$ and try
to put it in normal order using relation (\ref{uni1M_UNIQ.TEXdefpol}) ("to put in normal
order" means to rewrite as a linear combination of the monomials
$z^iz^{*j}$, $i,j\in{\mathbb Z}_+$). The procedure of normal ordering is
just the iterative process: at every subsequent step of this process we look
for all the monomials containing combination $z^*z$ at least once; then we
replace precisely one such combination with the right hand side of
(\ref{uni1M_UNIQ.TEXdefpol}) in each found monomial. In general, we can't put the
polynomial $f$ in normal order in a finite number of steps. But it is not
difficult to observe that for an arbitrarily large $N$ any polynomial $f$
can be written in a finite number of steps in the form
\begin{equation}\label{uni1M_UNIQ.TEXupor}
f=f_0+tf_1+t^2f_2+\ldots+t^Nf_N+t^{N+1}f_{N+1}
\end{equation}
with $f_0$, $f_1$,...,$f_N$ being normally ordered polynomials independent
of $t$. Thus, after "infinite number" of steps we should have a formal
series
$$
f_0+tf_1+t^2f_2+\ldots+t^Nf_N+\ldots\in\mathrm {Pol}({\mathbb C})_{q}[[t]]
$$ whose coefficients are
normally ordered polynomials. It can be proved that the series depends on
the polynomial $f$ only (i.e., it is independent of arbitrariness in choice
of the combination $z^*z$ at every step of the process). We shall call it
the asymptotic expansion of $f$.

Let us define a deformation $\star_t$ of the algebra $\mathrm {Pol}({\mathbb
C})_{q}$ as follows. Because of isomorphism (\ref{uni1M_UNIQ.TEXisom}), it is sufficient
to define $z^iz^{*j}\star_tz^kz^{*l}$. We set $$
z^iz^{*j}\star_tz^kz^{*l}=\mathrm{asymptotic}\quad\mathrm{expansion}\quad
\mathrm{of}\quad z^iz^{*j}z^kz^{*l}. $$

It is almost explicit from the algorithm of normal ordering that this
star-product is associative and defines a deformation with separation of
variables. Realness of this deformation is a less evident fact.
 \footnote{To prove the realness one actually should show that the procedure
of normal ordering "commutes" with application of the involution. It is easy
to observe that application of the involution to a polynomial $f$ followed
by normal ordering via relation (\ref{uni1M_UNIQ.TEXdefpol}) gives us the same series as
normal ordering via relation (\ref{uni1M_UNIQ.TEXdefpol1}) instead of (\ref{uni1M_UNIQ.TEXdefpol})
followed by application of the involution. Now realness of the deformation
is due to the following fact: the procedures of normal ordering via
(\ref{uni1M_UNIQ.TEXdefpol}) and (\ref{uni1M_UNIQ.TEXdefpol1}) lead to the same result.}

The deformation constructed in this way turns out to be covariant. It is a
consequence of the following statement \cite{uni1M_UNIQ.TEXSSV} (see also \cite{uni1M_UNIQ.TEXKL}).
\medskip
\begin{proposition}\label{uni1M_UNIQ.TEXpred1}
There exists a unique structure of $U_q \mathfrak{su}_{1,1}$-module algebra
in $\mathrm {Pol}({\mathbb C})_{q,t}$ satisfying (\ref{uni1M_UNIQ.TEXaction1}).
\end{proposition}
\medskip
Evidently, formulas (\ref{uni1M_UNIQ.TEXaction2}) hold in $\mathrm {Pol}({\mathbb
C})_{q,t}$ as well. To explain covariance of the deformation it suffices to
note that the $U_q \mathfrak{sl}_{2}$-action in $\mathrm {Pol}({\mathbb
C})_{q,t}$ is independent of $t$, commutes with multiplication by $t$, and
respects the normal order.

We have explained roughly the idea of constructing the Berezin deformation
of the algebra $\mathrm {Pol}({\mathbb C})_{q}$. Let us expound the above
'algorithm' rigorously.

To start with, let us agree about the following notation. If $F$ is a
${\mathbb C}[[t]]$-module then $F_n$ stands for the ${\mathbb
C}[[t]]$-module $F/t^nF$ and $\overline{F}$ for the completion of $F$ in the
$t$-adic topology. Note that if $F$ is an (involutive) ${\mathbb
C}[[t]]$-algebra then $F_n$ and $\overline{F}$ are naturally endowed with
structures of an (involutive) ${\mathbb C}[[t]]$-algebra and a topological
(involutive) ${\mathbb C}[[t]]$-algebra, respectively.

Consider the unital involutive ${\mathbb C}[[t]]$-algebra given by its
generator and relation (\ref{uni1M_UNIQ.TEXdefpol}). Denote it by ${\cal P}$. Evidently,
relation (\ref{uni1M_UNIQ.TEXdefpol1}) holds in this algebra as well. Using proposition
\ref{uni1M_UNIQ.TEXpred1}, one shows that there exists a unique structure of $U_q
\mathfrak{su}_{1,1}$-module algebra in ${\cal P}$ given by (\ref{uni1M_UNIQ.TEXaction1})
and (\ref{uni1M_UNIQ.TEXaction2}). Denote the ${\mathbb C}[[t]]$-algebra $\mathrm
{Pol}({\mathbb C})_{q}\otimes {\mathbb C}[[t]]$ by $P$. Endow the latter
algebra with a structure $U_q \mathfrak{su}_{1,1}$-module algebra via
${\mathbb C}[[t]]$-linearity.

Let $I_n$ stands for the unique ${\mathbb C}[[t]]$-linear map from $P_n$ to
${\cal P}_n$ such that $I_n:z^iz^{*j}\mapsto z^iz^{*j}$. Certainly, $I_n$ is
an embedding of ${\mathbb C}[[t]]$-modules and respects $U_q
\mathfrak{sl}_{2}$-actions. Moreover, $I_n(f^*)=(I_n(f))^*$ for any $f\in
P_n$. $I_n$ turns out to be an isomorphism of ${\mathbb C}[[t]]$-modules.
Let us prove this via induction on $n$.

The case $n=1$ is evident (moreover, due to (\ref{uni1M_UNIQ.TEXpol}) and (\ref{uni1M_UNIQ.TEXdefpol}),
$I_1$ is an isomorphism of $U_q \mathfrak{su}_{1,1}$-module algebras).
Suppose $I_k$ is an isomorphism of ${\mathbb C}[[t]]$-modules for any $k\leq
n$. It suffices to show that the embedding $I_{n+1}:P_{n+1}\rightarrow {\cal
P}_{n+1}$ is surjective.

Let $\varphi_{n+1}:P_{n+1}\rightarrow P_n$ and $\psi_{n+1}:{\cal
P}_{n+1}\rightarrow {\cal P}_n$ stand for the natural surjective ${\mathbb
C}[[t]]$-linear maps induced by the inclusions $t^{n+1}P\subset t^{n}P$ and
$t^{n+1}{\cal P}\subset t^{n}{\cal P}$. It is clear that
\begin{equation}\label{uni1M_UNIQ.TEXdia}
\psi_{n+1}\circ I_{n+1}=I_{n}\circ\varphi_{n+1}.
\end{equation}
Let $f\in{\cal P}_{n+1}$. In view of the inductive hypothesis, there exists
$g\in P_{n+1}$ such that $I_n\circ \varphi_{n+1}(g)=\psi_{n+1}(f)$.
(\ref{uni1M_UNIQ.TEXdia}) implies $\psi_{n+1}\circ I_{n+1}(g)=\psi_{n+1}(f)$, i.e.,
$I_{n+1}(g)-f\in\mathrm{Ker}\psi_{n+1}$. In other words,
$I_{n+1}(g)-f=t^ng_0$ for some $g_0\in{\cal P}_{n+1}$. Relation
(\ref{uni1M_UNIQ.TEXdefpol}) and the definition of the maps $I_n$ imply $t^n{\cal
P}_{n+1}\subset \mathrm{Im}I_{n+1}$. By (\ref{uni1M_UNIQ.TEXdia})
$f\in\mathrm{Im}I_{n+1}$. So, $I_{n}:P_{n}\rightarrow {\cal P}_{n}$ is an
isomorphism of ${\mathbb C}[[t]]$-modules.

The family of isomorphisms $I_n$, $n=1,2,\ldots$, satisfies (\ref{uni1M_UNIQ.TEXdia}). It
means that there exists the limit isomorphism of topological ${\mathbb
C}[[t]]$-modules $$ I_\infty:\overline{P}\rightarrow\overline{{\cal P}}. $$
Evidently, $\overline{P}$ and $\overline{{\cal P}}$ inherit structures of
$U_q \mathfrak{su}_{1,1}$-module ${\mathbb C}[[t]]$-algebras (moreover, $U_q
\mathfrak{su}_{1,1}$-module ${\mathbb C}[[t]]$-algebra $\overline{P}$ is
isomorphic to $\mathrm {Pol}({\mathbb C})_{q}[[t]]$). The map $I_\infty$ is
an isomorphism of $U_q \mathfrak{sl}_{2}$-modules and intertwines actions of
the involutions in $\overline{P}$ and $\overline{{\cal P}}$. But $I_\infty$
is not an {\it algebra} isomorphism. Define in $\overline{P}\simeq\mathrm
{Pol}({\mathbb C})_{q}[[t]]$ a new product by $$
f\star_tg=I^{-1}_\infty\left(I_\infty(f)\cdot I_\infty(g)\right). $$
Evidently, $\mathrm {Pol}({\mathbb C})_{q}[[t]]$ endowed with this new
product become a $U_q \mathfrak{su}_{1,1}$-module algebra.

This deformation is just the one described informally in the first part of
this section. We shall call it the Berezin deformation of $\mathrm
{Pol}({\mathbb C})_{q}$ (\cite{uni1M_UNIQ.TEXSSV}) and denote the attached star-product by
$\star_B$.

\medskip

{\bf Remark 2.} We present in \cite{uni1M_UNIQ.TEXSSV} explicit formulas for $\star_B$,
i.e., for the corresponding bilinear maps $m_i$ (see (\ref{uni1M_UNIQ.TEXstar})).

\bigskip
\section{Uniqueness of the RCW-deformation of the algebra $\mathrm
{Pol}({\mathbb C})_{q}$}

First of all, let us consider the following construction. Suppose $\star_t$
is a deformation of the algebra $\mathrm {Pol}({\mathbb C})_{q}$ and
$m_i:\mathrm {Pol}({\mathbb C})_{q}\times\mathrm {Pol}({\mathbb
C})_{q}\rightarrow \mathrm {Pol}({\mathbb C})_{q}$ are the attached bilinear
maps (see (\ref{uni1M_UNIQ.TEXstar})). Let $c(t)\in{\mathbb C}[[t]]$, $c(0)=0$. Using
$c(t)$, one can define a new star-product $\star'_t$ given by
\begin{equation}\label{uni1M_UNIQ.TEXstar'}
f\star'_tg=f\cdot g+c(t)\cdot m_1(f,g)+c(t)^2\cdot m_2(f,g)+\ldots.
\end{equation}
(The product of arbitrary elements $f,g\in\mathrm {Pol}({\mathbb C})_q[[t]]$
is defined via the conditions of ${\mathbb C}$-bilinearity and $t$-adic
continuity.) It is said that $\star'_t$ is obtained from $\star_t$ via
change of formal parameter. If the initial deformation $\star_t$ is a
RCW-deformation and the series $c(t)$ is real then $\star'_t$ is a
RCW-deformation as well.

The principal result of this paper is
\medskip
\begin{theorem}\label{uni1M_UNIQ.TEXmain}
Any RCW-deformation of the algebra $\mathrm {Pol}({\mathbb C})_q$ can be
obtained from the Berezin deformation via change of formal parameter.
\end{theorem}
\medskip
{\bf Proof of the Theorem.} Let us use the traditional notation $M[[t]]$ for
the ${\mathbb C}[[t]]$-module of formal series in $t$ over a vector space
$M$.

Let $\star_t$ denotes a RCW-deformation of the algebra $\mathrm
{Pol}({\mathbb C})_q$. Endow $\mathrm {Pol}({\mathbb C})_q[[t]]$ with the
algebra structure using the star-product $\star_t$. To prove the theorem we
shall show that the elements $z$ and $z^*$ of the algebra $\mathrm
{Pol}({\mathbb C})_q[[t]]$ satisfy a commutation relation which can be
obtained from (\ref{uni1M_UNIQ.TEXdefpol}) via change of formal parameter.

To start with, let us construct an embedding of the algebra $\mathrm
{Pol}({\mathbb C})_q[[t]]$ into the algebra of endomorphisms of a ${\mathbb
C}[[t]]$-module.

Let $J$ stands for the left ideal in $\mathrm {Pol}({\mathbb C})_q[[t]]$
generated by $z^*$: $$ J=\mathrm {Pol}({\mathbb C})_q[[t]]\star_tz^*. $$
Since $\star_t$ is a product with separation of variables, we have
$J=J_0[[t]]$ with $J_0=\mathrm {Pol}({\mathbb C})_qz^*$. Thus, $\mathrm
{Pol}({\mathbb C})_q[[t]]$ splits into the direct sum of ${\mathbb
C}[[t]]$-submodules
\begin{equation}\label{uni1M_UNIQ.TEXsplit}
\mathrm {Pol}({\mathbb C})_q[[t]]=J+{\mathbb C}[z]_q[[t]].
\end{equation}
From (\ref{uni1M_UNIQ.TEXsplit}), we have the natural isomorphism of ${\mathbb
C}[[t]]$-modules
\begin{equation}\label{uni1M_UNIQ.TEXfactor}
\varphi:\mathrm {Pol}({\mathbb C})_q[[t]]/J\rightarrow{\mathbb C}[z]_q[[t]].
\end{equation}
The ${\mathbb C}[[t]]$-module $\mathrm {Pol}({\mathbb C})_q[[t]]/J$ is
endowed with the evident $\mathrm {Pol}({\mathbb C})_q[[t]]$-module
structure. Using isomorphism (\ref{uni1M_UNIQ.TEXfactor}), one can "transfer" $\mathrm
{Pol}({\mathbb C})_q[[t]]$-module structure onto ${\mathbb C}[z]_q[[t]]$. As
a result, we get the algebra homomorphism $$ T:\mathrm {Pol}({\mathbb
C})_q[[t]]\rightarrow \mathrm{End}_{{\mathbb C}[[t]]}\left({\mathbb
C}[z]_q[[t]]\right). $$
\medskip
\begin{lemma}\label{uni1M_UNIQ.TEX1}
$T$ is an embedding.
\end{lemma}
\medskip
{\bf Proof of Lemma \ref{uni1M_UNIQ.TEX1}.} To prove injectivity of $T$, we shall describe
more or less explicitly the endomorphisms $T(z)$ and $T(z^*)$ of the
${\mathbb C}[[t]]$-module ${\mathbb C}[z]_q[[t]]$. It is easy to describe
the former endomorphism. Indeed, each direct summand in (\ref{uni1M_UNIQ.TEXsplit}) is
invariant under multiplication from the left by elements of the subalgebra
${\mathbb C}[z]_q[[t]]\subset \mathrm {Pol}({\mathbb C})_q[[t]]$. This
implies $$ T(z):z^n\mapsto z^{n+1}. $$

To describe explicitly $T(z^*)$, one should use an extra structure in the
$\mathrm {Pol}({\mathbb C})_q[[t]]$-module ${\mathbb C}[z]_q[[t]]$ arising
from the covariance of the deformation $\star_t$.

Let us consider the Hopf subalgebra $U_q\mathfrak{b}_{-}\subset
U_q\mathfrak{sl}_2$ generated by $K^{\pm1}$ and $F$. Certainly, $J$ is a
$U_q\mathfrak{b}_{-}$-invariant subspace in $\mathrm {Pol}({\mathbb
C})_q[[t]]$ (see the definition of $J$ and relations
(\ref{uni1M_UNIQ.TEXaction1}),(\ref{uni1M_UNIQ.TEXaction2})). Thus, the ${\mathbb C}[[t]]$-module
$\mathrm {Pol}({\mathbb C})_q[[t]]/J$ is endowed with a natural structure of
a $U_q\mathfrak{b}_{-}$-module. Moreover, since $\mathrm {Pol}({\mathbb
C})_q[[t]]$ is a $U_q\mathfrak{sl}_{2}$-module algebra, $\mathrm
{Pol}({\mathbb C})_q[[t]]/J$ is a $U_q\mathfrak{b}_{-}$-module $\mathrm
{Pol}({\mathbb C})_q[[t]]$-module, i.e., the multiplication map $$ \mathrm
{Pol}({\mathbb C})_q[[t]]\otimes (\mathrm {Pol}({\mathbb
C})_q[[t]]/J)\rightarrow\mathrm ({Pol}({\mathbb C})_q[[t]]/J) $$ is a
morphism of $U_q\mathfrak{b}_{-}$-modules. Furthermore, the
$U_q\mathfrak{b}_{-}$-action respects decomposition (\ref{uni1M_UNIQ.TEXsplit}) (this is a
direct consequence of (\ref{uni1M_UNIQ.TEXaction1}) and (\ref{uni1M_UNIQ.TEXaction2})). Therefore, the
structure of $\mathrm {Pol}({\mathbb C})_q[[t]]$-module in ${\mathbb
C}[z]_q[[t]]$ (defined above via isomorphism (\ref{uni1M_UNIQ.TEXfactor})) is compatible
with the natural $U_q\mathfrak{b}_{-}$-module structure. In other words,
${\mathbb C}[z]_q[[t]]$ is a $U_q\mathfrak{b}_{-}$-module $\mathrm
{Pol}({\mathbb C})_q[[t]]$-module. This observation allows to calculate
$T(z^*)$.

Suppose that $T(z^*):z^n\mapsto \phi_n(z)$ for some $\phi_n(z)\in{\mathbb
C}[z]_q[[t]]$. Formulas (\ref{uni1M_UNIQ.TEXaction1}), (\ref{uni1M_UNIQ.TEXaction2}) and
$U_q\mathfrak{b}_{-}$-moduleness of the $\mathrm {Pol}({\mathbb
C})_q[[t]]$-module ${\mathbb C}[z]_q[[t]]$-module imply
$K(\phi_n(z))=q^{2n-2}\phi_n(z)$. It is clear that the latter equality holds
only when $\phi_n(z)=c_n(t)z^{n-1}$, $c_n(t)\in{\mathbb C}[[t]]$. It is
possible to calculate $c(0)$: since $$ z^*\star_tz\equiv q^2zz^*+1-q^2 \quad
(\mathrm{mod}\quad t) $$ (see (\ref{uni1M_UNIQ.TEXpol})), we get $$ z^*\star_tz^n\equiv
q^{2n}z^nz^*+(1-q^{2n})z^{n-1} \quad (\mathrm{mod}\quad t). $$ Thus,
$c_n(t)\equiv 1-q^{2n} \quad (\mathrm{mod}\quad t)$. As a result, we have
\begin{equation}\label{uni1M_UNIQ.TEXz}
T(z):z^n\mapsto z^{n+1},
\end{equation}
\begin{equation}\label{uni1M_UNIQ.TEXz*}
T(z^*):z^n\mapsto c_n(t)z^{n-1},
\end{equation}
with
\begin{equation}\label{uni1M_UNIQ.TEXc}
c_n(t)\equiv 1-q^{2n} \quad (\mathrm{mod}\quad t).
\end{equation}

We are now in a position to prove injectivity of $T$. Let $f=\sum_{i,j}
a_{ij}z^iz^{*j}\in\mathrm {Pol}({\mathbb C})_q[[t]]$ be a non-zero element
such that $T(f)=0$. Suppose $j_0$ is the least non-negative integer such
that $a_{ij_0}\ne0$ for some $i$. From (\ref{uni1M_UNIQ.TEXz}) and (\ref{uni1M_UNIQ.TEXz*}), $$
T(f):z^{j_0}\mapsto\sum_ia_{ij_0}z^i. $$ Thus,
$a_{ij_0}(t)c_{j_0}(t)c_{j_0-1}(t)\ldots c_{1}(t)=0$ for any $i$. In view of
(\ref{uni1M_UNIQ.TEXc}), $a_{ij_0}(t)=0$ for any $i$. It contradicts our choice of $j_0$.
Q.E.D.

\medskip

\begin{lemma}\label{uni1M_UNIQ.TEX2}
There exists a formal series $c(t)$ with real coefficients such that
$c(0)=0$ and
\begin{equation}\label{uni1M_UNIQ.TEXrel}
T(z^*)T(z)=
q^2T(z)T(z^*)+1-q^2+\frac{(1-q^2)t}{1-q^2t}(1-T(z)T(z^*))(1-T(z^*)T(z)).
\end{equation}
\end{lemma}
{\bf Proof of Lemma \ref{uni1M_UNIQ.TEX2}.} Let us invoke $U_q\mathfrak{b}_{-}$-moduleness
of the $\mathrm {Pol}({\mathbb C})_q[[t]]$-module ${\mathbb C}[z]_q$ once
more. Apply the element $F\in U_q\mathfrak{b}_{-}$ to both hand sides of the
equality $T(z^*)(z^n)=c_n(t)z^{n-1}$ (see (\ref{uni1M_UNIQ.TEXz*})):
\begin{equation}\label{uni1M_UNIQ.TEX3}
T(Fz^*)(K^{-1}z^n)+T(z^*)(Fz^n)=c_n(t)Fz^{n-1}.
\end{equation}

Equalities (\ref{uni1M_UNIQ.TEXaction1}) and easy induction on $n$ give
\begin{equation}\label{uni1M_UNIQ.TEXaction3}
K^{-1}z^n=q^{-2n}z^n, \qquad Fz^n=q^{1/2}\frac{1-q^{-2n}}{1-q^{-2}}z^{n-1}.
\end{equation}
According to (\ref{uni1M_UNIQ.TEXaction2}) and (\ref{uni1M_UNIQ.TEXaction3}), (\ref{uni1M_UNIQ.TEX3}) can be rewritten
as follows
\begin{equation}\label{uni1M_UNIQ.TEX4}
-q^{5/2}T(z^{*2})(q^{-2n}z^n)+q^{1/2}
\frac{1-q^{-2n}}{1-q^{-2}}T(z^*)(z^{n-1})
=q^{1/2}\frac{1-q^{2-2n}}{1-q^{-2}}c_n(t)z^{n-2}.
\end{equation}
From (\ref{uni1M_UNIQ.TEXz*}),
\begin{equation}\label{uni1M_UNIQ.TEX5}
-q^{5/2-2n}c_n(t)c_{n-1}(t)z^{n-2}+q^{1/2}\frac{1-q^{-2n}}{1-q^{-2}}
c_{n-1}(t)z^{n-2} =q^{1/2}\frac{1-q^{2-2n}}{1-q^{-2}}c_n(t)z^{n-2}.
\end{equation}

For $n=0$ or $n=1$, the latter equality holds trivially. For $n>1$, it is
equivalent to the following 'recurrence' relation
\begin{equation}\label{uni1M_UNIQ.TEX6}
-q^{2-2n}(1-q^{-2})c_n(t)c_{n-1}(t)+(1-q^{-2n}) c_{n-1}(t)
=(1-q^{2-2n})c_n(t).
\end{equation}

Let us define the series $c(t)$ by
\begin{equation}\label{uni1M_UNIQ.TEX7}
c_1(t)=\frac{1-q^2}{1-q^2c(t)}.
\end{equation}

From (\ref{uni1M_UNIQ.TEXc}), we have $c(0)=0$. In view of (\ref{uni1M_UNIQ.TEX6}),
\begin{equation}\label{uni1M_UNIQ.TEX8} c_n(t)=\frac{1-q^{2n}}{1-q^{2n}c(t)}
\end{equation}
(it suffices to proceed by induction on $n$). So, there exists a series
$c(t)$
 such that $c(0)=0$ and
\begin{equation}\label{uni1M_UNIQ.TEXz1}
T(z)(z^n)= z^{n+1},
\end{equation}
\begin{equation}\label{uni1M_UNIQ.TEXz*1}
T(z^*)(z^n)=\frac{1-q^{2n}}{1-q^{2n}c(t)}z^{n-1}.
\end{equation}
 The latter equalities imply (\ref{uni1M_UNIQ.TEXrel}). Now, to finish proof of
 Lemma \ref{uni1M_UNIQ.TEX2}, it remains to show that the series $c(t)$ given by (\ref{uni1M_UNIQ.TEX8})
 has real coefficients. Because of the definition of $c_1(t)$,
\begin{equation}\label{uni1M_UNIQ.TEXc1}
z^*\star_t z=c_1(t)+f
\end{equation}
for some $f\in J$. Apply the element $K\in U_q\mathfrak{sl}_2$ to the both
hand sides of (\ref{uni1M_UNIQ.TEXc1}). Since the deformation $\star_t$ is covariant, we
have $Kf=f$, i.e., $f=\sum_ja_j(t)z^jz^{*j}$. Apply the involution to the
both hand sides of the equality $$ z^*\star_t
z=c_1(t)+\sum_ja_j(t)z^jz^{*j}. $$ Now realness of coefficients of $c_1(t)$
is due to realness of the deformation. It follows that coefficients of
$c(t)$ are real as well (see
 (\ref{uni1M_UNIQ.TEX8})). Lemma \ref{uni1M_UNIQ.TEX2} is proved.
\medskip

Lemmas \ref{uni1M_UNIQ.TEX1} and \ref{uni1M_UNIQ.TEX2} imply existence of a series $c(t)$ with real
coefficient such that $c(0)=0$ and
\begin{equation}\label{uni1M_UNIQ.TEXrel*}
z^*\star_tz=q^2z\star_tz^*+1-q^2+
\frac{(1-q^2)c(t)}{1-q^2c(t)}(1-z\star_tz^*)\star_t(1-z^*\star_tz).
\end{equation}
It is not difficult to finish the proof of Theorem \ref{uni1M_UNIQ.TEXmain}. The Berezin
deformation $\star_B$ certainly satisfies the condition
$$
z^*\star_Bz=q^2z\star_Bz^*+1-q^2+\frac{(1-q^2)t}{1-q^2t}(1-z\star_Bz^*)
\star_B(1-z^*\star_Bz)
$$
(see the previous section). Moreover, the series $z^*\star_Bz$ can be
reconstructed from the latter relation in a unique way. Similarly, there
are no two different deformations with separation of variables satisfying
(\ref{uni1M_UNIQ.TEXrel*}) (see Remark 1). It is clear how to construct one
such deformation: for $f,g\in \mathrm {Pol}({\mathbb C})_q$ one should set
$$ f\star_t g=f\cdot g+c(t)\cdot m_1(f,g)+c(t)^2\cdot m_2(f,g)+\ldots, $$
with $m_i$ given by $$ f\star_B g=f\cdot g+t\cdot m_1(f,g)+t^2\cdot
m_2(f,g)+\ldots. $$ It is straightforward now that the deformation
$\star_t$ can be obtained from the Berezin one via change of formal
parameter. Theorem \ref{uni1M_UNIQ.TEXmain} is proved.

\medskip

{\bf Remark 3.} At first sight, we didn't invoke the full "quantum
symmetry" in the proof of Theorem \ref{uni1M_UNIQ.TEXmain} (i. e.,
$U_q\mathfrak{su}_{1,1}$-covariance of the algebra $\mathrm {Pol}({\mathbb
C})_q[[t]]$). What we realy used was the $U_q\mathfrak{b}_{-}$-covariance.
But it is simple to show that realness and $U_q\mathfrak{b}_{-}$-covariance
of a deformation with separation of variables imply its
$U_q\mathfrak{su}_{1,1}$-covariance.





\title{\bf HARISH-CHANDRA EMBEDDING AND q-ANALOGUES OF BOUNDED SYMMETRIC
DOMAINS}

\author{\sl S. Sinel'shchikov \thanks{Partially supported by ISF grant
U2B200 and grant DKNT-1.4/12}\and \hspace{-3em} and \ \sl L. Vaksman
\thanks{Partially supported by the grant INTAS-94-4720, ISF grant U21200
and grant DKNT-1.4/12}}

\date{\tt Institute for Low Temperature Physics \& Engineering \\
National Academy of Sciences of Ukraine}

\newpage
\setcounter{section}{0}
\large

\thispagestyle{empty} \ \vfill \begin{center}\LARGE \bf PART II \\ UNITARY
REPRESENTATIONS AND NON-COMMUTATIVE HARMONIC ANALYSIS
\end{center}
\vfill \addcontentsline{toc}{chapter}{Part II \ \ UNITARY REPRESENTATIONS
AND NON-COMMUTATIVE HARMONIC ANALYSIS}
\newpage

\makeatletter
\renewcommand{\@oddhead}{HARISH-CHANDRA EMBEDDING \hfill \thepage}
\renewcommand{\@evenhead}{\thepage \hfill S. Sinel'shchikov and L. Vaksman}
\let\@thefnmark\relax
\@footnotetext{This lecture has been delivered at D. Volkov Memorial Seminar
Held in Kharkov, January 1997 and published in Supersymmetry and Quantum
Field Theory, J. Wess and P. Akulov (eds). Springer, 1998, 312--316.}
\addcontentsline{toc}{chapter}{\@title \\ {\sl S. Sinel'shchikov and L.
Vaksman}\dotfill} \makeatother

\maketitle

{\bf 1.} This work is devoted to study of a very restricted class of
homogeneous spaces associated to quantum groups \cite{uni2_CQABSD.TEXD,uni2_CQABSD.TEXJ}. We follow
\cite{uni2_CQABSD.TEXSiV2} in describing here the construction of algebras of functions and
differential forms on these quantum homogeneous spaces.

We hope to extend to the above context a great deal of the results of
function theory and harmonic analysis in bounded symmetric domains
\cite{uni2_CQABSD.TEXHu}. This is shown here to be available for the simplest one among
such domains, the quantum disc \cite{uni2_CQABSD.TEXSSV}.

Our subsequent constructions are q-analogues of the corresponding
Harish-Chandra's constructions which allow one to embed a Hermitian
symmetric space of non-compact type into $\mathbb{C}^N$ \cite{uni2_CQABSD.TEXHe}.

Let $A$ be a Hopf algebra, $\varepsilon$ its counit, and $S$ its antipode.
Consider an algebra $F$ equipped also with a structure of $A$-module. $F$ is
said to be an $A$-module (covariant) algebra if

 i) the multiplication $m:\:F \otimes F \to F,\; m:\,f_1 \otimes f_2
\mapsto f_1 \cdot f_2;\;f_1,f_2 \in F$ is a morphism of $A$-modules;

 ii) the unit $1 \in F$ is an invariant: $\xi 1=\varepsilon(\xi)1,\;\xi
\in A$.

\medskip

If $A$ is a Hopf $*$-algebra, and $F$ is also equipped with an involution,
then the definition of covariance should include the following compatibility
condition for involutions:
$$\forall \xi \in A,\,f \in F \quad (\xi f)^*=(S(\xi))^*f^*.$$

In the sequel all the algebras of "functions"
($\mathbb{C}[\mathfrak{g}_{-1}]_q,\,
\mathbb{C}[\overline{\mathfrak{g}}_{-1}]_q,
\,\mathrm{Pol}(\mathfrak{g}_{-1})_q$) and "differential forms" are covariant
algebras.

\bigskip

{\bf 2.} Let $\mathfrak{g}$ be a simple complex Lie algebra, $\mathfrak{h}
\subset \mathfrak{g}$ a Cartan subalgebra, $\alpha_j \in \mathfrak{h}^*,
j=1,\ldots,l$, a system of simple roots with $\alpha_{j_0}$ being one of
those roots. Consider the $\mathbb{Z}$-grading $\mathfrak{g}=\bigoplus
\limits_m \mathfrak{g}_m$ given by
$$\mathfrak{g}_m=\{\xi \in \mathfrak{g}|\,[H_0,\xi]=2m \xi \},$$
with $H_0 \in \mathfrak{h}$ such that
$$\alpha_j(H_0)=0,\,j \ne j_0;\;\alpha_{j_0}(H_0)=2.$$

If this grading terminates,
$$\mathfrak{g}=\mathfrak{g}_{-1}+\mathfrak{g}_0+\mathfrak{g}_{+1},$$
then clearly $\mathfrak{g}_{\pm 1}$ are Abelian Lie subalgebras. This is
just the case when Harish-Chandra's construction presents a bounded
symmetric domain $U$ in the vector space $\mathfrak{g}_{-1}$.

Note also that in the case
$\mathfrak{g}=\mathfrak{sl}_{m+n},\,\alpha_{j_0}=\alpha_m$, $U$ is the
matrix ball in the space of $m \times n$ matrices:
$$
\mathfrak{g}_{-1}\simeq \mathrm{Mat}(m,n);\quad U=\{T \in
\mathrm{Mat}(m,n)|\, \| T \|<1 \}.
$$

\bigskip

{\bf 3.} Turn to the quantum case and fix $q \in (0,1)$.

Remind that the Hopf algebra $U_q \mathfrak{sl}_2$ is determined by its
generators $K^{\pm 1},E,F$ and the relations
\begin{gather*}
K \cdot K^{-1}=K^{-1}\cdot K=1,\quad K^{\pm 1} \cdot E=q^{\pm 2}EK^{\pm 1},
\\
K^{\pm 1} \cdot F=q^{\mp 2}FK^{\pm 1},\quad EF-FE=(K-K^{-1})/(q-q^{-1}).
\end{gather*}
Comultiplication $\Delta:U_q \mathfrak{sl}_2 \to U_q \mathfrak{sl}_2 \otimes
U_q \mathfrak{sl}_2$ is given by $\Delta(E)=E \otimes 1+K \otimes E,\quad
\Delta(F)=F \otimes K^{-1}+1 \otimes F,\quad \Delta(K^{\pm 1})=K^\pm \otimes
K^{\pm 1}.$

This Hopf algebra was introduced by E. Sklyanin, and its generalization $U_q
\mathfrak{g}$ to the case of an arbitrary simple Lie algebra $\mathfrak{g}$
in the works of V. Drinfeld and M. Jimbo \cite{uni2_CQABSD.TEXD}. $U_q \mathfrak{g}$ is
determined by the generators $\{K_j^{\pm 1},E_j,F_j \}_{j=1,\ldots,l}$ and
the well known relations \cite{uni2_CQABSD.TEXJ}. In this setting, every simple root
$\alpha_j, j=1,\ldots,l$, generates an embedding
$\varphi_j:U_{q_j}\mathfrak{sl}_2 \to U_q \mathfrak{g}$ given by
$$
\varphi_j:K^{\pm 1}\mapsto K_j^{\pm 1},\;\varphi_j:E \mapsto
E_j,\;\varphi_j:F \mapsto F_j.
$$
Here $q_j=q^{d_j}$ with $d_j>0$ such that $d_ia_{ij}= a_{ji}d_j$ for all
$i,j$.

\bigskip

{\bf 4.} Equip $U_q \mathfrak{g}$ with a structure of graded algebra:
\begin{gather*}
\deg K_j \:=\:\deg E_j \:=\:\deg F_j \:=\:0,\quad j \ne j_0
\\ \deg K_{j_0}\:=\:0,\quad \deg E_{j_0}\:=\:1,\quad \deg F_{j_0}\:=\:-1.
\end{gather*}

The embedding $\mathfrak{g}_{-1}\subset U \mathfrak{g}$ has no good
q-analog. This forces us to use the generalized Verma modules instead of $U
\mathfrak{g}$.

Let $V$ be a graded $U_q \mathfrak{g}$-module determined by its generator $v
\in V$ and the relations
\begin{gather*}
E_iv=0,\quad K_i^{\pm 1}v=v,\quad i=1,\ldots,l,
\\ F_jv=0,\quad j \ne j_0.
\end{gather*}
In the classic limit $q \to 1$ there is an embedding
$\mathfrak{g}_{-1}\hookrightarrow V,\;\xi \mapsto \xi v$. This allows one to
treat the homogeneous component $V_{-1}\:=\:\{v \in V|\,\deg v \:=\:-1 \}$
as a q-analog of the vector space $\mathfrak{g}_{-1}$.

We need also the graded $U_q \mathfrak{g}$-module $V'$ given by its
generator $v'$ and the relations
\begin{gather*}
E_iv'=0,\quad K_i^{\pm 1}v'=q^{\mp a_{ij_0}}v',\quad i=1,\ldots,l,
\\ F_j^{-a_{ij_0}+1}v'=0,\quad j \ne j_0;\quad \deg v'\,=\,-1.
\end{gather*}

\bigskip

{\bf 5.} Introduce the notation $\mathbb{C}[\mathfrak{g}_{-1}]_q=\bigoplus
\limits_m(V_m)^*,\;\bigwedge^1(\mathfrak{g}_{-1})_q=\bigoplus
\limits_m(V_m')^*$ for the dual to the $U_q \mathfrak{g}$-modules $V$ and
$V'$ respectively graded $U_q \mathfrak{g}$-modules. The elements $f \in
\mathbb{C}[\mathfrak{g}_{-1}]_q$ are be called holomorphic polynomials, and
the elements $\omega \in \bigwedge^1(\mathfrak{g}_{-1})_q$ differential
1-forms. The linear operator $d:\mathbb{C}[\mathfrak{g}_{-1}]_q \to
\bigwedge^1(\mathfrak{g}_{-1})_q$ is defined via the adjoint operator
$d^*:V'\to V$. In turn, $d^*$ is defined as the unique $U_q
\mathfrak{g}$-module morphism with $d^*:v'\mapsto F_{j_0}v$. Evidently, the
differential $d$ is a morphism of $U_q \mathfrak{g}^{op}$-modules.

\bigskip

{\bf 6.} Besides the comultiplication $\Delta:U_q \mathfrak{g}\to U_q
\mathfrak{g}\otimes U_q \mathfrak{g}^{op}$ we need also the opposite
comultiplication $\Delta^{op}$. It is used to equip the vector spaces $V
\otimes V,\:V \otimes V',\: V'\otimes V$ with a structure of $U_q
\mathfrak{g}$-modules.

The maps $v \mapsto v \otimes v,\; v'\mapsto v \otimes v',\; v'\mapsto
v'\otimes v$ admit the unique extensions to morphisms of $U_q
\mathfrak{g}$-modules
$$V \to V \otimes V,\quad V'\to V \otimes V',\quad V'\to V'\otimes V.$$
The adjoint operator to the comultiplication $V \to V \otimes V$ equips
$\mathbb{C}[\mathfrak{g}_{-1}]_q$ with a structure of associative algebra.
Similarly, the operators dual to the above morphisms $V'\to V \otimes
V',\;V'\to V'\otimes V$ equip $\bigwedge^1(\mathfrak{g}_{-1})_q$ with a
structure of a bimodule over $\mathbb{C}[\mathfrak{g}_{-1}]_q$.

It is easy to show that $d(f_1f_2)=df_1 \cdot f_2+f_1 \cdot df_2$ for all
$f_1,f_2 \in \mathbb{C}[\mathfrak{g}_{-1}]_q$. This allows one to pass from
the 1-forms to the higher differential forms (see, for instance, the
construction of G. Maltsiniotis \cite{uni2_CQABSD.TEXM}).

\bigskip

{\bf 7.} It is possible to describe the above algebras by their generators
and relations. Even in the simplest case $\mathfrak{g}=\mathfrak{sl}_N$ our
approach yields the profound results \cite{uni2_CQABSD.TEXSiV1}.

It should be noted that our approach to the construction of algebras of
differential forms is completely analogous to that of V. G. Drinfeld to the
construction of the algebra of functions on a formal quantum group \cite{uni2_CQABSD.TEXD}.

\bigskip

{\bf 8.} If we replace in the above construction the $U_q
\mathfrak{g}$-module $V$ with a highest weight vector by the $U_q
\mathfrak{g}$-module with a lowest weight vector, we obtain the algebra
$\mathbb{C}[\overline{\mathfrak{g}}_{-1}]_q$ of antiholomorphic polynomials
on the quantum vector space $\mathfrak{g}_{-1}$.

The tensor product
$$
\mathrm{Pol}(\mathfrak{g}_{-1})_q=\mathbb{C}[\mathfrak{g}_{-1}]_q \otimes
\mathbb{C}[\overline{\mathfrak{g}}_{-1}]_q
$$
is equipped with a structure of algebra by means of the universal R-matrix
together with the corresponding "commutativity morphism" \cite{uni2_CQABSD.TEXD}:
$$
\check{R}:\mathbb{C}[\overline{\mathfrak{g}}_{-1}]_q \otimes
\mathbb{C}[\mathfrak{g}_{-1}]_q \to \mathbb{C}[\mathfrak{g}_{-1}]_q \otimes
\mathbb{C}[\overline{\mathfrak{g}}_{-1}]_q.
$$

\bigskip

{\bf 9.} For the sake of passage from 'complex quantum Lie groups to real
ones' equip the Hopf algebra $U_q \mathfrak{g}$ with an involution:
$$
E_j^*=\left \{\begin{array}{r|c}K_jF_j & j \ne j_0 \\ -K_jF_j & j=j_0
\end{array}\right.,\qquad F_j^*=\left \{\begin{array}{r|c}E_jK_j^{-1} & j
\ne j_0 \\ -E_jK_j^{-1} & j=j_0 \end{array}\right.,
$$
$(K_j^{\pm 1})^*=K_j^{\pm 1},\quad i,j \in \{1,\ldots,l \}$.

The involution in $\mathrm{Pol}(\mathfrak{g}_{-1})_q$ presented in
\cite{uni2_CQABSD.TEXSiV2} possess the property $\forall \xi \in U_q \mathfrak{g},\,f \in
\mathrm{Pol}(\mathfrak{g}_{-1})_q$ $\qquad(\xi f)^*=(S(\xi))^*f^*$.

In the simplest case $\mathfrak{g}=\mathfrak{sl}_2$ one obtains the well
known $*$-algebra given by the generators $z,z^*$ and the relation
$z^*z-q^2zz^*=1-q^2$ \cite{uni2_CQABSD.TEXSiV2}.

The passage from the polynomial algebra to the algebra of continuous
functions in the closure of a bounded symmetric domain is made by means of a
C$^*$-completion. In the special case of quantum disc this argument was
used, in particular, in \cite{uni2_CQABSD.TEXNN}.

In \cite{uni2_CQABSD.TEXSSV} the q-analogs for the basic integral representations of the
function theory in the quantum disc were obtained. Besides, there are
several results for the quantum ball \cite{uni2_CQABSD.TEXV}.

\bigskip

{\bf 10.} Finally, we express our gratitude to V. P. Akulov for helpful
discussions of the results of this work.




\title{\bf q-ANALOGUES OF SOME BOUNDED SYMMETRIC DOMAINS}

\author{D. L. Shklyarov \and S. D. Sinel'shchikov \and L. L. Vaksman}
\date{\tt Institute for Low Temperature Physics \& Engineering \\
National Academy of Sciences of Ukraine,\\ Kharkov 61103, Ukraine}

\renewcommand{\theequation}{\thesection.\arabic{equation}}

\newpage
\setcounter{section}{0}
\large

\makeatletter
\renewcommand{\@oddhead}{q-ANALOGUES OF SOME BOUNDED SYMMETRIC DOMAINS
\hfill \thepage}
\renewcommand{\@evenhead}{\thepage \hfill D. Shklyarov, S. Sinel'shchikov,
and L. Vaksman}
\let\@thefnmark\relax
\@footnotetext{This lecture has been delivered at the 8-th International
Colloquium 'Quantum groups and Integrable Systems' held in Prague, June
1999; published in Czechoslovak Journal of Physics {\bf 50} (2000), No 1,
175 -- 180} \addcontentsline{toc}{chapter}{\@title \\ {\sl D. L. Shklyarov,
S. D. Sinel'shchikov, and L. L. Vaksman}\dotfill} \makeatother

\maketitle

\begin{quotation}\small{\bf Abstract.}
We study q-analogues of matrix balls. A description of algebras of finite
functions in the quantum matrix balls, an explicit form for the invariant
integral in the space of finite functions, q-analogues for the weighted
Bergman spaces, together with an explicit formula for the corresponding
Bergman kernel, are presented.
\end{quotation}

\section{Introduction}

Hermitian symmetric spaces of non-compact type constitute one of the most
important classes of homogeneous symmetric spaces. A well known result by
Harish-Chandra claims that any such space can be realized as a bounded
symmetric domain in a complex vector space $V$ (via the so-called
Harish-Chandra embedding).

Irreducible bounded symmetric domains were classified by E. Cartan. They are
among the important subjects in Lie theory, geometry and function theory.

The first step in studying q-analogues of irreducible bounded symmetric
domains was made in \cite{uni3cont.texSV}. This work provides a q-analogue for the
Harish-Chandra embedding and, in particular, a construction for q-analogues
of the polynomial algebra and the differential calculus on $V$.

Our subject is the simplest class among those q-analogues, the quantum
matrix balls. In the classical case $q=1$ the corresponding vector space $V$
is the space $\mathrm{ Mat}_{m,n}$ of rectangle complex matrices, and the
matrix ball is defined as
$$\mathbb{U}=\{\mathbf{z}\in \mathrm{Mat}_{m,n}|\mathbf{zz}^*<1 \}.$$

This ball is a homogeneous space of the group $SU_{n,m}$.

Proofs of all presented results and also some results concerning the special
case of quantum disc ($n=m=1$) can be found in our electronic preprints
(http://xxx.lanl.gov/).

\bigskip

\section{ Polynomials and finite functions in the quantum matrix ball}

Everywhere in the sequel $q \in (0,1)$, $m,n \in \mathbb{N}$, $m \le n$,
$N=m+n$. We use the standard notation ${\it {su}_{n,m}}$ for the Lie algebra
of the group $SU_{n,m}$.

The Hopf algebra $U_q \mathfrak{sl}_N$ is determined by its generators
$\{E_i,F_i,K_i^{\pm 1}\}_{i=1,\ldots,N-1}$ and the well known Drinfeld-Jimbo
relations \cite{uni3cont.texR}. Equip $U_q \mathfrak{sl}_N$ with the involution defined
on the generators $K_j^{\pm 1}$, $E_j$, $F_j$, $j=1,\ldots,N-1$ by
$$
\left(K_j^{\pm 1}\right)^*=K_j^{\pm 1},\qquad E_j^*=\left \{
\begin{array}{rr}K_jF_j,& j \ne n \\ -K_jF_j,&j=n \end{array}\right., \qquad
F_j^*=\left \{\begin{array}{rr}E_jK_j^{-1},& j \ne n \\ -E_jK_j^{-1},&j=n
\end{array}\right..
$$
The Hopf $*$-algebra $U_q \mathfrak{su}_{n,m}=(U_q \mathfrak{sl}_N,*)$
arising this way is a q-analogue of the Hopf algebra $U
\mathfrak{su}_{n,m}$.

Remind some well known definitions. An algebra $F$ is said to be an
$A$-module algebra if it is a module over a Hopf algebra $A$, the unit of
$F$ is an invariant and the multiplication $F \otimes F \to F$, $f_1 \otimes
f_2 \mapsto f_1 \cdot f_2$, is a morphism of $A$-modules. In the case of a
$*$-algebra $F$ and a Hopf-$*$-algebra $A$, there is an additional
requirement that the involutions agree as follows:
\begin{equation}\label{uni3cont.texadreq}
(af)^*=(S(a))^*f^*,\qquad a \in A,\;f \in F,
\end{equation}
with $S:A \to A$ being the antipode of $A$.

In \cite{uni3cont.texSV} a $U_q \mathfrak{su}_{n,m}$-module algebra
$\mathrm{Pol}(\mathrm{Mat}_{m,n})_q$ and its $U_q \mathfrak{sl}_N$-module
subalgebra $\mathbb{C}[\mathrm{Mat}_{m,n}]_q$ were introduced (the notation
$\mathfrak{g}_{-1}$ was used in \cite{uni3cont.texSV} instead of $\mathrm{Mat}_{m,n}$).
These algebras are q-analogues of polynomial algebras in the vector spaces
$\mathrm{Mat}_{m,n}$. We present below a description of these algebras in
terms of generators and relations, together with explicit formulae for the
$U_q \mathfrak{sl}_N$-action in $\mathbb{C}[\mathrm{Mat}_{m,n}]_q$
(corresponding explicit formulae for $U_q \mathfrak{su}_{n,m}$-action in
$\mathrm{Pol}(\mathrm{Mat}_{m,n})_q$ can be produced via (\ref{uni3cont.texadreq})).

With the definitions of \cite{uni3cont.texSV} as a background, one can prove the
following two propositions.

\medskip

\begin{proposition}
There exists a unique family $\{z_a^\alpha
\}_{a=1,\ldots,n;\alpha=1,\ldots,m}$ of elements of the $U_q
\mathfrak{sl}_N$-module algebra $\mathbb{C}[\mathrm{ Mat}_{m,n}]_q$ such
that for all ${a=1,\ldots,n;\alpha=1,\ldots,m}$
\begin{align}\label{uni3cont.texh}
H_nz_a^\alpha&=\left \{\begin{array}{ccl}2z_a^\alpha &,&a=n \;\&\;\alpha=m
\\ z_a^\alpha &,&a=n \;\&\;\alpha \ne m \quad \mathrm{or}\quad a \ne n
\;\&\; \alpha=m \\ 0 &,&\mathrm{otherwise}\end{array}\right.,&
\\ \label{uni3cont.texf}
F_nz_a^\alpha&=q^{1/2}\cdot \left \{\begin{array}{ccl}1 &,& a=n \;\&
\;\alpha=m \\ 0 &,&\mathrm{otherwise}\end{array}\right.,&
\\ \label{uni3cont.texe}
E_nz_a^\alpha&=-q^{1/2}\cdot \left \{\begin{array}{ccl}q^{-1}z_a^mz_n^\alpha
&,&a \ne n \;\&\;\alpha \ne m \\ (z_n^m)^2 &,& a=n \;\&\;\alpha=m \\
z_n^mz_a^{\alpha} &,&\mathrm{otherwise}\end{array}\right.,&
\end{align}
and with $k \ne n$
\begin{align}\label{uni3cont.texh's}
H_kz_a^\alpha&= \left \{\begin{array}{ccl}z_a^\alpha &,& k<n \;\&\;a=k \quad
\mathrm{or}\quad k>n \;\&\;\alpha=N-k \\-z_a^\alpha &,& k<n \;\&\;a=k+1
\quad \mathrm{or}\quad k>n \;\&\;\alpha=N-k+1 \\ 0 &,&\mathrm{
otherwise}\end{array}\right.,
\\ \label{uni3cont.texf's}
F_kz_a^\alpha&=q^{1/2}\cdot \left \{\begin{array}{ccl}z_{a+1}^\alpha &,& k<n
\;\&\;a=k \\ z_a^{\alpha+1} &,& k>n \;\&\;\alpha=N-k \\ 0
&,&\mathrm{otherwise}\end{array}\right.,
\\ \label{uni3cont.texe's}
E_kz_a^\alpha&=q^{-1/2}\cdot \left \{\begin{array}{ccl}z_{a-1}^\alpha &,&
k<n \;\&\; a=k+1 \\ z_a^{\alpha-1}&,& k>n \;\&\;\alpha=N-k+1 \\ 0 &,&
\mathrm{otherwise}\end{array}\right..
\end{align}
\end{proposition}

\medskip

{\sc Remark.} The elements $\{H_j \}_{j=1,..N-1}$ and $\{K_j \}_{j=1,..N-1}$
are related as follows:
$$K_j=q^{H_j}$$
(the exact definition of $\{H_j \}$ one can find in {\cite{uni3cont.texSV}}).

\medskip

\begin{proposition}
$\{z_a^\alpha \}_{a=1,\ldots,n;\alpha=1,\ldots,m}$ generate
$\mathbb{C}[\mathrm{Mat}_{m,n}]_q$ as an algebra and
$\mathrm{Pol}(\mathrm{Mat}_{m,n})_q$ as a $*$-algebra. The complete list of
relations is as follows:
\begin{align}\label{uni3cont.texz}
z_a^\alpha z_b^\beta&=\left \{\begin{array}{ccl}qz_b^\beta z_a^\alpha &,&
a=b \;\&\;\alpha<\beta \quad \mathrm{or}\quad a<b \;\&\;\alpha=\beta \\
z_b^\beta z_a^\alpha &,& a<b \;\&\;\alpha>\beta \\ z_b^\beta
z_a^\alpha+(q-q^{-1})z_a^\beta z_b^\alpha &,& a<b \;\&\;\alpha<\beta
\end{array} \right.,
\\ \label{uni3cont.texz,z*}
(z_b^\beta)^*z_a^\alpha&=q^2 \cdot \sum_{a',b'=1}^n
\sum_{\alpha',\beta'=1}^m R(b,a,b',a')R(\beta,\alpha,\beta',\alpha) \cdot
z_{a'}^{\alpha'}\left(z_{b'}^{\beta'}\right)^*+
(1-q^2)\delta_{ab}\delta^{\alpha \beta},
\end{align}
with $\delta_{ab}$, $\delta^{\alpha \beta}$ being the Kronecker symbols and
\begin{equation}\label{uni3cont.texR}
R(b,a,b',a')=\left \{\begin{array}{ccl}q^{-1} &,& a \ne b
\;\&\;b=b'\;\&\;a=a' \\ 1 &,& a=b=a'=b' \\ -(q^{-2}-1) &,& a=b
\;\&\;a^\prime=b' \;\&\;a'>a \\ 0 &,&\mathrm{otherwise}\end{array}\right..
\end{equation}
\end{proposition}

\medskip

{\sc Example.} In the simplest case $m=n=1$ the relations presented above
describe a very well known $U_q \mathfrak{su}_{1,1}$-module algebra
\begin{equation}\label{uni3cont.texd} z^*z=q^2 zz^*+1-q^2. \end{equation}

Consider the $*$-algebra $\mathrm{Fun}(\mathbb{U})_q \supset
\mathrm{Pol}(\mathrm{Mat}_{m,n})_q$ derived from
$\mathrm{Pol}(\mathrm{Mat}_{m,n})_q$ by adding a generator $f_0$ such that
\begin{equation}\label{uni3cont.texf0} f_0=f_0^2=f_0^*,\quad \left(z_a^\alpha
\right)^*f_0=f_0z_a^\alpha=0, \quad a=1,\ldots,n;\;\alpha=1,\ldots,m.
\end{equation}
(Relations (\ref{uni3cont.texf0}) allow one to treat $f_0$ as a q-analogue of the
function that equal to 1 in the center of the ball and equal to 0 in other
points.)

\medskip

\begin{proposition}
There exists a unique extension of the structure of a $U_q
\mathfrak{su}_{n,m}$-module algebra from
$\mathrm{Pol}(\mathrm{Mat}_{m,n})_q$ onto $\mathrm{Fun}(\mathbb{U})_q$ such
that
\begin{equation}\label{uni3cont.texhfe}
H_nf_0=0,\qquad F_nf_0=-\frac{q^{1/2}}{q^{-2}-1}f_0 \cdot \left(z_n^m
\right)^*,\qquad E_nf_0=-\frac{q^{1/2}}{1-q^2}z_n^m \cdot f_0
\end{equation}
and with $k \ne n$
\begin{equation}H_kf_0=F_kf_0=E_kf_0=0.\end{equation}
\end{proposition}

\medskip

The two-sided ideal $D(\mathbb{U})_q \stackrel{\mathrm{def}}{=}\mathrm{
Fun}(\mathbb{U})_qf_0 \mathrm{Fun}(\mathbb{U})_q$ is a $U_q
\mathfrak{su}_{n,m}$-module algebra. Its elements will be called the finite
functions in the quantum matrix ball.

\bigskip

\section{Invariant integral}

It is well known that in the classical case $q=1$ the positive
$SU_{n,m}$-invariant integral could not be defined on the polynomial
algebra. However, it is well defined on the space of finite smooth
functions. These observations are still applicable in the quantum case.

Consider the representation $T$ of $\mathrm{Fun}(\mathbb{U})_q$ in the space
$\mathcal{H}=\mathrm{Fun}(\mathbb{U})_qf_0=
\mathrm{Pol}(\mathrm{Mat}_{m,n})_qf_0$:
\begin{equation}\label{uni3cont.texT}
T(f) \psi=f \psi,\qquad f \in \mathrm{Fun}(\mathbb{U})_q,\quad \psi \in
\mathcal{H}.
\end{equation}

\medskip

{\sc Remark.} It can be shown that there exists a unique positive scalar
product in ${\cal H}$ such that $(f_0,f_0)=1$, and
\begin{equation}\label{uni3cont.texscpr}
(T(f)\psi_1,\psi_2)=(\psi_1,T(f^*)\psi_2),\qquad f \in
\mathrm{Fun}(\mathbb{U})_q,\quad \psi_1,\psi_2 \in \mathcal{H}.
\end{equation}
Moreover one can prove that the $*$-algebra
$\mathrm{Pol}(\mathrm{Mat}_{m,n})_q$ admits a unique up to unitary
equivalence {\sl faithful} irreducible $*$-representation by {\sl bounded}
operators in a Hilbert space. This $*$-representation can be produced via
extending the operators $T(f)$, $f \in \mathrm{Pol}(\mathrm{Mat}_{m,n})_q$,
onto the completion of the pre-Hilbert space $\mathcal{H}$.

Remind the notation $U_q \mathfrak{b}_+$ for the subalgebra of $U_q
\mathfrak{sl}_N$ generated by the elements $\{E_i,K_i^{\pm
1}\}_{i=1,\ldots,N-1}$. Obviously,
$$U_q \mathfrak{b}_+\mathcal{H}\subset \mathcal{H},$$
and thus we obtain the representation $\Gamma$ of the algebra $U_q
\mathfrak{b}_+$ in $\mathcal{H}$. Let also
\begin{equation}\label{uni3cont.texro}
{\check{\rho}=\frac12 \sum_{j=1}^{N-1}j(N-j)H_j}.
\end{equation}

\medskip

\begin{proposition}
The linear functional
\begin{equation}\label{uni3cont.texii}
\int \limits_{\mathbb{U}_q}fd \nu=\mathrm{tr}(T(f)\Gamma(q^{-2
\check{\rho}})),\qquad f \in D(\mathbb{U})_q,
\end{equation}
is well defined, $U_q \mathfrak{su}_{n,m}$-invariant and positive (i.e.
$\int \limits_{\mathbb{U}_q}f^*fd \nu>0$ for $f \ne 0$).
\end{proposition}

\bigskip

\section{Weighted Bergman spaces and Bergman kernels}

Our intention is to produce q-analogues of weighted Bergman spaces. In the
case $q=1$ one has
\begin{equation}\label{uni3cont.texde}
\det(1-\mathbf{zz}^*)=1+\sum_{k=1}^m(-1)^k \mathbf{z}^{\wedge
k}\mathbf{z}^{*\wedge k},
\end{equation}
with $\mathbf{z}^{\wedge k}$, $\mathbf{z}^{*\wedge k}$ being the "exterior
powers" of the matrices $\mathbf{z}$, $\mathbf{z}^*$, that is, matrices
formed by the minors of order $k$.

Let $1 \le \alpha_1<\alpha_2<\ldots<\alpha_k \le m$, $1 \le
a_1<a_2<\ldots<a_k \le n$. Introduce q-analogues of minors for the matrix
$\mathbf{z}$:
\begin{equation}\label{uni3cont.texminor}
{\mathbf{z}^{\wedge k \,}}_{\{a_1,a_2,\ldots,a_k \}}
^{\{\alpha_1,\alpha_2,\ldots,\alpha_k \}}=\sum_{s \in
S_k}(-q)^{l(s)}z_{a_1}^{\alpha_{s(1)}}z_{a_2}^{\alpha_{s(2)}}\ldots
z_{a_k}^{\alpha_{s(k)}},
\end{equation}
with $l(s)=\mathrm{card}\{(i,j)|\;i<j \quad \& \quad s(i)>s(j)\}$ being the
length of the permutation $s$.

The q-analogue $y \in \mathrm{Pol}(\mathrm{Mat}_{m,n})_q$ for the polynomial
$\det(1-\mathbf{zz}^*)$ is defined by
\begin{equation}\label{uni3cont.texy}
y=1+\sum_{k=1}^m(-1)^k \sum_{\{J'|\;\mathrm{card}(J')=k
\}}\sum_{\{J''|\;\mathrm{card}(J'')=k \}}{\mathbf{z}^{\wedge k
\,}}_{J''}^{J'}\cdot \left({\mathbf{ z}^{\wedge k \,}}_{J''}^{J'}\right)^*.
\end{equation}
Let $\lambda>m+n-1$. Now one can define the integral with weight $y^\lambda$
as follows:
\begin{equation}\label{uni3cont.texwi}
\int \limits_{\mathbb{U}_q}fd \nu_ \lambda \stackrel{\mathrm{
def}}{=}C(\lambda)\int \limits_{\mathbb{U}_q}fy^\lambda d \nu ,\qquad f \in
D(\mathbb{U})_q,
\end{equation}
where $C(\lambda)=\displaystyle
\prod_{j=0}^{n-1}\prod_{k=0}^{m-1}\left(1-q^{2(\lambda+1-N)}q^
{2(j+k)}\right)$ provides $\displaystyle \int \limits_{{\bf U}_q}1d \nu_
\lambda=1.$

The Hilbert space $L^2(d \nu_ \lambda)_q$ is defined as a completion of the
space $D(\mathbb{U})_q$ of finite functions with respect to the norm $\|f
\|_ \lambda=\left(\displaystyle \int \limits_{\mathbb{U}_q}f^*fd \nu_
\lambda \right)^{1/2}$. The closure $L_a^2(d \nu_ \lambda)_q$ in $L^2(d \nu_
\lambda)_q$ of the algebra $\mathbb{C}[\mathrm{Mat}_{m,n}]_q$ will be called
a weighted Bergman space.

Consider the orthogonal projection $P_ \lambda$ in $L^2(d \nu_ \lambda)_q$
onto the weighted Bergman space $L_a^2(d \nu_ \lambda)_q$. Our goal here is
to show that $P_ \lambda$ could be written as an integral operator
\begin{equation}\label{uni3cont.texpl}
P_ \lambda f=\int \limits_{\mathbb{U}_q}K_
\lambda(\mathbf{z},\boldsymbol{\zeta}^*)f(\boldsymbol{\zeta})d \nu_
\lambda(\boldsymbol{\zeta}),\qquad f \in D(\mathbb{U})_q.
\end{equation}

The main intention of this section is to introduce the algebra
$\mathbb{C}[[\mathrm{Mat}_{m,n}\times \overline{\mathrm{Mat}}_{m,n}]]_q$ of
kernels of integral operators and to determine an explicit form of the
Bergman kernel $K_ \lambda \in \mathbb{C}[[\mathrm{Mat}_{m,n}\times
\overline{\mathrm{Mat}}_{m,n}]]_q$.

Introduce the notation
\begin{equation}\label{uni3cont.texpk}
h_i=\sum_{\genfrac{}{}{0pt}{2}{J'\subset \{1,2,\ldots,m
\}}{\mathrm{card}(J')=i}} \sum_{\genfrac{}{}{0pt}{0}{J''\subset
\{1,2,\ldots,n \}}{\mathrm{card}(J'')=i}} {\mathbf{z}^{\wedge i
\,}}_{J''}^{J'}\otimes \left({\mathbf{z}^{\wedge i \,}}_{J''}^{J'}\right)^*.
\end{equation}

Let $\mathbb{C}[\overline{\mathrm{Mat}}_{m,n}]_q \subset
\mathrm{Pol}(\mathrm{Mat}_{m,n})_q$ be the unital subalgebra generated by
$(z_a^\alpha)^*$, $a=1,2,\ldots,n$, $\alpha=1,2,\ldots,m$, and
$\mathbb{C}[\mathrm{Mat}_{m,n}]_q^\mathrm{op}$ the algebra which {\bf
differs from} $\mathbb{C}[\mathrm{Mat}_{m,n}]_q$ {\bf by a replacement of
its multiplication law to the opposite one}. The tensor product algebra
$\mathbb{C}[\mathrm{Mat}_{m,n}]_q^\mathrm{op}\otimes
\mathbb{C}[\overline{\mathrm{Mat}}_{m,n}]_q$ will be called an algebra of
polynomial kernels. It is possible to show that in this algebra
 $h_ih_j=h_jh_i$ for all $i,j=1,2,\ldots,m$.

We follow \cite{uni3cont.texSV} in equipping $\mathrm{Pol}(\mathrm{Mat}_{m,n})_q$ with a
$\mathbb{Z}$-gradation: $\deg(z_a^\alpha)=1$, $\deg((z_a^\alpha)^*)=-1$,
$a=1,2,\ldots,n$, $\alpha=1,2,\ldots,m$. In this context one has:
\begin{equation}\label{uni3cont.texdecomp1}
\mathbb{C}[\mathrm{Mat}_{m,n}]_q^\mathrm{op}= \bigoplus_{i=0}^\infty
\mathbb{C}[\mathrm{Mat}_{m,n}]_{q,i}^\mathrm{op},\qquad
\mathbb{C}[\overline{\mathrm{Mat}}_{m,n}]_q=\bigoplus_{j=0}^\infty
\mathbb{C}[\overline{\mathrm{Mat}}_{m,n}]_{q,-j},
\end{equation}
\begin{equation}\label{uni3cont.texpka}
\mathbb{C}[\mathrm{Mat}_{m,n}]_q^\mathrm{op}\otimes
\mathbb{C}[\overline{\mathrm{Mat}}_{m,n}]_q=\bigoplus_{i,j=0}^\infty
\mathbb{C}[\mathrm{Mat}_{m,n}]_{q,i}^\mathrm{op}\otimes
\mathbb{C}[\overline{\mathrm{Mat}}_{m,n}]_{q,-j}.
\end{equation}

The kernel algebra $\mathbb{C}[[\mathrm{Mat}_{m,n}\times
\overline{\mathrm{Mat}}_{m,n}]]_q$ will stand for a completion of
$\mathbb{C}[\mathrm{Mat}_{m,n}]_q^\mathrm{op}\otimes
\mathbb{C}[\overline{\mathrm{Mat}}_{m,n}]_q$ in the topology associated to
the gradation in (\ref{uni3cont.texpka}).

\begin{proposition}
Let $K_ \lambda$ be an element of the algebra
$\mathbb{C}[[\mathrm{Mat}_{m,n}\times \overline{\mathrm{Mat}}_{m,n}]]_q$
defined by
\begin{equation}\label{uni3cont.texwBk}
K_ \lambda=\prod_{j=0}^\infty \left(1+\sum_{i=1}^m(-q^{2(\lambda+j)})^ih_i
\right)\cdot \prod_{j=0}^\infty \left(1+\sum_{i=1}^m(-q^{2j})^i h_i
\right)^{-1}
\end{equation}
Then (\ref{uni3cont.texpl}) holds.
\end{proposition}

\medskip

{\sc Remark.} A q-analogue of the ordinary Bergman kernel for the matrix
ball is derivable from (\ref{uni3cont.texwBk}) by a substitution $\lambda=m+n$:
$$
K=\prod_{j=0}^{m+n-1}\left(1+\sum_{i=1}^m(-q^{2j})^ih_i
\right)^{-1}\underset{q \to 1}{\longrightarrow} \left(\det(1-\mathbf{z}\cdot
\boldsymbol{\zeta}^*)\right)^{-(m+n)}.
$$

\bigskip

\section*{Notes of the Editor}

The proofs of the results announced in this work can be found in \cite{uni3cont.texSSV1,uni3cont.texSSV2}. These results acquired a further development in \cite{uni3cont.texV}.




\title{\bf ON A q-ANALOGUE OF THE FOCK INNER PRODUCT}
\author{ D.Shklyarov}
\date{\tt Institute for Low Temperature Physics \& Engineering \\
National Academy of Sciences of Ukraine,\\ Kharkov 61103, Ukraine}

\renewcommand{\theequation}{\thesection.\arabic{equation}}

\newpage
\setcounter{section}{0}
\large

\makeatletter
\renewcommand{\@oddhead}{ON A Q-ANALOG OF THE FOCK INNER PRODUCT \hfill
\thepage}
\renewcommand{\@evenhead}{\thepage \hfill D. Shklyarov}
\addcontentsline{toc}{chapter}{\@title \\ {\sl D. L. Shklyarov}\dotfill}
\makeatother

\maketitle

\bigskip
\section{Introduction}

Let $\mathrm{Mat}_{m,n}$ be the space of complex $n\times m$ matrices. We
denote by $z_a^\alpha$, $a=1,\ldots,n$, $\alpha=1,\ldots,m$ the standard
coordinate functions on $\mathrm{Mat}_{m,n}$ given by the matrix entries.
Let $\mathbb{C}[\mathrm{Mat}_{m,n}]$ stands for the space of polynomials on
$\mathrm{Mat}_{m,n}$ (i. e. polynomials in $z_a^\alpha$'s).

We recall that the Fock inner product in the space
$\mathbb{C}[\mathrm{Mat}_{m,n}]$ is defined by

\begin{equation}\label{uni3fock.texclfip}
  \left(P\,,\,Q\right)_F=\int\limits_{\mathrm{Mat}_{m,n}}P(\mathbb{Z})
\overline{Q(\mathbb{Z})}e^{-\mathrm{tr}(\mathbb{Z}\mathbb{Z}^*)}d\mathbb{Z}
\end{equation}
with $\mathbb{Z}=(z_a^\alpha)$ and $d\mathbb{Z}$ being the
Lebesgue measure on $\mathrm{Mat}_{m,n}$ such that
$$\left(1\,,\,1\right)_F=1.$$

The inner product possesses the following remarkable property

\begin{equation}\label{uni3fock.tex1pro}
  \left(\frac{\partial P}{\partial z_a^\alpha}\,,\,Q\right)_F=
\left(P\,,\,z_a^\alpha Q\right)_F \quad \forall a,\alpha,
\end{equation}
which allows us to rewrite the product in the differential form
\begin{equation}\label{uni3fock.texdifform}
  \left(P\,,\,Q\right)_F=\partial_{P}(\overline{Q})(0)
\end{equation}
where $\partial_{P}$ stands for the differential operator with
constant coefficients derived from the polynomial $P$ by
substituting $z_a^\alpha\to\frac{\partial }{\partial z_a^\alpha}$,
$a=1,\ldots,n$, $\alpha=1,\ldots,m$, and $\overline{Q}$ stands for
the polynomial derived from $Q$ by changing the coefficients in
the monomial basis to the complex conjugate ones.


At first glance, there is no necessity to involve the matrix space
$\mathrm{Mat}_{m,n}$ into the definition of the Fock inner product since the
definition depends on the dimension $m\cdot n$ of $\mathrm{Mat}_{m,n}$ only
(indeed, $\mathrm{tr}(\mathbb{Z}\mathbb{Z}^*)$ coincides with the usual
Euclidean norm under the natural identification
$\mathrm{Mat}_{m,n}\cong\mathbb{C}^{mn}$). Our motivation comes from the
theory of bounded symmetric domains (see, for instance, \cite{uni3fock.texArazy}).
Suppose $\mathrm{D}\subset\mathbb{C}^d$ is a bounded symmetric domain. One
associates to $\mathrm{D}$ a Euclidean norm $\|\cdot\|$ in $\mathbb{C}^d$
and then the Fock inner product in the space of polynomials on
$\mathbb{C}^d$:

$$ \left(P\,,\,Q\right)_F=\int\limits_{\mathbb{C}^d}P(\mathbf{z})
\overline{Q(\mathbf{z})}e^{-\|\mathbf{z}\|^2}d\mathbf{z}$$ with
$\mathbf{z}=(z_1,z_2,\ldots,z_d)$ and $d\mathbf{z}$ being the
Lebesgue measure on $\mathbb{C}^d$ normalized by the same
condition as in (\ref{uni3fock.texclfip}). Let $G$ be the group of
biholomorphic automorphisms of $\mathrm{D}$ and $K\subset G$ the
stationary subgroup of the center $0\in \mathrm{D}$. Then $K$ is a
compact subgroup in $GL(d)$, and the Fock inner product is
invariant under the induced action of $K$ in the space of
polynomials on $\mathbb{C}^d$. The same $\mathbb{C}^d$ may contain
different bounded symmetric domains, and even if the associated
Fock inner products coincide, the choice of a particular bounded
symmetric domain indicates the symmetry group of the inner
product. For example, by using the notation $\mathrm{Mat}_{m,n}$
we indicate the connection of the inner product (\ref{uni3fock.texclfip}) to
the matrix unit ball
$\mathrm{D}_{m,n}=\{\mathbf{Z}\,|\,\mathbf{Z}\mathbf{Z}^*<\mathrm{I}\}\subset
\mathrm{Mat}_{m,n}$. The ball, we recall, is a homogeneous space
of the group $SU(n,m)$, and the stationary subgroup of the zero
matrix $0$ is $S(U(n)\times U(m))$. The action of the latter group
in $\mathbb{C}[\mathrm{Mat}_{m,n}]$ is described explicitly by $$
P^g(\mathbf{Z})=P(\overline{g}_2^{\,t}\mathbf{Z}g_1),\quad
g=(g_1,g_2)\in S(U(n)\times U(m)),$$ and the invariance of the
Fock inner product is written as follows

\begin{equation}\label{uni3fock.tex2pro}
  \left(P^g\,,\,Q^g\right)_F=\left(P\,,\,Q\right)_F,\quad\forall
g\in S(U(n)\times U(m)).
\end{equation}

It worth noting that easy computability and the invariance of the
Fock inner product make it extremely useful in function theory on
bounded symmetric domains.

We turn now to the subject of the present paper. It is known that the
algebra $\mathbb{C}[\mathrm{Mat}_{m,n}]$ has a non-commutative counterpart
$\mathbb{C}[\mathrm{Mat}_{m,n}]_q$ studied in quantum group theory. The
latter algebra is the unital algebra given by its generators $z_a^\alpha$,
$a=1,\ldots n$, $\alpha=1,\ldots m$, and the following relations
\begin{equation}\label{uni3fock.texz}
z_a^\alpha z_b^\beta=\left \{\begin{array}{ccl}qz_b^\beta z_a^\alpha &,& a=b
\;\&\;\alpha<\beta \quad{\rm or}\quad a<b \;\&\;\alpha=\beta \\ z_b^\beta
z_a^\alpha &,& a<b \;\&\;\alpha>\beta \\ z_b^\beta
z_a^\alpha+(q-q^{-1})z_a^\beta z_b^\alpha &,& a<b \;\&\;\alpha<\beta
\end{array} \right.,
\end{equation}
{\it Throughout the paper $q$ is supposed to be a number from the
interval $(0;1)$ and the ground field is $\mathbb{C}$.}

Around five years ago L.Vaksman and S.Sinel'shchikov \cite{uni3fock.texSV2}
constructed certain analogs of bounded symmetric domains in
framework of quantum group theory. Namely, they associated to a
bounded symmetric domain $\mathrm{D}$ in a complex vector space
$V$ certain non-commutative algebras $\mathbb{C}[V]_q$,
$\mathrm{Pol}(V)_q$, which they treated as the algebras of
holomorphic resp. arbitrary polynomials on the {\it quantum}
vector space $V$. The algebras of continuous or finite functions
on the {\it quantum} bounded symmetric domain $\mathrm{D}$ are
derived then from $\mathrm{Pol}(V)_q$ via some completion
procedure. Further investigation has shown \cite{uni3fock.texSSV1,uni3fock.texSSV3} that
in the case of the matrix unit ball (i.e. in the case
$V=\mathrm{Mat}_{m,n}$) the corresponding algebra of holomorphic
polynomials from \cite{uni3fock.texSV2} is just the algebra
$\mathbb{C}[\mathrm{Mat}_{m,n}]_q$.

It turned out \cite{uni3fock.texSSV1,uni3fock.texSSV2} that many constructions and
problems of classical function theory in bounded symmetric domains
admit natural generalization to the quantum setting. Thus, it is
reasonable to expect that there should be an appropriate
$q$-analog of the Fock inner product and that it will serve the
same purposes as its classical counterpart does. The aim of the
present paper is to present a candidate for that $q$-analog in the
case of the quantum matrix space $\mathrm{Mat}_{m,n}$.

The structure of the paper is as follows. In the next section we
formulate our main results. Sections 3, 4, and 5 contain an
auxiliary material concerning quantum group symmetry.
Specifically, in section 3 we discuss the notions of module
algebras and modules over module algebras. In sections 4,5 we
describe an action of a quantum group in algebras and spaces we
deal with. The proof of the main results is given in section 6.

\bigskip

The author is grateful to L. Vaksman for numerous interesting
discussions of the results.

\bigskip
\section{Statement of main results}

The equality (\ref{uni3fock.tex1pro}) served for us as a guide when we were looking for
a $q$-analog of the Fock inner product. The point is that there are natural
quantum analogs of the partial derivatives $\frac{\partial}{\partial
z_a^\alpha}$. They are constructed by using a first order differential
calculus over the algebra $\mathbb{C}[\mathrm{Mat}_{m,n}]_q$. Let us recall
a definition of that differential calculus. Let
$\Lambda^1(\mathrm{Mat}_{m,n})_q$ be the
$\mathbb{C}[\mathrm{Mat}_{m,n}]_q$-bimodule given by its generators
$dz_a^\alpha$, $a=1,\ldots n$, $\alpha=1,\ldots m$, and the relations

\begin{equation}\label{uni3fock.texzdz}
z_b^\beta dz_a^\alpha=\sum_{\alpha',\beta'=1}^m \sum_{a',b'=1}^n R_{\beta
\alpha}^{\beta'\alpha'}R^{b'a'}_{ba}dz_{a'}^{\alpha'}\cdot z_{b'}^{\beta'},
\end{equation}
with
\begin{equation}\label{uni3fock.texR}
R^{b'a'}_{ba}=\left \{\begin{array}{ccl}q^{-1} &,& a=b=a'=b' \\ 1 &,& a \ne
b \quad \&\quad a=a'\quad \&\quad b=b'\\ q^{-1}-q &,& a<b \quad \&\quad a=b'
\quad \&\quad b=a'\\ 0 &,& \mathrm{otherwise}\end{array}\right..
\end{equation}
The map $d:z_a^\alpha \mapsto dz_a^\alpha$ can be extended up to a
linear operator $d:\mathbb{C}[\mathrm{Mat}_{m,n}]_q \rightarrow
\Lambda^1(\mathrm{Mat}_{m,n})_q$ satisfying the Leibnitz rule. The
pair $\left(\Lambda^1(\mathrm{Mat}_{m,n})_q, d \right)$ is the
first order differential calculus over
$\mathbb{C}[\mathrm{Mat}_{m,n}]_q$. It worth noting that this
first order differential calculus coincides with certain
'canonical' one, defined in \cite{uni3fock.texSV2}.

The $q$-analogs of the partial derivatives may be defined now via
the differential $d$ as follows: $$
df=\sum_{a=1}^{n}\sum_{\alpha=1}^{m} \frac{\partial f}{\partial
z_a^{\alpha}}dz_a^{\alpha},\qquad f\in
\mathbb{C}[\mathrm{Mat}_{m,n}]_q.$$

It is reasonable to rise the question about existence of an inner product in
$\mathbb{C}[\mathrm{Mat}_{m,n}]_q$ with the property (\ref{uni3fock.tex1pro}) where the
classical partial derivatives are replaced by the quantum ones. If such an
inner product existed we would have the algebra {\it antihomomorphism}
\begin{equation}\label{uni3fock.texhom}
\mathbb{C}[\mathrm{Mat}_{m,n}]_q\to\mathbb{C}[\partial]_q, \qquad
z_a^\alpha\mapsto\frac{\partial }{\partial z_a^\alpha},\quad
\forall a,\alpha
\end{equation}
 where $\mathbb{C}[\partial]_q$ stands for the algebra
of quantum differential operators with constant coefficients (i.e., the
unital algebra of linear operators in $\mathbb{C}[\mathrm{Mat}_{m,n}]_q$
generated by the quantum partial derivatives). However, in reality
(\ref{uni3fock.texhom}) is the algebra {\it homomorphism}, that is, the partial
derivatives satisfy the same relations as the generators of
$\mathbb{C}[\mathrm{Mat}_{m,n}]_q$ do. To prove this, one should use the
higher order differential calculus over $\mathbb{C}[\mathrm{Mat}_{m,n}]_q$
associated to the first order one \cite{uni3fock.texSV1,uni3fock.texSV2,uni3fock.texSSV1}. Specifically, let
$\Lambda(\mathrm{Mat}_{m,n})_q$ be the unital algebra given by its
generators $z_a^\alpha$, $dz_a^\alpha$, $a=1,\ldots n$, $\alpha=1,\ldots m$,
satisfying (\ref{uni3fock.texz}), (\ref{uni3fock.texzdz}), and the relations

\begin{equation}\label{uni3fock.texdzdz}
dz_b^\beta dz_a^\alpha=-\sum_{\alpha',\beta'=1}^m \sum_{a',b'=1}^n
R_{\beta
\alpha}^{\beta'\alpha'}R^{b'a'}_{ba}dz_{a'}^{\alpha'}\cdot
dz_{b'}^{\beta'},
\end{equation}
with $R$ defined by (\ref{uni3fock.texR}). The algebra admits the natural
grading by degrees of differential forms. There exists a unique
extension of the differential
$d:\mathbb{C}[\mathrm{Mat}_{m,n}]_q\to\Lambda^1(\mathrm{Mat}_{m,n})_q$
to a linear operator
$d:\Lambda(\mathrm{Mat}_{m,n})_q\to\Lambda(\mathrm{Mat}_{m,n})_q$
which satisfies the (graded) Leibnitz rule and the property
$d^2=0$, which implies, in particular, $d:dz_a^\alpha\mapsto0$
(for all $a,\alpha$). The property means

$$ \sum_{a,a'=1}^{n}\sum_{\alpha,\alpha'=1}^{m}
\frac{\partial}{\partial z_a^{\alpha}}\frac{\partial}{\partial
z_{a'}^{\alpha'}}(f)dz_a^{\alpha}dz_{a'}^{\alpha'}=0,\qquad \forall f\in
\mathbb{C}[\mathrm{Mat}_{m,n}]_q.$$ This equality and the relations
(\ref{uni3fock.texdzdz}) imply the desired commutation relations between the quantum
partial derivatives:

\begin{equation}\label{uni3fock.texd/dz}
\frac{\partial}{\partial z_a^\alpha}\frac{\partial}{\partial
z_b^\beta}=\left \{\begin{array}{ccl}q\frac{\partial}{\partial
z_b^\beta}\frac{\partial}{\partial  z_a^\alpha} &,& a=b
\;\&\;\alpha<\beta \quad{\rm or}\quad a<b \;\&\;\alpha=\beta \\
\frac{\partial}{\partial z_b^\beta}\frac{\partial}{\partial
z_a^\alpha} &,& a<b \;\&\;\alpha>\beta \\ \frac{\partial}{\partial
z_b^\beta}\frac{\partial}{\partial
z_a^\alpha}+(q-q^{-1})\frac{\partial}{\partial z_a^\beta}
\frac{\partial}{\partial z_b^\alpha} &,& a<b \;\&\;\alpha<\beta
\end{array} \right..
\end{equation}

Thus we can't expect existence of an inner product in
$\mathbb{C}[\mathrm{Mat}_{m,n}]_q$ satisfying the property (\ref{uni3fock.tex1pro}).
However, the relations (\ref{uni3fock.texd/dz}) suggest a way to overcome the problem:
one can try to look for such an inner product which makes the operator
$\frac{\partial}{\partial z_a^\alpha}$ conjugate to the operator of {\it
right multiplication} by $z_a^\alpha$ in $\mathbb{C}[\mathrm{Mat}_{m,n}]_q$.
Such an inner product turned out to exist, and this observation is one of
the main results of the paper:

\bigskip
\begin{theorem}\label{uni3fock.tex1}

There exists a unique inner product $(\,\cdot\,,\,\cdot\,)_{F}$ in
$\mathbb{C}[\mathrm{Mat}_{m,n}]_q$ satisfying the properties
\begin{equation}\label{uni3fock.texproo}
  \left(1\,,\,1\right)_{F}=1,
\end{equation}
\begin{equation}\label{uni3fock.texpro}
  \left(\frac{\partial P}{\partial z_a^\alpha}\,,\,Q\right)_{F}=
\left(P\,,\, Q\cdot z_a^\alpha\right)_{F} \quad \forall a,\alpha.
\end{equation}
\end{theorem}
\bigskip
Note that uniqueness follows immediately from the two properties
since any partial derivative is an operator of degree $-1$ with
respect to the natural $\mathbb{Z}_+$-grading in
$\mathbb{C}[\mathrm{Mat}_{m,n}]_q$ by powers of monomials.

Our next result may be formulated as follows: the inner product
$(\,\cdot\,,\,\cdot\,)_{F}$ is invariant with respect to an action
of the {\it quantum} group $S(U(n)\times U(m))$ in
$\mathbb{C}[\mathrm{Mat}_{m,n}]_q$. To present a precise
formulation, we need some preparation.

Let us recall the notion of the Drinfeld-Jimdo quantized universal
enveloping algebra of $\mathfrak{sl}(k)$.  Let $(a_{ij})$ be the
Cartan matrix for $\mathfrak{sl}(k)$. The Hopf algebra $U_q
\mathfrak{sl}(k)$ is determined by the generators $E_i$, $F_i$,
$K_i$, $K_i^{-1}$, $i=1,\ldots,k-1$, and the relations

$$ K_iK_j=K_jK_i,\qquad K_iK_i^{-1}=K_i^{-1}K_i=1,\qquad
K_iE_j=q^{a_{ij}}E_jK_i, $$ $$ K_iF_j=q^{-a_{ij}}F_jK_i, \qquad
E_iF_j-F_jE_i=\delta_{ij}(K_i-K_i^{-1})/(q-q^{-1}) $$
\begin{equation}
E_i^2E_j-(q+q^{-1})E_iE_jE_i+E_jE_i^2=0,\qquad |i-j|=1
\end{equation}
$$F_i^2F_j-(q+q^{-1})F_iF_jF_i+F_jF_i^2=0,\qquad |i-j|=1$$
$$[E_i,E_j]=[F_i,F_j]=0,\qquad |i-j|\ne 1.$$

The comultiplication $\Delta$, the antipode $S$, and the counit
$\varepsilon$ are determined by
\begin{equation}
\Delta(E_i)=E_i \otimes 1+K_i \otimes E_i,\;\Delta(F_i)=F_i
\otimes K_i^{-1}+1 \otimes F_i,\;\Delta(K_i)=K_i \otimes K_i,
\end{equation}
\begin{equation}
S(E_i)=-K_i^{-1}E_i,\qquad S(F_i)=-F_iK_i,\qquad S(K_i)=K_i^{-1},
\end{equation}
$$\varepsilon(E_i)=\varepsilon(F_i)=0,\qquad \varepsilon(K_i)=1.$$

Let $U_q (\mathfrak{sl}(n)\oplus\mathfrak{sl}(m))$ stands for the
Hopf algebra $U_q \mathfrak{sl}(n)\otimes U_q\mathfrak{sl}(m)$.

Recall the standard terminology. Let $A$ be a Hopf algebra. An
algebra $F$ is said to be an $A$-module algebra if $F$ carries a
structure of $A$-module and multiplication in $F$ agrees with the
$A$-action (i. e. the multiplication $F\otimes F\rightarrow F$ is
a morphism of $A$-modules).

The algebra $\mathbb{C}[\mathrm{Mat}_{m,n}]_q$ possesses the well
known structure of
$U_q(\mathfrak{sl}(n)\oplus\mathfrak{sl}(m))$-module algebra:
\begin{equation}\label{uni3fock.texh's1}
K_i\otimes1(z_a^\alpha)=\left \{\begin{array}{ccl}qz_a^\alpha,
\quad a=i\\ q^{-1}z_a^\alpha,\quad a=i+1 \\ z_a^\alpha
,\quad\mathrm{otherwise}\end{array}\right., \quad 1\otimes
K_j(z_a^\alpha)=\left \{\begin{array}{ccl}qz_a^\alpha,\quad
\alpha=m-j
\\ q^{-1}z_a^\alpha,\quad \alpha=m-j+1 \\ z_a^\alpha
,\quad\mathrm{otherwise}\end{array}\right.,
\end{equation}

\begin{equation}\label{uni3fock.texf's1}
F_i\otimes1(z_a^\alpha)= \left
\{\begin{array}{ccl}q^{1/2}z_{a+1}^\alpha,\quad a=i \\ 0,\quad{\rm
otherwise}\end{array}\right.,\quad 1\otimes F_j(z_a^\alpha)= \left
\{\begin{array}{ccl}q^{1/2}z_a^{\alpha+1},\quad \alpha=m-j
\\ 0,\quad{\rm otherwise}\end{array}\right.,
\end{equation}

\begin{equation}\label{uni3fock.texe's1}
E_i\otimes1(z_a^\alpha)= \left
\{\begin{array}{ccl}q^{-1/2}z_{a-1}^\alpha,\quad  a=i+1 \\
 0,\quad\mathrm{otherwise}\end{array}\right.,\quad 1\otimes E_j(z_a^\alpha)= \left \{\begin{array}{ccl}
q^{-1/2}z_a^{\alpha-1},\quad \alpha=m-j+1
\\ 0,\quad\mathrm{otherwise}\end{array}\right..
\end{equation}

The Hopf algebra $U_q\mathfrak{s}(
\mathfrak{gl}(n)\oplus\mathfrak{gl}(m))$ is derived from
$U_q(\mathfrak{sl}(n)\oplus\mathfrak{sl}(m))$ by adding the
generator $K_0$, commuting with the other generators and satisfies
the properties
\begin{equation}\label{uni3fock.tex0}
\Delta(K_0)=K_0 \otimes K_0, \quad S(K_0)=K_0^{-1}, \quad
\varepsilon(K_0)=1.
\end{equation}
Let us extend the
$U_q(\mathfrak{sl}(n)\oplus\mathfrak{sl}(m))$-module algebra
structure in $\mathbb{C}[\mathrm{Mat}_{m,n}]_q$ to a
$U_q\mathfrak{s}(\mathfrak{gl}(n)\oplus\mathfrak{gl}(m))$-module
algebra structure as follows:

\begin{equation}\label{uni3fock.texh0}
K_0(z_a^\alpha)=q^{n+m}z_a^\alpha\quad\forall a,\alpha.
\end{equation}

The $*$-Hopf algebra $U_q(\mathfrak{su}(n)\oplus\mathfrak{su}(m))$
is the pair $(U_q(\mathfrak{sl}(n)\oplus\mathfrak{sl}(m)), *)$
with $*$ being the involution in
$U_q(\mathfrak{sl}(n)\oplus\mathfrak{sl}(m))$ given by

\begin{align*}
(E_i\otimes1)^*&=K_iF_i\otimes1,&(F_i\otimes1)^*&=E_iK_i^{-1}\otimes1,&(K_i^{\pm
1}\otimes1)^*&=K_i^{\pm 1}\otimes1,&
\\ (1\otimes E_j)^*&=1\otimes K_jF_j,&(1\otimes F_j)^*&=1\otimes E_jK_j^{-1},
&(1\otimes K_j^{\pm 1})^*&=1\otimes K_j^{\pm1}.
\end{align*}
This involution is extended to an involution in
$U_q\mathfrak{s}(\mathfrak{gl}(n)\oplus\mathfrak{gl}(m))$ by
setting $K_0^*=K_0.$ The resulting $*$-Hopf algebra
$(U_q\mathfrak{s}(\mathfrak{gl}(n)\oplus\mathfrak{gl}(m)),*)$ is
denoted by
$U_q\mathfrak{s}(\mathfrak{u}(n)\oplus\mathfrak{u}(m))$.

Now we are ready to formulate the invariance property of the
$q$-Fock inner product.

\begin{theorem}\label{uni3fock.tex2}

The inner product $(\,\cdot\,,\,\cdot\,)_{F}$ is
$U_q\mathfrak{s}(\mathfrak{u}(n)\oplus\mathfrak{u}(m))$-invariant,
i.e.
\begin{equation}\label{uni3fock.texinv}
  \left(\xi(P)\,,\,Q\right)_{F}=
\left(P\,,\, \xi^*(Q)\right)_{F}
\end{equation}
for all $P,Q\in\mathbb{C}[\mathrm{Mat}_{m,n}]_q$ and $\xi\in
U_q\mathfrak{s}(\mathfrak{gl}(n)\oplus\mathfrak{gl}(m))$.
\end{theorem}
\bigskip

Our proof of the above two theorems is based on the following more
or less explicit description of the $q$-Fock inner product.

Consider the unital involutive algebra $\mathcal{P}(m,n)_q$ with the
generators $z_a^\alpha$, $a=1,\ldots n$, $\alpha=1,\ldots m$, satisfying
(\ref{uni3fock.texz}) and the relations
\begin{equation}\label{uni3fock.texz,z*}
(z_b^\beta)^\ast z_a^\alpha=\sum_{a^\prime,b^\prime=1}^n
\sum_{\alpha^\prime,\beta^\prime=1}^m \widehat{R}_{b a}^{b^\prime
a^\prime}\widehat{R}_{\beta \alpha}^{ \beta^\prime \alpha^\prime}
z_{a^\prime}^{\alpha^\prime}\left(z_{b^\prime}^{\beta^\prime}\right)^\ast+
\delta_{ab}\delta^{\alpha \beta},
\end{equation}
with $\delta_{ab}$, $\delta^{\alpha \beta}$ being the Kronecker
symbols and
\begin{equation}\label{uni3fock.texRR}
\widehat{R}_{b a}^{b^\prime a^\prime}=\left \{\begin{array}{ccl}1
&,& a \ne b \;\&\;b=b^\prime \;\&\;a=a^\prime \\ q &,&
a=b=a^\prime=b^\prime \\ q-q^{-1} &,& a=b \;\&\;a^\prime=b^\prime
\;\&\;a^\prime>a \\ 0 &,&{\rm otherwise}\end{array}\right..
\end{equation}
Clearly, the algebra $\mathbb{C}[\mathrm{Mat}_{m,n}]_q$ is
embedded into $\mathcal{P}(m,n)_q$.

It worth noting that the algebra $\mathcal{P}(1,n)_q$ is
isomorphic to the well known twisted CCR-algebra introduced by
W.Pusz and S.Woronovicz \cite{uni3fock.texPW}. For arbitrary $m$ and $n$ the
algebra is isomorphic to the algebra
$\mathrm{Pol}(\mathrm{Mat}_{m,n})_q$ of polynomials on the quantum
matrix space \cite{uni3fock.texSSV1}. A precise definition of
$\mathrm{Pol}(\mathrm{Mat}_{m,n})_q$ and a description of the
isomorphism are to be found in Section 4.

Let us consider the $\mathcal{P}(m,n)_q$-module $H$ given by its unique
generator $e_{vac}$ and the relations $(z_a^\alpha)^*e_{vac}=0$ for all $a$
and $\alpha$. The aforementioned isomorphism
$\mathcal{P}(m,n)_q\cong\mathrm{Pol}(\mathrm{Mat}_{m,n})_q$ allows us to use
results from \cite{uni3fock.texSSV1,uni3fock.texSSV3} to derive the following statements concerning
$\mathcal{P}(m,n)_q$. First of all, the multiplication map
$m:\mathcal{P}(m,n)_q\otimes\mathcal{P}(m,n)_q\to\mathcal{P}(m,n)_q$ induces
the isomorphism of vector spaces
$\mathbb{C}[\mathrm{Mat}_{m,n}]_q\otimes\mathbb{C}[\overline{\mathrm{Mat}}_{m,n}]_q
\to\mathcal{P}(m,n)_q$ with
$\mathbb{C}[\overline{\mathrm{Mat}}_{m,n}]_q=\{f^*\,|\,f\in\mathbb{C}[\mathrm{Mat}_{m,n}]_q\}$
(note that surjectivity is a simple consequence of the relations
(\ref{uni3fock.texz,z*})). Thus

\begin{equation}\label{uni3fock.texisom}
 H=\mathbb{C}[\mathrm{Mat}_{m,n}]_q e_{vac}.
\end{equation}

Further, there exists a unique inner product in $H$ so that
\begin{equation}\label{uni3fock.teximpo}
(1\,,\,1)=1,\qquad (f e_1, e_2)=(e_1, f^* e_2)
\end{equation}
 for any $f\in\mathcal{P}(m,n)_q$ and $e_1,e_2\in H$. The equality
(\ref{uni3fock.texisom}) allows us to regard the inner product as the one on
$\mathbb{C}[\mathrm{Mat}_{m,n}]_q$.

\bigskip
\begin{theorem}\label{uni3fock.tex3}

The inner product $(\,\cdot\,,\,\cdot\,)$
 in $\mathbb{C}[\mathrm{Mat}_{m,n}]_q$ coincides with
 $(\,\cdot\,,\,\cdot\,)_{F}$, that is
\begin{equation}\label{uni3fock.texprooo}
\left(\frac{\partial P}{\partial z_a^\alpha}\,,\,Q\right)=( P\,,\,
Q\cdot z_a^\alpha) \quad \forall a,\alpha.
\end{equation}
\end{theorem}
\bigskip

Though we'll present a complete proof of the theorem, it still
seems to be somewhat mysterious since it relates two, at first
glance, different quantum analogs of the Weyl algebra:
$\mathcal{P}(m,n)_q$ and the algebra generated by quantum partial
derivatives and the operators of right multiplication by
$z_a^\alpha$'s (we'll denote the latter algebra by
$A_{m,n}(\mathbb{C})_q$).

Let us explain very briefly a logic of the proof of the three theorems.
Obviously, theorem 3 implies theorem 1. In section 6 we express the inner
product $(\,\cdot\,,\,\cdot\,)$ through another one denoted by
$\langle\,\cdot\,,\,\cdot\,\rangle$. The latter inner product is obviously
invariant and this gives us theorem 2. Thus, the relation (\ref{uni3fock.texprooo}) is,
in a sense, central in the paper, and most of section 6 is devoted to its
proof.

\bigskip
\section{Quantum symmetry}

The aim of this section is to remind some general notions from
quantum group theory.

In this section $A$ denotes a Hopf algebra with the
comultiplication $\Delta$, the antipode $S$, and the counit
$\varepsilon$. An algebra $F$ is said to be an $A$-module algebra
if $F$ carries a structure of $A$-module and multiplication in $F$
agrees with the $A$-action:
\begin{equation}\label{uni3fock.texcov}
\xi(f_1\cdot f_2)=\sum_j\xi'_j(f_1)\cdot\xi''_j(f_2), \quad f_1,
f_2\in F, \quad \xi\in A,
\quad\sum_j\xi'_j\otimes\xi''_j=\Delta(\xi)
\end{equation}
(i. e. the multiplication $F\otimes F\rightarrow F$ is a morphism
of $A$-modules). If $A$ or $F$ have some additional structures,
this definition includes some extra requirements. For example, if
$F$ is unital, one requires $A$-invariance of the unit:

$$\xi(1)=\varepsilon(\xi)\cdot1, \quad \xi\in A.$$

In the case of a $*$-algebra $F$ and a $*$-Hopf algebra $A$ one
imposes the requirement of agreement of the involutions:
\begin{equation}\label{uni3fock.texagree}
(\xi(f))^*=S(\xi)^*(f^*),\quad\xi\in A, f\in F.
\end{equation}

Some examples of module algebras naturally appear in
 representation theory and harmonic analysis.

 Suppose a smooth manifold $X$ is acted by a Lie group $G$.
 This induces an action of the Lie algebra $\mathfrak{g}$ of $G$ in the space $C^\infty(X)$
 by means of
vector fields. In turn, the $\mathfrak{g}$-action induces an
action of the universal enveloping algebra $U\mathfrak{g}$ in
$C^\infty(X)$ by means of differential operators. The usual
Leibnitz rule means that $C^\infty(X)$ is a $U\mathfrak{g}$-module
algebra.

Another important example of a module algebra is the algebra of
linear endomorphisms of a vector space, acted by a Hopf algebra.
Let $V$ be an $A$-module. Endow the space
$\mathrm{End}_{\mathbb{C}}(V)$ with a structure of $A$-module as
follows:
\begin{equation}\label{uni3fock.texend}
\mathrm{ad}\xi(T)=\sum_j \xi'_j \cdot T \cdot S(\xi''_j),
\end{equation}
where $\xi \in A$, $T \in \mathrm{End}_{\mathbb{C}}(V)$,
$\Delta(\xi)=\sum \limits_j \xi'_j \otimes \xi''_j$, $S$ is the
antipode of $A$, and the elements in the right-hand side are
multiplied within the algebra $\mathrm{End}_{\mathbb{C}}(V)$. It
follows from elementary properties of Hopf algebras that this
action of $A$ makes $\mathrm{End}_{\mathbb{C}}(V)$ into an
$A$-module algebra.

Now turn to another important notion, namely, that of modules over
module algebra. Let $F$ be an $A$-module algebra and $M$ be a left
(right) $F$-module. Then $M$ is said to be an $A$-module left
(resp. right) $F$-module, if $M$ carries a structure of $A$-module
and the multiplication $F\otimes M\to M$ (resp. $M\otimes F\to M$)
agrees with the $A$-action, that is, $F\otimes M\to M$ (resp.
$M\otimes F\to M$) is a morphism of $A$-modules.

 Let again $X$ be a smooth manifold acted by a Lie group $G$, and
 $T^*X\to X$ be the cotangent bundle. $T^*X$ inherits a natural
 $G$-action. Then the space of smooth sections of the bundle (i.e.
the space of $1$-forms) is acted by the universal enveloping
algebra $U\mathfrak{g}$, and the action agrees with multiplication
by smooth functions. Thus, $1$-forms constitute a
$U\mathfrak{g}$-module left (as well as right)
$C^\infty(X)$-module.

Let $V$ be an $A$-module, and $\mathrm{End}_{\mathbb{C}}(V)$ is given with a
structure of $A$-module as in (\ref{uni3fock.texend}). It is straightforward that the
natural action $T\otimes v\mapsto T(v)$ makes $V$ into a $A$-module left
$\mathrm{End}_{\mathbb{C}}(V)$-module.

Observe that the latter example is quite general. Indeed, let $F$ be an
$A$-module algebra and $M$ an $A$-module left $F$-module. The multiplication
$F\otimes M\to M$ induces the natural algebra homomorphism
$F\to\mathrm{End}_{\mathbb{C}}(M)$. It is easy to verify that this is a
morphism of $A$-module algebras (here $\mathrm{End}_{\mathbb{C}}(M)$ is
viewed with the $A$-module structure given by (\ref{uni3fock.texend})).

Finally, let us agree about the following useful notation. Suppose $F$ is an
algebra and an $A$-module. If, instead of (\ref{uni3fock.texcov}), we have
\begin{equation}\label{uni3fock.texcovv}
  \xi(f_1\cdot f_2)=\sum_j\xi''_j(f_1)\cdot\xi'_j(f_2)
\quad\forall f_1,f_2\in F,\quad \xi\in A
\end{equation}
 (with $ \sum_j\xi'_j\otimes\xi''_j=\Delta(\xi)$) then we call
$F$ an $A^{op}$-module algebra. Similarly, suppose $F$ is an
$A^{op}$-module algebra and $M$ is a left (right) $F$-module and
an $A$-module. Then $M$ is said to be an $A^{op}$-module left
(resp. right) $F$-module, if the multiplication $F\otimes M\to M$
(resp. $M\otimes F\to M$) satisfies the property

$$  \xi(f\cdot m)=\sum_j\xi''_j(f)\cdot\xi'_j(m) \quad\forall f\in
F,\quad m\in M,\quad \xi\in A$$ (resp. $  \xi(m\cdot
f)=\sum_j\xi''_j(m)\cdot\xi'_j(f)).$

An important example of $A^{op}$-module algebras is constructed as
follows. For an $A$-module $V$ the space
$\mathrm{End}_{\mathbb{C}}(V)$ admits the following alternative
structure of $A$-module:
\begin{equation}\label{uni3fock.texend'}
\mathrm{ad}'\xi(T)=\sum_j \xi''_j \cdot T \cdot S^{-1}(\xi'_j).
\end{equation}
This $A$-action makes $\mathrm{End}_{\mathbb{C}}(V)$ into an
$A^{op}$-module algebra, and the natural action $T\otimes v\mapsto
T(v)$ makes $V$ into an $A^{op}$-module left
$\mathrm{End}_{\mathbb{C}}(V)$-module.

Some important examples to the aforementioned notions appear in
the next sections.

\bigskip
\section{Examples of quantum symmetry: functions in the quantum matrix ball}

In this section we describe some non-trivial examples of a quantum
group action in algebras. The algebras we deal with in this
section are treated as the algebras of functions on certain
quantum $G$-spaces.

The first example of a module algebra already appeared in section 2. There
we described a well known structure
$U_q\mathfrak{s}(\mathfrak{gl}(n)\oplus\mathfrak{gl}(m))$-module algebra
structure in $\mathbb{C}[\mathrm{Mat}_{m,n}]_q$. It is observed in
\cite{uni3fock.texSSV1} that this
$U_q\mathfrak{s}(\mathfrak{gl}(n)\oplus\mathfrak{gl}(m))$-module algebra
structure may be extended to a $U_q\mathfrak{sl}(n+m)$-module one. First, we
have to explain in what sence the latter structure extends the former one.
The point is that there is a natural embedding of Hopf algebras
$U_q\mathfrak{s}(\mathfrak{gl}(n)\oplus\mathfrak{gl}(m))\hookrightarrow
U_q\mathfrak{sl}(n+m)$ determined by

$$ E_i\otimes1\mapsto E_i,\quad F_i\otimes1\mapsto F_i,\quad
K^{\pm1}_i\otimes1\mapsto K^{\pm1}_i,$$

$$ 1\otimes E_j\mapsto E_{n+j},\quad 1\otimes F_j\mapsto
F_{n+j},\quad 1\otimes K^{\pm1}_j\mapsto K^{\pm1}_{n+j},$$

$$K_0\mapsto
K^{nm}_n\cdot\prod_{i=1}^{n-1}K^{mi}_i\cdot\prod_{j=1}^{m-1}K^{n(m-j)}_{n+j}.$$

We also recall one the notation $U_q\mathfrak{su}(n,m)$ for the
$*$-Hopf algebra $(U_q\mathfrak{sl}(n+m),*)$ where the involution
is defined as follows $$\left(K_j^{\pm 1}\right)^*=K_j^{\pm
1},\qquad E_j^*=\left \{
\begin{array}{rr}K_jF_j,& j \ne n \\ -K_jF_j,&j=n \end{array}\right.,
\qquad F_j^*=\left \{\begin{array}{rr}E_jK_j^{-1},& j \ne n \\
-E_jK_j^{-1},&j=n \end{array}\right..$$ One verifies easily that
the above embedding respects the involutions in
$U_q\mathfrak{s}(\mathfrak{u}(n)\oplus\mathfrak{u}(m))$ and
$U_q\mathfrak{su}(n,m)$. Let us agree to use the notation
$U_q\mathfrak{s}(\mathfrak{gl}(n)\oplus\mathfrak{gl}(m))$
($U_q\mathfrak{s}(\mathfrak{u}(n)\oplus\mathfrak{u}(m))$) to
denote the image of the above embedding, i.e. the corresponding
Hopf (resp. $*$-Hopf) subalgebra in $U_q\mathfrak{sl}(n+m)$ (resp.
$U_q\mathfrak{su}(n,m)$).

It follows from (\ref{uni3fock.texh's1}) and (\ref{uni3fock.texh0}) that
\begin{equation}\label{uni3fock.texvak-h}
K_nz_a^\alpha=\left \{\begin{array}{ccl}q^2z_a^\alpha &,&a=n
\;\&\;\alpha=m
\\ qz_a^\alpha &,&a=n \;\&\;\alpha \ne m \quad{\rm or}\quad a \ne n \;\&\;
\alpha=m \\ z_a^\alpha &,&{\rm otherwise}\end{array}\right..
\end{equation}

To describe the $U_q\mathfrak{sl}(n+m)$-module algebra structure
in $\mathbb{C}[\mathrm{Mat}_{m,n}]_q$ completely, we have to add
to (\ref{uni3fock.texh's1}), (\ref{uni3fock.texf's1}), (\ref{uni3fock.texe's1}), and (\ref{uni3fock.texvak-h})
formulae for the action of the generators $E_n$ and $F_n$:
\begin{equation}\label{uni3fock.texvak-f}
F_nz_a^\alpha= \left \{\begin{array}{ccl}q^{1/2} &,& a=n \;\&
\;\alpha=m
\\ 0 &,&{\rm otherwise}\end{array}\right.,
\end{equation}
\begin{equation}\label{uni3fock.texvak-e}
E_nz_a^\alpha=\left \{\begin{array}{ccl}-q^{-1/2}z_a^mz_n^\alpha
 &,&a \ne n \;\&\;\alpha \ne m \\
-q^{1/2}(z_n^m)^2 &,& a=n \;\&\;\alpha=m \\
-q^{1/2}z_n^mz_a^{\alpha} &,&{\rm otherwise}\end{array}\right..
\end{equation}

Let us say few words about the nature of the above
$U_q\mathfrak{sl}(n+m)$-action in
$\mathbb{C}[\mathrm{Mat}_{m,n}]_q$. Its classical counterpart
admits the following description.
 As it is noted in Introduction, the vector space
$\mathrm{Mat}_{m,n}$ contains the so called matrix ball
$\mathrm{D}_{m,n}=\{\mathbf{Z}\,|\,\mathbf{Z}\mathbf{Z}^*<\mathrm{I}\}.
$
 The group
$SU(n,m)$ acts in $\mathrm{D}_{m,n}$ via biholomorphic
automorphisms. Thus, elements of the universal enveloping algebra
$U\mathfrak {su}(n,m)$ (and hence of its complexification
$U\mathfrak{sl}(n+m)$) act in the space of holomorphic functions
in $\mathrm{D}_{m,n}$ by means of differential operators. The
differential operators have polynomial coefficients and,
therefore, preserve the subspace of polynomials. This
$U\mathfrak{sl}(n+m)$-action is the 'classical limit' of the above
$U_q\mathfrak{sl}(n+m)$-action in
$\mathbb{C}[\mathrm{Mat}_{m,n}]_q$.

The algebra $\mathbb{C}[\mathrm{Mat}_{m,n}]_q$ is a particular
case of some algebras constructed in \cite{uni3fock.texSV2} (see
Introduction). Let us recall here the description of
$\mathbb{C}[\mathrm{Mat}_{m,n}]_q$ given in \cite{uni3fock.texSV2} (see also
\cite{uni3fock.texSSV1}, \cite{uni3fock.texSSV3}).

Let us consider the generalized Verma module $V(0)$ over
$U_q\mathfrak{sl}(n+m)$, given by its generator $v(0)$ and the
relations
\begin{equation}\label{uni3fock.texv1}
E_iv(0)=0, \quad K_iv(0)=v(0),\quad \forall i,
\end{equation}
\begin{equation}\label{uni3fock.texv2}
F_iv(0)=0, \quad i\ne n.
\end{equation}
Viewed as a
$U_q\mathfrak{s}(\mathfrak{gl}(n)\oplus\mathfrak{gl}(m))$-module,
the space $V(0)$ splits into direct sum of its finite dimensional
$U_q\mathfrak{s}(\mathfrak{gl}(n)\oplus\mathfrak{gl}(m))$-submodules
$V(0)_{-k}$, $k\in\mathbb{Z}_+$, with
\begin{equation}\label{uni3fock.texgrad}
V(0)_{-k}=\{v\in V(0)|\quad K_0v=q^{(n+m)k}v\}
\end{equation}
(here, we recall,
$K_0=K^{nm}_n\cdot\prod_{i=1}^{n-1}K^{mi}_i\cdot\prod_{j=1}^{m-1}K^{n(m-j)}_{n+j}$).
 Denote by
$V(0)^*$ the {\it graded} dual $U_q\mathfrak{sl}(n+m)$-module: $
V(0)^*= \bigoplus_{k\in\mathbb{Z}_+}\left(V(0)_{-k}\right)^*. $

Let us equip the tensor product $V(0)\otimes V(0)$ with a
$U_q\mathfrak{sl}(n+m)$-module structure via the opposite
comultiplication
\begin{equation}\label{uni3fock.texmodstr}
\xi(v_1\otimes v_2)=\sum_j\xi''_j(v_1)\otimes\xi'_j(v_2), \quad
\xi\in U_q\mathfrak{sl}(n+m),\quad v_1,v_2\in V(0).
\end{equation}
The relations (\ref{uni3fock.texv1}), (\ref{uni3fock.texv2}) imply that the maps $v(0)\mapsto
v(0)\otimes v(0)$, $ v(0)\mapsto1$ are extendable up to morphisms of
$U_q\mathfrak{sl}(n+m)$-modules $ \Delta_-:V(0)\rightarrow V(0)\otimes
V(0)$, $\varepsilon_-:V(0)\rightarrow\mathbb{C}. $ It can be shown that
$\Delta_-$ and $\varepsilon_-$ make $V(0)$ into a coassociative coalgebra
with a counit. Thus, the dual maps
$m=(\Delta_-)^*: V(0)^*\otimes V(0)^* \rightarrow V(0)^*$,
${1}=(\varepsilon_-)^*:\mathbb{C}\rightarrow V(0)^*$
make $V(0)^*$ into an associative unital algebra. Moreover, the product map
$m$ is a morphism of $U_q\mathfrak{sl}(n+m)$-modules and the unit $1$ is
$U_q\mathfrak{sl}(n+m)$-invariant, i.e. $V(0)^*$ is a
$U_q\mathfrak{sl}(n+m)$-module algebra. It turns out to be isomorphic to
$\mathbb{C}[\mathrm{Mat}_{m,n}]_q$ \cite{uni3fock.texSSV1}.

Remind the notation $\mathrm{Pol}(\mathrm{Mat}_{m,n})_q$ for the algebra of
(not necessary holomorphic) polynomials on the quantum matrix space (see
section 2). This algebra is the unital involutive algebra with the
generators $z_a^\alpha$, $a=1,\ldots n$, $\alpha=1,\ldots m$, satisfying
(\ref{uni3fock.texz}) and the relations
\begin{equation}\label{uni3fock.texz,z*1}
(z_b^\beta)^\ast z_a^\alpha=\sum_{a^\prime,b^\prime=1}^n
\sum_{\alpha^\prime,\beta^\prime=1}^m
\widehat{R}_{b,a}^{b^\prime,a^\prime}\widehat{R}_{\beta, \alpha}^{
\beta^\prime, \alpha^\prime}
z_{a^\prime}^{\alpha^\prime}\left(z_{b^\prime}^{\beta^\prime}\right)^\ast+
\delta_{ab}\delta^{\alpha \beta}\cdot(1-q^2)
\end{equation}
(the matrix $R$ is given by (\ref{uni3fock.texRR})). As we stated in section 2,
$\mathrm{Pol}(\mathrm{Mat}_{m,n})_q$ is isomorphic to the $*$-algebra
$\mathcal{P}(m,n)_q$: the isomorphism
$J:\mathrm{Pol}(\mathrm{Mat}_{m,n})_q\to\mathcal{P}(m,n)_q$ is determined by
$J:z_a^\alpha\mapsto (1-q^2)^{\frac12}z_a^\alpha$ for all $a$ and $\alpha$.

The algebra $\mathrm{Pol}(\mathrm{Mat}_{m,n})_q$ is a particular case of
involutive algebras introduced in \cite{uni3fock.texSV2}. It follows from a general
result of \cite{uni3fock.texSV2} that there exists a unique structure of
$U_q\mathfrak{su}(n,m)$-module algebra in
$\mathrm{Pol}(\mathrm{Mat}_{m,n})_q$ so that (\ref{uni3fock.texh's1}), (\ref{uni3fock.texf's1}),
(\ref{uni3fock.texe's1}), (\ref{uni3fock.texvak-h}), (\ref{uni3fock.texvak-f}), (\ref{uni3fock.texvak-e}) hold. This means
that we may use the relation (\ref{uni3fock.texagree}) to 'transfer' the
$U_q\mathfrak{sl}(n+m)$-action from the subalgebra
$\mathbb{C}[\mathrm{Mat}_{m,n}]_q\subset\mathrm{Pol}(\mathrm{Mat}_{m,n})_q$
to the subalgebra
$\mathbb{C}[\overline{\mathrm{Mat}}_{m,n}]_q=\{f^*\,|\,f\in\mathbb{C}[\mathrm{Mat}_{m,n}]_q\}
\subset\mathrm{Pol}(\mathrm{Mat}_{m,n})_q$ of 'antiholomorphic' polynomials,
and the resulting $U_q\mathfrak{sl}(n+m)$-action respects the commutation
relations (\ref{uni3fock.texz,z*1}).

The following example of a $U_q\mathfrak{su}(n,m)$-module algebra appeared
for the first time in \cite{uni3fock.texSSV1}. Let us add to $\mathrm{Pol}(\mathrm
{Mat}_{m,n})_q$ one more generator $f_0$ such that
\begin{equation}\label{uni3fock.texf0} f_0=f_0^2=f_0^*,\quad \left(z_a^\alpha
\right)^*f_0=f_0z_a^\alpha=0, \quad \forall a,\alpha.
\end{equation}
The relations allow us to treat $f_0$ as a q-analogue of the
function, which is equal to 1 in the center of the matrix ball
$\mathrm{D}_{m,n}$ and to 0 in other points. We denote the
$*$-algebra by $\mathrm{Fun}(\mathrm{D}_{m,n})_q$. There exists a
unique extension of the structure of a $U_q
\mathfrak{su}(n,m)$-module algebra from $\mathrm{Pol}(\mathrm
{Mat}_{m,n})_q$ to $\mathrm{Fun}(\mathrm{D}_{m,n})_q$ such that
\begin{equation}\label{uni3fock.texhfe} K_nf_0=f_0,\qquad F_nf_0=-{q^{1/2}\over
q^{-2}-1}f_0 \cdot \left(z_n^m \right)^*,\qquad
E_nf_0=-{q^{1/2}\over 1-q^2}z_n^m \cdot f_0,
\end{equation}

\begin{equation}\label{uni3fock.texcomp}
(K^{\pm1}_j-1)f_0=F_jf_0=E_jf_0=0,\quad j \ne n.
\end{equation}

Obviously, the two-sided ideal
$\mathcal{D}(\mathrm{D}_{m,n})_q{=}\mathrm{Fun}(\mathrm{D}_{m,n})_q\cdot
f_0\cdot\mathrm{Fun}(\mathrm{D}_{m,n})_q$ is a $U_q
\mathfrak{su}(n,m)$-module algebra.  We call its elements finite
functions in the quantum matrix ball $\mathrm{D}_{m,n}$.

In this section we dealt with module algebras of functions on
quantum $G$-spaces. The next section is devoted to a module
algebra of a completely different nature.

\bigskip
\section{Example of quantum symmetry: $q$-Weyl algebra}

In this section we present a remarkable structure of $U_q
\mathfrak{sl}(n+m)^{op}$-module algebra in the $q$-Weyl algebra
$A_{m,n}(\mathbb{C})_q$ (see section 2). The contents of this
section is closely related to that of paper \cite{uni3fock.texSSV4}. However,
the $q$-analog of the Weyl algebra, treated here, differs from
that considered in \cite{uni3fock.texSSV4}.

Remind that $A_{m,n}(\mathbb{C})_q$ is the unital subalgebra in
$\mathrm{End}_{\mathbb{C}}(\mathbb{C}[\mathrm{Mat}_{m,n}]_q)$
generated by the $q$-partial derivatives $\frac{\partial}{\partial
z_a^\alpha}$ and the operators $\hat{z}_a^\alpha$ of right
multiplication by $z_a^\alpha$.

Consider the structure of a $U_q \mathfrak{sl}(n+m)$-module in
$\mathrm{End}_{\mathbb{C}}(\mathbb{C}[\mathrm{Mat}_{m,n}]_q)$ given by
(\ref{uni3fock.texend'}). The following observation plays an important role in our proof
of the main results.

\bigskip
\begin{proposition}\label{uni3fock.tex4}
$A_{m,n}(\mathbb{C})_q$ is a $U_q \mathfrak{sl}(n+m)$-submodule in
$\mathrm{End}_{\mathbb{C}}(\mathbb{C}[\mathrm{Mat}_{m,n}]_q)$. The
induced $U_q \mathfrak{sl}(n+m)$-action makes
$A_{m,n}(\mathbb{C})_q$ into a $U_q
\mathfrak{sl}(n+m)^{op}$-module algebra.
\end{proposition}
\medskip
{\bf Proof.} The second statement is a straightforward consequence
of the first one since any invariant subalgebra in a module
algebra is automatically a module algebra. Let us proof the first
statement.

We have to explain that for arbitrary $\xi\in
U_q\mathfrak{sl}(n+m)$ and $a$, $\alpha$
\begin{equation}\label{uni3fock.texvak-in1}
\mathrm{ad}'\xi(\hat{z}_a^\alpha)\in A_{m,n}(\mathbb{C})_q,
\end{equation}
\begin{equation}\label{uni3fock.texvak-in2}
\mathrm{ad}'\xi(\frac{\partial}{\partial z_a^\alpha})\in
A_{m,n}(\mathbb{C})_q.
\end{equation}
Note that (\ref{uni3fock.texvak-in1}) is a simple consequence of $U_q
\mathfrak{sl}(n+m)$-moduleness of the algebra
$\mathbb{C}[\mathrm{Mat}_{m,n}]_q$. Let us prove (\ref{uni3fock.texvak-in2}).

The crucial role in the proof plays $U_q
\mathfrak{sl}(n+m)$-covariance of the first order differential
calculus $(\Lambda^1(\mathrm{Mat}_{m,n})_q, d)$ (see section 2)
observed for the first time in \cite{uni3fock.texSV1}. The covariance means
that there exists a unique structure of $U_q
\mathfrak{sl}(n+m)$-module in $\Lambda^1(\mathrm{Mat}_{m,n})_q$
such that the differential $d$ is a morphism of the $U_q
\mathfrak{sl}(n+m)$-modules and $\Lambda^1(\mathrm{Mat}_{m,n})_q$
is a $U_q \mathfrak{sl}(n+m)$-module
$\mathbb{C}[\mathrm{Mat}_{m,n}]_q$-bimodule.

That
$d:\mathbb{C}[\mathrm{Mat}_{m,n}]_q\to\Lambda^1(\mathrm{Mat}_{m,n})_q$
intertwines the $U_q \mathfrak{sl}(n+m)$-actions may be written as
follows

 $$ \xi\left(\sum_{b,\beta}\frac{\partial f}{\partial
z_b^{\beta}}dz_b^{\beta}\right)=\sum_{c,\gamma}\frac{\partial
\xi(f)}{\partial z_c^{\gamma}}dz_c^{\gamma},\qquad f\in
\mathbb{C}[\mathrm{Mat}_{m,n}]_q,\quad\xi\in U_q
\mathfrak{sl}(n+m).$$

Since $\Lambda^1(\mathrm{Mat}_{m,n})_q$ is, in particular, a $U_q
\mathfrak{sl}(n+m)$-module {\it left}
$\mathbb{C}[\mathrm{Mat}_{m,n}]_q$-module, we may rewrite the
latter equality:
\begin{equation}\label{uni3fock.texint}
  \sum_{b,\beta}\xi_{(1)}\left(\frac{\partial f}{\partial
z_b^{\beta}}\right)\xi_{(2)}(dz_b^{\beta})=\sum_{c,\gamma}\frac{\partial
\xi(f)}{\partial z_c^{\gamma}}dz_c^{\gamma}
\end{equation}
(here and further in the proof $\xi_{(1)}\otimes \xi_{(2)}$ stands
for $\Delta(\xi)$; similarly, $\xi_{(1)}\otimes
\xi_{(2)}\otimes\xi_{(3)}$ stands for
$1\otimes\Delta(\Delta(\xi))=\Delta\otimes1(\Delta(\xi))$ and so
on).

It is noted in \cite{uni3fock.texSSV1} that $\Lambda^1(\mathrm{Mat}_{m,n})_q$
is a free left $\mathbb{C}[\mathrm{Mat}_{m,n}]_q$-module with the
generators $dz_a^\alpha$, $a=1,\ldots n$, $\alpha=1,\ldots m$.
Thus, for any $\eta\in U_q\mathfrak{sl}(n+m)$ there exist elements
$f_{\beta,a}^{b,\alpha}(\eta)\in\mathbb{C}[\mathrm{Mat}_{m,n}]_q$,
$a, b=1,\ldots n$, $\alpha, \beta=1,\ldots m$, such that $$\eta(
dz_b^\beta)=\sum_{a,\alpha}
f_{\beta,a}^{b,\alpha}(\eta)dz_a^\alpha.$$ Thus, we may rewrite
(\ref{uni3fock.texint}) as follows:

$$  \sum_{b,\beta}\xi_{(1)}\left(\frac{\partial f}{\partial
z_b^{\beta}}\right)f_{\beta,a}^{b,\alpha}(\xi_{(2)})dz_a^{\alpha}=\sum_{c,\gamma}\frac{\partial
\xi(f)}{\partial z_c^{\gamma}}dz_c^{\gamma},$$
or, in terms of operators in $\mathbb{C}[\mathrm{Mat}_{m,n}]_q$,

\begin{equation}\label{uni3fock.texint1}
  \sum_{b,\beta}\mathcal{R}_{f_{\beta,a}^{b,\alpha}(\xi_{(2)})}\cdot\xi_{(1)}\cdot\frac{\partial }{\partial
z_b^{\beta}}=\frac{\partial }{\partial z_a^{\alpha}}\cdot\xi
\end{equation}
where $\mathcal{R}_{f}$ stands for the operator in
$\mathbb{C}[\mathrm{Mat}_{m,n}]_q$ of right multiplication by
$f\in\mathbb{C}[\mathrm{Mat}_{m,n}]_q$.

Now we are ready to prove (\ref{uni3fock.texvak-in2}). By the definition

$$\mathrm{ad}'\xi\left(\frac{\partial }{\partial
z_a^{\alpha}}\right)=\xi_{(2)}\cdot\frac{\partial }{\partial
z_a^{\alpha}}\cdot S^{-1}(\xi_{(1)}). $$ (\ref{uni3fock.texint1}) implies $$
\xi_{(2)}\cdot\frac{\partial }{\partial z_a^{\alpha}}\cdot
S^{-1}(\xi_{(1)})=
\sum_{b,\beta}\xi_{(3)}\cdot\mathcal{R}_{f_{\beta,a}^{b,\alpha}(S^{-1}(\xi_{(1)}))}\cdot
S^{-1}(\xi_{(2)})\cdot\frac{\partial }{\partial z_b^{\beta}}.$$ Observe that
for any $f\in\mathbb{C}[\mathrm{Mat}_{m,n}]_q$ and any $\eta\in
U_q\mathfrak{sl}(n+m)$

$$\eta\cdot\mathcal{R}_{f}=\mathcal{R}_{\eta_{(2)}(f)}\cdot\eta_{(1)}.$$
Thus

$$
\sum_{b,\beta}\xi_{(3)}\cdot\mathcal{R}_{f_{\beta,a}^{b,\alpha}(S^{-1}(\xi_{(1)}))}\cdot
S^{-1}(\xi_{(2)})\cdot\frac{\partial }{\partial
z_b^{\beta}}=\sum_{b,\beta}\mathcal{R}_{\xi_{(4)}(f_{\beta,a}^{b,\alpha}(S^{-1}(\xi_{(1)})))}\cdot\xi_{(3)}
\cdot S^{-1}(\xi_{(2)})\cdot\frac{\partial }{\partial
z_b^{\beta}}.$$ What remains is to take into account the equality

$$\xi_{(1)}\otimes\xi_{(3)}
S^{-1}(\xi_{(2)})\otimes\xi_{(4)}=\xi_{(1)}\otimes\xi_{(2)}.$$ We
finally get

$$
\sum_{b,\beta}\mathcal{R}_{\xi_{(4)}(f_{\beta,a}^{b,\alpha}(S^{-1}(\xi_{(1)})))}\cdot\xi_{(3)}
\cdot S^{-1}(\xi_{(2)})\cdot\frac{\partial }{\partial
z_b^{\beta}}=\sum_{b,\beta}\mathcal{R}_{\xi_{(2)}(f_{\beta,a}^{b,\alpha}(S^{-1}(\xi_{(1)})))}\cdot\frac{\partial
}{\partial z_b^{\beta}}.$$ This finishes the proof of the
proposition.\hfill $\blacksquare$

\bigskip

It is not very difficult to describe the  $U_q
\mathfrak{sl}(n+m)$-action in $A_{m,n}(\mathbb{C})_q$ explicitly.
But, for the purposes of the present paper, we need only the
following partial result.

Impose the notation $$L=\text{\rm linear}\,\text{\rm
span}\,\text{\rm of}\,\{z_a^\alpha\}_{a,\alpha}\subset
\mathbb{C}[\mathrm{Mat}_{m,n}]_q,$$ $$L'=\text{\rm
linear}\,\text{\rm span}\,\text{\rm
of}\,\{\hat{z}_a^\alpha\}_{a,\alpha}\subset
\mathrm{End}_{\mathbb{C}}(\mathbb{C}[\mathrm{Mat}_{m,n}]_q),$$
$$L''=\text{\rm linear}\,\text{\rm span}\,\text{\rm
of}\,\left\{\frac{\partial}{\partial
z_a^\alpha}\right\}_{a,\alpha}\subset
\mathrm{End}_{\mathbb{C}}(\mathbb{C}[\mathrm{Mat}_{m,n}]_q).$$

Note that $L$ is a (finite dimensional)
$U_q\mathfrak{s}(\mathfrak{gl}(n)\oplus\mathfrak{gl}(m))$-submodule in
$\mathbb{C}[\mathrm{Mat}_{m,n}]_q$ due to (\ref{uni3fock.texh's1}), (\ref{uni3fock.texf's1}),
(\ref{uni3fock.texe's1}), (\ref{uni3fock.texvak-h}).

\bigskip
\begin{proposition}\label{uni3fock.tex5}\hfill
\medskip

 i) The map $z_a^\alpha\mapsto\hat{z}_a^\alpha$, $a=1,\ldots,n$,
 $\alpha=1,\ldots,m$ is extended to an isomorphism of the
 $U_q\mathfrak{s}(\mathfrak{gl}(n)\oplus\mathfrak{gl}(m))$-modules
 $j\,':L\to L'$;

\medskip

 ii) the map $z_a^\alpha\mapsto\frac{\partial}{\partial z_a^\alpha}$, $a=1,\ldots,n$,
 $\alpha=1,\ldots,m$ is extended to an isomorphism of the vector
 spaces $j\,'':L\to L''$, satisfying the
 following intertwining property:
\begin{equation}\label{uni3fock.texintpro}
 j\,''(\xi(f))=\mathrm{ad}'\omega(\xi)(j\,''(f)),\quad \xi\in
 U_q\mathfrak{s}(\mathfrak{gl}(n)\oplus\mathfrak{gl}(m)),
 f\in\mathbb{C}[\mathrm{Mat}_{m,n}]_q.
\end{equation}
 Here $\omega$ is the
 automorphism of
 $U_q\mathfrak{s}(\mathfrak{gl}(n)\oplus\mathfrak{gl}(m))$ (the
 Chevalley involution) given by $$ \omega(E_i)=-F_i,\quad
 \omega(F_i)=-E_i, \quad \omega(K^{\pm1}_i)=K^{\mp1}_i;$$

\medskip

 iii) $\frac{\partial}{\partial z_1^1}$ is a primitive
 vector in the $U_q\mathfrak{sl}(n+m)$-module
 $\mathrm{End}_{\mathbb{C}}(\mathbb{C}[\mathrm{Mat}_{m,n}]_q)$:

  $$\mathrm{ad}'F_j\left(\frac{\partial}{\partial
  z_1^1}\right)=0 \quad\forall j,$$

  $$\mathrm{ad}'K_j\left(\frac{\partial}{\partial
  z_1^1}\right)=\left\{\begin{array}{ccl}q^{-1}\frac{\partial}{\partial
  z_1^1}&,&j=1\;{\rm or}\; j=n+m-1 \\ \frac{\partial}{\partial
  z_1^1}&,&{\rm otherwise}\end{array}\right..
  $$

\end{proposition}
\medskip
{\bf Proof.} Statement i) follows immediately from the $U_q
\mathfrak{sl}(n+m)$-moduleness of the algebra
$\mathbb{C}[\mathrm{Mat}_{m,n}]_q$.

Let us prove statement ii). First of all, let us prove that $L''$
is a (finite dimensional)
$U_q\mathfrak{s}(\mathfrak{gl}(n)\oplus\mathfrak{gl}(m))$-submodule
in $\mathrm{End}_{\mathbb{C}}(\mathbb{C}[\mathrm{Mat}_{m,n}]_q)$.

 We have derived the following
formula for the action of $U_q \mathfrak{sl}(n+m)$ on $q$-partial
derivatives (see the proof of the previous proposition):
\begin{equation}\label{uni3fock.texexpli}
  \mathrm{ad}'\xi\left(\frac{\partial }{\partial
z_a^{\alpha}}\right)=\sum_{b,\beta}\mathcal{R}_{\xi_{(2)}(f_{\beta,a}^{b,\alpha}
(S^{-1}(\xi_{(1)})))} \cdot\frac{\partial }{\partial z_b^{\beta}}
\end{equation}
 where $\mathcal{R}_{f}$ stands for the operator in
$\mathbb{C}[\mathrm{Mat}_{m,n}]_q$ of right multiplication by $f$,
and the functions $f_{\beta,a}^{b,\alpha} (\eta)$ are defined by
means of the equality

$$\eta( dz_a^\alpha)=\sum_{b,\beta}
f_{\beta,a}^{b,\alpha}(\eta)dz_b^\beta,\quad\eta\in U_q
\mathfrak{sl}(n+m).$$ Observe that the elements
$f_{\beta,a}^{b,\alpha}(\eta)\in\mathbb{C}[\mathrm{Mat}_{m,n}]_q$
are constants for $\eta\in
U_q\mathfrak{s}(\mathfrak{gl}(n)\oplus\mathfrak{gl}(m))$. This is
a consequence of (\ref{uni3fock.texh's1}), (\ref{uni3fock.texf's1}), (\ref{uni3fock.texe's1}), and
(\ref{uni3fock.texvak-h}). Thus, if $\xi\in
U_q\mathfrak{s}(\mathfrak{gl}(n)\oplus\mathfrak{gl}(m))$, the
elements $\mathcal{R}_{\xi_{(2)}(f_{\beta,a}^{b,\alpha}
(S^{-1}(\xi_{(1)})))}$ are constants as well. The formula
(\ref{uni3fock.texexpli}) then implies that $L''$ is a
$U_q\mathfrak{s}(\mathfrak{gl}(n)\oplus\mathfrak{gl}(m))$-invariant
subspace in
$\mathrm{End}_{\mathbb{C}}(\mathbb{C}[\mathrm{Mat}_{m,n}]_q)$.

Remind (see section 3) that the map
\begin{equation}\label{uni3fock.texopmod}
\mathrm{End}_{\mathbb{C}}(\mathbb{C}[\mathrm{Mat}_{m,n}]_q)\otimes
\mathbb{C}[\mathrm{Mat}_{m,n}]_q\to\mathbb{C}[\mathrm{Mat}_{m,n}]_q,
\quad T\otimes f\mapsto T(f)
\end{equation}
 makes $\mathbb{C}[\mathrm{Mat}_{m,n}]_q$ into a
$U_q\mathfrak{sl}(n+m)^{op}$-module left
$\mathrm{End}_{\mathbb{C}}(\mathbb{C}[\mathrm{Mat}_{m,n}]_q)$-module.
Consequently, the restriction of the map onto the subspace
$L''\otimes L$ defines the map $\tau: L''\otimes
L\to\mathbb{C}[\mathrm{Mat}_{m,n}]_q$, satisfying

$$\xi(\tau(T\otimes
f))=\tau(\mathrm{ad}'\xi_{(2)}(T)\otimes\xi_{(1)}(f)), \quad\xi\in
U_q\mathfrak{s}(\mathfrak{gl}(n)\oplus\mathfrak{gl}(m)).$$
 Clearly, the image of $\tau$ is the subspace of
constants $\mathbb{C}\subset\mathbb{C}[\mathrm{Mat}_{m,n}]_q$.
Thus, $\tau$ defines a
$U_q\mathfrak{s}(\mathfrak{gl}(n)\oplus\mathfrak{gl}(m))$-invariant
pairing $\tau^{op}=\tau\circ\sigma: L\otimes L''\to\mathbb{C}$
with $\sigma$ being the flip of tensor multipliers. Note that the
basis $\{z_a^\alpha\}_{a,\alpha}$ is dual to the basis
$\left\{\frac{\partial}{\partial z_a^\alpha}\right\}_{a,\alpha}$
with respect to the pairing. This observation, together with the
invariance of the pairing, implies, in particular, that the matrix
of the operator $\mathrm{ad}'E_j$ in the basis
$\left\{\frac{\partial}{\partial z_a^\alpha}\right\}_{a,\alpha}$
coincides with the transposed to the matrix of the operator
$-E_jK^{-1}_j$ in the basis $\{z_a^\alpha\}_{a,\alpha}$. But the
latter matrix is easily seen to be the transposed to the matrix of
the operator $-F_i$ in the basis $\{z_a^\alpha\}_{a,\alpha}$. This
proves (\ref{uni3fock.texintpro}) in the case $\xi=E_j$, $j\ne n$. The cases
$\xi=F_j$, $\xi=K^{\pm1}_j$ may be proved in the same way.

Turn to statement iii) of the proposition. Since $K_i$ for any $i$ and $F_j$
for $j\ne n$ belong to
$U_q\mathfrak{s}(\mathfrak{gl}(n)\oplus\mathfrak{gl}(m))$, a part of the
statement follows from statement ii) and formulae (\ref{uni3fock.texh's1}),
(\ref{uni3fock.texe's1}), and (\ref{uni3fock.texvak-h}). What we have to prove is the equality
$$\mathrm{ad}'F_n\left(\frac{\partial}{\partial
  z_1^1}\right)=0.$$

  Use the formula (\ref{uni3fock.texexpli}):

$$  \mathrm{ad}'F_n\left(\frac{\partial }{\partial
z_1^{1}}\right)=\sum_{b,\beta}\mathcal{R}_{F_n(f_{1,a}^{1,\alpha} (1))}
\cdot\frac{\partial }{\partial
z_b^{\beta}}-\sum_{b,\beta}\mathcal{R}_{K^{-1}_n(f_{1,a}^{1,\alpha}
(K_nF_n))} \cdot\frac{\partial }{\partial z_b^{\beta}}. $$ Clearly,
$f_{1,a}^{1,\alpha} (1)=\delta_{a,1}\delta_{\alpha,1}$, and
$F_n(f_{1,a}^{1,\alpha} (1))=0$ for any $a$ and $\alpha$. By using
(\ref{uni3fock.texvak-f}), we also get $\mathcal{R}_{K^{-1}_n(f_{1,a}^{1,\alpha}
(K_nF_n))}=0$. That is

$$  \mathrm{ad}'F_n\left(\frac{\partial }{\partial
z_1^{1}}\right)=0.$$ The proposition is
proved.\hfill$\blacksquare$

\bigskip
\section{Proof of the main results}

In this section we prove the theorems formulated in section 2.

Remind the notation $(\,\cdot\,,\,\cdot\,)$ for the inner product in
$\mathbb{C}[\mathrm{Mat}_{m,n}]_q$ introduced in section 2 (see theorem
\ref{uni3fock.tex3}). First of all, we are going to express this inner product via
another one, which is constructed by means of the algebra
$\mathcal{D}(\mathrm{D}_{m,n})_q$ of finite functions in the quantum matrix
ball (see section 4).

Let $P,Q$ be elements of $\mathbb{C}[\mathrm{Mat}_{m,n}]_q$.
Consider the finite function $f_0\cdot Q^*\cdot P\cdot
f_0\in\mathcal{D}(\mathrm{D}_{m,n})_q$ (here $f_0$ is the
'distinguished' finite function which, along with
$\mathrm{Pol}(\mathrm{Mat}_{m,n})_q$, generates the algebra
$\mathrm{Fun}(\mathrm{D}_{m,n})_q$; see section 4). The properties
(\ref{uni3fock.texf0}) imply that the finite function $f_0\cdot Q^*\cdot
P\cdot f_0$ differs from $f_0$ by a constant. Obviously, the
constant depends linearly on $P$ and conjugate linearly on $Q$.
Thus, we get the sesquilinear form
$\langle\,\cdot\,,\,\cdot\,\rangle$ on
$\mathbb{C}[\mathrm{Mat}_{m,n}]_q$:

\begin{equation}\label{uni3fock.texsesq}
f_0\cdot Q^*\cdot P\cdot f_0=\langle P\,,\,Q\rangle\cdot f_0.
\end{equation}

It turns out to be closely related to the inner product
$(\,\cdot\,,\,\cdot\,)$. Namely, recall one the notation $J$ for
the algebra isomorphism
$\mathrm{Pol}(\mathrm{Mat}_{m,n})_q\to\mathcal{P}(m,n)_q$ (see
section 4). It is explicitly given by $J:z_a^\alpha\mapsto
(1-q^2)^{\frac12}\cdot z_a^\alpha$ for all $a$ and $\alpha$. It
follows almost immediately from the definition of the inner
product $(\,\cdot\,,\,\cdot\,)$ that

\begin{equation}\label{uni3fock.texrela}
\langle P\,,\,Q\rangle=(J(P)\,,\,J(Q)),\qquad
P,Q\in\mathbb{C}[\mathrm{Mat}_{m,n}]_q.
\end{equation}
Indeed, the subspace
$\mathbb{C}[\mathrm{Mat}_{m,n}]_qf_0\subset\mathcal{D}(\mathrm{D}_{m,n})_q$
is invariant under left multiplication by elements of
$\mathrm{Pol}(\mathrm{Mat}_{m,n})_q$ (see (\ref{uni3fock.texf0})). This allows us to
define a structure of $\mathcal{P}(m,n)_q$-module in the subspace
$\mathbb{C}[\mathrm{Mat}_{m,n}]_qf_0\subset\mathcal{D}(\mathrm{D}_{m,n})_q$
by means of the isomorphism $J$:

$$ F\otimes (P\cdot f_0)\mapsto J^{-1}(F)\cdot P\cdot f_0,\quad
F\in \mathcal{P}(m,n)_q, P\in\mathbb{C}[\mathrm{Mat}_{m,n}]_q.$$ Obviously,
the $\mathcal{P}(m,n)_q$-module $\mathbb{C}[\mathrm{Mat}_{m,n}]_qf_0$ is
isomorphic to $H$ (see section 2), and the inner product
$\langle\,\cdot\,,\,\cdot\,\rangle$, regarded as an inner product on
$\mathbb{C}[\mathrm{Mat}_{m,n}]_qf_0$, satisfies the properties
(\ref{uni3fock.teximpo}). Since such an inner product is unique, we have (\ref{uni3fock.texrela}).

In what follows, we denote by $J$ the automorphism of
$\mathbb{C}[\mathrm{Mat}_{m,n}]_q$ given by $J:z_a^\alpha\mapsto
(1-q^2)^{\frac12}\cdot z_a^\alpha$ for all $a$ and $\alpha$. Suppose we have
proved theorem \ref{uni3fock.tex3}. Then theorem \ref{uni3fock.tex2} is a straightforward
consequence of (\ref{uni3fock.texrela}). Indeed, $J$ commutes with the
$U_q\mathfrak{s}(\mathfrak{gl}(n)\oplus\mathfrak{gl}(m))$-action since $J$
is just a constant operator on any homogeneous component of
$\mathbb{C}[\mathrm{Mat}_{m,n}]_q$, and the homogeneous components are
$U_q\mathfrak{s}(\mathfrak{gl}(n)\oplus\mathfrak{gl}(m))$-submodules in
$\mathbb{C}[\mathrm{Mat}_{m,n}]_q$. So we have to prove
$U_q\mathfrak{s}(\mathfrak{u}(n)\oplus\mathfrak{u}(m))$-invariance of the
inner product $\langle\,\cdot\,,\,\cdot\,\rangle$. This may be derived from
the $U_q\mathfrak{s}(\mathfrak{u}(n)\oplus\mathfrak{u}(m))$-moduleness of
the algebra $\mathrm{Fun}(\mathrm{D}_{m,n})_q$ (see section 4) and the
relations (\ref{uni3fock.texcomp}) as follows. The relations (\ref{uni3fock.texcomp}) mean
$U_q\mathfrak{s}(\mathfrak{gl}(n)\oplus\mathfrak{gl}(m))$-invariance of the
element $f_0$:

$$\xi(f_0)=\varepsilon(\xi)\cdot f_0, \quad \xi\in
U_q\mathfrak{s}(\mathfrak{u}(n)\oplus\mathfrak{u}(m)).$$ This and the
equality (\ref{uni3fock.texagree}) imply $$\xi(f_0\cdot Q^*\cdot P\cdot
f_0)=f_0\cdot\xi_{(1)}(Q^*)\cdot\xi_{(2)}(P)\cdot
f_0=f_0\cdot(S(\xi_{(1)})^*(Q))^* \cdot\xi_{(2)}(P)\cdot f_0, \quad
P,Q\in\mathbb{C}[\mathrm{Mat}_{m,n}]_q $$ (here $\xi_{(1)}\otimes \xi_{(2)}$
stands for $\Delta(\xi)$). On the other hand, by (\ref{uni3fock.texsesq}) $f_0\cdot
Q^*\cdot P\cdot f_0=\langle P\,,\,Q\rangle\cdot f_0$, i.e.
$$f_0\cdot(S(\xi_{(1)})^*(Q))^* \cdot\xi_{(2)}(P)\cdot
f_0=\xi(\langle P\,,\,Q\rangle\cdot f_0)=\varepsilon(\xi)\cdot\langle
P\,,\,Q\rangle\cdot f_0.$$ Take into account (\ref{uni3fock.texsesq}) once again:

$$\langle\xi_{(2)}(P)\,,\,S(\xi_{(1)})^*(Q)\rangle
=\varepsilon(\xi)\cdot\langle P\,,\,Q\rangle, \quad \xi\in
U_q\mathfrak{s}(\mathfrak{u}(n)\oplus\mathfrak{u}(m)).$$ The
latter property is obviously equivalent to
$U_q\mathfrak{s}(\mathfrak{u}(n)\oplus\mathfrak{u}(m))$-invariance
of the inner product $\langle\,\cdot\,,\,\cdot\,\rangle$.

In view of the above arguments, the first two theorems, stated in section 2,
follow from the third one. The remaining part of this section is devoted to
a proof of theorem \ref{uni3fock.tex3}. It will be convenient for us to prove the
following statement instead of theorem \ref{uni3fock.tex3} itself. It is equivalent to
the statement of the theorem due to (\ref{uni3fock.texrela}).

\bigskip
\begin{theorem}\label{uni3fock.tex3'}
\begin{equation}\label{uni3fock.texprooo'}
\left\langle\frac{\partial P}{\partial
z_a^\alpha}\,,\,Q\right\rangle=\frac1{1-q^2}\cdot\langle P\,,\, Q\cdot
z_a^\alpha\rangle \quad \forall a,\alpha.
\end{equation}
\end{theorem}
\bigskip
{\bf Proof.} Let us agree about the following notation: if
$T\in\mathrm{End}_{\mathbb{C}}(\mathbb{C}[\mathrm{Mat}_{m,n}]_q)$
then
$T^\dag\in\mathrm{End}_{\mathbb{C}}(\mathbb{C}[\mathrm{Mat}_{m,n}]_q)$
stands for the conjugate operator to $T$ with respect to the inner
product $\langle\,\cdot\,,\,\cdot\,\rangle$. In this notation
(\ref{uni3fock.texprooo'}) says

\begin{equation}\label{uni3fock.tex33'}
  (\hat{z}_a^\alpha)^\dag=(1-q^2)\cdot\frac{\partial}{\partial
z_a^\alpha} \quad \forall a, \alpha.
\end{equation}

Let us explain how
$U_q\mathfrak{s}(\mathfrak{u}(n)\oplus\mathfrak{u}(m))$-invariance
of the inner product $\langle\,\cdot\,,\,\cdot\,\rangle$ allows
one to reduce the general case to the case $a=1$, $\alpha=1$.
Observe that the invariance of the inner product
$\langle\,\cdot\,,\,\cdot\,\rangle$ means
\begin{equation}\label{uni3fock.texdag*}
\xi^\dag=\xi^*, \quad \xi\in
U_q\mathfrak{s}(\mathfrak{u}(n)\oplus\mathfrak{u}(m)).
\end{equation}

Suppose we have already proved (\ref{uni3fock.tex33'}) in the case $a=1$, $\alpha=1$.
Proposition \ref{uni3fock.tex5} (i) and the formulae (\ref{uni3fock.texf's1}) imply
$$\hat{z}_2^1=q^{-1/2}\cdot\mathrm{ad}'F_1(\hat{z}_1^1).$$ Thus
$$(\hat{z}_2^1)^\dag=q^{-1/2}\cdot(\mathrm{ad}'F_1(\hat{z}_1^1))^\dag=
q^{-1/2}\cdot(F_1\cdot\hat{z}_1^1-K_1^{-1}\cdot\hat{z}_1^1\cdot K_1\cdot
F_1)^\dag=$$ $$=q^{-1/2}((\hat{z}_1^1)^\dag \cdot F_1^\dag-F_1^\dag\cdot
K_1^{\dag}\cdot(\hat{z}_1^1)^\dag\cdot
(K_1^{-1})^\dag)=q^{-1/2}((\hat{z}_1^1)^\dag \cdot
E_1K_1^{-1}-E_1K_1^{-1}\cdot K_1\cdot(\hat{z}_1^1)^\dag\cdot K_1^{-1})=$$
$$=-q^{-1/2}(1-q^2)(E_1\cdot\frac{\partial}{\partial z_1^1}\cdot
K_1^{-1}-\frac{\partial}{\partial z_1^1} \cdot E_1K_1^{-1})=
-q^{-1/2}(1-q^2)\mathrm{ad}'E_1\left(\frac{\partial}{\partial
z_1^1}\right)=(1-q^2)\frac{\partial}{\partial z_2^1}$$ (the latter equality
follows from proposition \ref{uni3fock.tex5} (ii)). The other cases may be proved in a
completely analogous way. Thus, it remains to prove (\ref{uni3fock.tex33'}) in the case
$a=1$, $\alpha=1$.

\bigskip
\begin{lemma}\label{uni3fock.texlem1}
Suppose $T\in\mathrm{End}_{\mathbb{C}}(\mathbb{C}[\mathrm{Mat}_{m,n}]_q)$
satisfies the properties

\begin{equation}\label{uni3fock.texprop1}
i)\quad \mathrm{ad}'K_j(T)=\left\{\begin{array}{ccl}q^{-1}T&,&j=1\;{\rm
or}\; j=n+m-1 \\ T&,&{\rm otherwise}\end{array}\right.;
\end{equation}

\begin{equation}\label{uni3fock.texprop2}
ii)\quad\mathrm{ad}'F_j(T)=0 \quad\forall j;
\end{equation}

$$iii)\quad T(z_a^\alpha)=0 \quad \forall a, \alpha.$$

Then $T\equiv0.$
\end{lemma}
\bigskip
{\bf Proof of the lemma.} Let us denote by
$\mathbb{C}[\mathrm{Mat}_{m,n}]_{q,k}$ the $k$-th homogeneous component in
$\mathbb{C}[\mathrm{Mat}_{m,n}]_q$ (remind that the latter algebra admits
the natural $\mathbb{Z}_{+}$-grading by powers of monomials). As we have
already noted, each subspace $\mathbb{C}[\mathrm{Mat}_{m,n}]_{q,k}$ is a
(finite dimensional)
$U_q\mathfrak{s}(\mathfrak{gl}(n)\oplus\mathfrak{gl}(m))$-submodule in
$\mathbb{C}[\mathrm{Mat}_{m,n}]_q$
($U_q\mathfrak{s}(\mathfrak{gl}(n)\oplus\mathfrak{gl}(m))$-invariance of the
subspaces is a straightforward consequence of the formulae (\ref{uni3fock.texh's1}),
(\ref{uni3fock.texf's1}), (\ref{uni3fock.texe's1}), (\ref{uni3fock.texvak-h})). The subspaces
$\mathbb{C}[\mathrm{Mat}_{m,n}]_{q,k}$ admit the following 'coordinateless'
description (see (\ref{uni3fock.texh0}))

\begin{equation}\label{uni3fock.texcoordless}
\mathbb{C}[\mathrm{Mat}_{m,n}]_{q,k}=\{f\in\mathbb{C}
[\mathrm{Mat}_{m,n}]_{q}\,|\,K_0(f)=q^{(n+m)k}f\}
\end{equation}
(remind that $K_0$ corresponds to the element
$K^{nm}_n\cdot\prod_{i=1}^{n-1}K^{mi}_i\cdot\prod_{j=1}^{m-1}K^{n(m-j)}_{n+j}$
under the embedding
$U_q\mathfrak{s}(\mathfrak{gl}(n)\oplus\mathfrak{gl}(m))\subset
U_q\mathfrak{sl}(n+m)$; see section 4). This description, together
with (\ref{uni3fock.texprop1}), implies

$$T(\mathbb{C}[\mathrm{Mat}_{m,n}]_{q,k})\subset\mathbb{C}[\mathrm{Mat}_{m,n}]_{q,k-1}.$$
We are going to prove the equality
$T|_{\mathbb{C}[\mathrm{Mat}_{m,n}]_{q,k}}=0$ by induction in $k$.
The property iii) of $T$ is the induction base.

Suppose we have already proved that
$T|_{\mathbb{C}[\mathrm{Mat}_{m,n}]_{q,k}}=0$ for any $k\leq M$
($M\geq1$), and let
$T|_{\mathbb{C}[\mathrm{Mat}_{m,n}]_{q,M+1}}\neq0$. Let
$f\in\mathbb{C}[\mathrm{Mat}_{m,n}]_{q,M+1}$ be such an element
that $T(f)\neq0$. We may assume that $f$ is a weight vector of the
$U_q \mathfrak{sl}(n+m)$-module
$\mathbb{C}[\mathrm{Mat}_{m,n}]_{q}$ (i.e. an eigenvector of each
$K_j$, $j=1,\ldots,n+m-1$). In this case the element $T(f)$ is a
weight vector as well (this is a consequence of (\ref{uni3fock.texprop1}) and
$U_q\mathfrak{sl}(n+m)^{op}$-moduleness of the left
$\mathrm{End}_{\mathbb{C}}(\mathbb{C}[\mathrm{Mat}_{m,n}]_q)$-module
$\mathbb{C}[\mathrm{Mat}_{m,n}]_q$).  We intend to construct an
element $\tilde{f}\in\mathbb{C}[\mathrm{Mat}_{m,n}]_{q,M+1}$ so
that

$$a)\quad T(\tilde{f})\neq0;$$ $$b)\quad F_j(T(\tilde{f}))=0\quad
\forall j;$$ $$c)\quad T(\tilde{f})\quad \mathrm{is}\quad
\mathrm{a}\quad \mathrm{weight}\quad \mathrm{vector}.$$

To start with, let us note that (\ref{uni3fock.texprop2}) means
\begin{equation}\label{uni3fock.texprop3}
F_jT=\left\{\begin{array}{ccl}qTF_j&,&j=1\;{\rm or}\; j=n+m-1
\\ TF_j&,&{\rm otherwise}\end{array}\right..
\end{equation}
Let some $j_1\neq n$ satisfies $F_{j_1}(T(f))\neq0$. If there is
no such $j_1$ then $\tilde{f}=f$. Indeed, $f$ satisfies the
properties a) and c), and the property b) for $j\neq n$.
$F_n(T(f))=T(F_n(f))$ due to (\ref{uni3fock.texprop3}). But the formula
(\ref{uni3fock.texcoordless}) implies
$F_n(f)\in\mathbb{C}[\mathrm{Mat}_{m,n}]_{q,M}$, and thus
$T(F_n(f))=0$ (the induction assumption). So $\tilde{f}$ has been
built.

If $j_1$ with the above property exists, we set $f_1=F_{j_1}(f)$.
Let $j_2\neq n$ satisfies $F_{j_2}(T(f_1))\neq0$. If there is no
such $j_2$ then we set $\tilde{f}=f_1$ and so on. Since
$\mathbb{C}[\mathrm{Mat}_{m,n}]_{q,M}$ is finite dimensional, this
process will give us an element $\tilde{f}$ satisfying the
properties a), b), and c). Thus we get a non-zero primitive weight
vector (namely, $T(\tilde{f})$) in the $U_q
\mathfrak{sl}(n+m)$-module $\mathbb{C}[\mathrm{Mat}_{m,n}]_{q}$,
which belongs to $\mathbb{C}[\mathrm{Mat}_{m,n}]_{q,M}$ with
$M\geq1$. Let us show that this is a contradiction.

Remind that the $U_q \mathfrak{sl}(n+m)$-module
$\mathbb{C}[\mathrm{Mat}_{m,n}]_{q}$ is the graded dual to a
generalized Verma module $V(0)$ (see section 4). Let us denote by
$(\,\cdot\,,\,\cdot\,)_0$ the pairing
$\mathbb{C}[\mathrm{Mat}_{m,n}]_{q}\times V(0)\to\mathbb{C}$. The
term 'dual $U_q \mathfrak{sl}(n+m)$-module' means
\begin{equation}\label{uni3fock.texinvpai}
(\xi(f)\,,\,v)_0=(f\,,\,S(\xi)(v))_0, \quad \forall\xi\in U_q
\mathfrak{sl}(n+m), f\in\mathbb{C}[\mathrm{Mat}_{m,n}]_{q}, v\in
V(0).
\end{equation}

 Any element from $V(0)$ has the form $\xi(v(0))$, where $v(0)$
 is the generator of $V(0)$
 (see section 4) and $\xi$ is an element from the unital
 subalgebra $U_q\mathfrak{n}_-\subset U_q \mathfrak{sl}(n+m)$ generated by all $F_j$'s.
 This observation is a consequence of the equalities
 (\ref{uni3fock.texv1}) and a PBW-type theorem for the
 quantized universal enveloping algebra $U_q \mathfrak{sl}(n+m)$ \cite{uni3fock.texJ}.

 Let us return now to the non-zero element
 $T(\tilde{f})$. It satisfies
(\ref{uni3fock.texprop2}). One then has

$$(T(\tilde{f})\,,\,\xi(v(0)))_0=
(S^{-1}(\xi)(T(\tilde{f}))\,,\,v(0))_0=\varepsilon(\xi)\cdot(T(\tilde{f})\,,\,v(0))_0$$
for any $\xi\in U_q\mathfrak{n}_-$. This means that the functional
$T(\tilde{f})\in V(0)^*$ is equal to $0$ on any graded component $V(0)_{-k}$
(see (\ref{uni3fock.texgrad})). That is
$T(\tilde{f})\in\mathbb{C}[\mathrm{Mat}_{m,n}]_{q,0}=\mathbb{C}$. Since
$T(\tilde{f})\in\mathbb{C}[\mathrm{Mat}_{m,n}]_{q,M}$ ($M\geq1$), we get
$T(\tilde{f})=0$.

 We see that the assumption
$T|_{\mathbb{C}[\mathrm{Mat}_{m,n}]_{q,M+1}}\neq0$ leads to a
contradiction. The lemma is proved.\hfill$\blacksquare$

\bigskip

\begin{corollary}\label{uni3fock.texcor}
The subspace in
$\mathrm{End}_{\mathbb{C}}(\mathbb{C}[\mathrm{Mat}_{m,n}]_q)$ of
 operators, satisfying the properties i) and ii) of the
previous lemma, is one dimensional.
\end{corollary}
\medskip
{\bf Proof.} Let $T'$ and $T''$ be two non-zero linear operators,
satisfying the properties i) and ii). One has
\begin{equation}\label{uni3fock.tex000}
   T'(z_a^\alpha)= T''(z_a^\alpha)=0,\quad a\neq1\quad
   \mathrm{or}\quad
\alpha\neq1.
\end{equation}
Indeed, by the formulae (\ref{uni3fock.texf's1})

$$z_a^\alpha=c_a^\alpha\cdot F_{a-1}F_{a-2}\ldots
F_{1}F_{n+m+1-\alpha}\ldots F_{n+m-2}F_{n+m-1}(z_1^1)$$ for a
non-zero constant $c_a^\alpha$. Then

$$ T'(z_a^\alpha)=c_a^\alpha\cdot T'(F_{a-1}F_{a-2}\ldots
F_{1}F_{n+m+1-\alpha}\ldots
F_{n+m-2}F_{n+m-1}(z_1^1))=$$$$=c_a^\alpha\cdot
F_{a-1}F_{a-2}\ldots F_{1}F_{n+m+1-\alpha}\ldots
F_{n+m-2}F_{n+m-1}(T'(z_1^1))$$ (the same is true for $T''$). What
remains is to use the fact that $T'(z_1^1),
T''(z_1^1)\in\mathbb{C}$ (we pointed out in the proof of the above
lemma that any linear operator $T$ in
$\mathbb{C}[\mathrm{Mat}_{m,n}]_q$, satisfying i) and ii),
possesses the property
$T(\mathbb{C}[\mathrm{Mat}_{m,n}]_{q,k})\subset\mathbb{C}[\mathrm{Mat}_{m,n}]_{q,k-1}.$)

Suppose $$T'(z_1^1)=c_1,\quad T''(z_1^1)=c_2$$ for certain
constants $c_1$ and $c_2$. If, for example, $c_1=0$ then $T'=0$
due to the above lemma. Thus both constants are non-zero. Then, by
using the lemma once again, we get $$c_2\cdot T'-c_1\cdot
T''=0.$$\hfill$\blacksquare$

\bigskip

To complete the proof of theorem \ref{uni3fock.tex3'}, it suffices to establish that

1. $(\hat{z}_1^1)^\dag$ satisfies the conditions i), ii) of lemma
\ref{uni3fock.texlem1};

2. $\frac{\partial}{\partial z_1^1}$ satisfies the conditions i), ii) of
lemma \ref{uni3fock.texlem1};

3. $(\hat{z}_1^1)^\dag(z_1^1)=(1-q^2)\frac{\partial}{\partial
z_1^1}(z_1^1).$

Note that point 2 is just the statement iii) of proposition \ref{uni3fock.tex5}. Point 3
is a simple consequence of definitions: clearly,
$(1-q^2)\cdot\frac{\partial}{\partial z_1^1}(z_1^1)=1-q^2$; on the other
hand, $(\hat{z}_1^1)^\dag(z_1^1)$ is easily seen to be a constant, so
$$(\hat{z}_1^1)^\dag(z_1^1)=\langle
(\hat{z}_1^1)^\dag(z_1^1)\,,\,1\rangle=\langle
z_1^1\,,\,z_1^1\rangle=1-q^2.$$

What remains is to establish point 1. Observe, that the equalities
$$\mathrm{ad}'K_j((\hat{z}_1^1)^\dag)=\left\{\begin{array}{ccl}q^{-1}(\hat{z}_1^1)^\dag&,&j=1\;{\rm
or}\; j=n+m-1 \\ (\hat{z}_1^1)^\dag&,&{\rm
otherwise}\end{array}\right.,$$
\begin{equation}\label{uni3fock.texpolupro}
\mathrm{ad}'F_j((\hat{z}_1^1)^\dag)=0 \quad j\neq n
\end{equation}
 may be derived from (\ref{uni3fock.texdag*}) and the statement i) of
proposition \ref{uni3fock.tex5}. For example, we prove the equality
$\mathrm{ad}'F_1((\hat{z}_1^1)^\dag)=0$:

$$\mathrm{ad}'F_1((\hat{z}_1^1)^\dag)=
F_1\cdot(\hat{z}_1^1)^\dag-K_1^{-1}\cdot(\hat{z}_1^1)^\dag\cdot K_1\cdot
F_1=$$ $$=(\hat{z}_1^1\cdot F_1^\dag-F_1^\dag\cdot
K_1^{\dag}\cdot\hat{z}_1^1\cdot (K_1^{-1})^\dag)^\dag=(\hat{z}_1^1 \cdot
E_1K_1^{-1}-E_1K_1^{-1}\cdot K_1\cdot\hat{z}_1^1\cdot K_1^{-1})^\dag=$$
$$=-(E_1\cdot\hat{z}_1^1\cdot K_1^{-1}-\hat{ z}_1^1 \cdot E_1K_1^{-1})=
-(\mathrm{ad}'E_1(\hat{z}_1^1))^\dag=0.$$ The other cases in (\ref{uni3fock.texpolupro})
are proved just as this one.

Finally, we have to show that $\mathrm{ad}'F_n((\hat{z}_1^1)^\dag)=0$ or,
equivalently, $F_n\cdot(\hat{z}_1^1)^\dag=(\hat{z}_1^1)^\dag\cdot F_n$ (see
(\ref{uni3fock.texprop3})). The latter equality, in turn, is equivalent to
\begin{equation}\label{uni3fock.texcommut}
  F_n^\dag\cdot \hat{z}_1^1=\hat{z}_1^1\cdot F_n^\dag.
\end{equation}

\bigskip

\begin{lemma}\label{uni3fock.texfinal}
$F_n^\dag=-E_nK_n^{-1}+\frac{q^{1/2}}{1-q^2}\hat{z}_n^m.$
\end{lemma}

\medskip
{\bf Proof of the lemma} to be found in Appendix.
\hfill$\blacksquare$

\bigskip

Now we are ready to prove (\ref{uni3fock.texcommut}). Let
$P\in\mathbb{C}[\mathrm{Mat}_{m,n}]_q$. Then
$$F_n^\dag\hat{z}_1^1(P)-\hat{z}_1^1F_n^\dag(P)=
F_n^\dag(P\cdot z_1^1)-F_n^\dag(P)\cdot z_1^1.$$ By the previous lemma
$$F_n^\dag(P\cdot z_1^1)-F_n^\dag(P)\cdot z_1^1=$$
$$=-E_nK_n^{-1}(P\cdot z_1^1)+\frac{q^{1/2}}{1-q^{2}}P\cdot z_1^1\cdot z_n^m+
E_nK_n^{-1}(P)\cdot z_1^1-\frac{q^{1/2}}{1-q^{2}}P\cdot z_n^m\cdot z_1^1=$$
$$=-E_n(K_n^{-1}(P)\cdot K_n^{-1}(z_1^1))+E_nK_n^{-1}(P)\cdot z_1^1+
\frac{q^{1/2}}{1-q^{2}}P\cdot (z_1^1\cdot z_n^m-z_n^m\cdot
z_1^1)\stackrel{\mathrm{see}(\ref{uni3fock.texvak-h})}{=}$$
$$=-E_nK_n^{-1}(f)\cdot z_1^1-P\cdot E_n(z_1^1)+E_nK_n^{-1}(P)\cdot z_1^1+
\frac{q^{1/2}}{1-q^{2}}P\cdot (z_1^1\cdot z_n^m-z_n^m\cdot z_1^1)=$$
$$=-P\cdot E_n(z_1^1)+
\frac{q^{1/2}}{1-q^{2}}P\cdot (z_1^1\cdot z_n^m-z_n^m\cdot
z_1^1)\stackrel{\mathrm{see}(\ref{uni3fock.texvak-e})}{=}$$
$$=q^{-1/2}P\cdot z_1^m\cdot z_n^1+\frac{q^{1/2}}{1-q^{2}}P\cdot
(z_1^1\cdot z_n^m-z_n^m\cdot z_1^1)\stackrel{\mathrm{see}(\ref{uni3fock.texz})}{=}0.$$

Theorem \ref{uni3fock.tex3'} is proved completely.\hfill$\blacksquare$

\bigskip
\section{Appendix. Proof of lemma \ref{uni3fock.texfinal}} Let
$P,Q\in\mathbb{C}[\mathrm{Mat}_{m,n}]_q$. By the definition

$$\langle F_n(P)\,,\,Q\rangle=f_0\cdot Q^*\cdot F_n(P)\cdot f_0.$$
The equalities (\ref{uni3fock.texhfe}) and $U_q\mathfrak{su}(n,m)$-moduleness of the
algebra $\mathrm{Fun}(\mathrm{D}_{m,n})_q$ (see section 4) imply
$$F_n(P)\cdot f_0= F_n(P)\cdot K_n^{-1}(f_0)=F_n(P\cdot
f_0)-P\cdot F_n(f_0) =F_n(P\cdot f_0)+\frac{q^{1/2}}
{q^{-2}-1}\cdot P\cdot f_0\cdot (z_n^m)^*.$$ Thus $$\langle
F_n(P)\,,\,Q\rangle=f_0\cdot Q^*\cdot F_n(P\cdot
f_0)+\frac{q^{1/2}} {q^{-2}-1}\cdot f_0\cdot Q^*\cdot P\cdot
f_0\cdot (z_n^m)^*=$$ $$=F_n(f_0\cdot Q^*\cdot P\cdot
f_0)-F_n(f_0\cdot Q^*)\cdot K_n^{-1}(P\cdot f_0) +\frac{q^{1/2}}
{q^{-2}-1}\cdot f_0\cdot Q^*\cdot P\cdot f_0\cdot (z_n^m)^*=$$
$$=F_n(\langle P\,,\,Q\rangle\cdot  f_0)-F_n(f_0\cdot Q^*)\cdot
K_n^{-1}(P)\cdot f_0 +\frac{q^{1/2}} {q^{-2}-1}\cdot \langle
P\,,\,Q\rangle\cdot f_0\cdot
(z_n^m)^*\stackrel{\mathrm{see}(\ref{uni3fock.texhfe})}{=}$$$$=\frac{-q^{1/2}}
{q^{-2}-1}\cdot \langle P\,,\,Q\rangle \cdot f_0\cdot (z_n^m)^*-F_n(f_0\cdot
Q^*)\cdot K_n^{-1}(P)\cdot f_0+\frac{q^{1/2}} {q^{-2}-1}\cdot \langle
P\,,\,Q\rangle\cdot f_0\cdot (z_n^m)^*=$$$$= -F_n(f_0\cdot Q^*)\cdot
K_n^{-1}(P)\cdot f_0.$$

Evidently, $$F_n(f_0\cdot Q^*)=F_n(f_0)\cdot K_n^{-1}(Q^*)+f_0\cdot F_n(Q^*)
\stackrel{\mathrm{see}(\ref{uni3fock.texhfe})}{=}$$
$$=-\frac{q^{1/2}} {q^{-2}-1}\cdot f_0\cdot (z_n^m)^*\cdot
K_n^{-1}(Q^*)+f_0\cdot F_n(Q^*) \stackrel{\mathrm{see}(\ref{uni3fock.texagree})}=$$
$$=-\frac{q^{1/2}} {q^{-2}-1}\cdot f_0\cdot (z_n^m)^*\cdot
(K_n(Q))^*+q^2\cdot f_0\cdot (E_n(Q))^*.$$ Finally we get
$$\langle F_n(P)\,,\,Q\rangle=\frac{q^{1/2}} {q^{-2}-1}\cdot
f_0\cdot (z_n^m)^*\cdot (K_n(Q))^*\cdot K_n^{-1}(P)\cdot
f_0-q^2\cdot f_0\cdot (E_n(Q))^*\cdot K_n^{-1}(P)\cdot f_0=$$
$$=\frac{q^{1/2}} {q^{-2}-1}\cdot f_0\cdot (K_n(Q)\cdot
z_n^m)^*\cdot K_n^{-1}(P)\cdot f_0- q^2\cdot f_0\cdot
(E_n(Q))^*\cdot K_n^{-1}(P)\cdot f_0=$$ $$=\frac{q^{1/2}}
{q^{-2}-1}\cdot \langle K_n^{-1}(P)\,,\,K_n(Q)\cdot z_n^m\rangle-
q^2\cdot \langle K_n^{-1}(P)\,,\,E_n(Q)\rangle=$$
$$=\frac{q^{1/2}} {q^{-2}-1}\cdot \langle
P\,,\,K_n^{-1}(K_n(Q)\cdot z_n^m)\rangle- q^2\cdot \langle
P\,,\,K_n^{-1}E_n(Q)\rangle\stackrel{\mathrm{see}(\ref{uni3fock.texvak-h})}=$$
$$=\frac{q^{-3/2}} {q^{-2}-1}\cdot \langle P\,,\,Q\cdot
z_n^m\rangle-\langle P\,,\,E_nK_n^{-1}(Q)\rangle=$$ $$=\langle
P\,,\,-E_nK_n^{-1}(Q)+\frac{q^{-3/2}} {q^{-2}-1}\cdot Q\cdot
z_n^m\rangle.$$

 \pagebreak




\makeatletter \@addtoreset{equation}{section}\makeatother
\renewcommand{\theequation}{\thesection.\arabic{equation}}

\title{\bf GEOMETRIC REALIZATIONS FOR SOME SERIES OF REPRESENTATIONS OF THE
QUANTUM GROUP $\mathbf{SU_{2,2}}$}

\author{D. Shklyarov \and S. Sinel'shchikov \and L. Vaksman}
\date{\tt Institute for Low Temperature Physics \& Engineering\\
47 Lenin Avenue, 61103 Kharkov, Ukraine}

\newpage
\setcounter{section}{0}
\large

\makeatletter
\renewcommand{\@oddhead}{GEOMETRIC REALIZATIONS FOR REPRESENTATIONS OF THE
QUANTUM $\mathbf{SU_{2,2}}$ \hfill \thepage}
\renewcommand{\@evenhead}{\thepage \hfill D. Shklyarov, S. Sinel'shchikov,
and L. Vaksman}
\let\@thefnmark\relax
\@footnotetext{This research was partially supported by Award No UM1-2091 of
the Civilian Research \& Development Foundation \newline \indent This
lecture has been delivered at the 13-th International Hutsulian Workshop,
Kiev, September 2000; published in Matematicheskaya Fizika. Analiz. Geometriya,
{\bf 8} (2001), No 1, 90 -- 110}
\addcontentsline{toc}{chapter}{\@title \\ {\sl D.  Shklyarov, S.
Sinel'shchikov, and L. Vaksman}\dotfill} \makeatother

\maketitle

\section{Introduction}

We consider some series of modules over the quantum universal enveloping
Drinfeld-Jimbo algebra $U_q \mathfrak{g}$ in the special case
$\mathrm{dim}\,\mathfrak{g}<\infty$, $0<q<1$. The finite dimensional $U_q
\mathfrak{g}$-modules are closely related to compact quantum groups; those
were investigated well enough \cite{uni4_rqsu22.texJ,uni4_rqsu22.texR}. Infinite dimensional $U_q
\mathfrak{g}$-modules we deal with in this work originate from our earlier
paper \cite{uni4_rqsu22.texSV}, together with some applications therein to the theory of
q-Cartan domains. To make the exposition more transparent, we restrict
ourselves to a q-analogue of the ball in the space of all complex $2 \times
2$ matrices $\mathbb{U}=\{{\boldsymbol z}\in
\mathrm{Mat}_2|\:\boldsymbol{zz}^*<1 \}$, which is among the simplest Cartan
domains.

The classes of infinite dimensional $U_q \mathfrak{g}$-modules in question
differ from those considered by Letzter \cite{uni4_rqsu22.texL}. The problem of producing
and investigating of the principal series of quantum Harish-Chandra modules
in our case appears to be essentially more complicated.

It is worthwhile to note that some properties of the ladder representation
of the quantum $SU_{2,2}$ described below are already well known \cite{uni4_rqsu22.texDo}.

Everywhere in the sequel $0<q<1$, the ground field is $\mathbb{C}$, and all
the algebras are assumed to be unital, unless the contrary is stated
explicitly.

Consider the Hopf algebra $U_q \mathfrak{g}=U_q \mathfrak{sl}_4$ determined
by the standard lists of generators $E_j$, $F_j$, $K_j^{\pm 1}$, $j=1,2,3$,
and relations \cite{uni4_rqsu22.texJ,uni4_rqsu22.texR}. The coproduct $\triangle$, the counit
$\varepsilon$, and the antipode $S$ are given as follows:
\begin{align*}
\triangle{E_j}&=E_j \otimes 1+K_j \otimes E_j,& \varepsilon(E_j)&=0,&
S(E_j)&=-K_j^{-1}E_j,\\ \triangle{F_j}&=F_j \otimes K_j^{-1}+1 \otimes F_j,&
\varepsilon(F_j)&=0,& S(F_j)&=-F_jK_j,\\ \triangle{K_j}&=K_j \otimes K_j,&
\varepsilon(K_j)&=1,& S(K_j)&=K_j^{-1}.
\end{align*}

We call a $U_q \mathfrak{g}$-module $V$ $\mathbb{R}^3$-weight module if
$V=\bigoplus \limits_{\boldsymbol \mu}V_{\boldsymbol \mu}$ with
${\boldsymbol \mu}=(\mu_1,\mu_2,\mu_3)\in \mathbb{R}^3$, $V_{\boldsymbol
\mu}=\{v \in V|\:K_j^{\pm 1}v=q^{\pm \mu_j}v,\,j=1,2,3 \}$. Let $U_q
\mathfrak{k}\subset U_q \mathfrak{g}$ be the Hopf subalgebra generated by
$K_2^{\pm 1}$, $E_j$, $F_j$, $K_j^{\pm 1}$, $j=1,3$. Every $U_q
\mathfrak{g}$-module inherits a structure of $U_q \mathfrak{k}$-module. We
are interested in quantum $(\mathfrak{g},\mathfrak{k})$-modules, i.e.
$\mathbb{R}^3$-weight $U_q \mathfrak{g}$-modules which are direct sums of
finite dimensional $U_q \mathfrak{k}$-modules.

Equip the Hopf algebra $U_q \mathfrak{g}$ with an involution:
$$E_2^*=-K_2F_2,\qquad F_2^*=-E_2K_2^{-1},\qquad K_2^*=K_2,$$
$$E_j^*=K_jF_j,\qquad F_j^*=E_jK_j^{-1},\qquad K_j^*=K_j,\qquad j=1,3.$$
 We thus get a $*$-Hopf algebra $(U_q \mathfrak{g},*)$ which is a q-analogue
of $U \mathfrak{su}_{2,2}$ and its subalgebra $(U_q \mathfrak{k},*)$ is a
q-analogue of $U \mathfrak{s}(\mathfrak{u}_2 \times \mathfrak{u}_2)$.

A quantum $(\mathfrak{g},\mathfrak{k})$-module $V$ is said to be
unitarizable if $(\xi v_1,v_2)=(v_1,\xi^*v_2)$ for some Hermitian scalar
product in $V$ and all $v_1,v_2 \in V$, $\xi \in U_q \mathfrak{g}$. Our
purpose here is to produce some series of unitarizable quantum
$(\mathfrak{g},\mathfrak{k})$-modules by means of non-commutative geometry
and non-commutative function theory in q-Cartan domains \cite{uni4_rqsu22.texSV,uni4_rqsu22.texSSV1,uni4_rqsu22.texSSV2,uni4_rqsu22.texSSV3,uni4_rqsu22.texSSV4}.

The third named author would like to express his gratitude to H. P.
Jakobsen, \\A. Klimyk, A. Stolin and L. Turowska for helpful discussions.

\bigskip

\section{The $\mathbf{U_q \mathfrak{su}_{2,2}}$-module algebra
$\mathbf{Pol(Pl_{2,4})_{q,x}}$}

Let $e_1,e_2,e_3,e_4$ be the standard basis in $\mathbb{C}^4$. Associate to
every linear operator in $\mathbb{C}^2$ its graph, a two-dimensional
subspace in $\mathbb{C}^4=\mathbb{C}^2 \times \mathbb{C}^2$, which has
trivial intersection with the linear span of $e_1,e_2$. We are interested in
the pairs $(L,\omega)$, with $L$ a subspace as above and $\omega$ its
non-zero volume form (an skew-symmetric bilinear form) in $L$. We need a
q-analogue of this algebraic variety which we call the Pl\"ucker manifold
$\mathrm{Pl}_{2,4}$. The matrix elements $\begin{pmatrix} \alpha & \beta
\\ \gamma & \delta \end{pmatrix}$ of the linear operator $L$, together
with $t^{\pm 1}$ related to the volume element $\omega$, work as
'coordinates' on $\mathrm{Pl}_{2,4}$.

An algebra $F$ is called a $U_q \mathfrak{g}$-module algebra if the
multiplication $m:F \otimes F \to F$ is a morphism of $U_q
\mathfrak{g}$-modules, and the unit $1 \in F$ is a $U_q
\mathfrak{g}$-invariant. To rephrase, one can say that for all $f_1,f_2 \in
F$, $j=1,2,3$,
\begin{align*}
E_j(f_1f_2)&=E_j(f_1)f_2+(K_jf_1)(E_jf_2),& E_j1&=0,\\
F_j(f_1f_2)&=(F_jf_1)(K_j^{-1}f_2)+f_1(F_jf_2),& F_j1&=0,\\ K_j^{\pm
1}(f_1f_2)&=(K_j^{\pm 1}f_1)(K_j^{\pm 1}f_2),& K_j^{\pm 1}1&=1.&
\end{align*}

In the case of a $*$-algebra $F$ one should impose an additional
compatibility requirement for involutions:
$$(\xi f)^*=(S(\xi))^*f^*,\quad \xi \in U_q \mathfrak{g},\;f \in F.$$

Once the $*$-algebra $F$ is given by the list of its generators and
relations, the $U_q \mathfrak{g}$-module structure in $F$ is determined
unambiguously by the action of the generators $E_j$, $F_j$, $K_j^{\pm 1}$,
$j=1,2,3$, on the generators of $F$.

Consider the $*$-algebra $\mathrm{Pol}(\mathrm{Mat}_2)_q$ given by its
generators $\alpha$, $\beta$, $\gamma$, $\delta$ and the following
commutation relations (the initial six of those are well known and the rest
was obtained in \cite{uni4_rqsu22.texSSV2}):
$$\left \{\begin{array}{ccc}\alpha \beta&=&q \beta
\alpha \\ \gamma \delta&=&q \delta \gamma \end{array}\right.\qquad \left
\{\begin{array}{ccc}\alpha \gamma&=&q \gamma \alpha \\ \beta \delta&=&q
\delta \beta \end{array}\right.\qquad \left \{\begin{array}{ccl}\beta
\gamma&=&\gamma \beta \\ \alpha \delta&=&\delta \alpha+(q-q^{-1})\beta
\gamma \end{array}\right.$$
$$\left \{\begin{array}{ccl}\delta^*\alpha&=&\alpha \delta^*\\
\delta^*\beta&=& q \beta \delta^* \\ \delta^*\gamma&=&q\gamma \delta^*\\
\delta^*\delta&=&q^2 \delta \delta^*+1-q^2 \end{array}\right.\qquad \left \{
\begin{array}{ccl}\gamma^*\alpha&=&q \alpha \gamma^*-(q^{-1}-q)\beta \delta^*
\\ \gamma^*\beta&=&\beta \gamma^*\\ \gamma^*\gamma&=&q^2 \gamma
\gamma^*-(1-q^2)\delta \delta^*+1-q^2\end{array}\right.$$
$$\left \{\begin{array}{ccl}\beta^*\alpha&=&q \alpha \beta^*-(q^{-1}-q)\gamma
\delta^*\\ \beta^*\beta&=&q^2 \beta \beta^*-(1-q^2)\delta
\delta^*+1-q^2\end{array}\right.$$
$$\alpha^*\alpha=q^2\alpha \alpha^*-(1-q^2)(\beta \beta^*+\gamma
\gamma^*)+(q^{-1}-q)^2 \delta \delta^*+1-q^2.$$

The $*$-algebra $\mathrm{Pol}(\mathrm{Pl}_{2,4})_{q,x}$ is given by the
generators $\alpha$, $\beta$, $\gamma$, $\delta$, $t$, $t^{-1}$, the
commutation relations as in the above definition of
$\mathrm{Pol}(\mathrm{Mat}_2)_q$, and the additional relations
$tt^{-1}=t^{-1}t=1$, $tt^*=t^*t$, $zt=qtz$, $zt^*=qt^*z$, with $z \in \{
\alpha,\beta,\gamma,\delta \}$\footnote{The notation $x=tt^*$ and
$\mathrm{Pol}(\mathrm{Pl}_{2,4})_{q,x}$ are justified by the fact that the
algebra $\mathrm{Pol}(\mathrm{Pl}_{2,4})_{q,x}$ in question can be derived
as a localization of another useful algebra
$\mathrm{Pol}(\mathrm{Pl}_{2,4})_q$ with respect to the multiplicative
system $x^{\mathbb{N}}$.}.

An application of a q-analogue for the above geometric interpretation of the
Pl\"ucker manifold allows one to prove

\medskip

\begin{proposition} i) There exists a unique structure of $U_q
\mathfrak{su}_{2,2}$-module algebra in $\mathrm{Pol}(\mathrm{Mat}_2)_q$ such
that
\begin{align*}\begin{pmatrix}E_1 \alpha & E_1 \beta \\ E_1 \gamma & E_1
\delta \end{pmatrix}&=q^{-1/2}\begin{pmatrix}0 & \alpha \\ 0 & \gamma
\end{pmatrix},& \begin{pmatrix}E_3 \alpha & E_3 \beta \\ E_3 \gamma &
E_3 \delta \end{pmatrix}&=q^{-1/2}\begin{pmatrix}0 & 0 \\ \alpha & \beta
\end{pmatrix}\\ \begin{pmatrix}F_1 \alpha & F_1 \beta \\ F_1 \gamma & F_1
\delta \end{pmatrix}&=q^{1/2}\begin{pmatrix}\beta & 0\\ \delta & 0
\end{pmatrix},& \begin{pmatrix}F_3 \alpha & F_3 \beta \\ F_3 \gamma & F_3
\delta \end{pmatrix}&=q^{1/2}\begin{pmatrix}\gamma & \delta \\ 0 & 0
\end{pmatrix}\\ \begin{pmatrix}K_1 \alpha & K_1 \beta \\ K_1 \gamma & K_1
\delta \end{pmatrix}&=\begin{pmatrix}q \alpha & q^{-1}\beta \\ q \gamma &
q^{-1}\delta \end{pmatrix},& \begin{pmatrix}K_3 \alpha & K_3 \beta \\ K_3
\gamma & K_3 \delta \end{pmatrix}&=\begin{pmatrix}q \alpha & q \beta \\
q^{-1}\gamma & q^{-1}\delta \end{pmatrix}\\ \begin{pmatrix}E_2 \alpha & E_2
\beta \\ E_2 \gamma & E_2 \delta
\end{pmatrix}&=-q^{1/2}\begin{pmatrix}q^{-1}\beta \gamma & \delta \beta
\\ \delta \gamma & \delta^2 \end{pmatrix},& \begin{pmatrix}F_2 \alpha &
F_2 \beta \\ F_2 \gamma & F_2 \delta \end{pmatrix}&=q^{1/2}\begin{pmatrix}0
& 0 \\ 0 & 1 \end{pmatrix}\\ \begin{pmatrix}K_2 \alpha & K_2 \beta \\ K_2
\gamma & K_2 \delta \end{pmatrix}&=\begin{pmatrix}\alpha & q \beta \\ q
\gamma & q^2 \delta \end{pmatrix}.
\end{align*}
ii) There exists a unique structure of $U_q \mathfrak{su}_{2,2}$-module
algebra in $\mathrm{Pol}(\mathrm{Pl}_{2,4})_{q,x}$ such that the action of
$E_j$, $F_j$, $K_j^{\pm 1}$ on $\alpha$, $\beta$, $\gamma$, $\delta$ is
given by the above equations and
$$\left \{\begin{array}{ccl}E_jt &=& 0 \\ F_jt &=& 0 \\ K_jt &=&
t \end{array}\right.,\quad j=1,3;\qquad \left \{\begin{array}{ccl}E_2t &=&
q^{-1/2}t \delta \\ F_2t &=& 0 \\ K_2t &=& q^{-1}t \end{array}\right..$$
\end{proposition}

\medskip

Note that a much more general result is obtained in \cite{uni4_rqsu22.texSSV2}.

To produce the series of quantum $(\mathfrak{g},\mathfrak{k})$-modules
considered in the sequel we use essentially the specific dependencies of the
elements $E_2t^\lambda$, $F_2t^\lambda$, $K_2^{\pm 1}t^\lambda$,
$E_2((\alpha \delta-q \beta \gamma)^\lambda)$, $F_2((\alpha \delta-q \beta
\gamma)^\lambda)$, $K_2^{\pm 1}((\alpha \delta-q \beta \gamma)^\lambda)$ on
$q^\lambda$. These are easily deducible from the definitions that for all
$\lambda \in \mathbb{Z}_+$
$$\shoveleft{E_2t^\lambda=q^{-3/2}\frac{q^{-2\lambda}-1}{q^{-2}-1}\delta
t^\lambda,\qquad F_2t^\lambda=0,\qquad K_2^{\pm 1}t^\lambda=q^{\mp
\lambda}t^\lambda,}$$
$$E_2((\alpha \delta-q \beta \gamma)^\lambda)=-q^{1/2}\frac{1-q^{2
\lambda}}{1-q^2}\delta(\alpha \delta-q \beta \gamma)^\lambda,$$
$$F_2((\alpha \delta-q \beta \gamma)^\lambda)=q^{1/2}\frac{q^{-2
\lambda}-1}{q^{-2}-1}\alpha(\alpha \delta-q \beta \gamma)^{\lambda-1},\qquad
\lambda \ne 0,$$
$$K_2^{\pm 1}((\alpha \delta-q \beta \gamma)^\lambda)=q^{\pm 2
\lambda}(\alpha \delta-q \beta \gamma)^\lambda.$$

For instance, the first relation is obvious for $\lambda=0$, and the general
case is accessible via an induction argument:
\begin{multline*}E_2(t^{\lambda+1})=(E_2t)t^\lambda+(K_2t)(E_2t^\lambda)=
q^{-1/2}t \delta t^\lambda+q^{-1}tq^{-3/2}\frac{q^{-2
\lambda}-1}{q^{-2}-1}\delta t^\lambda=\\=\left(q^{-1/2}+q^{-5/2}\frac{q^{-2
\lambda}-1}{q^{-2}-1}\right)q^{-1}\delta
t^{\lambda+1}=q^{-3/2}\frac{q^{-2(\lambda+1)}-1}{q^{-2}-1}\delta
t^{\lambda+1}.
\end{multline*}

\bigskip

\section{The analytic continuation of the holomorphic discrete series: step
one}

Consider the subalgebra $\mathbb{C}[\mathrm{Pl}_{2,4}]_{q,t}\subset
\mathrm{Pol}(\mathrm{Pl}_{2,4})_{q,x}$ generated by $\alpha$, $\beta$,
$\gamma$, $\delta$, $t$, $t^{-1}$. Equip it with a $\mathbb{Z}$-grading:
$\deg \alpha=\deg \beta=\deg \gamma=\deg \delta=0$, $\deg(t^{\pm 1})=\pm 1$.
The homogeneous components of this algebra are quantum
$(\mathfrak{g},\mathfrak{k})$-modules\footnote{The notation
$\mathbb{C}[\mathrm{Pl}_{2,4}]_{q,t}$ can be justified in the same way as
the notation $\mathrm{Pol}(\mathrm{Pl}_{2,4})_{q,x}$ introduced in the
previous section.}.

Consider the subalgebra $\mathbb{C}[\mathrm{Mat}_2]_q \subset
\mathrm{Pol}(\mathrm{Mat}_2)_q$ generated by $\alpha$, $\beta$, $\gamma$,
$\delta$. This algebra constitutes a famous subject of a research in the
quantum group theory. Associate to each $\lambda \in \mathbb{Z}$ a linear
operator $i_\lambda:\mathbb{C}[\mathrm{Mat}_2]_q \to
\mathbb{C}[\mathrm{Pl}_{2,4}]_{q,t}$, $i_\lambda:f \mapsto ft^{-\lambda}$.
This isomorphism between the vector space $\mathbb{C}[\mathrm{Mat}_2]_q$ and
a homogeneous component of $\mathbb{C}[\mathrm{Pl}_{2,4}]_{q,t}$ allows one
to transfer the structure of $U_q \mathfrak{sl}_4$-module from
$\mathbb{C}[\mathrm{Pl}_{2,4}]_{q,t}$ to $\mathbb{C}[\mathrm{Mat}_2]_q$.
Thus we obtain a representation of $U_q \mathfrak{sl}_4$ in
$\mathbb{C}[\mathrm{Mat}_2]_q$, to be denoted by $\pi_{q^\lambda}$. For all
$\xi \in U_q \mathfrak{sl}_4$, $f \in \mathbb{C}[\mathrm{Mat}_2]_q$, the
vector valued function $\pi_{q^\lambda}(\xi)f$ appears to be a Laurent
polynomial of an indeterminate $\zeta=q^\lambda$. This leads to the
canonical analytic continuation of the operator valued function
$\pi_{q^\lambda}$. The term 'analytic continuation of the holomorphic
discrete series' stands for the above family $\pi_{q^\lambda}$ of
representations of $U_q \mathfrak{sl}_4$.

The results of the work by H. P. Jakobsen \cite{uni4_rqsu22.texJak} imply that the quantum
$(\mathfrak{g},\mathfrak{k})$-modules $\pi_{q^\lambda}$ are unitarizable for
all $\lambda>1$. We follow \cite{uni4_rqsu22.texSSV3} in finding an explicit form for the
related scalar product.

Consider the $\mathrm{Pol}(\mathrm{Mat}_2)_q$-module given by a single
generator $v$ and the relations $\alpha^*v=\beta^*v=\gamma^*v=\delta^*v=0$.
The associated representation $T$ of $\mathrm{Pol}(\mathrm{Mat}_2)_q$ in the
vector space $H=\mathbb{C}[\mathrm{Mat}_2]_qv$ is faithful; it is called the
vacuum representation.

Let $\check{\rho}$ be the linear operator in $H$ that realizes the action of
the 'half-sum of positive coroots':
$$
\check{\rho}(\alpha^a \beta^b \gamma^c
\delta^dv)=(3a+2b+2c+d)\alpha^a \beta^b \gamma^c \delta^dv,
$$
with $a,b,c,d \in \mathbb{Z}_+$. We need also the element
$$
y=1-(\alpha\alpha^*+\beta\beta^*+\gamma\gamma^*+\delta\delta^*)+
(\alpha\delta-q \beta\gamma)(\alpha\delta-q \beta\gamma)^*,
$$
which is a q-analogue of the determinant $\det(1-\mathbf{zz}^*)$, with
$\mathbf{z}=\begin{pmatrix}\alpha&\beta \\ \gamma&\delta \end{pmatrix}$.

As a consequence of the results of \cite{uni4_rqsu22.texSSV3} we have

\medskip

\begin{proposition} i) For $\lambda>3$ the linear functional
$$\int \limits_{\mathbb{U}_q}fd \nu_\lambda \stackrel{\mathrm{def}}{=}
\frac{\mathrm{tr}(T(fy^\lambda)q^{-2
\check{\rho}})}{\mathrm{tr}(T(y^\lambda)q^{-2 \check{\rho}})}$$
 is well defined and positive on $\mathrm{Pol}(\mathrm{Mat}_2)_q$.

ii) For $\lambda>3$ the scalar product $(f_1,f_2)_{q^{2 \lambda}}=\int
\limits_{\mathbb{U}_q}f_2^*f_1d \nu_\lambda$ in
$\mathbb{C}[\mathrm{Mat}_2]_q$ is well defined, positive, and
$$(\pi_{q^\lambda}(\xi)f_1,f_2)_{q^{2
\lambda}}=(f_1,\pi_{q^\lambda}(\xi^*)f_2)_{q^{2 \lambda}},\qquad \xi \in U_q
\mathfrak{g},\,f_1,f_2 \in \mathbb{C}[\mathrm{Mat}_2]_q.$$
\end{proposition}

\medskip

The representations $\pi_{q^\lambda}$, $\lambda=3,4,5,\ldots$, are
q-analogues of the holomorphic discrete series representations, and the
completions of $\mathbb{C}[\mathrm{Mat}_2]_q$ with respect to the norms $\|f
\|_{q^{2\lambda}}=(f,f)_{q^{2\lambda}}^{1/2}$ are q-analogues of the
weighted Bergman spaces. Our intention in what follows is to present
explicit formulae for the analytic continuation of the scalar product
$(f_1,f_2)_{q^{2 \lambda}}$ with respect to the parameter $q^{2 \lambda}$,
and to prove the positivity of this scalar product for $\lambda>1$.

To conclude, consider the $U_q \mathfrak{k}$-invariants
$$
y_1=\alpha \alpha^*+\beta \beta^*+\gamma \gamma^*+\delta \delta^*,\qquad
y_2=(\alpha\delta-q \beta\gamma)(\alpha\delta-q \beta\gamma)^*.
$$
Prove that $T(y_1)T(y_2)=T(y_2)T(y_1)$, or equivalently, $y_1y_2=y_2y_1$. In
fact, observe that $H$ admits a structure of $U_q \mathfrak{k}$-module
($\xi(fv)=(\xi f)v$, $\xi \in U_q \mathfrak{k}$, $f \in
\mathbb{C}[\mathrm{Mat}_2]_q$) and splits into a sum of pairwise
non-isomorphic simple $U_q \mathfrak{k}$-modules $H=\bigoplus \limits_{k_1
\ge k_2 \ge 0}H^{(k_1,k_2)}$, $H^{(k_1,k_2)}=U_q \mathfrak{k}
\delta^{k_1-k_2}(\alpha\delta-q \beta\gamma)^{k_2}v$. What remains is to
take into account that the restrictions of $T(y_1)$, $T(y_2)$ onto
$H^{(k_1,k_2)}$ are scalar operators by the 'Schur lemma'. Those scalars are
easily deducible:
$$T(y_1)|_{H^{(k_1,k_2)}}=1-q^{2k_1}+q^{-2}(1-q^{2k_2}),$$
$$T(y_2)|_{H^{(k_1,k_2)}}=q^{-2}(1-q^{2k_2})(1-q^{2(k_1+1)}).$$
Just as one could expect, the joint spectrum of the operators $T(y_1)$,
$T(y_2)$ tends to
$$
\{(\mathrm{tr}(\boldsymbol{zz}^*),\det(\boldsymbol{zz}^*)|\:\boldsymbol{z}
\in \mathbb{U}\}=\{(y_1,y_2)|\:0 \le y_1 \le 2 \quad \& \quad 0 \le y_2 \le
y_1^2/4 \}
$$
as $q$ goes to 1.

\bigskip

\section{An invariant integral on the Shilov boundary}

Let $c=\alpha \delta-q \beta \gamma$ and $\mathbb{C}[GL_2]_q$ be the
localization of $\mathbb{C}[\mathrm{Mat}_2]_q$ with respect to the
multiplicative system $c^{\mathbb{N}}$. It is easy to prove the existence
and uniqueness of an extension of the $U_q \mathfrak{g}$-module structure
from $\mathbb{C}[\mathrm{Mat}_2]_q$ onto $\mathbb{C}[GL_2]_q$. Equip the
$U_q \mathfrak{g}$-module algebra $\mathbb{C}[GL_2]_q$ with an involution:
\begin{align*}
\alpha^*&=q^{-2}(\alpha\delta-q \beta\gamma)^{-1}\delta,&
\beta^*&=-q^{-1}(\alpha\delta-q \beta\gamma)^{-1}\gamma,\\
\gamma^*&=-q^{-1}(\alpha\delta-q \beta\gamma)^{-1}\beta, &
\delta^*&=(\alpha\delta-q \beta\gamma)^{-1}\alpha.
\end{align*}
and introduce the notation
$\mathrm{Pol}(S(\mathbb{U}))_q=(\mathbb{C}[GL_2]_q,*)$.

The following propositions justifies our choice of the involution.

\medskip

\begin{proposition} For all $f \in \mathrm{Pol}(S(\mathbb{U}))_q$, $\xi \in
U_q \mathfrak{g}$ one has
$$(\xi f)^*=(S(\xi))^*f^*.$$
\end{proposition}

\medskip

\begin{proposition} There exists a unique homomorphism of $U_q
\mathfrak{g}$-module $*$-algebras $j:~\mathrm{Pol}(\mathrm{Mat}_2)_q \to
\mathrm{Pol}(S(\mathbb{U}))_q$ such that $j(\alpha)=\alpha$,
$j(\beta)=\beta$, $j(\gamma)=\gamma$, $j(\delta)=\delta$.
\end{proposition}

\medskip

These statements are proved in an essentially more general form in
\cite{uni4_rqsu22.texSSV4}. It also follows from the results of that work that the $U_q
\mathfrak{g}$-module $*$-algebra $\mathrm{Pol}(S(\mathbb{U}))_q$ is a
q-analogue of the polynomial algebra on the Shilov boundary $S(\mathbb{U})$
of the unit ball $\mathbb{U}$ in the space $\mathrm{Mat}_2$ of complex $2
\times 2$ matrices.

The $U_q \mathfrak{k}$-module $\mathrm{Pol}(S(\mathbb{U}))_q$ splits into a
sum of pairwise non-isomorphic simple finite dimensional submodules. In
particular, the trivial $U_q \mathfrak{k}$-module appears in
$\mathrm{Pol}(S(\mathbb{U}))_q$ with multiplicity 1 and there exists a
unique $U_q \mathfrak{k}$-invariant integral
$\mu:~\mathrm{Pol}(S(\mathbb{U}))_q \to \mathbb{C}$, $f \mapsto \int
\limits_{S(\mathbb{U})_q}fd \mu$, with $\int \limits_{S(\mathbb{U})_q}1d
\mu=1$.

\medskip

\begin{proposition} The above $U_q \mathfrak{k}$-invariant integral is
positive definite.
\end{proposition}

\smallskip

{\bf Proof.} Consider the $*$-algebra $\mathrm{Pol}(U_2)_q$ of regular
functions on the quantum $U_2$ \cite{uni4_rqsu22.texKoe}, together with the
$*$-homomorphism of algebras $i:\mathrm{Pol}(S(\mathbb{U}))_q \to
\mathrm{Pol}(U_2)_q$ given by
\begin{align*}i(\alpha)&=q^{-1}\alpha,& i(\beta)&=q^{-1}\beta,\\
i(\gamma)&=\gamma,& i(\delta)&=\delta.
\end{align*}

The positivity of an invariant integral on the quantum group $U_2$
constitutes a well known fact. So, what remains is to prove the invariance
of the integral
$$\mathrm{Pol}(U_2)_q \to \mathbb{C},\qquad f \mapsto \int
\limits_{S(\mathbb{U})_q}i^{-1}(f)d \mu$$
 with respect to the action of $U_q \mathfrak{u}_2$ by 'right translations'
on the quantum $U_2$. This is a consequence of the invariance of $\mu$ with
respect to the action of the subalgebra in $U_q \mathfrak{k}$ generated by
$E_1$, $F_1$, $K_1^{\pm 1}$, $(K_1K_2^2K_3)^{\pm 1}$.\hfill $\square$

\medskip

Now introduce an auxiliary $U_q \mathfrak{g}$-module $*$-algebra
$\mathrm{Pol}(\widehat{S}(\mathbb{U}))_q$, to be used in a construction of
the principal degenerate series of quantum
$(\mathfrak{g},\mathfrak{k})$-modules.

 The
$*$-algebra $\mathrm{Pol}(\widehat{S}(\mathbb{U}))_q$ is defined by adding
$t$, $t^{-1}$ to the list $\alpha$, $\beta$, $\gamma$, $\delta$, $c^{-1}$ of
generators of $\mathrm{Pol}(S(\mathbb{U}))_q$ and
$$tt^{-1}=t^{-1}t=1,\qquad tt^*=t^*t,$$
$$zt=qtz,\qquad zt^*=qt^*z,\qquad \text{with}\quad z \in
\{\alpha,\beta,\gamma,\delta \}$$
 to the list of relations.

The next two statements follow from the results of \cite{uni4_rqsu22.texSSV4}.

\medskip

\begin{proposition} i) There exists a unique extension of the structure of
$U_q \mathfrak{g}$-module \hbox{$*$-algebra} from
$\mathrm{Pol}(S(\mathbb{U}))_q$ onto
$\mathrm{Pol}(\widehat{S}(\mathbb{U}))_q$ such that
$$\left \{\begin{array}{ccc}E_jt&=&0 \\ F_jt&=&0 \\ K_jt&=&0
\end{array}\right.,\quad j=1,3, \qquad \left
\{\begin{array}{ccc}E_2t&=&q^{-1/2}t\delta \\ F_2t&=&0 \\ K_2t&=&q^{-1}t
\end{array}\right..$$

ii) There exists a unique homomorphism
$\widehat{j}:\mathrm{Pol}(\mathrm{Pl}_{2,4})_{q,x}\to
\mathrm{Pol}(\widehat{S}(\mathbb{U}))_q$ of $U_q \mathfrak{g}$-module
$*$-algebras such that
$$\widehat{j}(\alpha)=\alpha,\qquad \widehat{j}(\beta)=\beta,\qquad
\widehat{j}(\gamma)=\gamma,\qquad \widehat{j}(\delta)=\delta,\qquad
\widehat{j}(t^{\pm 1})=t^{\pm 1}.$$
\end{proposition}

\medskip

\begin{proposition}\label{uni4_rqsu22.texii0}
The subspace $t^{*-2}\mathrm{Pol}(S(\mathbb{U}))_qt^{-2}$ is a submodule of
the $U_q \mathfrak{g}$-module $\mathrm{Pol}(\widehat{S}(\mathbb{U}))_q$, and
the linear functional
$$t^{*-2}\mathrm{Pol}(S(\mathbb{U}))_qt^{-2}\to \mathbb{C},\qquad
t^{*-2}ft^{-2}\mapsto \int \limits_{S(\mathbb{U})_q}fd \mu$$
 is an invariant integral (i.e. a morphism of $U_q \mathfrak{g}$-modules).
\end{proposition}

\bigskip

\section{An analytic continuation of the holomorphic discrete series: step
two}

Just as in the classical case $q=1$, one has
$$\mathbb{C}[\mathrm{Mat}_2]_q=\bigoplus_{k_1 \ge k_2 \ge
0}\mathbb{C}[\mathrm{Mat}_2]_q^{(k_1,k_2)}=U_q
\mathfrak{k}\delta^{k_1-k_2}(\alpha \delta-q \beta \gamma)^{k_2},$$
 with $\mathbb{C}[\mathrm{Mat}_2]_q^{(k_1,k_2)}$ being simple pairwise
non-isomorphic $U_q \mathfrak{k}$-submodules of the $U_q
\mathfrak{k}$-module $\mathbb{C}[\mathrm{Mat}_2]_q$. Introduce the notation
$f^{(k_1,k_2)}$ for a projection of $f$ onto the $U_q \mathfrak{k}$-isotypic
component $\mathbb{C}[\mathrm{Mat}_2]_q^{(k_1,k_2)}$ parallel to the sum of
all other $U_q \mathfrak{k}$-isotypic components.

By the 'Schur lemma', every $U_q \mathfrak{k}$-invariant Hermitian form
$(f_1,f_2)$ on $\mathbb{C}[\mathrm{Mat}_2]_q$ is given by
$$(f_1,f_2)=\sum_{k_1 \ge k_2 \ge 0}c(k_1,k_2)\int
\limits_{S(\mathbb{U})_q}(f_2^{(k_1,k_2)})^*f_1^{(k_1,k_2)}d \mu.$$

We are going to obtain this decomposition for $(f_1,f_2)_{q^{2 \lambda}}$,
$\lambda>3$. Recall the notation $(a;q^2)_m=\prod
\limits_{j=0}^{m-1}(1-aq^{2j})$.

\medskip

\begin{proposition}\label{uni4_rqsu22.texqii} For all $\lambda>3$, $f_1,f_2 \in
\mathbb{C}[\mathrm{Mat}_2]_q$,
$$\int \limits_{\mathbb{U}_q}f_2^*f_1d \nu_\lambda=\sum_{k_1 \ge k_2
\ge 0}c(k_1,k_2,q^{2 \lambda})\int
\limits_{S(\mathbb{U})_q}(f_2^{(k_1,k_2)})^*f_1^{(k_1,k_2)}d \mu,$$
 with
\begin{equation}\label{uni4_rqsu22.texc}c(k_1,k_2,q^{2 \lambda})=
\frac{(q^4;q^2)_{k_1}(q^2;q^2)_{k_2}}{(q^{2
\lambda};q^2)_{k_1}(q^{2(\lambda-1)};q^2)_{k_2}}.
\end{equation}
\end{proposition}

\smallskip

{\bf Proof.} In the case $q=1$ a similar result was obtained by Faraut and
Koranyi \cite{uni4_rqsu22.texFK} in a very big generality. Our proof here imitates that of
\cite{uni4_rqsu22.texFK}.

First introduce the subalgebra $\mathbb{C}[\overline{\mathrm{Mat}}_2]_q
\subset \mathrm{Pol}(\mathrm{Mat}_2)_q$ generated by $\alpha^*$, $\beta^*$,
$\gamma^*$, $\delta^*$, and the algebra
$\mathbb{C}[\overline{\mathrm{Mat}}_2]_q^{\mathrm{op}}$ which differs from
$\mathbb{C}[\overline{\mathrm{Mat}}_2]_q$ by replacement of the
multiplication law with an opposite one. We use the algebra
$\mathbb{C}[\mathrm{Mat}_2 \times \overline{\mathrm{Mat}}_2]_q
\stackrel{\mathrm{def}}{=}\mathbb{C}[\mathrm{Mat}_2]_q \otimes
\mathbb{C}[\overline{\mathrm{Mat}}_2]_q^{\mathrm{op}}$ as a q-analogue for
the algebra of (degenerate) kernels of integral operators.

Equip $\mathbb{C}[\mathrm{Mat}_2 \times \overline{\mathrm{Mat}}_2]_q$ with a
bigrading
$$
\deg(\alpha \otimes 1)=\deg(\beta \otimes 1)=\deg(\gamma \otimes
1)=\deg(\delta \otimes 1)=(1,0),
$$
$$
\deg(1 \otimes \alpha^*)=\deg(1 \otimes \beta^*)=\deg(1 \otimes
\gamma^*)=\deg(1 \otimes \delta^*)=(0,1)
$$
and the associated topology. The completed algebra
$\mathbb{C}[[\mathrm{Mat}_2 \times \overline{\mathrm{Mat}}_2]]_q$ will work
as the algebra of generalized kernels of integral operators \cite{uni4_rqsu22.texSSV3}.

Just as in the case $q=1$ one can deduce proposition \ref{uni4_rqsu22.texqii} from the
following three lemmas.

\medskip

\begin{lemma}\label{uni4_rqsu22.texP} Given $k_1,k_2 \in \mathbb{Z}$, $k_1 \ge k_2 \ge 0$,
denote by $P_{k_1,k_2}$ the projection in $\mathbb{C}[\mathrm{Mat}_2]_q$
onto the component $\mathbb{C}[\mathrm{Mat}_2]_q^{(k_1,k_2)}$ parallel to
the sum of all other $U_q \mathfrak{k}$-isotypic components. There exists a
unique element $p_{k_1,k_2}\in \mathbb{C}[\mathrm{Mat}_2 \times
\overline{\mathrm{Mat}}_2]_q$ such that
$$P_{k_1,k_2}f({\boldsymbol{z}})=\int
\limits_{S(\mathbb{U})_q}p_{k_1,k_2}({\boldsymbol{z}},{\boldsymbol
\zeta})f({\boldsymbol \zeta})d \mu({\boldsymbol \zeta})$$
 for all $f \in \mathbb{C}[\mathrm{Mat}_2]_q$.
\end{lemma}

\medskip

Introduce the notation $L^2(d \nu_\lambda)_q$, $L_a^2(d \nu_\lambda)_q$ for
completions of vector spaces $\mathrm{Pol}(\mathrm{Mat}_2)_q$,
$\mathbb{C}[\mathrm{Mat}_2]_q$ respectively, with respect to the norm $\|f
\|_{q^{2\lambda}}=\left(\int \limits_{\mathbb{U}_q}f^*fd \nu_\lambda
\right)^{1/2}$. These are well defined for $\lambda>3$, and certainly
$L_a^2(d \nu_\lambda)_q \subset L^2(d \nu_\lambda)_q$.

\medskip

\begin{lemma}\label{uni4_rqsu22.texK} Given $\lambda>3$, denote by $P_\lambda$ the
orthogonal projection in $L^2(d \nu_\lambda)_q$ onto $L_a^2(d
\nu_\lambda)_q$. There exists a unique $K_\lambda \in
\mathbb{C}[[\mathrm{Mat}_2 \times \overline{\mathrm{Mat}}_2]]_q$ such that
$$
P_\lambda f(\boldsymbol{z})=\int
\limits_{\mathbb{U}_q}K_\lambda(\boldsymbol{z},\boldsymbol{\zeta})
f({\boldsymbol \zeta})d \nu_\lambda({\boldsymbol \zeta}).$$
 for all $f \in \mathrm{Pol}(\mathrm{Mat}_2)_q$.
\end{lemma}

\medskip

\begin{lemma}\label{uni4_rqsu22.texKp} In $\mathbb{C}[[\mathrm{Mat}_2 \times
\overline{\mathrm{Mat}}_2]]_q$ one has
$$K_\lambda=\sum_{k_1 \ge k_2 \ge 0}\frac1{c(k_1,k_2,\lambda)}p_{k_1,k_2},$$
 with $c(k_1,k_2,\lambda)$ being given by (\ref{uni4_rqsu22.texc}).
\end{lemma}

\medskip

Lemmas \ref{uni4_rqsu22.texP}, \ref{uni4_rqsu22.texK} can be proved in the same way as in the case $q=1$.
Turn to the proof of lemma \ref{uni4_rqsu22.texKp}.

\smallskip

We are going to use the Schur polynomials
$$s_{k_1k_2}(x_1,x_2)=(x_1x_2)^{k_2}\cdot
\frac{x_1^{k_1-k_2+1}-x_2^{k_1-k_2+1}}{x_1-x_2}.$$
 These are expressible in terms of elementary symmetric polynomials:
$$s_{k_1k_2}(x_1,x_2)=u_{k_1k_2}(x_1+x_2,x_1x_2).$$
 (The polynomials $u_{k_1k_2}$ are closely related to the well known
Chebyshev polynomials of second kind $U_{k_1-k_2}(x)$).

Recall the notation $[j]_q=\dfrac{q^j-q^{-j}}{q-q^{-1}}$,
$(a;q^2)_\infty=\prod \limits_{j=0}^\infty(1-aq^{2j})$ and consider the
kernels $\chi_1=\alpha \otimes \alpha^*+\beta \otimes \beta^*+\gamma \otimes
\gamma^*+\delta \otimes \delta^*$, \ $\chi_2=c \otimes c^*$ with $c=(\alpha
\delta-q \beta \gamma)\in \mathrm{Pol}(\mathrm{Mat}_2)_q$.

\medskip

\begin{lemma}\label{uni4_rqsu22.texPK}i) $p_{k_1,k_2}=q^{k_1+k_2}\cdot[k_1-k_2+1]_q \cdot
u_{k_1k_2}(\chi_1,\chi_2)$,

ii) $K_\lambda=\prod \limits_{j=0}^\infty \left
(1-q^{2(\lambda+j)}\chi_1+q^{4(\lambda+j)}\chi_2 \right) \left(\prod
\limits_{j=0}^\infty \left(1-q^{2j}\chi_1+q^{4j}\chi_2 \right)\right)^{-1}$.
\end{lemma}

\medskip

The first statement of lemma \ref{uni4_rqsu22.texPK} are easily deducible from the
orthogonality relations for matrix elements of representations of the
quantum group $U_2$. The second statement follows from the results of
\cite{uni4_rqsu22.texSSV3}.

Lemma \ref{uni4_rqsu22.texKp} is a consequence of lemma \ref{uni4_rqsu22.texPK} and the following well
known relation in the theory of Schur polynomials \cite{uni4_rqsu22.texM}:
\begin{multline*}\frac{(q^{2 \lambda}x_1;q^2)_\infty}{(x_1;q^2)_\infty}\cdot
\frac{(q^{2 \lambda}x_2;q^2)_\infty}{(x_2;q^2)_\infty}=\\ =\sum_{k_1 \ge k_2
\ge 0}\frac{(q^{2
\lambda};q^2)_{k_1}(q^{2(\lambda-1)};q^2)_{k_2}}{(q^4;q^2)_{k_1}
(q^2;q^2)_{k_2}}[k_1-k_2+1]_q \cdot q^{(k_1+k_2)}s_{k_1k_2}(x_1,x_2).
\end{multline*}

The above proof of proposition \ref{uni4_rqsu22.texqii} is transferable quite literally
onto the case of quantum $SU_{n,n}$ and a q-analogue of the unit ball in the
space of $n \times n$ matrices.

\bigskip

\section{Analytic continuation of the holomorphic discrete series: ladder
representation of the quantum group $\mathbf{SU_{2,2}}$}

It is explained in \cite{uni4_rqsu22.texFK} that the results like our proposition \ref{uni4_rqsu22.texqii}
allow one to solve the problems of irreducibility, unitarizability, and
composition series of the representations $\pi_{q^\lambda}$. We restrict
ourselves to some simplest corollaries from proposition \ref{uni4_rqsu22.texqii}.

\medskip

\begin{proposition}\label{uni4_rqsu22.texPKL}
Suppose that either $\lambda>1$ or $\mathrm{Im}\,\lambda \in \frac{\pi}{\lg
q}\mathbb{Z}$. Then the sesquilinear form $(f_1,f_2)_{q^{2 \lambda}}$ is
positive definite, and for all $f_1,f_2 \in \mathbb{C}[\mathrm{Mat}_2]_q$,
$\xi \in U_q \mathfrak{g}$ one has
\begin{equation}\label{uni4_rqsu22.texpkl}
(\pi_{q^\lambda}(\xi)f_1,f_2)_{q^{2
\lambda}}=(f_1,\pi_{q^\lambda}(\xi^*)f_2)_{q^{2 \lambda}}.
\end{equation}
\end{proposition}

\smallskip

{\bf Proof.} The positivity follows from proposition \ref{uni4_rqsu22.texqii}. Let
$\zeta=q^\lambda$. If $\mathrm{Im}\,\zeta=0$, both sides of (\ref{uni4_rqsu22.texpkl}) are
rational functions of $\zeta$. So, what remains is to use the fact that this
equality is true for $0<\zeta<q^3$. \hfill $\square$

\medskip

Turn to the case $\lambda=1$. It follows from proposition \ref{uni4_rqsu22.texPKL} that the
kernel of the sesquilinear form $\langle f_1,f_2 \rangle=\lim
\limits_{\lambda \to 1+0}(1-q^{2 \lambda-2})(f_1,f_2)_{q^{2 \lambda}}$ is a
common invariant subspace for all the operators $\pi_1(\xi)$, $\xi \in U_q
\mathfrak{g}$. Explicitly, this kernel is
$$L=\bigoplus_{k=0}^\infty \mathbb{C}[\mathrm{Mat}_2]_q^{(k,0)}.$$
On $L$ one has a well defined Hermitian form $(f_1,f_2)=\lim
\limits_{\lambda \to 1+0}(f_1,f_2)_{q^{2 \lambda}}$, and hence the quantum
$(\mathfrak{g},\mathfrak{k})$-module associated to the restriction
$\pi_q|_L$ is unitarizable. The representation $\pi_q|_L$ is a q-analogue of
the well known {\sl ladder representation}.

In the case $q=1$ the subspace $\bigoplus_{k=0}^\infty
\mathbb{C}[\mathrm{Mat}_2]^{(k,0)}$ coincides with the kernel of the
covariant differential operator $\square=\dfrac{\partial}{\partial
\alpha}\dfrac{\partial}{\partial \delta}-\dfrac{\partial}{\partial
\beta}\dfrac{\partial}{\partial \gamma}$. Our intention is to obtain a
q-analogue of this result\footnote{A similar result was obtained by V.
Dobrev \cite{uni4_rqsu22.texDo} and H. P. Jakobsen \cite{uni4_rqsu22.texJak1} in a different context.}.

We use a notion of the first order differential calculus over an algebra $A$
and a covariant first order differential calculus as in \cite{uni4_rqsu22.texKS}.

Among well known $U_q \mathfrak{k}$-invariant first order differential
calculi over $\mathbb{C}[\mathrm{Mat}_2]_q$ one has to distinguish a unique
$U_q \mathfrak{g}$-invariant calculus. A general method of producing such
differential calculi (with hidden symmetry) for q-Cartan domains is
described in \cite{uni4_rqsu22.texSV}.

The first order differential calculus we need is determined by the following
'commutation relations between coordinates and differentials' (these are
written in \cite{uni4_rqsu22.texSSV2} in R-matrix form):
\begin{flalign*}
d \alpha \cdot \alpha&=q^2 \alpha d \alpha;& d \alpha \cdot \beta&=q \beta d
\alpha-(1-q^2)\alpha d \beta;\\ d \alpha \cdot \gamma&=q \gamma d
\alpha-(1-q^2)\alpha d \gamma;& d \alpha \cdot \delta&=\delta d
\alpha-(q^{-1}-q)(\gamma d \beta+\beta d \gamma)+(q^{-1}-q)^2 \alpha d
\delta;\\ d \beta \cdot \alpha&=q \alpha \cdot d \beta;& d \beta \cdot
\beta&=q^2 \beta d \beta;\\ d \beta \cdot \gamma&=\gamma d
\beta-(q^{-1}-q)\alpha d \delta;& d \beta \cdot \delta&=q \delta d
\beta-(1-q^2)\beta d \delta;\\ d \gamma \cdot \alpha&=q \alpha d \gamma;& d
\gamma \cdot \gamma&=q^2 \gamma d \gamma;\\ d \gamma \cdot \beta&=\beta d
\gamma-(q^{-1}-q)\alpha d \delta;& d \gamma \cdot \delta&=q \delta d
\gamma-(1-q^2)\gamma d \delta;\\ d \delta \cdot \alpha&=\alpha d \delta;& d
\delta \cdot \gamma&=q \gamma d \delta;\\ d \delta \cdot \beta&=q \beta d
\delta;& d \delta \cdot \delta&=q^2 \delta d \delta.
\end{flalign*}

It is worthwhile to note that it admits an extension up to a $U_q
\mathfrak{g}$-module first order differential calculus over $U_q
\mathfrak{g}$-module algebra $\mathbb{C}[\mathrm{Pl}_{2,4}]_{q,t}$: $dt
\cdot t=q^{-2}tdt$,
$$
dz \cdot t=q^{-1}tdz, \qquad dt \cdot z=q^{-1}zdt+(q^{-2}-1)tdz \qquad
\text{for all}\quad z \in \{\alpha,\beta,\gamma,\delta \}.
$$

Turn back to $\mathbb{C}[\mathrm{Mat}_2]_q$. The operator $d$ is given on
the generators of this algebra in an obvious way and is extended onto the
entire algebra via the Leibnitz rule. The operators
$\dfrac{\partial}{\partial \alpha}$, $\dfrac{\partial}{\partial \beta}$,
$\dfrac{\partial}{\partial \gamma}$, $\dfrac{\partial}{\partial \delta}$ in
$\mathbb{C}[\mathrm{Mat}_2]_q$ are imposed in a standard way:
$$
df=\frac{\partial f}{\partial \alpha}d \alpha+\frac{\partial f}{\partial
\beta}d \beta+\frac{\partial f}{\partial \gamma}d \gamma+\frac{\partial
f}{\partial \delta}d \delta.
$$

As an easy consequence of the definitions one has

\medskip

\begin{proposition}
Let $\square_q=\dfrac{\partial}{\partial \alpha}\dfrac{\partial}{\partial
\delta}-q \dfrac{\partial}{\partial \beta}\dfrac{\partial}{\partial
\gamma}$.

i) $\square_q$ intertwines $\pi_q$ and $\pi_{q^3}$:
$$\pi_{q^3}(\xi)\square_q=\square_q \pi_q(\xi),\qquad \xi \in U_q
\mathfrak{g}.$$

ii) $L=\mathrm{Ker}\,\square_q$.

iii) $(\alpha \delta-q \beta
\gamma)\square_q|_{\mathbb{C}[\mathrm{Mat}_2]_q^{(k_1,k_2)}}=q^{-2}\cdot
\dfrac{1-q^{2k_2}}{1-q^2}\cdot \dfrac{1-q^{2(k_1+1)}}{1-q^2}$.
\end{proposition}

\medskip

\begin{corollary}
For all $s \in \mathbb{N}$
$$
\square_q(\alpha \delta-q \beta \gamma)^s=b_q(s)(\alpha \delta-q \beta
\gamma)^{s-1},
$$
$$
b_q(s)=q^{-2}\cdot \frac{1-q^{2s}}{1-q^2}\cdot \frac{1-q^{2(s+1)}}{1-q^2}.
$$
\end{corollary}

\medskip

$b_q(s)$ is a q-analogue of the Sato-Bernstein polynomial $b(s)=s(s+1)$ for
the prehomogeneous vector space $\mathrm{Mat}_2$. In a recent preprint
\cite{uni4_rqsu22.texKa} and the works cited therein, another approach to q-analogues for
algebras of differential operators was used to produce q-analogues of the
Bernstein polynomials.

Consider the vector space $\mathbb{C}^4$ (with its standard coordinate
system $t_1$, $t_2$, $t_3$, $t_4$), together with the associated projective
space $\mathbb{CP}^3$. Let $\mathbb{L}\subset \mathbb{CP}^3$ be a
projectivization of the plane $t_3=t_4=0$. It is well known that in the case
$q=1$ the ladder representation is isomorphic to the natural representation
of $U \mathfrak{g}$ in the cohomologies $H^1(\mathbb{CP}^3 \setminus
\mathbb{L},{\cal O}(-2))$. A computation of these cohomologies by the
\v{C}ech method leads to the Laurent polynomials:
$$
H^1(\mathbb{CP}^3 \setminus \mathbb{L},{\cal O}(-2))=\left
\{\sum_{(j_1,j_2,j_3,j_4)\in
J}c_{j_1,j_2,j_3,j_4}t_1^{j_1}t_2^{j_2}t_3^{j_3}t_4^{j_4} \right \},
$$
with $J=\{(j_1,j_2,j_3,j_4)\in \mathbb{Z}^4|\:j_1 \ge 0,\,j_2 \ge
0,\,j_3<0,\,j_4<0,\,j_1+j_2+j_3+j_4=-2 \}$.

So, one has two geometric realizations of the ladder representation of
$SU_{2,2}$ (those in $H^1(\mathbb{CP}^3 \setminus \mathbb{L},{\cal O}(-2))$
and in $\mathrm{Ker}\,\square$).

The lowest weight subspace in $H^1(\mathbb{CP}^3\setminus \mathbb{L},{\cal
O}(-2))$ is generated by the Laurent polynomial $t_3^{-1}t_4^{-1}$, and in
the kernel of $\square=\dfrac{\partial}{\partial
\alpha}\dfrac{\partial}{\partial \delta}-\dfrac{\partial}{\partial
\beta}\dfrac{\partial}{\partial \gamma}$ by the constant function $1$. There
exists a unique isomorphism between the two realizations of the ladder
representation which takes $t_3^{-1}t_4^{-1}$ to $1$. This operator is very
essential in the mathematical physics and is called the Penrose transform
\cite{uni4_rqsu22.texBE}. A replacement of the commutation relation $t_it_j=t_jt_i$ by
$t_it_j=qt_jt_i$, $i<j$, allows one to transfer easily the above
observations onto the case $0<q<1$ (more precisely, everything but the
notion of cohomologies for quasi-coherent sheaves). It is just the way of on
which another realization of the ladder representation and the quantum
Penrose transform appear.

\bigskip

\section{The principal degenerate series of quantum Harish-Chandra modules}

In the classical theory the principal series of Harish-Chandra modules are
associated to parabolic subgroups $P$. Our purpose is to produce a
q-analogue of the principal series of Harish-Chandra modules associated to a
stability group $P$ for a point of the Shilov boundary $p \in
S(\mathbb{U})$.

We call a $U_q \mathfrak{g}$-module $V$ $\mathbb{Z}^3$-weight module if
$V=\bigoplus \limits_{\boldsymbol \mu}V_{\boldsymbol \mu}$ with
${\boldsymbol \mu}=(\mu_1,\mu_2,\mu_3)\in \mathbb{Z}^3$, $V_{\boldsymbol
\mu}=\{v \in V|\:K_j^{\pm 1}v=q^{\pm \mu_j}v,\,j=1,2,3 \}$.

A quantum Harish-Chandra module is a finitely generated
$\mathbb{Z}^3$-weight $U_q \mathfrak{g}$-module $V$ such that

i) $V$ is a sum of finite dimensional simple $U_q \mathfrak{k}$-modules,

ii) each simple finite dimensional $U_q \mathfrak{k}$-module $W$ occurs in
$V$ with finite multiplicity ($\dim \mathrm{Hom}_{U_q
\mathfrak{k}}(W,V)<\infty$).

Quantum Harish-Chandra modules are quantum
$(\mathfrak{g},\mathfrak{k})$-modules, and the notion of unitarizability is
applicable here. The rest of this section is devoted to producing the
principal degenerate series of the unitarizable quantum Harish-Chandra
modules. Note that producing and classification of simple unitarizable
quantum Harish-Chandra modules still constitute an open problem even in our
special case of quantum $SU_{2,2}$.

In the case $\lambda \in -2 \mathbb{Z}_+$ one has a well defined linear
operator $\mathrm{Pol}(S(\mathbb{U}))_q \to
\mathrm{Pol}(\widehat{S}(\mathbb{U}))_q$, $f \mapsto f\cdot(\alpha \delta-q
\beta \gamma)^{-\lambda/2}t^{-\lambda}$. The same argument as that applied
in section 3 to produce $\pi_{q^\lambda}$, yields

\medskip

\begin{proposition}\label{uni4_rqsu22.textl}
There exists a unique one-parameter family $\tau_{q^\lambda}$ of
representations of $U_q \mathfrak{g}$ in the space
$\mathrm{Pol}(S(\mathbb{U}))_q$ of polynomials on the Shilov boundary of the
quantum matrix ball such that

i) for all $\lambda \in -2 \mathbb{Z}_+$, $\xi \in U_q \mathfrak{g}$, $f \in
\mathrm{Pol}(S(\mathbb{U}))_q$ one has
$$(\tau_{q^\lambda}(\xi)f)(\alpha \delta-q \beta
\gamma)^{-\lambda/2}t^{-\lambda}=\xi(f(\alpha \delta-q \beta
\gamma)^{-\lambda/2}t^{-\lambda});$$

ii) for all $\xi \in U_q \mathfrak{g}$, $f \in
\mathrm{Pol}(S(\mathbb{U}))_q$, the vector function
$\tau_{q^\lambda}(\xi)f$is a Laurent polynomial of the indeterminate
$\zeta=q^\lambda$.
\end{proposition}

\medskip

Note that the multiple $(\alpha \delta-q \beta \gamma)^{-\lambda/2}$
provides the integral nature for weight of $\tau_{q^\lambda}$. We are to
produce a q-analogue of the principal degenerate series of Harish-Chandra
modules associated to the Shilov boundary $S(\mathbb{U})$.

\medskip

{\sc remark.} For $q=1$ the construction of degenerate discrete series
involves a finite dimensional irreducible representation $\tau$ of a
reductive subgroup $M$. (This subgroup is determined in a standard way:
$$P=MAN,\qquad S(\mathbb{U})\approx P \setminus SU_{2,2}).$$
We have produced q-analogues of those representations of degenerate discrete
series where $\tau$ is trivial, i.e. a q-analogue of the spherical principal
degenerate series:
$$
\tau_{q^\lambda}(\xi)1=\varepsilon(\xi)1,\qquad \xi \in U_q \mathfrak{k}.
$$

\medskip

Turn to a construction of the corresponding principal unitary series.

\medskip

\begin{proposition}
In the case $\mathrm{Re}\,\lambda=2$ the quantum Harish-Chandra module
associated to $\tau_{q^\lambda}$ is unitarizable:
\begin{equation}\label{uni4_rqsu22.texuni}
\int \limits_{S(\mathbb{U})_q}f_2^*(\tau_{q^\lambda}(\xi)f_1)d \mu=\int
\limits_{S(\mathbb{U})_q}(\tau_{q^\lambda}(\xi^*)f_2)^*f_1d \mu
\end{equation}
 for all $f_1,f_2 \in \mathrm{Pol}(S(\mathbb{U}))_q$, $\xi \in U_q
\mathfrak{g}$.
\end{proposition}

\smallskip

{\bf Proof.} The representation $\tau_{q^\lambda}$ can be defined in a
different way, as one can extend the $U_q \mathfrak{g}$-module algebra
$\mathrm{Pol}(\widehat{S}(\mathbb{U}))_q$ via adding to the list of
generators the powers $t^{\lambda}$, $(t^*)^{\lambda}$, $(\alpha \delta-q
\beta \gamma)^{\lambda}$ for any $\lambda\in\mathbb{C}$. The relations
between the generators of the extended algebra as well as the action of
$E_j$, $F_j$, $K_j^{\pm 1}$, $j=1,2,3$, on them are derived from the
corresponding formulae for integral powers of $t$, $t^*$, and $\alpha
\delta-q \beta \gamma$ via the analytic continuation which uses Laurent
polynomials of the indeterminate $\zeta=q^\lambda$. Moreover, this new
algebra may be endowed with an involution as follows
$$
(t^{\lambda})^*=(t^*)^{\overline{\lambda}}, \qquad ((\alpha \delta-q \beta
\gamma)^{\lambda})^*=q^{-2\overline{\lambda}}\cdot(\alpha \delta-q \beta
\gamma)^{-\overline{\lambda}}
$$
(where bar denotes the complex conjugation), and thus it is made a $U_q
\mathfrak{g}$-module $*$-algebra.

Now the relation (\ref{uni4_rqsu22.texuni}) follows from

\medskip

\begin{lemma}\label{uni4_rqsu22.texii}
Let $\mathrm{Re}\,\lambda=2$. The linear subspace
$$
((\alpha \delta-q \beta \gamma)^{-{\lambda}/2}\cdot t^{-\lambda})^*\cdot
\mathrm{Pol}(S(\mathbb{U}))_q \cdot (\alpha \delta-q \beta
\gamma)^{-\lambda/2}\cdot t^{-\lambda}$$ is a $U_q \mathfrak{g}$-module, and
the linear functional
$$
((\alpha \delta-q \beta \gamma)^{-\lambda/2}\cdot t^{-\lambda})^*\cdot f
\cdot(\alpha \delta-q \beta \gamma)^{-\lambda/2}\cdot t^{-\lambda}\mapsto
\int \limits_{S(\mathbb{U})_q}fd \mu
$$
is a $U_q \mathfrak{g}$-invariant integral.
\end{lemma}

\smallskip

{\bf Proof} \ of lemma \ref{uni4_rqsu22.texii}. Suppose that $\lambda=2+i\rho$ with $\rho
\in \mathbb{R}$. Then, by definitions,
$$
((\alpha \delta-q \beta \gamma)^{-\lambda/2}\cdot t^{-\lambda})^*\cdot f
\cdot(\alpha \delta-q \beta \gamma)^{-\lambda/2}\cdot t^{-\lambda}=
$$
$$
(t^*)^{-\overline{\lambda}}\cdot((\alpha \delta-q \beta
\gamma)^*)^{-\overline{\lambda}/2}\cdot f\cdot(\alpha \delta-q \beta
\gamma)^{-\lambda/2}\cdot t^{-\lambda}=$$
$$
={\rm const}(\rho)\cdot(t^*)^{-2}\cdot(t^*)^{i\rho}\cdot(\alpha \delta-q
\beta \gamma)^{-{i\rho}}\cdot f \cdot t^{-i\rho}\cdot t^{-2}.$$
 Now it suffices to apply proposition \ref{uni4_rqsu22.texii0}, the equality
$$
\int \limits_{S(\mathbb{U})_q}fd \mu=\int
\limits_{S(\mathbb{U})_q}f^{(0,0)}d \mu,
$$
and the observation that the element $t^*t^{-1}(\alpha \delta-q \beta
\gamma)^{-1}\in \mathrm{Pol}(\widehat{S}(\mathbb{U}))_q$ commutes with any
element of the subalgebra $\mathrm{Pol}(S(\mathbb{U}))_q$ and is a $U_q
\mathfrak{g}$-invariant.

\medskip

A construction of the second part $\tau_{q^\lambda}'$ of the principal
degenerate series of quantum Harish-Chandra modules we are interested in is
described in the following proposition. Its proof is just the same as that
of proposition \ref{uni4_rqsu22.textl}.

\medskip

\begin{proposition}
There exists a unique one-parameter family $\tau_{q^\lambda}'$ of
representations of $U_q \mathfrak{g}$ in the space
$\mathrm{Pol}(S(\mathbb{U}))_q$ such that

i) for all $\lambda \in -2 \mathbb{Z}_+$, $\xi \in U_q \mathfrak{g}$, $f \in
\mathrm{Pol}(S(\mathbb{U}))_q$, one has
$$
(\tau_{q^\lambda}'(\xi)f)(\alpha \delta-q \beta
\gamma)^{-\lambda/2}t^{-\lambda-1}=\xi(f(\alpha \delta-q \beta
\gamma)^{-\lambda/2}t^{-\lambda-1})
$$

ii) for all $\xi \in U_q \mathfrak{g}$, $f \in
\mathrm{Pol}(S(\mathbb{U}))_q$, the vector function
$\tau_{q^\lambda}'(\xi)f$ is a Laurent polynomial of the indeterminate
$\zeta=q^\lambda$.
\end{proposition}

\medskip

{\sc Remark.} Both parts $\tau_{q^\lambda}$, $\tau_{q^\lambda}'$ of the
series of quantum Harish-Chandra modules in question could be also derived
via embeddings of vector spaces $\mathrm{Pol}(S(\mathbb{U}))_q \to
\mathrm{Pol}(\widehat{S}(\mathbb{U}))_q$, $f \mapsto ft^{l_1}t^{*l_2}$. For
that, with $l_1-l_2 \in \mathbb{Z}$ being fixed, one should arrange 'an
analytic continuation in $\zeta=q^{l_1+l_2}$'. An equivalence of the two
above approaches to producing the principal degenerate series follows from
properties of the element $t^*t^{-1}(\alpha \delta-q \beta \gamma)^{-1}$
(see proof of lemma \ref{uni4_rqsu22.texii}).

\bigskip

\section{The principal non-degenerate series of quantum Harish-Chandra
modules}

The finite dimensional simple weight $U_q \mathfrak{g}$-modules allow a
plausible description in terms of generators and relations when the highest
weight vectors are chosen as generators. In the infinite dimensional case
the capability of this approach is much lower. The well known method of
inducing from a parabolic subgroup in our case is also inapplicable due to
the absence of a valuable q-analogue of the Iwasawa decomposition.

Fortunately, there exists one more approach to a description of
Harish-Chandra modules, that of Beilinson and Bernstein \cite{uni4_rqsu22.texS}. Within the
framework of this approach simple Harish-Chandra modules are produced in
cohomologies with supports on $K$-orbits in the space of full flags $X=G/B$
(in our case $G=SL_4$, $K=S(GL_2 \times GL_2)$, and $B$ a standard Borel
subgroup). The principal non-degenerate series is related to an open orbit,
which is an affine algebraic variety. This fact sharply simplifies the
problem of producing the principal non-degenerate series, and makes it
possible to solve the problem for the quantum case.

An application of the results of Kostant \cite{uni4_rqsu22.texKo} allows one to obtain an
analogue of proposition \ref{uni4_rqsu22.texii0} for full flags and to prove that, together
with every quantum Harish-Chandra module, the principal non-degenerate
series contains its dual.

\bigskip

\section*{Appendix 1. A complete list of irreducible $*$-representations of
\boldmath $\mathrm{Pol}(\mathrm{Mat}_{2})_q$}

This appendix presents an outline of the results of L. Turowska \cite{uni4_rqsu22.texT} on
classification of irreducible $*$-representations of
$\mathrm{Pol}(\mathrm{Mat}_{2})_q$.

To forestall the exposition, note that every irreducible representation from
the list of L. Turowska possesses a distinguished vector $v$ (determined up
to a scalar multiple) and is a completion of the
$\mathrm{Pol}(\mathrm{Mat}_{2})_q$-module
$V=\mathrm{Pol}(\mathrm{Mat}_{2})_qv$ with respect to a suitable topology.
Our intention is to produce the list of relations which determine the above
$\mathrm{Pol}(\mathrm{Mat}_{2})_q$-modules. As one can observe from the
results of L. Turowska, the non-negative linear functionals
$$
l_q:\mathrm{Pol}(\mathrm{Mat}_{2})_q \to \mathbb{C},\qquad l_q:f \mapsto
(fv,v)
$$
lead in the classical limit $q \to 1$ to non-negative linear functionals on
the polynomial algebra $\mathrm{Pol}(\mathrm{Mat}_2)$. The limit functionals
are just the delta-functions in some points of the closure of the unit ball
$\mathbb{U}$.

We list below those points, together with the lists of determining relations
for the associated
$\mathrm{Pol}(\mathrm{Mat}_2)_q$-modules\footnote{Consider the Poisson
bracket $\{f_1,f_2 \}=\lim \limits_{h \to 0}\frac{f_1f_2-f_2f_1}{ih}$,
$h=2\log(q^{-1})$, and associate to each of those points a bounded
symplectic leaf containing this point. An important invariant of the
irreducible $*$-representation is the dimension of the associated symplectic
leaf.}.

\medskip

\underline{0-dimensional leaves}
\begin{flalign*}\begin{pmatrix}e^{i \varphi_1}&0 \\ 0&e^{i
\varphi_2}\end{pmatrix}&&
\begin{array}{lllll}\alpha v=e^{i\varphi_1}v,& \beta v=0,& \gamma v=0,&
\delta v=e^{i \varphi_2}v,& \\ \alpha^*v=e^{-i \varphi_1}v,& \beta^*v=0,&
\gamma^*v=0,& \delta^*v=e^{-i \varphi_2}v,& \varphi_1,\varphi_2 \in
\mathbb{R}/2 \pi \mathbb{Z}.\end{array}
\end{flalign*}

\medskip

\underline{2-dimensional leaves}
\begin{flalign*}&\begin{pmatrix}0&0 \\ 0&e^{i \varphi}\end{pmatrix}&
\begin{array}{lllll}& \beta v=0,& \gamma v=0,& \delta v=e^{i \varphi}v,& \\
\alpha^*v=0,& \beta^* v=0,& \gamma^*v=0,& \delta^*v=e^{-i \varphi}v,&
\varphi \in \mathbb{R}/2 \pi \mathbb{Z}.\end{array}&
\end{flalign*}

\smallskip

\begin{flalign*}&\begin{pmatrix}0&e^{i \varphi_1} \\ e^{i \varphi_2}&0
\end{pmatrix}&
\begin{array}{llll}\alpha v=0,& \beta v=e^{i \varphi_1}v,& \gamma v=e^{i
\varphi_2}v,& \\ \alpha^*v=-q^{-1}e^{-i(\varphi_1+\varphi_2)}\delta v,&
\beta^* v=e^{-i \varphi_1}v,& \gamma^*v=e^{-i \varphi_2}v,& \delta^*v=0,\\
\varphi_1,\varphi_2 \in \mathbb{R}/2 \pi \mathbb{Z}.&&&\end{array}&
\end{flalign*}

\smallskip

\underline{4-dimensional leaves}
\begin{flalign*}&\begin{pmatrix}0&0 \\ e^{i \varphi}&0\end{pmatrix}&
\begin{array}{lllll}\alpha v=0, & & \gamma v=e^{i \varphi}v,&&
\\ \alpha^*v=0,&
\beta^*v=0,& \gamma^*v=e^{-i \varphi}v,& \delta^*v=0,& \varphi \in
\mathbb{R}/2 \pi \mathbb{Z}.\end{array}&
\end{flalign*}

\smallskip

\begin{flalign*}&\begin{pmatrix}0&e^{i \varphi_1} \\ 0&0\end{pmatrix}&
\begin{array}{lllll}\alpha v=0, & \beta v=e^{i \varphi}v,& & &
\\ \alpha^*v=0,&
\beta^*v=e^{-i \varphi}v,& \gamma^*v=0,& \delta^*v=0,& \varphi \in
\mathbb{R}/2 \pi \mathbb{Z}.\end{array}&
\end{flalign*}

\medskip

\underline{6-dimensional leaves}
\begin{flalign*}&\begin{pmatrix}e^{i \varphi}&0 \\ 0&0\end{pmatrix}&
\begin{array}{lll}\alpha v=e^{i \varphi}v,&& \\ \alpha^*v=e^{-i \varphi}v,&
\beta^*v=\gamma^*v=\delta^*v=0,& \varphi \in \mathbb{R}/2 \pi
\mathbb{Z}.\end{array}&
\end{flalign*}

\medskip

\underline{8-dimensional leaf}
\begin{flalign*}&\begin{pmatrix}0&0 \\ 0&0 \end{pmatrix}&
\alpha^*v=\beta^*v=\gamma^*v=\delta^*v=0.&&
\end{flalign*}

\medskip

It follows from the results of L. Turowska that every of the above
$\mathrm{Pol}(\mathrm{Mat}_{2})_q$-modules can be equipped with a structure
of pre-Hilbert space in such a way that the
$\mathrm{Pol}(\mathrm{Mat}_{2})_q$-action is extendable onto the associated
Hilbert space, and this procedure provides a complete list of irreducible
$*$-representations of $\mathrm{Pol}(\mathrm{Mat}_{2})_q$\footnote{More
precisely, the work by L. Turowska \cite{uni4_rqsu22.texT} presents explicit formulae that
describe the action of the operators $\alpha,\beta,\gamma,\delta$ in the
Hilbert space $l^2(\mathbb{Z}_+)^{\otimes d/2}$, with $d$ being the
dimension of the corresponding symplectic leaf.}.

Note that the $*$-representation associated to the 8-dimensional symplectic
leaf is faithful; it is unique (up to a unitary equivalence) faithful
irreducible $*$-representation. The uniqueness is easily deducible from the
commutation relations between $\alpha$, $\beta$, $\gamma$, $\delta$,
$\alpha^*$, $\beta^*$, $\gamma^*$, $\delta^*$, $y$ (the later element is
defined in section 3).

Another two series of $*$-representations are related to the leaves that
contain unitary matrices
$$\begin{pmatrix}e^{i \varphi_1}&0 \\ 0&e^{i \varphi_2}\end{pmatrix},\qquad
\begin{pmatrix}0&e^{i \varphi_1} \\ e^{i \varphi_2}&0\end{pmatrix},\qquad
\varphi_1,\varphi_2 \in \mathbb{R}/2 \pi \mathbb{Z}.$$
 These two series are due to the $*$-homomorphism
$\mathrm{Pol}(\mathrm{Mat}_2)_q \to \mathbb{C}[U_2]_q$ described in the main
sections of this work. They could be obtained within the theory of
$*$-representations of the algebra $\mathbb{C}[U_2]_q$ of regular functions
on the quantum $U_2$.

\bigskip

\section*{Notes of the Editor}

The results of L. Turowska mentioned in this work were applied recently in
\cite{uni4_rqsu22.texPT}.

A quantum analogue of the Shilov boundary in a more general context was
produced in \cite{uni4_rqsu22.texV}.

The correspondence of the notion of Shilov boundary used in this work and
its well known counterpart by W.~Arveson is discussed in \cite{uni4_rqsu22.texV1}.





\title{\bf NON-COMPACT QUANTUM GROUPS AND HARISH-CHANDRA MODULES}

\author{D.Shklyarov$^\dagger$ \and S.Sinel'shchikov$^\dagger$ \and
A.Stolin$^\ddagger$ \and L.Vaksman$^\dagger$} \date{}

\newpage
\setcounter{section}{0}
\large \thispagestyle{empty} \ \vfill \begin{center}\LARGE \bf PART III \\
QUANTUM HARISH-CHANDRA MODULES ASSOCIATED TO q-CARTAN DOMAINS
\end{center}
\vfill \addcontentsline{toc}{chapter}{Part III \ \ QUANTUM HARISH-CHANDRA
MODULES ASSOCIATED TO q-CARTAN DOMAINS}
\newpage

\makeatletter
\renewcommand{\@oddhead}{NON-COMPACT QUANTUM GROUPS AND HARISH-CHANDRA
MODULES \hfill \thepage}
\renewcommand{\@evenhead}{\thepage \hfill D. Shklyarov, S. Sinel'shchikov,
A. Stolin, and L. Vaksman}
\let\@thefnmark\relax
\@footnotetext{This research was supported in part by Award No UM1-2091 of
the US Civilian Research \& Development Foundation \newline \indent This
lecture has been delivered at the Conference 'Supersymmetry and Quantum
Field Theory'; published in Nucl. Physics B (Proc. Suppl.) 102\& 103 (2001),
334 -- 337.}
\addcontentsline{toc}{chapter}{\@title \\ {\sl D. Shklyarov, S.
Sinel'shchikov, A. Stolin, and L.Vaksman}\dotfill} \makeatother

\maketitle

\centerline{$^\dagger$Institute for Low temperature Physics \& Engineering}
\centerline {47 Lenin Ave., 61103 Kharkov, Ukraine}

\medskip

\centerline{$^\ddagger$Chalmers Tekniska H\"ogskola, Mathematik}
\centerline{412 96, G\"oteborg, Sweden}

\medskip

\begin{quotation}\small {\bf Abstract.} An important problem of the quantum
group theory is to construct and classify the Harish-Chandra modules; it is
discussed in this work. The way of producing the principal non-degenerate
series representations of the quantum $SU_{n,n}$ is sketched. A q-analogue
for the Penrose transform is described.
\end{quotation}

\bigskip

A general theory of non-compact quantum groups which could include, for
instance, quantum $SU_{2,2}$, does not exist. However, during the recent
years, a number of problems on non-commutative geometry and harmonic
analysis on homogeneous spaces of such 'groups' was solved. In these
researches, the absent notion of non-compact quantum group was replaced by
Harish-Chandra modules over quantum universal enveloping algebra $U_q
\mathfrak{g}$. This work approaches an important and still open problem in
the theory of quantum groups, the problem of constructing and classifying
quantum Harish-Chandra modules. A construction of the principal
non-degenerate series of quantum Harish-Chandra modules is described in the
special case of the quantum $SU_{2,2}$. A q-analogue of the Penrose
transform is investigated.

The last named author is grateful to V. Akulov for numerous discussions of
geometric aspects of the quantum group theory.

Everywhere in the sequel $\mathfrak{g}$ stands for a simple complex Lie
algebra and $\{\alpha_1,\alpha_2,\ldots,\alpha_l \}$ for its system of
simple roots with the standard ordering. The field $\mathbb{C}(q)$ of
rational functions of the deformation parameter $q$ normally works as a
ground field (when solving the problems of harmonic analysis, it is more
convenient to assume $q \in(0,1)$ and to set $\mathbb{C}$ as a ground
field).

A background in quantum universal enveloping algebras was made up by V.
Drinfeld and M.~Jimbo in mid-80-ies. The principal results of this theory at
its early years are expounded in the review \cite{uni5nqghcm.texR} and in the lectures
\cite{uni5nqghcm.texJan}. We inherit the notation of these texts; in particular, we use
the standard generators $\{E_j,F_j,K_j^{\pm 1}\}_{j=1,2,\ldots,l}$ of the
Hopf algebra $U_q \mathfrak{g}$ and relatively prime integers $d_j$,
$j=1,2,\ldots,l$, which symmetrize the Cartan matrix of $\mathfrak{g}$ (note
that $d_j=1$, $j=1,2,\ldots,l$, in the case
$\mathfrak{g}=\mathfrak{sl}_{l+1}$). We restrict ourselves to considering
$\mathbb{Z}^l$-{\sl weight} $U_q \mathfrak{g}$-modules $V$ i.e. those
admitting a decomposition into a sum of weight subspaces
$$V=\bigoplus_{\boldsymbol \mu} V_{\boldsymbol \mu},\;
\boldsymbol{\mu}=(\mu_1,\mu_2,\ldots,\mu_l)\in \mathbb{Z}^l,\quad
V_{\boldsymbol \mu}=\{v \in V|\:K_j^{\pm 1}v=q^{\pm
d_j\mu_j}v,\,j=1,2,\ldots,l \}.$$

Recall that some of the simple roots $\alpha_1,\alpha_2,\ldots,\alpha_l$
determine Hermitian symmetric spaces of non-compact type \cite{uni5nqghcm.texH}.(The
coefficients of such simple roots in an expansion of the highest root of
$\mathfrak{g}$ is 1.) For example, in the case
$\mathfrak{g}=\mathfrak{sl}_{l+1}$ all simple roots possess this property.
Choose one such root $\alpha_{j_0}$ and introduce the notation $U_q
\mathfrak{k}$ for the Hopf subalgebra generated by $K_{j_0}^{\pm 1}$, $E_j$,
$F_j$, $K_j^{\pm 1}$, $j \ne j_0$. Of course, every $U_q
\mathfrak{g}$-module $V$ is also a $U_q \mathfrak{k}$-module.

A finitely generated $\mathbb{Z}^l$-weight $U_q \mathfrak{g}$-module $V$ is
called quantum Harish-Chandra module if

\begin{enumerate}

\item $U_q \mathfrak{k}$-module $V$ is a sum of finite dimensional simple
$U_q \mathfrak{k}$-modules,

\item every finite dimensional simple $U_q \mathfrak{k}$-module $W$ occurs
with finite multiplicity \\ ($\mathrm{dim}\,\mathrm{Hom}_{U_q
\mathfrak{k}}(W,V)<\infty$).

\end{enumerate}

In the classical theory there are several methods of construction and
classification of Harish-Chandra modules \cite{uni5nqghcm.texSch}. A similar problem for
quantum Harish-Chandra modules is still open. In our opinion, it is among
the most important of the quantum group theory.

To describe obstacles that appear in solving these problems, consider the
Hermitian symmetric space $SU_{2,2}/S(U_2 \times U_2)$. It is determined by
the Lie algebra $\mathfrak{g}=\mathfrak{sl}_4$ and the specified simple root
$\alpha_{j_0}=\alpha_2$ of this Lie algebra.

Our primary desire is to construct a q-analogue for the principal
non-degenerate series of Harish-Chandra modules. This interest is partially
inspired by Casselman's theorem \cite{uni5nqghcm.texBB} which claims that every classical
simple Harish-Chandra module admits an embedding into a principal
non-degenerate series module. A well known method of producing this series
is just the induction from a parabolic subalgebra. Regretfully, this
subalgebra has no q-analogue (this obstacle does not appear if one
substitutes the subalgebra $U \mathfrak{k}$, thus substituting a subject of
research, cf. \cite{uni5nqghcm.texL}). Fortunately, a quantization is available for
another less known method of producing the principal non-degenerate series.
We describe this method in a simple special case of quantum $SU_{2,2}$.

Let $G=SL_4(\mathbb{C})$, $K=S(GL_2 \times GL_2)$, $B \subset G$ be the
standard Borel subgroup of upper triangular matrices, and $X=G/B$ the
variety of complete flags. It is known that there exists an open $K$-orbit
in $X$. It is well known that this orbit is an affine algebraic variety
\cite{uni5nqghcm.texSch}. The regular functions on this orbit constitute a Harish-Chandra
module of the principal non-degenerate series. The regular differential
forms of the highest degree form another module of this series. The general
case is {\sl essentially} approached by considering a generic homogeneous
algebraic line bundle on $X$ and subsequent restricting it to the open
$K$-orbit. The above construction procedure can be transferred to the
quantum case and leads to the principal non-degenerate series of quantum
Harish-Chandra modules.

In the classical representation theory, the above interplay between
$K$-orbits on the variety and the theory of Harish-Chandra modules
constitutes a generic phenomenon. In the theory of Beilinson-Bernstein,
simple Harish-Chandra modules are derived from the so called standard
Harish-Chandra modules. Furthermore, every series of standard Harish-Chandra
modules is associated to a $K$-orbit in $X=G/B$. However, given an orbit $Q$
of a codimension $s>0$, one should consider the local cohomology
$H_Q^s(X,{\cal F})$ instead of functions on $Q$, with ${\cal F}$ being a
sheaf of sections of a homogeneous algebraic line bundle (here our
description of a standard module is somewhat naive but hopefully more
plausible; the precise construction is expounded in \cite{uni5nqghcm.texSch}).

Our conjecture is that the standard quantum Harish-Chandra modules can be
produced via some q-analogue of local cohomology $H_Q^s(X,{\cal F})$. An
immediate obstacle that appears this way is in a lack of critical background
in non-commutative algebraic geometry.

Probably the case of a closed $K$-orbit $Q \subset X$ is the simplest and
most important one. Note that closed $K$-orbits are related to discrete
series of Harish-Chandra modules, which are of an essential independent
interest \cite{uni5nqghcm.texZ}.

To conclude, consider a very simple example of a closed $K$-orbit in the
space of {\sl incomplete} flags which leads to the well known 'ladder
representation' of the quantum $SU_{2,2}$ and a q-analogue of the Penrose
transform \cite{uni5nqghcm.texBE}.

We start with a purely algebraic description of the corresponding simple
quantum Harish-Chandra module, to be succeeded with its two geometric
realizations. Quantum Penrose transform intertwines these geometric
realizations and is given by an explicit integral formula.

Consider the generalized Verma module $M$ over $U_q \mathfrak{sl}_4$ given
by its single generator $v$ and the relations
$$F_jv=E_jv=(K_j^{\pm 1}-1)v=0,\qquad j=1,3,$$
$$F_2v=0,\qquad K_2^{\pm 1}v=q^{\pm 1}v.$$

Consider its largest $U_q \mathfrak{sl}_4$-submodule $J \varsubsetneqq M$
and the associated quotient module $L=M/J$. Thus we get a simple quantum
Harish-Chandra module which is a q-analogue of the $U_q
\mathfrak{sl}_4$-module related to the ladder representation of $SU_{2,2}$.

The first geometric realization is related to the $K$-action in the variety
$\mathbb{CP}^3$ of lines in $\mathbb{C}^4$. We use the coordinate system
$(u_1,u_2,u_3,u_4)$ in $\mathbb{C}^4$. Let $Q \subset \mathbb{CP}^3$ be the
subvariety of lines which are inside the plane $u_3=u_4=0$. Then, as it was
demonstrated in \cite{uni5nqghcm.texBE}, the local cohomology $H_Q^2(\mathbb{CP}^3,{\cal
O}(-2))$ can be described in terms of Laurent polynomials with complex
coefficients \multlinegap=0mm
\begin{multline*}H_{[\mathbb{L}]}^2(\mathbb{CP}^3,{\cal O}(-2))\simeq
\left\{
\sum_{j_1j_2j_3j_4}c_{j_1j_2j_3j_4}u_1^{j_1}u_2^{j_2}u_3^{j_3}u_4^{j_4}
\right |\\ \left. j_1 \ge 0\quad \&\quad j_2 \ge 0 \quad \& \quad j_3 \le -1
\quad \&\quad j_4 \le -1 \quad \&\quad \sum_{j_i}=-2 \right \}.
\end{multline*}
This space of Laurent polynomials is a $U \mathfrak{sl}_4$-module since it
is a quotient of the $U \mathfrak{sl}_4$-module
$$\left \{\left.
\sum_{j_1j_2j_3j_4}c_{j_1j_2j_3j_4}u_1^{j_1}u_2^{j_2}u_3^{j_3}u_4^{j_4}
\right|\;j_1 \ge 0,\quad j_2 \ge 0 \right \}.$$

A similar geometric realization of the quantum Harish-Chandra module $L$ is
accessible via replacing the classical vector space $\mathbb{C}^4$ with the
quantum one, which means just replacement of the commutation relations
$u_ju_j=u_ju_i$ with $u_iu_j=qu_ju_i$, $1 \le i<j \le 4$.

The second geometric realization is in considering the vector space of
polynomial solutions of the 'wave equation' in the space
$\mathrm{Mat}_{2,2}$ of $2 \times 2$ matrices:
$$\left \{\left.\psi \left(\begin{pmatrix}\alpha & \beta \\ \gamma &
\delta \end{pmatrix}\right)\in
\mathbb{C}[\mathrm{Mat}_{2,2}]\right|\;\frac{\partial^2 \psi}{\partial
\alpha
\partial \delta}-\frac{\partial^2 \psi}{\partial \beta \partial
\gamma}=0 \right \}.$$
 The structure of a $U \mathfrak{sl}_4$-module in the space
$\mathbb{C}[\mathrm{Mat}_{2,2}]$ of polynomials on the space of matrices is
introduced via embedding into the space of rational functions on the
Grassmann variety of two-dimensional subspaces in $\mathbb{C}^4$:
\begin{equation}\label{uni5nqghcm.texpsi}\psi
\left(\begin{pmatrix}\alpha & \beta \\ \gamma & \delta
\end{pmatrix}\right)\mapsto \psi
\begin{pmatrix}t_{\{3,4 \}}^{-1}t_{\{1,3 \}}& t_{\{3,4 \}}^{-1}t_{\{2,3 \}}\\
t_{\{3,4 \}}^{-1}t_{\{1,4 \}}& t_{\{3,4 \}}^{-1}t_{\{2,4
\}}\end{pmatrix}t_{\{3,4 \}}^{-1}.
\end{equation}
 Here $t_{\{i,j \}}=t_{1i}t_{2j}-t_{1j}t_{2i}$, $i<j$, and $t_{ij}$ are the
generators of the algebra $\mathbb{C}[\mathrm{Mat}_{2,4}]$ of polynomials on
$2 \times 4$ matrices. This realization is can be transferred onto the
quantum case, with the ordinary wave equation being replaced by its
q-analogue
$$\frac{\partial^2 \psi}{\partial \alpha \partial \delta}-q \frac{\partial^2
\psi}{\partial \beta \partial \gamma}=0,$$ and ordinary matrices by quantum
matrices (cf. \cite{uni5nqghcm.texD}). It is worthwhile to note that we use $U_q
\mathfrak{sl}_4$-invariant differential calculus on the quantum space of $2
\times 2$ matrices.

The isomorphism of the two geometric realizations of the ladder
representation is unique up to a constant multiple. In the case $q=1$ it is
given by the Penrose transform. It can be defined explicitly by the
following integral formula:
$$f(u_1,u_2,u_3,u_4)\mapsto \int f \left(\sum_{i=1}^2
\zeta_it_{i1},\sum_{i=1}^2 \zeta_it_{i2},\sum_{i=1}^2
\zeta_it_{i3},\sum_{i=1}^2 \zeta_it_{i4}\right)d \nu(\zeta),$$
 where it is implicit that
$$(\zeta_1t_{13}+\zeta_2t_{23})^{-1}=\frac{1}{\zeta_2t_{23}}\cdot
\sum_{k=0}^\infty(-1)^k
\left(\frac{\zeta_1t_{13}}{\zeta_2t_{23}}\right)^k,$$
$$(\zeta_1t_{14}+\zeta_2t_{24})^{-1}=\sum_{k=0}^\infty(-1)^k
\left(\frac{\zeta_2t_{24}}{\zeta_1t_{14}}\right)^k\cdot
\frac{1}{\zeta_1t_{14}},$$
 and the $U \mathfrak{sl}_2$-invariant integral is given by
$$\int \left(\sum c_{i_1i_2}\zeta_1^{i_1}\zeta_2^{i_2}\right)d
\nu \stackrel{\rm def}{=}c_{-1,-1}$$
 (a passage from the Pl\"ucker coordinates $t_{\{i,j \}}$ to polynomials on
$\mathrm{Mat}_{2,2}$ is described by (\ref{uni5nqghcm.texpsi})). For example, in the case
of a lowest weight vector $\dfrac{1}{u_3u_4}$ we have
$$\frac{1}{u_3u_4}\mapsto \frac{1}{t_{13}t_{24}-t_{14}t_{23}}\mapsto 1.$$
To pass from the classical case to the quantum one it suffices to replace
the ordinary space $\mathbb{C}^2$ with its quantum analogue: $\zeta_1
\zeta_2=q\zeta_2 \zeta_1$, the ordinary product in $\sum \limits_ i
\zeta_it_{ij}$ with the tensor product and to order multiples in the above
formulae in a proper way.

We thus get the quantum Penrose transform, which is an isomorphism of the
two geometric realizations of the quantum Harish-Chandra module $L$.

It is well known \cite{uni5nqghcm.texJac} that the quantum Harish-Chandra module $L$ is
unitarizable. The second geometric realization of this module allows one to
find the corresponding scalar product as in \cite{uni5nqghcm.texO} (using an analytic
continuation of the scalar product involved into the definition of the
holomorphic discrete series \cite{uni5nqghcm.texSSV2}).

The precise formulations and complete proofs of the results announced in
this work will be placed to the Eprint Archives
(http://www.arXiv.org/find/math).

\bigskip

\section*{Notes of the Editor}

An appreciable break-through in the theory of non-compact quantum groups was
advanced in the recent work \cite{uni5nqghcm.texKK}.




\title{\bf ON A q-ANALOGUE OF THE PENROSE TRANSFORM}

\author{D.Shklyarov$^\dagger$ \and S.Sinel'shchikov$^\dagger$ \and
A.Stolin$^\ddagger$ \and L.Vaksman$^\dagger$} \date{}

\makeatletter \@addtoreset{equation}{section}\makeatother
\renewcommand{\theequation}{\thesection.\arabic{equation}}

\newpage
\setcounter{section}{0}
\large

\makeatletter
\renewcommand{\@oddhead}{ON A q-ANALOGUE OF THE PENROSE TRANSFORM \hfill
\thepage}
\renewcommand{\@evenhead}{\thepage \hfill D. Shklyarov, S. Sinel'shchikov,
A. Stolin, and L. Vaksman}
\let\@thefnmark\relax \@footnotetext{Ukrainian Journal of Physics, {\bf 47}
(2002), No 3, 288 -- 292}
\addcontentsline{toc}{chapter}{\@title \\ {\sl D. Shklyarov, S.
Sinel'shchikov, A. Stolin, and L.Vaksman}\dotfill} \makeatother

\maketitle

\centerline{$^\dagger$Institute for Low temperature Physics \& Engineering}
\centerline {47 Lenin Ave., 61103 Kharkiv, Ukraine}

\medskip

\centerline{$^\ddagger$Chalmers Tekniska H\"ogskola, Mathematik}
\centerline{412 96, G\"oteborg, Sweden}

\medskip

\begin{center}
e-mail: sinelshchikov@ilt.kharkov.ua, vaksman@ilt.kharkov.ua,
astolin@math.chalmers.se
\end{center}

\makeatletter
\let\@thefnmark\relax
\@footnotetext{This research was partially supported by Award No UM1-2091 of
the Civilian Research \& Development Foundation} \makeatother

\section{Introduction}

In the framework of the theory of quantum groups and their homogeneous
spaces we consider two geometric realizations for the quantum ladder
representation, together with an intertwining linear transformation --- the
quantum Penrose transform.

In section 2 we supply a preliminary material on the classical Penrose
transform and prove (\ref{uni6qapt.texcpt}). The q-analogue of (\ref{uni6qapt.texcpt}) is to be used
in section 3 to produce a quantum Penrose transform.

Our results hint that a great deal of constructions specific for the theory
of quasi-coherent sheaves admit non-commutative analogues. This research is
motivated by a possibility to use the results of non-commutative algebraic
geometry for producing and studying Harish-Chandra modules over quantum
universal enveloping algebras.

There is a plenty of literature on the Penrose transform, quantum groups,
and non-commutative algebraic geometry. We restrict ourselves to mentioning
the monographs \cite{uni6qapt.texBE,uni6qapt.texJan,uni6qapt.texCP}, papers \cite{uni6qapt.texAZ,uni6qapt.texRo}, and the preprint
\cite{uni6qapt.texMT}.

Note that a noncommutative analogues for the Penrose transform and covariant
differential operators are also considered in the preprints \cite{uni6qapt.texKKO,uni6qapt.texZ}
and in the papers \cite{uni6qapt.texJak,uni6qapt.texDo,uni6qapt.texKa} respectively in a completely different
context.

\bigskip

\section{The classical case}

To recall the definition of the Penrose transform, we restrict ourselves to
a simplest substantial example. In this special case, the Penrose transform
intertwines the cohomology of the sheaf $\mathcal{O}(-2)$ on
$$
U'=\{(u_1:u_2:u_3:u_4)\in \mathbb{CP}^3|\,u_3 \ne 0 \quad \rm{or}\quad u_4
\ne 0 \}
$$
and sections of the sheaf $\mathcal{O}(-1)$ on some open affine submanifold
of the Grassmann manifold $\mathrm{Gr}_2(\mathbb{C}^4)\hookrightarrow
\mathbb{CP}^5$. Instead of the Grassmann manifold, we prefer to consider the
Stifel manifold of ordered linear independent pairs of vectors in
$\mathbb{C}^4$. In this context, the $GL_2$-covariant sections on the Stifel
manifold work as sections of the sheaf $\mathcal{O}(-1)$ on
$\mathrm{Gr}_2(\mathbb{C}^4)$.

Associate to each matrix $\mathbf{t}=(t_{ij})_{i=1,2;j=1,2,3,4}\in
\mathrm{Mat}_{2,4}$ the pairs of vectors in $\mathbb{C}^4$:
$$(t_{11},t_{12},t_{13},t_{14}),\qquad (t_{21},t_{22},t_{23},t_{24}).$$
Consider $U''=\{\mathbf{t}\in
\mathrm{Mat}_{2,4}|\,t_{13}t_{24}-t_{14}t_{23}\ne 0\}$. Every point
$\mathbf{u}=(u_1:u_2:u_3:u_4)\in U'$ determines a one-dimensional subspace
$L_\mathbf{u}\subset \mathbb{C}^4$, and every point $\mathbf{t}\in U''$
determines a two-dimensional subspace $L_\mathbf{t}$ generated by the
vectors of the corresponding pair. Let $U=\{(\mathbf{u},\mathbf{t})\in U'
\times U''|\,L_\mathbf{u}\subset L_\mathbf{t}\}$. We thus get a 'double
fibration' $U'\underset{\eta}{\leftarrow}U \underset{\tau}{\to}U''$, which
leads to the Penrose transform. It should be noted that every line
$L_\mathbf{t}$ is of the form
$$
L=\mathbb{C}(\zeta_1,\zeta_2)
\begin{pmatrix}t_{11}&t_{12}&t_{13}&t_{14}
\\ t_{21}&t_{22}&t_{23}&t_{24}\end{pmatrix},\qquad(\zeta_1,\zeta_2)\in
\mathbb{C}^2.
$$
Hence the above double fibration is isomorphic to the double fibration
$$
U'\underset{\pi}{\leftarrow}\mathbb{CP}^1 \times
U''\underset{\mathrm{pr}_2}{\to}U'',
$$
with
\begin{multline*}
\pi:\left((\zeta_1:\zeta_2),\begin{pmatrix}t_{11}&t_{12}&t_{13}&t_{14}
\\ t_{21}&t_{22}&t_{23}&t_{24}\end{pmatrix}\right)\mapsto
\\ ((\zeta_1t_{11}+\zeta_2t_{21}):(\zeta_1t_{12}+\zeta_2t_{22}):
(\zeta_1t_{13}+\zeta_2t_{23}):(\zeta_1t_{14}+\zeta_2t_{24})).
\end{multline*}
We thus get a coordinate description for the double fibration in question;
this coordinate description is going to be implicit in all subsequent
computations. Let us look at the cohomologies.

Consider an open affine cover $U'=U_1 \cup U_2$,
\begin{gather*}
U_1=\{(u_1:u_2:u_3:u_4)\in U'|\,u_3 \ne 0 \},
\\ U_2=\{(u_1:u_2:u_3:u_4)\in U'|\,u_4 \ne 0 \},
\end{gather*}
and compute the \v{C}ech cohomology $\check{H}^1(U',\mathcal{O}(-2))$. Let
$\mathbb{C}[u_1,u_2,u_3^{\pm 1},u_4^{\pm 1}]$ be the Laurent polynomials in
indeterminates $u_3$, $u_4$, with coefficients from $\mathbb{C}[u_1,u_2]$.
Introduce in a similar way $\mathbb{C}[u_1,u_2,u_3^{\pm 1},u_4]$,
$\mathbb{C}[u_1,u_2,u_3,u_4^{\pm 1}]$; of course, these appear to be $U
\mathfrak{sl}_4$-modules.

It follows from the definition of the \v{C}ech complex that there exists a
natural isomorphism of $U \mathfrak{sl}_4$-modules:
\begin{multline*}
\check{H}^1(U',\mathcal{O}(-2))=
\\ \{f \in \mathbb{C}[u_1,u_2,u_3^{\pm 1},u_4^{\pm
1}]\left/\left(\mathbb{C}[u_1,u_2,u_3^{\pm
1},u_4]+\mathbb{C}[u_1,u_2,u_3,u_4^{\pm 1}]\right)\right|
\\ \deg f=-2 \}.
\end{multline*}
Hence the Laurent polynomials
$$
\frac{u_1^{j_1}u_2^{j_2}}{u_3^{j_3}u_4^{j_4}},\qquad j_3 \ge 1 \quad \&\quad
j_4 \ge 1 \quad \&\quad j_3+j_4=j_1+j_2+2,
$$
form a basis of the vector space $H^1(U',\mathcal{O}(-2))$.

Consider the trivial bundle over $U''$ with fiber
$H^1(\mathbb{CP}^1,\mathcal{O}(-2))$. It is known that
$H^1(\mathbb{CP}^1,\mathcal{O}(-2))\simeq \mathbb{C}$, and the isomorphism
is available via choosing an open affine cover
$\mathbb{CP}^1=\{(\zeta_1:\zeta_2)|\,\zeta_1 \ne 0 \} \cup
\{(\zeta_1:\zeta_2)|\,\zeta_2 \ne 0 \}$. Specifically, $\sum
\limits_{j+k=-2}c_{jk}\zeta_1^j \zeta_2^k \mapsto c_{-1,-1}$. In a different
notation $f \mapsto \mathrm{CT}(\zeta_1\zeta_2f)$, with $\mathrm{CT}:\sum
\limits_{j,k}c_{jk}\zeta_1^j \zeta_2^k \mapsto c_{00}$ (the constant term of
a series). Now $\mathscr{P}f$ is defined as a higher direct image of the
cohomology class $\eta^*f$: $\mathscr{P}f=\tau_*^1 \eta^*f$. The linear map
$\tau_*^1$ is called the integration along the fibers of $\tau$. We restrict
ourselves to computing this 'integral' inside the infinitesimal neighborhood
of $\mathbf{t}_0=\begin{pmatrix}0&0&0&1 \\ 0&0&1&0 \end{pmatrix}$ by using
formal series in $t_{11}$, $t_{12}$, $t_{13}$, $t_{14}^{-1}$, $t_{21}$,
$t_{23}^{-1}$, $t_{24}$ with coefficients from $\mathbb{C}[\zeta_1^{\pm
1},\zeta_2^{\pm 1}]$. Of course, $\eta^*:f(\boldsymbol{u})\mapsto
f(\boldsymbol{\zeta t})$. So, in the coordinate description
\begin{equation}\label{uni6qapt.texcpt} \mathscr{P}:f(\boldsymbol{u})\mapsto
CT_{\boldsymbol{\zeta}}(\zeta_1\zeta_2f(\boldsymbol{\zeta t})),\qquad f \in
u_3^{-1}u_4^{-1}\mathbb{C}[u_1,u_2,u_3^{-1},u_4^{-1}], \end{equation} with
$\mathrm{CT}_{\boldsymbol{\zeta}}$ being the constant term in the
indeterminate $\boldsymbol{\zeta}$.

\medskip

{\sc Example.} Compute $\mathscr{P}(1/(u_3u_4))$. One has:
\begin{align*}
\frac1{\zeta_1t_{13}+\zeta_2t_{23}}&=
\frac1{\zeta_2t_{23}}\sum_{i=0}^\infty(-1)^i
\left(\frac{\zeta_1t_{13}}{\zeta_2t_{23}}\right)^i,
\\ \frac1{\zeta_1t_{14}+\zeta_2t_{24}}&=
\frac1{\zeta_1t_{14}}\sum_{j=0}^\infty(-1)^j
\left(\frac{\zeta_2t_{24}}{\zeta_1t_{14}}\right)^j.
\end{align*}
Hence,
\begin{multline*}
\mathscr{P}\left(\frac1{u_3u_4}\right)=
CT_{\boldsymbol{\zeta}}\left(\frac1{t_{23}t_{14}}\sum_{i,j=0}^\infty
(-1)^{i+j}\left(\frac{\zeta_1t_{13}}{\zeta_2t_{23}}\right)^i
\left(\frac{\zeta_2t_{24}}{\zeta_1t_{14}}\right)^j\right)=
\\ =\frac1{t_{23}t_{14}}\sum_{k=0}^\infty
\left(\frac{t_{13}t_{24}}{t_{23}t_{14}}\right)^k=\frac1{t_{23}t_{14}}\cdot
\frac1{1-\frac{t_{13}t_{24}}{t_{23}t_{14}}}=
-\frac1{t_{13}t_{24}-t_{14}t_{23}}.
\end{multline*}

\medskip

{\sc Remark.} It is known that the Penrose transform is an isomorphism
between the two realizations for the 'ladder' representation of
$\mathfrak{sl}_4$: the representation in $H^1(U',\mathcal{O}(-2))$ and the
representation in
$$
\{\psi(z_1^1,z_2^1,z_1^2,z_2^2)(t_{13}t_{24}-t_{14}t_{23})^{-1}\in
H^0(U'',\mathcal{O}(-1))|\,\square \psi=0 \},
$$
where $\square \psi \stackrel{\mathrm{def}}{=}\displaystyle \frac{\partial^2
\psi}{\partial z_1^1 \partial z_2^2}-\frac{\partial^2 \psi}{\partial z_2^1
\partial z_1^2}$, and
\begin{align*}
z_1^1&=\frac{t_{11}t_{23}-t_{13}t_{21}}{t_{13}t_{24}-t_{14}t_{23}},
&z_2^1&=\frac{t_{12}t_{23}-t_{13}t_{22}}{t_{13}t_{24}-t_{14}t_{23}}
\\ z_1^2&=\frac{t_{11}t_{24}-t_{14}t_{21}}{t_{13}t_{24}-t_{14}t_{23}},
&z_2^2&=\frac{t_{12}t_{24}-t_{14}t_{22}}{t_{13}t_{24}-t_{14}t_{23}}.
\end{align*}
The vectors $1/(u_3u_4)$ and $1/(t_{13}t_{24}-t_{14}t_{23})$ are lowest
weight vectors for the above representations of $\mathfrak{sl}_4$. Of
course, $z_1^1$, $z_2^1$, $z_1^2$, $z_2^2$ can be considered as the standard
coordinates on the big cell $t_{13}t_{24}-t_{14}t_{23}\ne 0$ of the
Grassmanian $\mathrm{Gr}_2(\mathbb{C}^4)$.

\bigskip

\section{The quantum case}

In the previous section we produced the formula (\ref{uni6qapt.texcpt}) which can be
treated as a definition of the Penrose transform in the classical case. Now
our intention is to produce a q-analogue of (\ref{uni6qapt.texcpt}). The principal
difference from the constructions of section 2 is in replacement of the
functors $\eta^*$, $\tau_*^1$ of the sheaf theory with the corresponding
morphisms of $U_q \mathfrak{sl}_4$-modules. (Here $U_q \mathfrak{sl}_4$ is a
quantum universal enveloping algebra. It is a Hopf algebra over the ground
field $\mathbb{C}(q)$ and is determined by the generators
$\{E_i,F_i,K_i^{\pm 1}\}_{i=1,2,3}$ and the well known Drinfeld-Jimbo
relations \cite{uni6qapt.texJan}.)

The quantum projective space $\mathbb{CP}_{\mathrm{quant}}^3$ is defined in
terms of a $\mathbb{Z}_+$-graded algebra $\mathbb{C}[u_1,u_2,u_3,u_4]_q$
whose generators $u_1$, $u_2$, $u_3$, $u_4$ are subject to the commutation
relations
$$u_iu_j=qu_ju_i,\qquad i<j.$$
Just as in the classical case, $\deg u_k=1$, $k=1,2,3,4$. The localization
$\mathbb{C}[u_1,u_2,u_3^{\pm 1},u_4^{\pm 1}]_q$ of
$\mathbb{C}[u_1,u_2,u_3,u_4]_q$ with respect to the multiplicative system
$(u_3u_4)^{\mathbb{N}}$ is equipped in a standard way with a structure of
$U_q \mathfrak{sl}_4$-module algebra. The subalgebras
$\mathbb{C}[u_1,u_2,u_3^{\pm 1},u_4]_q$, $\mathbb{C}[u_1,u_2,u_3,u_4^{\pm
1}]_q$ constitute $U_q \mathfrak{sl}_4$-submodules of the $U_q
\mathfrak{sl}_4$-module $\mathbb{C}[u_1,u_2,u_3^{\pm 1},u_4^{\pm 1}]_q$.
Thus we come to
\begin{equation}\label{uni6qapt.texV'}
V'=\mathbb{C}[u_1,u_2,u_3^{\pm 1},u_4^{\pm
1}]_q/\left(\mathbb{C}[u_1,u_2,u_3^{\pm
1},u_4]_q+\mathbb{C}[u_1,u_2,u_3,u_4^{\pm 1}]_q \right)
\end{equation}
as a q-analogue of the $U \mathfrak{sl}_4$-module $H^1(U',\mathcal{O}(-2))$.

We have produced a q-analogue for the first geometric realization of the
'ladder representation'. Turn to a construction of its second geometric
realization.

The algebra $\mathbb{C}[\mathrm{Mat}_{2,4}]_q$ of polynomials on the quantum
matrix space is determined by its generators $\{t_{ij}\}_{i=1,2;j=1,2,3,4}$
and the well known commutation relations
\begin{align*}
&t_{ik}t_{jk}=qt_{jk}t_{ik},\qquad t_{ki}t_{kj}=qt_{kj}t_{ki},&i<j,
\\ &t_{ij}t_{kl}=t_{kl}t_{ij},&i<k \quad \& \quad j>l,
\\ &t_{ij}t_{kl}-t_{kl}t_{ij}=(q-q^{-1})t_{ik}t_{jl},&i<k \quad \& \quad
j<l.
\end{align*}
The element $t=t_{13}t_{24}-qt_{14}t_{23}$ quasi-commutes with all the
generators $t_{ij}$, $i=1,2$, $j=1,2,3,4$. Let
$\mathbb{C}[\mathrm{Mat}_{2,4}]_{q,t}$ be a localization of
$\mathbb{C}[\mathrm{Mat}_{2,4}]_q$ with respect to the multiplicative system
$t^{\mathbb{N}}$ and $U_q \mathfrak{sl}_2$ the quantum universal enveloping
algebra (determined by the generators $E$, $F$, $K^{\pm 1}$ and the
Drinfeld-Jimbo relations).

$\mathbb{C}[\mathrm{Mat}_{2,4}]_{q,t}$ is equipped in a standard way with a
structure of $U_q \mathfrak{sl}_2 \otimes U_q \mathfrak{sl}_4$-module
algebra. In particular, $\mathbb{C}[\mathrm{Mat}_{2,4}]_{q,t}$ is a $U_q
\mathfrak{sl}_4$-module algebra.

Introduce the notation:
\begin{align*}
z_1^1&=t^{-1}(t_{11}t_{23}-qt_{13}t_{21}),
&z_2^1&=t^{-1}(t_{12}t_{23}-qt_{13}t_{22}),
\\ z_1^2&=t^{-1}(t_{11}t_{24}-qt_{14}t_{21}),
&z_2^2&=t^{-1}(t_{12}t_{24}-qt_{14}t_{22}).
\end{align*}
It is well known and easily deducible that
\begin{align*}
&z_k^iz_k^j=qz_k^jz_k^i,\qquad z_i^kz_j^k=qz_j^kz_i^k,&i<j,
\\ &z_j^iz_l^k=z_l^kz_j^i,&i<k \quad \& \quad j>l,
\\ &z_j^iz_l^k-z_l^kz_j^i=(q-q^{-1})z_k^iz_l^j,&i<k \quad \& \quad
j<l.
\end{align*}
It follows that the subalgebra generated by $z_1^1$, $z_2^1$, $z_1^2$,
$z_2^2$ is 'canonically' isomorphic to the algebra
$\mathbb{C}[\mathrm{Mat}_{2,2}]_q$ of 'polynomials on the quantum matrix
space'. It is easy to demonstrate that
$\mathbb{C}[\mathrm{Mat}_{2,2}]_qt^{-1}$ is a $U_q
\mathfrak{sl}_4$-submodule of the $U_q \mathfrak{sl}_4$-module
$\mathbb{C}[\mathrm{Mat}_{2,4}]_{q,t}$.

The simple submodule of the $U_q \mathfrak{sl}_4$-module
$\mathbb{C}[\mathrm{Mat}_{2,2}]_qt^{-1}$ we are interested in is
distinguished via a q-analogue $\square_q$ of the wave operator $\square$:
$$
\square_q=\frac{\partial}{\partial z_1^1}\frac{\partial}{\partial z_2^2}-q
\frac{\partial}{\partial z_2^1}\frac{\partial}{\partial z_1^2}.
$$
Specifically, $V''=\{\psi t^{-1}|\,\square_q \psi=0,\,\psi \in
\mathbb{C}[\mathrm{Mat}_{2,2}]_q \}$. It is worthwhile to note that the
operators $\dfrac{\partial}{\partial z_j^i}$ are defined in terms of a $U_q
\mathfrak{sl}_4$-invariant first order differential calculus in
$\mathbb{C}[\mathrm{Mat}_{2,2}]_q$:
$$df=\sum_{i,j}\frac{\partial f}{\partial z_j^i}dz_j^i.$$
In turn, this first order differential calculus is defined by the following
well known 'commutation' relations:
\begin{align*}
&
\begin{cases}
z_1^1dz_1^1=q^{-2}dz_1^1 \cdot z_1^1
\\ z_1^1dz_2^1=q^{-1}dz_2^1 \cdot z_1^1
\\ z_1^1dz_1^2=q^{-1}dz_1^2 \cdot z_1^1
\\ z_1^1dz_2^2=dz_2^2 \cdot z_1^1
\end{cases}
\\ &
\begin{cases}
z_2^1dz_1^1=q^{-1}dz_1^1 \cdot z_2^1+(q^{-2}-1)dz_2^1 \cdot z_1^1
\\ z_2^1dz_2^1=q^{-2}dz_2^1 \cdot z_2^1
\\ z_2^1dz_1^2=dz_1^2 \cdot z_2^1+(q^{-1}-q)dz_2^2 \cdot z_1^1
\\ z_2^1dz_2^2=q^{-1}dz_2^2 \cdot z_2^1
\end{cases}
\\ &
\begin{cases}
z_1^2dz_1^1=q^{-1}dz_1^1 \cdot z_1^2+(q^{-2}-1)dz_1^2 \cdot z_1^1
\\ z_1^2dz_2^1=dz_2^1 \cdot z_1^2+(q^{-1}-q)dz_2^2 \cdot z_1^1
\\ z_1^2dz_1^2=q^{-2}dz_1^2 \cdot z_1^2
\\ z_1^2dz_2^2=q^{-1}dz_2^2 \cdot z_1^2
\end{cases}
\\ &
\begin{cases}
z_2^2dz_1^1=dz_1^1 \cdot z_2^2+(q^{-1}-q)dz_2^1 \cdot z_1^2+(q^{-1}-q)dz_1^2
\cdot z_2^1+(q^{-1}-q)^2dz_2^2 \cdot z_1^1
\\ z_2^2dz_2^1=q^{-1}dz_2^1 \cdot z_2^2+(q^{-2}-1)dz_2^2 \cdot z_2^1
\\ z_2^2dz_1^2=q^{-1}dz_1^2 \cdot z_2^2+(q^{-2}-1)dz_2^2 \cdot z_1^2
\\ z_2^2dz_2^2=q^{-2}dz_2^2 \cdot z_2^2
\end{cases}
\end{align*}

We thus get the two $U_q \mathfrak{sl}_4$-modules $V'$, $V''$; our intention
is to find an explicit form of the linear map which provides an isomorphism
$\mathscr{P}:V'\to V''$.

We follow the ideas of classical constructions described in section 2 in
considering the quantum projective space $\mathbb{CP}_{\mathrm{quant}}^1$.
More precisely, let us consider a $\mathbb{Z}_+$-graded algebra
$\mathbb{C}[\zeta_1,\zeta_2]_q$:
$$\zeta_1 \zeta_2=q \zeta_2 \zeta_1,\qquad \deg(\zeta_1)=\deg(\zeta_2)=1,$$
together with its localization $\mathbb{C}[\zeta_1^{\pm 1},\zeta_2^{\pm
1}]_q$ with respect to the multiplicative system $(\zeta_1
\zeta_2)^{\mathbb{N}}$. The algebra $\mathbb{C}[\zeta_1^{\pm 1},\zeta_2^{\pm
1}]_q$ is equipped in a standard way with a structure of $U_q
\mathfrak{sl}_2$-module algebra. The following homomorphism of algebras will
work as the operator $f(\boldsymbol{u})\mapsto f(\boldsymbol{\zeta t})$:
\begin{gather*}
\eta^*:\mathbb{C}[u_1,u_2,u_3,u_4]_q \to \mathbb{C}[\zeta_1,\zeta_2]_q
\otimes \mathbb{C}[\mathrm{Mat}_{2,4}]_q,
\\ \eta^*:u_j \mapsto \zeta_1 \otimes t_{1j}+\zeta_2 \otimes t_{2j},\qquad
j=1,2,3,4.
\end{gather*}

To follow the constructions of section 2, we have to invert the elements
$\zeta_1 \otimes t_{13}+\zeta_2 \otimes t_{23}$, $\zeta_1 \otimes
t_{14}+\zeta_2 \otimes t_{24}$ in a suitable localization of
$\mathbb{C}[\zeta_1^{\pm 1},\zeta_2^{\pm 1}]_q \otimes
\mathbb{C}[\mathrm{Mat}_{2,4}]_{q,t}$. It is easy to verify that
\begin{gather*}
(t_{14}t_{23})^2 \cdot \mathbb{C}[\mathrm{Mat}_{2,4}]_{q,t}\subset
\mathbb{C}[\mathrm{Mat}_{2,4}]_{q,t}\cdot(t_{14}t_{23}),
\\ \mathbb{C}[\mathrm{Mat}_{2,4}]_{q,t}\cdot(t_{14}t_{23})^2 \subset
(t_{14}t_{23})\cdot \mathbb{C}[\mathrm{Mat}_{2,4}]_{q,t}.
\end{gather*}

Thus we have a well defined localization of
$\mathbb{C}[\mathrm{Mat}_{2,4}]_{q,t}$ with respect to the multiplicative
system $(t_{14}t_{23})^\mathbb{N}$. In an appropriate completion of this
algebra one has the following relations:
\begin{gather*}
(\zeta_1 \otimes t_{13}+\zeta_2 \otimes t_{23})^{-1}=(\zeta_2 \otimes
t_{23})^{-1}\sum_{i=0}^\infty(-1)^i \left(\zeta_1 \zeta_2^{-1}\right)^i
\otimes \left(t_{13}t_{23}^{-1}\right)^i,
\\ (\zeta_1 \otimes t_{14}+\zeta_2 \otimes
t_{24})^{-1}=\left(\sum_{j=0}^\infty(-1)^j \left(\zeta_1^{-1}\zeta_2
\right)^j \otimes \left(t_{14}^{-1}t_{24}\right)^j \right)(\zeta_1 \otimes
t_{14})^{-1}.
\end{gather*}

We define the quantum Penrose transform by
$$
\mathscr{P}_qf=(\mathrm{CT}\otimes \mathrm{id})(\zeta_1 \zeta_2 \otimes
1)(\eta^*f),
$$
where, just as above, $\mathrm{CT}:\sum \limits_{ij}c_{ij}\zeta_1^i
\zeta_2^j \mapsto c_{-1,-1}$, and $f$ belongs to the linear span of the
elements
\begin{equation}\label{uni6qapt.texmonom}
u_1^{j_1}u_2^{j_2}u_3^{-j_3}u_4^{-j_4},\qquad j_3 \ge 1 \quad \& \quad j_4
\ge 1 \quad \& \quad j_1+j_2-j_3-j_4=-2.
\end{equation}
Now (\ref{uni6qapt.texV'}) determines a $U_q \mathfrak{sl}_4$-module structure in this
linear span since the monomials (\ref{uni6qapt.texmonom}) form a basis in the vector
space $V'$.

\bigskip

\section{Appendix}

We sketch here the proof of the fact that $\mathscr{P}_q$ is an isomorphism
of $U_q \mathfrak{sl}_4$-modules $V'\widetilde{\longrightarrow}V''$.

It follows from the definition that $\mathscr{P}_q$ is a morphism of $U_q
\mathfrak{sl}_4$-modules. In view of the simplicity of $V'$ and $V''$, it
suffices to prove that $\mathscr{P}_q$ takes the (lowest weight) vector
$u_3^{-1}u_4^{-1}\in V'$ to the (lowest weight) vector
$-(t_{13}t_{24}-qt_{14}t_{23})^{-1}$. We start with an auxiliary statement:
\begin{equation}\label{uni6qapt.texaux}
\left(1-\left(t_{23}^{-1}t_{13}\right)
\left(t_{14}^{-1}t_{24}\right)\right)^{-1}=\sum_{k=0}^\infty
q^{-2k}\left(t_{23}^{-1}t_{13}\right)^k \left(t_{14}^{-1}t_{24}\right)^k.
\end{equation}
It follows from the commutation relation
$$
\left(t_{14}^{-1}t_{24}\right)\left(t_{23}^{-1}t_{13}\right)=
q^{-2}\left(t_{23}^{-1}t_{13}\right)\left(t_{14}^{-1}t_{24}\right)+1-q^{-2}
$$
and the relation (6.5) of \cite{uni6qapt.texSSV}.

An application of (\ref{uni6qapt.texaux}) allows one to prove that
$$
\mathscr{P}_q
\left(u_3^{-1}u_4^{-1}\right)=-(t_{13}t_{24}-qt_{14}t_{23})^{-1}.
$$
In fact,
\begin{multline*}
\mathscr{P}_q \left(u_3^{-1}u_4^{-1}\right)=\mathrm{CT}\otimes
\mathrm{id}\Bigg(\zeta_2^{-1} \otimes t_{23}^{-1}\cdot
\\ \left.\left(\sum_{i=0}^\infty(-1)^i(\zeta_1 \zeta_2^{-1})^i \otimes
q^{-i}\left(t_{23}^{-1}t_{13}\right)^i
\right)\left(\sum_{j=0}^\infty(-1)^j(\zeta_1 \zeta_2^{-1})^j \otimes
\left(t_{14}^{-1}t_{24}\right)^j \right)\zeta_1^{-1} \otimes
t_{14}^{-1}\right).
\end{multline*}
On the other hand,
$$
\zeta_2^{-1}\left(\zeta_1 \zeta_2^{-1}\right)^k \left(\zeta_1^{-1} \zeta_2
\right)^k
\zeta_1^{-1}=q^{-k}\zeta_2^{-1}\zeta_1^{-1}=q^{-k-1}\zeta_1^{-1}\zeta_2^{-1}.
$$
Hence,
\begin{multline*}
\mathscr{P}_q
\left(u_3^{-1}u_4^{-1}\right)=q^{-1}t_{23}^{-1}\left(\sum_{k=0}^\infty
q^{-2k}\left(t_{14}^{-1}t_{24}\right)^k
\left(t_{23}^{-1}t_{13}\right)^k\right)t_{14}^{-1}=
\\ =q^{-1}\left(t_{14}\left(1-
t_{23}^{-1}t_{13}t_{14}^{-1}t_{24}\right)t_{23}\right)^{-1}=
-(t_{13}t_{24}-qt_{14}t_{23})^{-1},
\end{multline*}
which completes the proof.




\makeatletter
\@addtoreset{equation}{section}
\makeatother

\renewcommand{\theequation}{\thesection.\arabic{equation}}

\title{\bf SPHERICAL PRINCIPAL NON-DEGENERATE SERIES OF REPRESENTATIONS FOR
THE QUANTUM GROUP $\mathbf{SU_{2,2}}$}
\author{S. Sinel'shchikov$^\dagger$ \and A. Stolin$^\ddagger$ \and L.
Vaksman$^\dagger$}
\date{}

\newpage
\setcounter{section}{0}
\large

\makeatletter
\renewcommand{\@oddhead}{SPHERICAL PRINCIPAL NON-DEGENERATE SERIES FOR
THE QUANTUM $\mathbf{SU_{2,2}}$ \hfill \thepage}
\renewcommand{\@evenhead}{\thepage \hfill S. Sinel'shchikov, A. Stolin,
and L. Vaksman}

\let\@thefnmark\relax \@footnotetext{This research was supported in part by
of the US Civilian Research \& Development Foundation (Award No UM1-2091)
and by Swedish Academy of Sciences (project No 11293562).

This lecture has been delivered at the 10th International Colloquium
'Quantum Groups and Integrable Systems' in Prague, June 2001. The text is
published in Czechoslovak Journal of Physics {\bf 51} (2001), No 12, 1431
-- 1440.}
\addcontentsline{toc}{chapter}{\@title \\ {\sl S. Sinel'shchikov, A. Stolin,
and L.Vaksman}\dotfill}
\makeatother

\maketitle

\centerline{$^\dagger$Institute for Low temperature Physics \& Engineering}
\centerline {47 Lenin Ave., 61103 Kharkov, Ukraine}

\medskip

\centerline{$^\ddagger$Chalmers Tekniska H\"ogskola, Mathematik}
\centerline{412 96, G\"oteborg, Sweden}

\section{Introduction}

The first step in studying quantum bounded symmetric domains was done in
\cite{uni7spnssu2.texSV}. Later on, an explicit form of the Plancherel measure was found in
\cite{uni7spnssu2.texSSV1} for the simplest among the above domains, the quantum disc. It
is still open to extend this result onto an arbitrary quantum bounded
symmetric domain. Presumably, the initial step in this direction should be
in producing and studying q-analogues of the principal non-degenerate series
representations for automorphism groups of bounded symmetric domains. (It is
certainly implicit in this setting that the quantum universal enveloping
Drinfeld-Jimbo algebras work as quantum groups and the Harish-Chandra
modules over those algebras work as representations of quantum groups
\cite{uni7spnssu2.texSSSV}.)

The present work deals with the simplest bounded symmetric domain of rank 2
$$
\mathscr{D}=\{\mathbf{z}\in \mathrm{Mat}_2(\mathbb{C})|\:\mathbf{zz^*}<I\},
$$
its q-analogue (the quantum matrix ball), together with the associated
Harish-Chandra modules. In section 2 we present a construction of the
principal non-degenerate series of such modules. In section 3 we prove that,
together with every Harish-Chandra module this series contains the dual
Harish-Chandra module.

The third named author is grateful to A. Rosenberg and Ya. Soibelman for a
discussion of the results expounded in section 2 of this paper.

\bigskip

\section{A construction of the spherical principal non-degenerate series}

From now on we assume $q$ to be transcendental and $\mathbb{C}$ to be the
ground field.

Consider the quantum universal enveloping algebra $U_q \mathfrak{sl}_4$. It
is determined by its generators $\{E_j,F_j,K_j^{\pm 1}\}_{j=1,2,3}$ and the
well known Drinfeld-Jimbo relations (see \cite{uni7spnssu2.texJa}). Recall that $U_q
\mathfrak{sl}_4$ is a Hopf algebra. (The comultiplication $\triangle:U_q
\mathfrak{sl}_4 \to U_q \mathfrak{sl}_4 \otimes U_q \mathfrak{sl}_4$ is
defined by $\triangle(K_j^{\pm 1})=K_j^{\pm 1}\otimes K_j^{\pm 1}$,
$\triangle(E_j)=E_j \otimes 1+K_j \otimes E_j$, $\triangle(F_j)=F_j \otimes
K_j^{-1}+1 \otimes F_j$, $j=1,2,3$, the counit $\varepsilon:U_q
\mathfrak{sl}_4 \to \mathbb{C}$ by
$\varepsilon(E_j)=\varepsilon(F_j)=\varepsilon(K_j^{\pm 1}-1)=0$, $j=1,2,3$,
and the antipode $S:U_q \mathfrak{sl}_4 \to U_q \mathfrak{sl}_4$ by
$S(K_j^{\pm 1})=K_j^{\mp 1}$, $S(E_j)=-K_j^{-1}E_j$, $S(F_j)=-F_jK_j$,
$j=1,2,3$). An ordinary universal enveloping algebra $U \mathfrak{sl}_4$ is
derivable from $U_q \mathfrak{sl}_4$ via substituting $q=e^{-h/2}$,
$K_j^{\pm 1}=e^{\mp hH_j/2}$ and a formal passage to the limit $h \to 0$.

Introduce the notation $U_q \mathfrak{g}$ for the Hopf algebra $U_q
\mathfrak{sl}_4$ and $U_q \mathfrak{k}$ for its Hopf subalgebra generated by
$K_2^{\pm 1}$, $E_j$, $F_j$, $K_j^{\pm 1}$, $j=1,3$. A $U_q
\mathfrak{g}$-module $V$ is said to be $\mathbb{Z}^3$-weight module if
$V=\bigoplus \limits_{\boldsymbol \mu}V_{\boldsymbol \mu}$,
$\boldsymbol{\mu}=(\mu_1,\mu_2,\mu_3)\in \mathbb{Z}^3$, $V_{\boldsymbol
\mu}=\{v \in V|\,K_j^{\pm 1}v=q^{\pm \mu_j}v,\;j=1,2,3 \}$.

In what follows we are going to consider only $U_q \mathfrak{g}$-modules $V$
of the above form, which allows one to introduce the linear operators $H_1$,
$H_2$, $H_3$, in $V$ with $K_i^{\pm 1}=q^{\pm H_i}$, $i=1,2,3$.
Specifically, $H_i|_{V_{\boldsymbol{\mu}}}=\mu_iI$, $i=1,2,3$.

\medskip

{\sc Remark.} For every simple Lie algebra and every lattice $L$ which is
between the lattices of roots and weights, a class of $L$-weight $U_q
\mathfrak{g}$-modules can be introduced in the same way.

\medskip

A $\mathbb{Z}^3$-weight $U_q \mathfrak{g}$-module $V$ is called a
Harish-Chandra module if

\begin{enumerate}

\item $U_q \mathfrak{k}$-module $V$ is a sum of finite dimensional simple
$U_q \mathfrak{k}$-modules,

\item every finite dimensional simple $U_q \mathfrak{k}$-module $W$ occurs
with finite multiplicity \\ ($\mathrm{dim}\,\mathrm{Hom}_{U_q
\mathfrak{k}}(W,V)<\infty$)\footnote{The class of $U_q \mathfrak{g}$-modules
we are interested in differs essentially from that considered by G. Letzter
\cite{uni7spnssu2.texL1}. Both classes are accessible via the general approach of
\cite[Definition 2.2]{uni7spnssu2.texL}.}.

\end{enumerate}

The initial step in studying such $U_q \mathfrak{g}$-modules was done by the
authors in \cite{uni7spnssu2.texSSSV}. To produce the principal non-degenerate series of
Harish-Chandra modules, it was suggested to use the approach of
Beilinson-Bernstein \cite{uni7spnssu2.texBB}. Within this approach, in the classical case
($q=1$) the principal non-degenerate series admits a geometric realization
on the open $S(GL_2 \times GL_2)$-orbit $U$ in the projective variety $X$ of
flags in $\mathbb{C}^4$.

A passage to the quantum case should be started with producing a q-analogue
for the open orbit $U$. We need a well known q-analogue of the the affine
algebraic variety $\widehat{X}$ associated to $X$.

We are about to introduce a q-analogue $\mathbb{C}[\widehat{X}]_q$ of the
algebra $\mathbb{C}[\widehat{X}]$ of regular functions on $\widehat{X}$.
Start with the algebra $\mathbb{C}[SL_4]_q$ of 'regular functions' on the
quantum group $SL_4$. Its description in terms of generators $t_{ij}$,
$i,j=1,2,3$, and relations is well known \cite{uni7spnssu2.texKl-Sch}. It admits a
canonical embedding $\mathbb{C}[SL_4]_q \hookrightarrow(U_q \mathfrak{g})^*$
which sends $t_{ij}$ into matrix elements of the vector representation of
$U_q \mathfrak{g}$ in the standard basis \cite{uni7spnssu2.texD}. Thus $\mathbb{C}[SL_4]_q$
is equipped with a structure of $(U_q \mathfrak{g})^\mathrm{op}\otimes U_q
\mathfrak{g}$-module algebra\footnote{For this and some other notions
widespread in the quantum group theory the reader is referred to
\cite{uni7spnssu2.texKl-Sch}. The Hopf algebra $(U_q \mathfrak{g})^\mathrm{op}$ is
derivable from $(U_q \mathfrak{g})$ via a replacement of its multiplication
law with the opposite one.}.

Assume $J \subset \{1,2,3,4 \}$ and $\mathrm{card}(J)\le 3$. Introduce the
notation $t_J$ for the quantum minor of the matrix $(t_{ij})$ formed by the
initial $\mathrm{card}(J)$ lines and the columns with indices from $J$. For
example, $t_{\{1,2 \}}=t_{11}t_{22}-qt_{12}t_{21}$, $t_{\{3,4
\}}=t_{13}t_{24}-qt_{14}t_{23}$.

A complete list of relations between $t_J$, $\mathrm{card}(J)\le 3$,
includes commutation relations and the generalized Pl\"ucker relations
\cite{uni7spnssu2.texTT,uni7spnssu2.texNYM}. $\mathbb{C}[\widehat{X}]_q$ is defined as a unital
subalgebra of $\mathbb{C}[SL_4]_q$ generated by $t_J$, $\mathrm{card}(J)\le
3$. As a subalgebra of $\mathbb{C}[SL_4]_q$, $\mathbb{C}[\widehat{X}]_q$ has
no divisors of zero. It inherits a structure of $U_q \mathfrak{g}$-module
algebra from $\mathbb{C}[SL_4]_q$. Furthermore, it admits a structure of
$\mathbb{Z}^3$-graded algebra :
$$
\deg t_J=
\begin{cases}
(1,0,0), & \mathrm{card}(J)=1
\\ (0,1,0), & \mathrm{card}(J)=2
\\ (0,0,1), & \mathrm{card}(J)=3.
\end{cases}
$$
We use the $(U_q \mathfrak{g})^\mathrm{op}$-module structure in
$\mathbb{C}[SL_4]_q^\mathrm{op}$ in the above definition.

Turn to the classical case $q=1$. Let $B,N$ stand for the standard Borel and
the unipotent subgroups of $SL_4$, respectively. One has $X \simeq B
\backslash SL_4$, $\widehat{X}\setminus \{0 \}\simeq N \backslash SL_4$,
which leads to a canonical morphism $\pi:\widehat{X}\setminus \{0 \}\to X$.
An inverse image of the open orbit $U \subset X$ with respect to $\pi$ is
the open set
$$
t_{\{1,2 \}}t_{\{3,4 \}}\left(t_{\{1 \}}t_{\{2,3,4 \}}-t_{\{2 \}}t_{\{1,3,4
\}}\right)\ne 0.
$$
Thus we obtain an embedding of the algebra of regular functions
$\mathbb{C}[U]$ into the localization of $\mathbb{C}[\widehat{X}]$ with
respect to the multiplicative subset $\left(t_{\{1,2 \}}t_{\{3,4
\}}\left(t_{\{1 \}}t_{\{2,3,4 \}}-t_{\{2 \}}t_{\{1,3,4
\}}\right)\right)^\mathbb{N}$. Of course, the above $\mathbb{Z}^3$-grading
admits a unique extension onto this localization, and the image
$\mathbb{C}[U]$ coincides with the subalgebra of all elements of degree
zero. This description of $\mathbb{C}[U]$ is to be used to produce its
q-analogue. The $U_q \mathfrak{sl}_2 \otimes U_q \mathfrak{sl}_2$-invariant
element
$$\eta=t_{\{1 \}}t_{\{2,3,4 \}}-qt_{\{2 \}}t_{\{1,3,4 \}}$$
of $\mathbb{C}[\widehat{X}]_q$ will work as the $SL_2 \times SL_2$-invariant
element $t_{\{1 \}}t_{\{2,3,4 \}}-t_{\{2 \}}t_{\{1,3,4 \}}\in
\mathbb{C}[\widehat{X}]$.

\medskip

\begin{proposition}
The multiplicative subset $\left(t_{\{1,2 \}}t_{\{3,4 \}}\eta
\right)^\mathbb{N}\subset \mathbb{C}[\widehat{X}]_q$ satisfies the Ore
condition.
\end{proposition}

\smallskip

{\bf Proof.} It is easy to deduce commutation relations between $t_{\{1,2
\}}$, $t_{\{3,4 \}}$, $\eta$, and the generators $t_J$ of
$\mathbb{C}[\widehat{X}]_q$.

Describe some of those. The general properties of the universal R-matrix and
the fact that all $t_J$ are matrix elements of simple $U_q
\mathfrak{g}$-modules, imply the following

\medskip

\begin{lemma}
The elements $t_{\{1,2 \}}$, $t_{\{3,4 \}}$ quasi-commute with all the
generators $t_J$ of $\mathbb{C}[\widehat{X}]_q$.
\end{lemma}

\medskip

All the generators $t_J$ are derivable from the elements
$$
t_{\{2 \}},\quad t_{\{4 \}},\quad t_{\{1,2 \}},\quad t_{\{3,4 \}},\quad
t_{\{2,4 \}},\quad t_{\{1,2,4 \}},\quad t_{\{2,3,4 \}},
$$
via an application of $E_1$, $E_3$, $E_1E_3$. Hence, in view of the
relations
$$E_j \eta=0,\qquad K_j^{\pm 1}\eta=\eta,\qquad j \ne 2,$$
it suffices to obtain commutation relations between $\eta$ and the
generators listed above. An application of the well known properties of the
quantum determinants allows one to deduce the following relations:

\medskip

\begin{lemma}\label{uni7spnssu2.texrs}
\begin{align*}
&\eta t_{\{2 \}}=q^{-1}t_{\{2 \}}\eta,&&\eta t_{\{4 \}}=qt_{\{4 \}}\eta,
\\ &\eta t_{\{1,2 \}}=q^{-1}t_{\{1,2 \}}\eta,&&\eta t_{\{3,4
\}}=qt_{\{3,4 \}}\eta,
\\ &\eta t_{\{2,4 \}}=qt_{\{2,4 \}}\eta-(1-q^2)t_{\{2 \}}t_{\{1,2,4
\}}t_{\{3,4 \}},
\\ &\eta t_{\{1,2,4 \}}=q^{-1}t_{\{1,2,4 \}}\eta,&&\eta t_{\{2,3,4
\}}=qt_{\{2,3,4 \}}\eta.
\end{align*}
\end{lemma}

\medskip

It follows from the above commutation relations that for all the generators
$t_J$ of $\mathbb{C}[\widehat{X}]_q$ one has
\begin{gather*}
\left(t_{\{1,2 \}}t_{\{3,4 \}}\eta \right)^2t_J \in
\mathbb{C}[\widehat{X}]_qt_{\{1,2 \}}t_{\{3,4 \}}\eta
\\ t_J \left(t_{\{1,2 \}}t_{\{3,4 \}}\eta \right)^2 \in t_{\{1,2 \}}t_{\{3,4
\}}\eta \mathbb{C}[\widehat{X}]_q,
\end{gather*}
that is, the Ore condition is satisfied. \hfill $\square$

\medskip

\begin{corollary}
For all $l \in \mathbb{N}$
$$
\eta^lt_{\{2,4 \}}=\left(qt_{\{2,4 \}}\eta-q^2(q^{-l}-q^l)t_{\{2
\}}t_{\{1,2,4 \}}t_{\{3,4 \}}\right)\eta^{l-1}.
$$
\end{corollary}

\medskip

Let $\widehat{U}=\left.\left \{x \in \widehat{X}\right|\:t_{\{1,2
\}}t_{\{3,4 \}}\left(t_{\{1 \}}t_{\{2,3,4 \}}-t_{\{2 \}}t_{\{1,3,4
\}}\right)\ne 0 \right \}$, and denote by $\mathbb{C}[\widehat{U}]_q$ the
localization of $\mathbb{C}[\widehat{X}]_q$ with respect to the
multiplicative subset $\left(t_{\{1,2 \}}t_{\{3,4 \}}\eta \right
)^\mathbb{N}$. Obviously, $\mathbb{C}[\widehat{X}]_q \hookrightarrow
\mathbb{C}[\widehat{U}]_q$, and $\mathbb{Z}^3$-grading is uniquely
extendable from $\mathbb{C}[\widehat{X}]_q$ onto
$\mathbb{C}[\widehat{U}]_q$.

\medskip

\begin{proposition}\label{uni7spnssu2.texext}
There exists a unique extension of the structure of $U_q
\mathfrak{g}$-module algebra from $\mathbb{C}[\widehat{X}]_q$ onto
$\mathbb{C}[\widehat{U}]_q$.
\end{proposition}

\smallskip

{\bf Proof.} The uniqueness of the extension is obvious. We are going to
construct such extension by applying the following statement.

\medskip

\begin{lemma} For every $\xi \in U_q \mathfrak{g}$, $f \in
\mathbb{C}[\widehat{X}]_q$, there exists a unique Laurent polynomial
$p_{\xi,f}(\lambda)$ with coefficients from $\mathbb{C}[\widehat{U}]_q$ such
that
$$
p_{\xi,f}(q^l)=\xi \left(f \cdot \left(t_{\{1,2 \}}t_{\{3,4 \}}\eta
\right)^l\right)\left(t_{\{1,2 \}}t_{\{3,4 \}}\eta \right)^{-l},\qquad l \in
\mathbb{Z}_+.
$$
\end{lemma}

\medskip

We can use the same Laurent polynomials for defining $\xi \left(f \cdot
\left(t_{\{1,2 \}}t_{\{3,4 \}}\eta \right)^l\right)$ for all integers $l$.
Of course, we need to verify that the action of $\xi \in U_q \mathfrak{g}$
in $\mathbb{C}[\widehat{U}]_q$ is well defined, and that we obtain this way
a structure of $U_q \mathfrak{g}$-module algebra.

For that, we have to prove some identities for $\xi \in U_q \mathfrak{g}$,
$\left(f_1 \cdot \left(t_{\{1,2 \}}t_{\{3,4 \}}\eta \right)^l\right)$,
$\left(f_2 \cdot \left(t_{\{1,2 \}}t_{\{3,4 \}}\eta \right)^l\right)$,
$f_1,f_2 \in \mathbb{C}[\widehat{U}]_q$, $l \in \mathbb{Z}$. Observe that
the left and right hand sides of those identities (up to multiplying by the
same powers of $t_{\{1,2 \}}t_{\{3,4 \}}\eta$) are just Laurent polynomials
of the indeterminate $\lambda=q^l$. So, it suffices to prove them for $l \in
\mathbb{Z}_+$, due to the well known uniqueness theorem for Laurent
polynomials. What remains is to use the fact that
$\mathbb{C}[\widehat{X}]_q$ is a $U_q \mathfrak{g}$-module algebra. \hfill
$\square$

\medskip

A more general but less elementary approach to proving statements like
proposition \ref{uni7spnssu2.texext} have been obtained in a recent work by Lunts and
Rosenberg \cite{uni7spnssu2.texLR}.

Consider the subalgebra $\mathbb{C}[U]_q=\left \{f \in
\mathbb{C}[\widehat{U}]_q|\:\deg f=(0,0,0)\right \}$. This $U_q
\mathfrak{g}$-module algebra is a q-analogue of the algebra $\mathbb{C}[U]$
of regular functions on the open $S(GL_2 \times GL_2)$-orbit $U \subset X$.

The following quite plausible statement allows one to embed the above $U_q
\mathfrak{g}$-module into a two-parameter family, the principal
non-degenerate series.

\medskip

\begin{proposition}
There exists a unique two-parameter family $\pi_{\lambda',\lambda''}$ of
representations of $U_q \mathfrak{g}$ in the vector space $\mathbb{C}[U]_q$
such that

i) for all $l',l''\in \mathbb{Z}$, $\xi \in U_q \mathfrak{g}$, $f \in
\mathbb{C}[U]_q$,
$$
\left(\pi_{q^{l'},q^{l''}}(\xi)f\right)\left(t_{\{1,2 \}}t_{\{3,4
\}}\right)^{l'}\eta^{l''}=\xi \left(f \cdot \left(t_{\{1,2 \}}\cdot t_{\{3,4
\}}\right)^{l'}\cdot \eta^{l''}\right)
$$

ii) for all $\xi \in U_q \mathfrak{g}$, $f \in \mathbb{C}[U]_q$, the vector
function $\pi_{\lambda',\lambda''}(\xi)f$ is a Laurent polynomial with
coefficients from $\mathbb{C}[U]_q$.
\end{proposition}

\medskip

Note that the method of 'analytic continuation' we use was suggested in
\cite{uni7spnssu2.texV}.

\medskip

\begin{proposition}
The $U_q \mathfrak{g}$-modules associated to $\pi_{q^{l'},q^{l''}}$, are
Harish-Chandra modules.
\end{proposition}

\smallskip

{\bf Proof.} It suffices to consider the $U_q \mathfrak{g}$-module
$\mathbb{C}[U]_q$, i.e. the case $l'=l''=0$. Obviously, this $U_q
\mathfrak{g}$-module is a $\mathbb{Z}^3$-weight module. Prove that the
multiplicity of occurrence of simple $U_q \mathfrak{k}$-modules in
$\mathbb{C}[U]_q$ is the same as in the classical case $q=1$. Equip the
vector space $\mathbb{C}[U]_q$ with a filtration: $\mathbb{C}[U]_q=\bigcup
\limits_j\mathbb{C}[U]_q^{(j)}$, where
$$
\mathbb{C}[U]_q^{(j)}=\left \{f \cdot(t_{\{1,2 \}}t_{\{3,4
\}}\eta)^{-j}\left|\:f \in \mathbb{C}[\widehat{X}]_q,\;\deg f=(j,2j,j)\right
\}\right..
$$
The $U_q \mathfrak{g}$-module
$$
\mathbb{C}[U]_q^{(j)}\cdot(t_{\{1,2 \}}t_{\{3,4 \}}\eta)^j=\left \{\left.f
\in \mathbb{C}[\widehat{X}]_q\right|\;\deg f=(j,2j,j)\right \}
$$
is simple, and $t_{\{1 \}}^jt_{\{1,2 \}}^{2j}t_{\{1,2,3 \}}^j$ is its
highest weight vector. A restriction of this representation onto $U_q
\mathfrak{k}$ splits into a direct sum of irreducible representations of
this subalgebra with the same multiplicities as in the classical case $q=1$.
Hence the multiplicities in $\mathbb{C}[U]_q^{(j)}/\mathbb{C}[U]_q^{(j-1)}$
are also classical. What remains is to observe that in the classical case
$\mathbb{C}[U]$ is a Harish-Chandra module. \hfill $\square$

\medskip

We thus obtain the spherical principal non-degenerate series:
$$
\pi_{q^{l'},q^{l''}}(\xi)1=\varepsilon(\xi)1,\qquad \xi \in U_q
\mathfrak{k}.
$$

\bigskip

\section{Duality}

Introduce the notation $V_{q^{l'},q^{l''}}$ for the Harish-Chandra modules
corresponding to the representations $\pi_{q^{l'},q^{l''}}$.

The outcome of the previous section was in producing the $U_q
\mathfrak{g}$-modules $V_{q^{l'},q^{l''}}$ via demonstrating their geometric
realization. This geometric realization is to be used in the present section
to produce a $U_q \mathfrak{g}$-invariant pairing $V_{q^{l'},q^{l''}}\times
V_{q^{-1-l'},q^{-2-l''}}\to \mathbb{C}$ (the associated linear map
$V_{q^{l'},q^{l''}}\otimes V_{q^{-1-l'},q^{-2-l''}}\to \mathbb{C}$ is going
to be a morphism of $U_q \mathfrak{g}$-modules). We restrict ourselves to
the special case $l',l''\in \mathbb{Z}$. A passage to the general case can
be performed using the analytic continuation in the parameters $u'=q^{l'}$,
$u''=q^{l''}$ as in the proof of proposition \ref{uni7spnssu2.texext}.

Recall that $\mathbb{C}[\widehat{U}]_q=\bigoplus
\limits_{j',j'',j'''}\mathbb{C}[\widehat{U}]_q^{(j',j'',j''')}$,
$$
\mathbb{C}[\widehat{U}]_q^{(j',j'',j''')}=\left.\left \{f \in
\mathbb{C}[\widehat{U}]_q \right|\:\deg f=(j',j'',j''')\right \},
$$
and one has well defined bilinear maps
$$
\mathbb{C}[\widehat{U}]_q^{(i',i'',i''')}\times
\mathbb{C}[\widehat{U}]_q^{(j',j'',j''')}\to
\mathbb{C}[\widehat{U}]_q^{(i'+j',i''+j'',i'''+j''')},\qquad f_1 \times f_2
\mapsto f_1f_2.
$$
Use the canonical isomorphisms
$$
V_{q^{l'},q^{l''}}\simeq \mathbb{C}[\widehat{U}]_q^{(l'',2l',l'')},\qquad
V_{q^{-1-l'},q^{-2-l''}}\simeq
\mathbb{C}[\widehat{U}]_q^{(-2-l'',-2-2l',-2-l'')}
$$
to obtain a bilinear map
$$
V_{q^{l'},q^{l''}}\times V_{q^{-1-l'},q^{-2-l''}}\to
\mathbb{C}[\widehat{U}]_q^{(-2,-2,-2)}.
$$
The non-degenerate $U_q \mathfrak{g}$-invariant pairing
$$
V_{q^{l'},q^{l''}}\times V_{q^{-1-l'},q^{-2-l''}}\to \mathbb{C}
$$
we are interested in can be obtained as a composition of the above bilinear
map and a $U_q \mathfrak{g}$-invariant integral
\begin{equation}\label{uni7spnssu2.texii}
\mathbb{C}[\widehat{U}]_q^{(-2,-2,-2)}\to \mathbb{C},\qquad f \mapsto \int f
d \nu,
\end{equation}
normalized as follows:
\begin{equation}\label{uni7spnssu2.texnorm}
\int\left(t_{\{1,2 \}}t_{\{3,4 \}}\eta^2 \right)^{-1}d \nu=1.
\end{equation}
What remains is to prove the existence and uniqueness of the above $U_q
\mathfrak{g}$-invariant integral $\nu$.

\medskip

\begin{proposition}\label{uni7spnssu2.texmg}
There exists a unique morphism of $U_q \mathfrak{g}$-modules (\ref{uni7spnssu2.texii})
which satisfies (\ref{uni7spnssu2.texnorm}).
\end{proposition}

\smallskip

{\bf Proof.} We start with proving the uniqueness. Consider the classical
case ($q=1$). Obviously, the $U \mathfrak{k}$-module $\mathbb{C}[U]$ of
regular functions on the orbit $U$ splits as a sum of finite dimensional
isotypical components. Hence, the multiplicity of a simple $U
\mathfrak{k}$-module $W$ in $\mathbb{C}[U]^{(j)}=\left \{f \cdot(t_{\{1,2
\}}t_{\{3,4 \}}\eta)^{-j}\left|\:f \in \mathbb{C}[\widehat{X}],\;\deg
f=(j,2j,j)\right \}\right.$ is constant for $j$ large enough. It is easy to
prove, using the standard methods of the quantum group theory \cite{uni7spnssu2.texKl-Sch},
the existence of a similar decomposition and stabilization of the
multiplicity for all the transcendental $q$. Hence, the multiplicity of
simple $U_q \mathfrak{k}$-modules in the $U_q \mathfrak{k}$-module
$\mathbb{C}[U]_q$ is finite. This multiplicity is the same as in the
classical case since it is determined by the multiplicity of weights of
finite dimensional $U_q \mathfrak{k}$-modules $\mathbb{C}[U]_q^{(j)}$. In
particular, the subspace of $U_q \mathfrak{k}$-invariants in
$\mathbb{C}[U]_q$ is one-dimensional, so we get the following

\medskip

\begin{lemma}
There exists a unique morphism of $U_q \mathfrak{k}$-modules (\ref{uni7spnssu2.texii})
which satisfies (\ref{uni7spnssu2.texnorm}).
\end{lemma}

\medskip

Turn to the proof of existence of a morphism of $U_q \mathfrak{g}$-modules
(\ref{uni7spnssu2.texii}) which satisfies (\ref{uni7spnssu2.texnorm}).

\medskip

\begin{lemma}\label{uni7spnssu2.texgen}
\hfill

i) The element $\omega=\left(t_{\{1,2 \}}t_{\{3,4 \}}\eta^2 \right)^{-1}$
generates the $U_q \mathfrak{g}$-module
$\mathbb{C}[\widehat{U}]_q^{(-2,-2,-2)}$.

ii) The following relations are valid: $\xi \omega=\varepsilon(\xi)\omega$
for $\xi \in U_q \mathfrak{k}$;\ $z \omega=\varepsilon(z)\omega$, for $z \in
Z(U_q \mathfrak{g})$. The standard notation $Z(U_q \mathfrak{g})$ is used
here to denote the center of $U_q \mathfrak{g}$.
\end{lemma}

\smallskip

{\bf Proof.} i) Since $q$ is transcendental, one can easily reduce the proof
of the first statement to the proof of a similar statement in the classical
case $q=1$. For this case the result is a consequence of the well known
Kostant theorem about cyclic vectors for principal series \cite{uni7spnssu2.texKo} (observe
that $\omega$ is the primitive weight vector with respect to the left action
of $U_q \mathfrak{g}$ and its weight is $-2 \rho$, with $\rho$ being the
half sum of positive roots).

ii) The relation $\xi \omega=\varepsilon(\xi)\omega$, $\xi \in U_q
\mathfrak{k}$, is obvious.

Consider an element $z \in Z(U_q \mathfrak{g})$, together with its image
$\psi_z$ under the Harish-Chandra isomorphism. One has for $m \in
\mathbb{Z}_+$:
\begin{equation}\label{uni7spnssu2.texz1}
z\left(\left(t_{\{1,2 \}}t_{\{3,4 \}}\eta^2 \right)^m\right)=\psi_z(2m
\rho)\left(t_{\{1,2 \}}t_{\{3,4 \}}\eta^2 \right)^m.
\end{equation}
That is,
\begin{equation}\label{uni7spnssu2.texz2}
\left(t_{\{1,2 \}}t_{\{3,4 \}}\eta^2 \right)^{-m}z \left(\left(t_{\{1,2
\}}t_{\{3,4 \}}\eta^2 \right)^m \right)=\psi_z(2m \rho).
\end{equation}
Both sides of (\ref{uni7spnssu2.texz2}) are Laurent polynomials of $q^m$. Hence (\ref{uni7spnssu2.texz1}),
(\ref{uni7spnssu2.texz2}) are valid for all $m \in \mathbb{Z}$. In particular,
$z\left(\left(t_{\{1,2 \}}t_{\{3,4 \}}\eta^2 \right)^{-1}\right)=\psi_z(-2
\rho)\left(t_{\{1,2 \}}t_{\{3,4 \}}\eta^2 \right)^{-1}$. What remains is to
prove that $\psi_z(-2 \rho)=\varepsilon(z)$. It suffices to prove the
equalities $\psi_z(-2 \rho)=\psi_z(0)$, $\psi_z(0)=\varepsilon(z)$. The
first of them is valid due to the invariance of $\psi_z$ with respect to the
standard action of the Weyl group. The second equality is obvious. \hfill
$\square$

\medskip

{\sc Remark.} It is convenient to pass from the Hopf algebras $U_q
\mathfrak{k} \subset U_q \mathfrak{g}$ to their 'extensions' $U_q
\mathfrak{k}^\mathrm{ext}\subset U_q \mathfrak{g}^\mathrm{ext}$ as in
\cite[chapter 8.5.3]{uni7spnssu2.texKl-Sch}. \footnote{The reason is that the well known
q-analogues of standard generators of $Z(U \mathfrak{g})$ belong to $U_q
\mathfrak{g}^\mathrm{ext}$ \cite{uni7spnssu2.texJ}} (These are derived via adding certain
products of $K_i^{\pm 1/4}$, $i=1,2,3$). Obviously,
$\mathbb{C}[\widehat{U}]_q^{(-2,-2,-2)}$ is a $U_q
\mathfrak{g}^\mathrm{ext}$-module, and lemma \ref{uni7spnssu2.texgen} can be proved for the
$U_q \mathfrak{g}^\mathrm{ext}$-module structure (in fact in the proof of
lemma \ref{uni7spnssu2.texgen} it is possible to use a quantum analogue of the
Harish-Chandra isomorphism for algebras over the formal series of $h$).

\medskip

\begin{corollary}
Consider the $U_q \mathfrak{g}^\mathrm{ext}$-module $V$ with a single
generator $v$ subject to the relations
$$
\xi v=\varepsilon(\xi)v,\quad \xi \in U_q \mathfrak{k}^\mathrm{ext};\qquad
zv=\varepsilon(z)v,\quad z \in Z(U_q \mathfrak{g}^\mathrm{ext}).
$$
The map $v \mapsto \omega$ is uniquely extendable up to an onto morphism of
$U_q \mathfrak{g}^\mathrm{ext}$-modules $\varphi:V \to
\mathbb{C}[\widehat{U}]_q^{(-2,-2,-2)}$.
\end{corollary}

\medskip

Turn back to the proof of proposition \ref{uni7spnssu2.texmg}. Evidently, there exists a
unique morphism of $U_q \mathfrak{g}^\mathrm{ext}$-modules $\mu:V \to
\mathbb{C}$ such that $\mu(1)=1$. So, the existence of a morphism
$\nu:\mathbb{C}[\widehat{U}]_q^{(-2,-2,-2)}\to \mathbb{C}$ follows from the
isomorphism $\mathbb{C}[\widehat{U}]_q^{(-2,-2,-2)}\simeq
V/\mathrm{Ker}\,\varphi$ and the following

\medskip

\begin{lemma}
$\mathrm{Ker}\,\mu \supset \mathrm{Ker}\,\varphi$.
\end{lemma}

\smallskip

{\bf Proof.} It suffices to prove firstly, that the $U_q
\mathfrak{k}^\mathrm{ext}$-module $V$ splits into a sum of finite
dimensional simple $U_q \mathfrak{k}^\mathrm{ext}$-modules and, secondly,
that $\mathrm{Ker}\,\varphi$ contains no non-zero $U_q
\mathfrak{k}^\mathrm{ext}$-invariant vectors. Consider the Hopf subalgebras
$U_q \mathfrak{p}_{\pm}^\mathrm{ext}$ generated respectively by
$\{E_j,F_j\}_{j \ne 2}$, $E_2$, $U_q \mathfrak{k}^\mathrm{ext}$; \
$\{E_j,F_j \}_{j \ne 2}$, $F_2$, $U_q \mathfrak{k}^\mathrm{ext}$. The first
statement follows from the fact that the $U_q
\mathfrak{k}^\mathrm{ext}$-module $U_q \mathfrak{g}^\mathrm{ext}/U_q
\mathfrak{k}^\mathrm{ext}\simeq(U_q \mathfrak{p}_+^\mathrm{ext}/U_q
\mathfrak{k}^\mathrm{ext})\otimes(U_q \mathfrak{p}_-^\mathrm{ext}/U_q
\mathfrak{k}^\mathrm{ext})$ splits as a sum of simple finite dimensional
$U_q \mathfrak{k}^\mathrm{ext}$-modules (to prove this, one can use suitable
bases in $U_q \mathfrak{g}^\mathrm{ext}$, cf. \cite{uni7spnssu2.texSSV2}). The second
statement is due to the following

\medskip

\begin{lemma}
$$
\dim V^{U_q \mathfrak{g}^\mathrm{ext}}=1.\;\footnote{$V^{U_q
\mathfrak{g}^\mathrm{ext}}=\{v' \in V|\:\xi v'=\varepsilon(\xi)v',\;\xi \in
U_q \mathfrak{g}^\mathrm{ext}\}$}
$$
\end{lemma}

\smallskip

{\bf Proof.} Consider a $U_q \mathfrak{g}^\mathrm{ext}$-module $U_q
\mathfrak{g}^\mathrm{ext}/U_q \mathfrak{k}^\mathrm{ext}$. It suffices to
prove that the canonical map
\begin{equation}\label{uni7spnssu2.texcm}
Z(U_q \mathfrak{g}^\mathrm{ext})\to \left(U_q \mathfrak{g}^\mathrm{ext}/U_q
\mathfrak{k}^\mathrm{ext}\right)^{U_q \mathfrak{k}^\mathrm{ext}}
\end{equation}
is onto.

Equip $U_q \mathfrak{g}^\mathrm{ext}$ with a $U_q
\mathfrak{k}^\mathrm{ext}$-invariant filtration via setting
$$
\deg(E_2)=\deg(F_2)=1;\quad \deg(\xi)=0,\;\xi \in U_q
\mathfrak{k}^\mathrm{ext}.
$$
Evidently, the subspaces $\{\xi \in U_q \mathfrak{g}^\mathrm{ext}/U_q
\mathfrak{k}^\mathrm{ext}|\:\deg \xi \le j\}$ are finite dimensional. It
suffices to prove that for all $j$ the canonical map
$$
\{\xi \in Z(U_q \mathfrak{g}^\mathrm{ext})|\:\deg \xi \le j\}\to \{\xi \in
(U_q \mathfrak{g}^\mathrm{ext}/U_q \mathfrak{k}^\mathrm{ext})^{U_q
\mathfrak{k}^\mathrm{ext}}|\:\deg \xi \le j\}
$$
is onto.

An application of the standard generators of $Z(U \mathfrak{g})$ allows one
to reduce easily the quantum case ($q$ transcendental) to the classical one
($q=1$). In the classical case the surjectivity of the map (\ref{uni7spnssu2.texcm}) is a
well known fact (see \cite[Th. 1.3.1]{uni7spnssu2.texKo}). \hfill $\square$

\bigskip

\section*{Notes of the Editor}

The methods of this work can be used to produce the principal {\bf unitary}
non-degenerate series, that is, to prove the unitarizability of the
Harish-Chandra modules $V_{q^{l'},q^{l''}}$ in the case
$\mathrm{Im}\,l'=-1/2$, $\mathrm{Im}\,l''=-1$.

Here is an outline of the proof. Let us consider the pseudo-Hermitian space
$\mathbb{C}^4$ endowed with the metric $-|t_1|^2-|t_2|^2+|t_3|^2+|t_4|^2$
and the projective algebraic variety of uncomplete 'Lagrange flags'
$$
X_0=\left\{\left.0 \varsubsetneqq L_1 \varsubsetneqq L_2 \subset
\mathbb{C}^4 \right|\;L_2=L_2^\perp \right \}\subset X.
$$
Let $\widehat{X}_0 \subset \widehat{X}$ be the {\bf real} affine algebraic
variety associated to $X_0$. Introduce a $U_q \mathfrak{su}_{2,2}$-module
algebra $\mathrm{Pol}(\widehat{X}_0)_q$, which is a q-analogue of the
algebra of regular functions on $\widehat{X}_0$.

In the category of non-involutive algebras $\mathrm{Pol}(\widehat{X}_0)_q$
can be defined as a subalgebra of $\mathbb{C}[\widehat{X}]_q$ generated by
\begin{gather*}
t_{\{i \}}t_{\{j_1,j_2,j_3 \}},\qquad 1 \le i \le 4,\quad 1 \le j_1<j_2<j_3
\le 4,
\\ t_{\{i_1,i_2 \}}t_{\{j_1,j_2 \}},\qquad 1 \le i_1<i_2 \le 4,\quad 1 \le
j_1<j_2 \le 4.
\end{gather*}
It is easy to prove the existence and uniqueness of the involution $*$ which
equips $\mathrm{Pol}(\widehat{X}_0)_q$ with a structure of $U_q
\mathfrak{su}_{2,2}$-module algebra and such that the elements
$$
\left(t_{\{1,2 \}}t_{\{3,4 \}}\right)^j \eta^k,\quad j,k \in \mathbb{Z}_+
$$
appear to be selfadjoint.\footnote{The uniqueness follows from the
selfadjointness of the monomials $(t_{\{1,2 \}}t_{\{3,4 \}})^j \eta^k$ and
from the fact that $\mathrm{Pol}(\widehat{X}_0)_q$ is a sum of simple $U_q
\mathfrak{su}_{2,2}$-modules generated by the above monomials:
$\mathrm{Pol}(\widehat{X}_0)_q=\bigoplus \limits_{j,k=0}^\infty L(k,2j,k)$.
Sketch the proof of the existence. Equip every subspace $L(k,2j,k)$ with an
antilinear map $*$ which satisfies all the necessary requirements, possibly
except
$$(f_1f_2)^*=f_2^*f_1^*,\qquad f_1,f_2 \in L(k,2j,k).$$
What remains is to extend it onto $\mathrm{Pol}(\widehat{X}_0)_q$ by
antilinearity and to prove the latter identity using the above claim of
uniqueness.}

The involution in $\mathrm{Pol}(\widehat{X_0})_q$ is canonically extendable
onto the localization of this $*$-algebra with respect to the
multiplicatively closed set $(t_{\{1,2 \}}t_{\{3,4 \}})^j \eta^k$, $j,k \in
\mathbb{Z}_+$.

Consider the $U_q \mathfrak{k}$-module subalgebra of zero degree homogeneous
elements. Its $U_q \mathfrak{k}$-isotypic components are finite dimensional
and are smooth functions of $q \in (0,1]$ (this could be well rephrased
rigorously). In particular, the dimensions of $U_q \mathfrak{k}$-isotypic
components are independent of $q \in(0,1]$. Hence there exists a unique $U_q
\mathfrak{k}$-invariant integral $\mu$ on this $*$-subalgebra with the
property $\mu(1)=1$. This integral is positive in the classical case $q=1$,
hence for all $0<q<1$ (otherwise for some $U_q \mathfrak{k}$-isotypic
component and some $0<q<1$ the $U_q \mathfrak{k}$-invariant scalar product
given by $(f_1,f_2)=\mu(f_2^*f_1)$ would appear to be non-negative but not
strictly positive, hence degenerate). In view of proposition \ref{uni7spnssu2.texmg} the
scalar product
$$
\left(f_1 \cdot (t_{\{1,2 \}}t_{\{3,4 \}})^{l'}\eta^{l''},f_2 \cdot
(t_{\{1,2 \}}t_{\{3,4
\}})^{l'}\eta^{l''}\right)\stackrel{\mathrm{def}}{=}\mu(f_2^*f_1)
$$
is $U_q \mathfrak{su}_{2,2}$-invariant in the case $\mathrm{Im}\,l'=-1/2$,
$\mathrm{Im}\,l''=-1$.




\title{\bf HIDDEN SYMMETRY OF THE DIFFERENTIAL CALCULUS ON THE QUANTUM
MATRIX SPACE}

\author{S. Sinel'shchikov\thanks{Partially supported by ISF grant
U2B200 \newline E-mail address: sinelshchikov@ilt.kharkov.ua}
\and L. Vaksman\thanks{Partially supported by ISF grant U21200 and
INTAS grant 4720 \newline E-mail address: vaksman@ilt.kharkov.ua}}
\date{\tt Mathematics Department, Institute for Low Temperature Physics and
Engineering, 47 Lenin Avenue, 61103 Kharkov, Ukraine\\
\bigskip\bigskip\raggedright PAX: 02.40-k geometry\\
PAX: 03.65-Fd algebraic methods, differential geometry, topology}

\newpage
\setcounter{section}{0}
\large
\thispagestyle{empty}
\ \vfill \begin{center}\LARGE \bf PART IV \\
ADDITONAL RESULTS ON SOME QUANTUM VECTOR SPACES
\end{center}
\vfill
\addcontentsline{toc}{chapter}{Part IV \ \
ADDITONAL RESULTS ON SOME QUANTUM VECTOR SPACES}
\newpage
\makeatletter

\makeatletter
\renewcommand{\@oddhead}{HIDDEN SYMMETRY OF THE DIFFERENTIAL CALCULUS \hfill
\thepage}
\renewcommand{\@evenhead}{\thepage \hfill S. Sinel'shchikov and L. Vaksman}
\let\@thefnmark\relax
\@footnotetext{This research was partially supported by Award No UM1-2091 of
the Civilian Research \& Development Foundation \newline \indent This
lecture has been delivered at the P. Vogel seminar, Paris, 1996; published
in J. Phys. A: Math. Gen., {\bf 30} (1997), L23 -- L26}
\addcontentsline{toc}{chapter}{\@title \\ {\sl S. Sinel'shchikov and L.
Vaksman}\dotfill} \makeatother

\maketitle

\begin{quotation}\small{\bf Abstract.}
A standard bicovariant differential calculus on the quantum matrix space
$\mathrm{Mat}(m,n)_q$ is considered. Our main result is in proving that the
$U_q\mathfrak{s}(\mathfrak{gl}_m \times \mathfrak{gl}_n)$-module
differential algebra $\Omega^*(\mathrm{Mat}(m,n))_q$ is in fact a $U_q
\mathfrak{sl}(m+n)$-module differential algebra.
\end{quotation}

{\bf 1.} This work solves a problem whose simple special case occurs in a
construction of a quantum unit ball of $\mathbb{C}^n$ (in the spirit of
\cite{uni8HSDCQMS.TEXR}). Within the framework of that theory, the action of the subgroup
$SU(n,1)\subset SL(n+1)$ by automorphisms of the unit matrix ball is
essential. The problem is that the Wess-Zumino differential calculus in
quantum $\mathbb{C}^n$ \cite{uni8HSDCQMS.TEXWZ} at a first glance seems to be only $U_q
\mathfrak{sl}_n$-invariant. In that particular case the lost $U_q
\mathfrak{sl}_{n+1}$-symmetry can be easily detected. The main result of
this work is in disclosing the hidden $U_q \mathfrak{sl}_n $-symmetry for
bicovariant differential calculus in the quantum matrix space
$\mathrm{Mat}(m,n)$. (Note that for $n=1$ we have the case of a ball).

The authors are grateful to V. Akulov and G. Maltsiniotis for a helpful
discussion of the results.

\bigskip

{\bf 2.} We start with recalling the definition of the Hopf algebra
$U_q \mathfrak{sl}_N$, $N>1$, over the field $\mathbb{C}(q)$ of rational
functions of an indeterminate $q$ \cite{uni8HSDCQMS.TEXD,uni8HSDCQMS.TEXJ}. (We follow the notations of
\cite{uni8HSDCQMS.TEXCK}).

For $i,j \in \{1,\ldots,N-1 \}$ let
$$
a_{ij}=\left \{\begin{array}{rl}{2,}&{i-j=0}\cr{-1,}&{|i-j|=1}\cr{0,}&
{|i-j|>1.}\end{array}\right.
$$

\smallskip

The algebra $U_q \mathfrak{sl}_N$ is defined by the generators
$\{E_i,\:F_i,\:K_i,\:K_i^{-1}\}$ and the relations

$$K_iK_j \,=\,K_jK_i,\quad K_iK_i^{-1}\,=\,K_i^{-1}K_i \,=\,1$$
$$K_iE_j \,=\,q^{a_{ij}}E_jK_i,\quad K_iF_j \,=\,q^{-a_{ij}}F_jK_i$$
$$E_iF_j \,-\,F_jE_i \:=\:\delta_{ij}(K_i-K_i^{-1})/(q-q^{-1})$$
$$
E_i^2E_j\,-\,(q+q^{-1})E_iE_jE_i \,+\,E_jE_i^2 \:=\:0,\quad |i-j| \,=\,1
$$
$$
F_i^2F_j \,-\,(q+q^{-1})F_iF_jF_i \,+\,F_jF_i^2\:=\:0, \quad |i-j|\,=\,1
$$
$$[E_i,E_j]\,=\,[F_i,F_j]\,=\,0, \quad |i-j|\,\ne \,1.$$

\smallskip

A comultiplication $\Delta$, an antipode $S$ and a counit $\varepsilon$ are
defined by
$$
\Delta E_i \:=\:E_i \otimes 1 \,+\,K_i \otimes E_i,\quad \Delta F_i \:=\:F_i
\otimes K_i^{-1}\,+\,1 \otimes F_i,
$$
$$\Delta K_i \,=\,K_i \otimes K_i,\quad S(E_i)\,=\,-K_i^{-1}E_i,$$
$$S(F_i)\,=\,-F_iK_i,\quad S(K_i)\,=\,K_i^{-1},$$
$$\varepsilon(E_i)=\varepsilon(F_i)=0,\quad \varepsilon(K_i)=1.$$

\bigskip

{\bf 3.} Remind a description of a differential algebra
$\Omega^*(\mathrm{Mat}(m,n))_q$ on a quantum matrix space \cite{uni8HSDCQMS.TEXCP,uni8HSDCQMS.TEXM1}.

Let $i,j,i',j'\in \{1,2,\ldots,m+n \}$, and
$$
\check{R}_{ij}^{i'j'}\;=\;\left
\{\begin{array}{rl}{q^{-1},}&{i=j=i'=j'}\cr{1,}& {i'=j \;\mathrm{and}\;j'=i
\;\mathrm{and}\;i \ne j}\cr{q^{-1}-q,}&{i=i'\;\mathrm{and}\;j=j'\;
\mathrm{and}\;i<j}\cr{0,}&{\mathrm{otherwise}}\end{array}\right.
$$
$\Omega^*(\hbox{Mat}(m,n))_q$ is given by the generators $\{t_a^\alpha \}$
and relations
$$
\sum_{\gamma,\delta}\check{R}_{\gamma \delta}^{\alpha \beta}t_a^\gamma
t_b^\delta \,=\,\sum_{c,d}\check{R}_{ab}^{cd}t_d^\beta t_c^\alpha
$$
$$
\sum_{a',b',\gamma',\delta'}\check{R}_{\gamma'\delta'}^{\alpha \beta}
\check{R}_{ab}^{a'b'}t_{a'}^{\gamma'}dt_{b'}^{\delta'}\;=\;dt_a^\alpha
t_b^\beta
$$
$$
\sum_{a',b',\gamma',\delta'}\check{R}_{\gamma'\delta'}^{\alpha \beta}
\check{R}_{ab}^{a'b'}dt_{a'}^{\gamma'}dt_{b'}^{\delta'}\;=\;-dt_a^\alpha
dt_b^\beta
$$

\noindent $(a,b,c,d,a',b'\,\in\{1,\ldots,n\};\quad
\alpha,\beta,\gamma,\delta,\gamma',\delta'\,\in\{1,\ldots,m\}).$

Let us define a grading by
$\deg(t_a^\alpha)\,=\,0,\;\deg(dt_a^\alpha)\,=\,1$. With that,
$\mathbb{C}[\mathrm{Mat}(m,n)]_q\,=\,\Omega^0(\mathrm{Mat}(m,n)))_q$ will
stand for a subalgebra of zero degree elements .

\bigskip

{\bf 4.} Let $A$ be a Hopf algebra and $F$ an algebra with unit and an
$A$-module the same time. $F$ is said to be a $A$-module algebra \cite{uni8HSDCQMS.TEXA} if
the multiplication $m:\:F \otimes F \,\to \,F$ is a morphism of $A$-modules,
and $1 \in F$ is an invariant (that is\\ $a(f_1f_2)\,=\,\sum \limits_j
a_j'f_1 \otimes a_j''f_2,\quad a1 \,=\,\varepsilon(a)1$ for all $a \in
A;\;f_1,f_2 \in F$, with $\Delta(a)\,=\,\sum \limits_ja_j'\otimes a_j''$).

An important example of an $A$-module algebra appears if one supplies
$A^*$ with the structure of an $A$-module: $\langle af,b \rangle \,=\,\langle
f,ba \rangle,\;a,b \in A,\:f \in A^*$.

\bigskip

{\bf 5.} Our immediate goal is to furnish $\mathbb{C}[\mathrm{Mat}(m,n)]_q$
with a structure of a $U_q \mathfrak{sl}_{m+n}$-module algebra via an
embedding $\mathbb{C}[\mathrm{Mat}(m,n)]_q \hookrightarrow (U_q
\mathfrak{sl}_{m+n})^*$.

Let $\{e_{ij}\}$ be a standard basis in $\mathrm{Mat}(m+n)$ and $\{f_{ij}\}$
the dual basis in $\mathrm{Mat}(m+n)^*$. Consider a natural representation
$\pi$ of $U_q \mathfrak{sl}_{m+n}$:
$$
\pi(E_i)\,=\,e_{i \:i+1},\quad \pi(F_i)\,=\,e_{i+1 \:i}, \quad
\pi(K_i)\:=\:qe_{ii}\,+\,q^{-1}e_{i+1 \:i+1}\,+\,\sum_{j \ne i,i+1}e_{jj}.
$$

The matrix elements $u_{ij}\:=\:f_{ij}\pi \in(U_q \mathfrak{sl}_{m+n})^*$ of
the natural representation may be treated as "coordinates" on the quantum
group $SL_{m+n}$ \cite{uni8HSDCQMS.TEXD}. To construct "coordinate" functions on a big cell
of the Grassmann manifold, we need the following elements of
$\mathbb{C}[\mathrm{Mat}(m,n)]_q$
$$
x(j_1,j_2,\ldots,j_m)\:=\:\sum \limits_{w \in
S_m}\,(-q)^{l(w)}u_{1j_{w(1)}}u_{2j_{w(2)}}\ldots u_{mj_{w(m)}},
$$
with $1 \le j_1<j_2<\ldots<j_m \le m+n$, and
$l(w)\;=\;\mathrm{card}\{(a,b)|\:a<b \;\mathrm{and}\;w(a)>w(b)\}$ being
the "length" of a permutation $w \in S_m$.

\medskip

\begin{proposition}\label{uni8HSDCQMS.TEX1}
$x(1,2,\ldots,m)$ is invertible in $(U_q \mathfrak{sl}_{m+n})^*$, and the
map
$$
t_a^\alpha \mapsto
x(1,2,\ldots,m)^{-1}x(1,\ldots,\widehat{m+1-\alpha},\ldots,m,m+a)
$$
can be extended up to an embedding
$$
i:\;\mathbb{C}[\mathrm{Mat}(m,n)]_q \hookrightarrow(U_q
\mathfrak{sl}_{m+n})^*.
$$
(The sign $\widehat{\;}$ here indicates the item in a list that should be
omitted).
\end{proposition}

\medskip

Proposition \ref{uni8HSDCQMS.TEX1} allows one to equip $\mathbb{C}[\mathrm{Mat}(m,n)]_q$
with the structure of a $U_q \mathfrak{sl}_{m+n}$-module algebra :\\
$$
i \xi t_a^\alpha \,=\,\xi it_a^\alpha, \qquad \xi \in U_q
\mathfrak{sl}_{m+n},\; a \in \{1,\ldots,n \},\;\alpha \in \{1,\ldots,m \}.
$$

\bigskip

{\bf 6.} The main result of our work is the following

\medskip

\begin{theorem}\label{uni8HSDCQMS.TEX2}
$\Omega^*(\mathrm{Mat}(m,n))_q$ admits a unique structure of a
$U_q \mathfrak{sl}_{m+n}$-module algebra such that the embedding
$$
i:\:\mathbb{C}[\mathrm{Mat}(m,n)]_q \hookrightarrow
\Omega^*(\mathrm{Mat}(m,n))_q
$$
and the differential
$$d:\;\Omega^*(\mathrm{Mat}(m,n))_q \,\to \,\Omega^*(\mathrm{Mat}(m,n))_q$$
are the morphisms of $U_q \mathfrak{sl}_{m+n}$-modules.
\end{theorem}

\medskip

{\sc Remark 1.} The bicovariance of the differential calculus on the quantum
matrix space allows one to equip the algebra $\Omega^*(\mathrm{Mat}(m,n))_q$
with a structure of $U_q \mathfrak{s}(\mathfrak{gl}_m \times
\mathfrak{gl}_n)$-module, which is compatible with multiplication in
$\Omega^*(\mathrm{Mat}(m,n))_q$ and differential $d$. Theorem \ref{uni8HSDCQMS.TEX2}
implies that $\Omega^*(\mathrm{Mat}(m,n))_q$ possess an additional hidden
symmetry since $U_q \mathfrak{sl}_{m+n}\,\supsetneqq \,U_q
\mathfrak{s}(\mathfrak{gl}_m \times \mathfrak{gl}_n)$.

\medskip

{\sc Remark 2.} Let $q_0 \in \mathbb{C}$ and $q_0$ is not a root of unity.
It follows from the explicit formulae for $E_mt_a^\alpha,\;F_mt_a^\alpha,
\;K_m^{\pm 1}t_a^\alpha,\quad a \in \{1,\ldots,n \},\quad \alpha \in \{1,
\ldots,m \}$, that the "specialization" $\Omega^*(\mathrm{Mat}(m,n))_{q_0}$
is a $U_{q_0}\mathfrak{sl}_{m+n}$-module algebra.

\bigskip

{\bf 7.} Supply the algebra $U_q \mathfrak{sl}_{m+n}$ with a grading as
follows:
$$
\deg(K_i)\,=\,\deg(E_i)\,=\,\deg(F_i)\,=\,0,\quad \mathrm{for}\quad i \ne m,
$$
$$\deg(K_m)\,=\,0,\quad \deg(E_m)\,=\,1,\quad \deg(F_m)\,=\,0.$$

The proofs of Proposition \ref{uni8HSDCQMS.TEX1} and Theorem \ref{uni8HSDCQMS.TEX2} reduce to the
construction of graded $U_q \mathfrak{sl}_{m+n}$-modules which are dual
respectively to the modules of functions $\Omega^0(\mathrm{Mat}(m,n))_q$ and
that of 1-forms $\Omega^1(\mathrm{Mat}(m,n))_q$. The dual modules are
defined by their generators and correlations. While proving the completeness
of the correlation list we implement the "limit specialization" $q_0=1$ (see
\cite[p. 416]{uni8HSDCQMS.TEXCK}).

The passage from the order one differential calculus
$\Omega^0(\mathrm{Mat}(m,n))_q \,\stackrel{d}{\to}\,
\Omega^1(\mathrm{Mat}(m,n))_q$ to $\Omega^*(\mathrm{Mat}(m,n))_q$ is done
via a universal argument described in a paper by G. Maltsiniotis \cite{uni8HSDCQMS.TEXM2}.
This argument doesn't break $U_q \mathfrak{sl}_{m+n}$-symmetry.

\bigskip

{\bf 8.} Our approach to the construction of order one differential calculus
is completely analogous to that of V. Drinfel'd \cite{uni8HSDCQMS.TEXD} used initially to
produce the algebra of functions on a quantum group by means of a universal
enveloping algebra.

\bigskip

{\bf 9.} The space of matrices is the simplest example of an irreducible
prehomogeneous vector space of parabolic type \cite{uni8HSDCQMS.TEXKi}. Such space can be
also associated to a pair constituted by a Dynkin diagram of a simple Lie
algebra $\mathcal{G}$ and a distinguished vertex of this diagram. Our method
can work as an efficient tool for producing $U_q \mathcal{G}$-invariant
differential calculi on the above prehomogeneous vector spaces.

Note that $U_q \mathcal{G}$-module algebras of polynomials on quantum
prehomogeneous spaces of parabolic type were considered in a recent work of
M. S. Kebe \cite{uni8HSDCQMS.TEXKe}.

\bigskip

\noindent{\bf Acknowledgement}

\smallskip

The authors would like to express their gratitude to Prof. A. Boutet de
Monvel at  University Paris VII  for the warm hospitality during the work on
this paper.

\bigskip

\section*{Notes of the Editor}

The proofs of the results announced in this work can be found in
\cite{uni8HSDCQMS.TEXSSV1}.




\makeatletter \@addtoreset{equation}{section}\makeatother
\renewcommand{\theequation}{\thesection.\arabic{equation}}

\title{\bf q-ANALOGS OF CERTAIN PREHOMOGENEOUS VECTOR SPACES: COMPARISON OF
SEVERAL APPROACHES}

\author{D. Shklyarov}
\date{}

\newpage
\setcounter{section}{0}
\large

\makeatletter
\renewcommand{\@oddhead}{q-ANALOGS OF PREHOMOGENEOUS VECTOR SPACES: SEVERAL
APPROACHES \hfill \thepage}
\renewcommand{\@evenhead}{\thepage \hfill D. Shklyarov}
\let\@thefnmark\relax
\@footnotetext{This research was supported in part by Award No UM1-2091 of
the US Civilian Research \& Development Foundation \newline \indent This
lecture has been delivered at the seminar 'Quantum groups', Kharkov, April
2000; published in Matematicheskaya Fizika. Analiz. Geometriya, {\bf 8}
(2001), No 3, 325 -- 345.}
\addcontentsline{toc}{chapter}{\@title \\ {\sl D. Shklyarov}\dotfill}
\makeatother

\maketitle

\section{Introduction: prehomogeneous vector spaces of commutative parabolic
type}

Let $\mathfrak{g}$ be a complex simple Lie algebra, $\mathfrak{h}$ its
Cartan subalgebra, $\{\alpha_i \}_{i=\overline{1,l}}$ the simple roots with
respect to $\mathfrak{h}$. Let us associate to each $\alpha_0 \in \{\alpha_i
\}_{i=\overline{1,l}}$ the element $H_0 \in \mathfrak{h}$ such that
$$
\alpha_i(H_0)=
\begin{cases}
2, & \alpha_i=\alpha_0,\\ 0, & \text{otherwise}.
\end{cases}
$$
Such $H_0$ can be used to equip $\mathfrak{g}$ with a $\mathbb{Z}$-grading
as follows
$$
\mathfrak{g}=\bigoplus_{i \in \mathbb{Z}}\mathfrak{g}_i,\qquad
\mathfrak{g}_i \stackrel{\mathrm{def}}{=}\{\xi \in
\mathfrak{g}|\quad[H_0,\xi] =2i \xi \}.
$$
If $\mathfrak{g}_i$ is nonzero only for $i \in \{-1,0,1 \}$ then the
subspace $\mathfrak{g}_{-1}$ is said to be a prehomogeneous vector space of
commutative parabolic type (see \cite{uni9pproach.texRu}).

\medskip

{\sc Remark 1.} Let us explain the adjective "prehomogeneous". Let $K
\subset{\mathrm{Aut}} \mathfrak{g}$ be the subgroup of those automorphisms
which preserve the decomposition
$$
\mathfrak{g}=\mathfrak{g}_{-1}\bigoplus \mathfrak{g}_{0}\bigoplus
\mathfrak{g}_{+1}.
$$
Then $\mathfrak{g}_{0}={\mathrm{Lie}} K$. The group $K$ acts in
$\mathfrak{g}_{-1}$ in such way that there exists a Zariski open $K$-orbit.
A pair $(G,V)$ ($G$ is an algebraic group acting in the vector space $V$)
which possesses this property is called a prehomogeneous vector space.

\medskip

Since $K$ acts in $\mathfrak{g}_{-1}$ and
$\mathfrak{g}_{0}=\mathrm{Lie}\,K$, one may consider the corresponding
representation of $U \mathfrak{g}_{0}$ in the space
$\mathbb{C}[\mathfrak{g}_{-1}]$ of holomorphic polynomials on
$\mathfrak{g}_{-1}$.

\medskip

{\sc Remark 2.} The Killing form of $\mathfrak{g}$ makes the vector spaces
$\mathfrak{g}_{-1}$ and $\mathfrak{g}_{+1}$ dual to each other. This allows
one to identify the algebras $\mathbb{C}[\mathfrak{g}_{-1}]$ and
$\mathrm{S}(\mathfrak{g}_{+1})$ (the symmetric algebra over
$\mathfrak{g}_{+1}$). The latter algebra is isomorphic to $U
\mathfrak{g}_{+1}$ for $\mathfrak{g}_{+1}$ is an Abelian Lie subalgebra in
$\mathfrak{g}$. The action of $U \mathfrak{g}_{0}$ in
$\mathbb{C}[\mathfrak{g}_{-1}]$ we deal with corresponds (under the
isomorphism $\mathbb{C}[\mathfrak{g}_{-1}]\simeq U \mathfrak{g}_{+1}$) to
the adjoint action of $U \mathfrak{g}_0$ in $U \mathfrak{g}_{+1}$.

\bigskip

There exist several approaches to constructing a $q$-analog of the algebra
$\mathbb{C}[\mathfrak{g}_{-1}]$. In the present paper we concern with those
developed in \cite{uni9pproach.texJ}, \cite{uni9pproach.texKMT}, \cite{uni9pproach.texSV2}. Within framework of each
approach a (noncommutative) analog of $\mathbb{C}[\mathfrak{g}_{-1}]$ is
endowed with an action of the quantum universal enveloping algebra $U_q
\mathfrak{g}_0$. We prove that $q$-analogs of
$\mathbb{C}[\mathfrak{g}_{-1}]$ constructed in \cite{uni9pproach.texJ}, \cite{uni9pproach.texKMT},
\cite{uni9pproach.texSV2} are isomorphic as $U_q \mathfrak{g}_0$-module algebras.

For the sake of simplicity we carry out this program for
$\mathfrak{g}=\mathfrak{sl}_4(\mathbb{C})$, $\alpha_0=\alpha_2$ (in this
case $\mathfrak{g}_{-1}$ is the space of $2 \times2$ complex matrices,
$K=S(GL_2(\mathbb{C})\times GL_2(\mathbb{C}))$). But we present such proofs
that are transferable to the case of an arbitrary prehomogeneous vector
space of commutative parabolic type.

The author thanks L. Vaksman for stating the problem and discussing the
results.

\bigskip

\section{Notation and auxiliary facts}

Let $\mathfrak{g}=\mathfrak{sl}_4(\mathbb{C})$. For convenience, we identify
$\mathfrak{g}$ with the Lie algebra of $4 \times4$ complex matrices with
zero trace. Let $\mathfrak{h}\subset \mathfrak{g}$ be the Cartan subalgebra
of diagonal matrices. Denote by $\Delta$ and $\mathrm{W}$ the root system of
$\mathfrak{g}$ with respect to $\mathfrak{h}$ and the Weyl group of this
system, respectively. Let also $\alpha_1$, $\alpha_2$, $\alpha_3$ be the
simple roots in $\Delta$ given by
$$\alpha_i(H)=a_i-a_{i+1}$$
with $H=\mathrm{diag}(a_1, a_2, a_3, a_4)\in \mathfrak{h}$.

There exists an isomorphism of the group $\mathrm{W}$ onto the group
$\mathrm{S}_4$ such that the simple reflections $s_{\alpha_1}$,
$s_{\alpha_2}$, $s_{\alpha_3}$ correspond to the transpositions $(1,2)$,
$(2,3)$, $(3,4)$. Let $\Delta_{+}\subset \Delta$ be the system of positive
roots:
$$
\Delta_{+}=\{\alpha_1, \alpha_2, \alpha_3, \alpha_1+\alpha_2,
\alpha_2+\alpha_3, \alpha_1+\alpha_2+\alpha_3 \}.
$$
Denote by $(\cdot|\cdot)$ the $\mathrm{W}$-invariant scalar product in
$\mathfrak{h}^*$ such that $(\alpha_i|\alpha_i)=2$.

The root $\alpha_2$ plays the role of the 'distinguished' root $\alpha_0$
(see Introduction). The associated element $H_0 \in \mathfrak{h}$ is given
by $H_0=\mathrm{diag}(1,1,-1,-1)$, or
\begin{equation}\label{uni9pproach.texh0}
H_0=H_1+2H_2+H_3
\end{equation}
with $H_1=\mathrm{diag}(1,-1,0,0)$, $H_2=\mathrm{diag}(0,1,-1,0)$,
$H_3=\mathrm{diag}(0,0,1,-1)$.

Let $\Delta_{\mathrm{c}}\stackrel{\mathrm{def}}=\{\alpha_1, \alpha_3,
-\alpha_1, -\alpha_3\}\subset \Delta$,
$\Delta_{\mathrm{n}}\stackrel{\mathrm{def}}= \Delta \setminus \Delta_
\mathrm{c}$. Then
$$
\mathfrak{g}_0=\mathfrak{h}\bigoplus \left (\bigoplus_{\alpha \in \Delta_
\mathrm{c}}\mathfrak{g}_{\alpha}\right),
$$
$$
\mathfrak{g}_{+1}=\bigoplus_{\alpha \in \Delta_+\cap \Delta_
\mathrm{n}}\mathfrak{g}_{\alpha},\qquad \mathfrak{g}_{-1}=\bigoplus_{-\alpha
\in \Delta_+\cap \Delta_ \mathrm{n}}\mathfrak{g}_{\alpha},
$$
with $\mathfrak{g}_{\alpha}$ being the root subspace in $\mathfrak{g}$
corresponding to $\alpha \in \Delta$.

Let $\mathrm{W}_\mathrm{c}$ be the subgroup in $\mathrm{W}$ generated by
$s_{\alpha_1}$, $s_{\alpha_3}$. Thus, $\mathrm{W}_\mathrm{c}\simeq
\mathrm{S}_2 \times \mathrm{S}_2$ is the Weyl group of the Lie subalgebra
$\mathfrak{sl}_2 \oplus \mathfrak{sl}_2 \in \mathfrak{g}$.

In the rest of this paper the ground field will be the field of rational
functions $\mathbb{C}(q^{1/4})$.

Let us recall one some definitions and facts of the quantum group theory (we
follow \cite{uni9pproach.texDC}).

The quantum universal enveloping algebra $U_q \mathfrak{g}$ is the algebra
with the generators $\{E_i,F_i,K^{\pm1}_i \}_{i=\overline{1,3}}$ satisfying
the following relations
$$K_iK_j=K_jK_i,\qquad K_iK_i^{-1}=K_i^{-1}K_i=1,$$
$$K_iE_j=q^{a_{ij}}E_jK_i,\qquad K_iF_j=q^{-a_{ij}}F_jK_i,$$
$$E_iF_j-F_jE_i=\delta_{ij}\frac{K_i-K_i^{-1}}{q-q^{-1}},$$
$$
\sum_{s=0}^{1-a_{ij}}(-1)^s
\genfrac{[}{]}{0pt}{0}{1-a_{ij}}{s}_qE_i^{1-a_{ij}-s}E_jE_i^s=0,\qquad i \ne
j,
$$
$$
\sum_{s=0}^{1-a_{ij}}(-1)^s
\genfrac{[}{]}{0pt}{0}{1-a_{ij}}{s}_qF_i^{1-a_{ij}-s}F_jF_i^s=0,\qquad i \ne
j,
$$
where $(a_{ij})$ is the Cartan matrix for $\mathfrak{g}$:
$$
a_{ij}=\left \{\begin{array}{ccl}2 &,& i-j=0 \\
                                 -1 &,& |i-j|=1 \\
                                 0 &,& \mathrm{otherwise} \end{array}\right.,
$$
$$
\genfrac{[}{]}{0pt}{0}{n}{j}_q
\stackrel{\mathrm{def}}=\frac{[n]_q!}{[n-j]_q![j]_q!},\qquad
[n]_q!\stackrel{\mathrm{def}}=[n]_q \cdot[n-1]_q \cdot \ldots
\cdot[1]_q,\qquad [n]_q \stackrel{\mathrm{def}}=\frac{q^n-q^{-n}}{q-q^{-1}}.
$$
The algebra $U_q \mathfrak{g}$ is endowed with a structure of a Hopf algebra
as follows
$$
\Delta(E_i)=E_i \otimes 1+K_i \otimes E_i,\quad \Delta(F_i)=F_i \otimes
K_i^{-1}+1 \otimes F_i,\quad \Delta(K_i)=K_i \otimes K_i,
$$
$$S(E_i)=-K_i^{-1}E_i,\qquad S(F_i)=-F_iK_i,\qquad S(K_i)=K_i^{-1},$$
$$\varepsilon(E_i)=\varepsilon(F_i)=0,\qquad \varepsilon(K_i)=1,$$
with $\Delta$, $S$, $\varepsilon$ being the comultiplication, the antipode,
and the counit, respectively.

Let us use the short notation $x_{(1)}\otimes x_{(2)}$ for the element
$\Delta(x)\in U_q \mathfrak{g}\bigotimes U_q \mathfrak{g}$ ($x \in U_q
\mathfrak{g}$). For example, coassociativity of the comultiplication
$\Delta:U_q \mathfrak{g}\rightarrow U_q \mathfrak{g}\bigotimes U_q
\mathfrak{g}$ looks in this notation as follows
\begin{equation}\label{uni9pproach.texcoas}
x_{(1)}\otimes x_{(2)(1)}\otimes x_{(2)(2)}= x_{(1)(1)}\otimes
x_{(1)(2)}\otimes x_{(2)}.
\end{equation}
Sometimes we use the notation $x_{(1)}\otimes x_{(2)}\otimes x_{(3)}$ for
the right (and left) hand side of (\ref{uni9pproach.texcoas}). Then the obvious meaning has
the notation $x_{(1)}\otimes x_{(2)}\otimes x_{(3)}\otimes x_{(4)}$ etc.

The adjoint representation of the algebra $U_q \mathfrak{g}$ is defined as
follows
$$\mathrm{ad}x(y)\stackrel{\mathrm{def}}=x_{(1)}\cdot y \cdot S(x_{(2)})$$
with $x,y \in U_q \mathfrak{g}$. This adjoint action makes $U_q
\mathfrak{g}$ a $U_q \mathfrak{g}$-module algebra. It means that the product
map $U_q \mathfrak{g}\bigotimes U_q \mathfrak{g}\rightarrow U_q
\mathfrak{g}$ is a morphism of $U_q \mathfrak{g}$-modules and the unit $1
\in U_q \mathfrak{g}$ is $U_q\mathfrak{g}$-invariant.

We fix the following notation for some subalgebras in $U_q \mathfrak{g}$:
$$
U_q^{\geq0}=\langle E_i,K^{\pm1}_i|\quad i=\overline{1,3}\rangle,\qquad
U_q^{\leq0}=\langle F_i,K^{\pm1}_i|\quad i=\overline{1,3}\rangle,
$$
$$
U_q^{+}=\langle E_i|\quad i=\overline{1,3}\rangle,\qquad U_q^{-}=\langle
F_i|\quad i=\overline{1,3}\rangle,
$$
$$
U_q^{0}=\langle K^{\pm1}_i|\quad i=\overline{1,3}\rangle, \qquad U_q
\mathfrak{g}_0=\langle K^{\pm1}_i, E_j, F_j|\quad i=\overline{1,3},\quad j
\ne2 \rangle.
$$

Let us recall a definition of the Lusztig automorphisms $T_i$,
$i=\overline{1,3}$, of the algebra $U_q \mathfrak{g}$. The action of $T_i$
on the subalgebra $U_q^{\geq0}$ is given by
$$T_i(K_j)=K_j \cdot K_i^{-a_{ij}},\qquad T_i(E_i)=-F_iK_i,$$
$$T_i(E_j)=(-\mathrm{ad}\,E_i)^{-a_{ij}}(E_j), \quad i \ne j.$$
To define $T_i$ completely one sets
$$T_i \circ \it{k}=\it{k}\circ T_i,$$
where $\it{k}$ is the conjugate linear antiautomorphism of the
$\mathbb{C}(q^{1/4})$-algebra $U_q \mathfrak{g}$ given by
$$
{\it k}(E_i)=F_i,\quad {\it k}(F_i)=E_i,\quad {\it k}(K_i)=K^{-1}_i,\quad
{\it k}(q^{1/4})=q^{-1/4}.
$$

Let $w \in \mathrm{W}$ and $w=s_{i_1}s_{i_2}\ldots s_{i_k}$ be a reduced
expression (we write '$s_{i}$' instead of $s_{\alpha_i}$). It is well known
that the automorphism $T_{w}\stackrel{\mathrm{def}}=T_{i_1}T_{i_2}\ldots
T_{i_k}$ does not depend on particular choice of a reduced expression of
$w$.

All $U_q\mathfrak{g}$-modules we consider possess the property
$$
V=\bigoplus_{\mu \in \mathbb{Z}^3} V_ \mu,\quad V_ \mu
\stackrel{\mathrm{def}}=\{v \in V|\quad K_iv=q^{\mu_i}v, i=\overline{1,3}\}
$$
with $\mu=(\mu_1, \mu_2, \mu_3)$. This allows one to introduce endomorphisms
$H_i$, $i=\overline{1,3}$, of any $U_q \mathfrak{g}$-module $V$ by
$$
H_iv=\mu_iv \qquad \Leftrightarrow \qquad v \in V_ \mu, \quad \mu=(\mu_1,
\mu_2, \mu_3).
$$
Formally this can be written as $K_i=q^{H_i}$.

Let $K_0 \stackrel{\mathrm{def}}=K_1 \cdot K_2^2 \cdot K_3$ (i.e.
$K_0=q^{H_0}$ with $H_0$ given by (\ref{uni9pproach.texh0})). It is an important
consequence of definitions that $K_0$ belongs to the centre of the algebra
$U_q \mathfrak{g}_0$:
\begin{equation}\label{uni9pproach.texho}
\mathrm{ad}K_0(\xi)=\xi,\qquad \xi \in U_q \mathfrak{g}_0.
\end{equation}

Let us recall one some facts concerning the universal $R$-matrix for $U_q
\mathfrak{g}$ (in context of the present paper the universal $R$-matrix have
to be understood just as in \cite{uni9pproach.texDC}).

$R$ satisfies some well known identities. We don't adduce a full list of
these identities but recall one those important for us:
\begin{equation}\label{uni9pproach.texprop1}
\mathrm{id}\otimes \Delta^{\mathrm{op}}(R)=R^{12}\cdot R^{13},
\end{equation}
\begin{equation}\label{uni9pproach.texprop2}
\Delta^{\mathrm{op}}(\eta)\cdot R=R \cdot \Delta(\eta),\quad \eta \in U_q
\mathfrak{g},
\end{equation}
with $\Delta^{\mathrm{op}}(x)\stackrel{\mathrm{def}}=x_{(2)}\otimes x_{(1)},
\quad R^{12}\stackrel{\mathrm{def}}=\sum_ia_i \otimes b_i \otimes1,\quad
R^{13}\stackrel{\mathrm{def}}=\sum_ia_i \otimes1 \otimes b_i$.

Remind an explicit formula for the $R$-matrix (the so-called multiplicative
formula). Let $w_0 \in \mathrm{W}$ be the maximal length element.
Identifying $\mathrm{W}$ with $\mathrm{S}_4$ we get
$$
w_0=\genfrac{(}{)}{0pt}{0}{1 \quad 2 \quad 3 \quad 4}{4 \quad 3 \quad 2
\quad 1}.
$$
The length of any reduced expression of $w_0$ is equal to $6$. To a reduced
expression \\$w_0=s_{i_1}s_{i_2}\ldots s_{i_6}$ one attaches the following
data:

i) the total order on the set $\Delta_+$ of positive roots:
$$
\beta_1=\alpha_{i_1},\quad \beta_2=s_{i_1}(\alpha_{i_2}),\quad \ldots \quad
\beta_6=s_{i_1}s_{i_2}\ldots s_{i_{5}}(\alpha_{i_6});
$$

ii) the set of elements $E_{\beta_1},E_{\beta_2},\ldots E_{\beta_6}\in
U_q^{+}$, $F_{\beta_1},F_{\beta_2},\ldots F_{\beta_6}\in U_q^{-}$ which are
$q$-analogs of root vectors in $\mathfrak{g}$:
$$
E_{\beta_1}=E_{i_1},\quad E_{\beta_2}=T_{i_1}(E_{i_2}),\quad \ldots \quad
E_{\beta_6}=T_{i_1}T_{i_2}\ldots T_{i_{5}}(E_{i_6}),
$$
$$
F_{\beta_1}=F_{i_1},\quad F_{\beta_2}=T_{i_1}(F_{i_2}),\quad \ldots \quad
F_{\beta_6}=T_{i_1}T_{i_2}\ldots T_{i_{5}}(F_{i_6});
$$

iii) the multiplicative formula for the $R$-matrix: \\$R=$
\begin{equation}\label{uni9pproach.texmult}
\mathrm{exp}_{q^2}((q^{-1}-q)E_{\beta_6}\otimes F_{\beta_6})\cdot \ldots
\cdot \mathrm{exp}_{q^2}((q^{-1}-q)E_{\beta_{2}}\otimes F_{\beta_{2}})\cdot
\mathrm{exp}_{q^2}((q^{-1}-q)E_{\beta_1}\otimes F_{\beta_1})\cdot
q^{\mathrm{t}}
\end{equation}
with $\mathrm{t}\stackrel{\mathrm{def}}=-\sum_{i,j}c_{ij}H_i \otimes H_j$,
the matrix $(c_{ij})$ being the inverse to the Cartan matrix, and
$$
\mathrm{exp}_{q^2}(t)\stackrel{\mathrm{def}}=\sum_{k=0}^{\infty}\frac{t^k}
{(k)_{q^2}!}, \quad (k)_{q^2}!\stackrel{\mathrm{def}}=\prod_{j=1}^k
\frac{1-q^{2j}}{1-q^2}.
$$

\bigskip
\section{Three approaches to quantization of prehomogeneous vector spaces of
commutative parabolic type} \subsection{First approach} In this subsection
we describe an approach to constructing $q$-analogs of prehomogeneous vector
spaces of commutative parabolic type developed in \cite{uni9pproach.texSV2}.

Let us consider the generalized Verma module $V(0)$ over $U_q\mathfrak{g}$
given by its generator $v(0)$ and the relations
\begin{equation}\label{uni9pproach.texv1}
E_iv(0)=0, \quad K_iv(0)=v(0),\quad i=\overline{1,3},
\end{equation}
\begin{equation}\label{uni9pproach.texv2}
F_iv(0)=0, \quad i\ne2.
\end{equation}
$V(0)$ splits into direct sum of its finite dimensional subspaces
$V(0)_{k}$, $-k\in\mathbb{Z}_+$, with
$$V(0)_{k}\stackrel{\mathrm{def}}=\{v \in V(0)|\quad H_0v=2kv \}.$$
Consider the {\sl graded} dual $U_q \mathfrak{g}$-module:
$$
\mathbb{C}[\mathfrak{g}_{-1}]_q \stackrel{\mathrm{def}}= \bigoplus_{-k \in
\mathbb{Z}_+}\left(V(0)_{k}\right)^*.
$$

Let us equip the tensor product $V(0)\bigotimes V(0)$ with a $U_q
\mathfrak{g}$-module structure as follows
\begin{equation}\label{uni9pproach.texmodstr}
\xi(v_1 \otimes v_2)=\xi_{(2)}(v_1)\otimes \xi_{(1)}(v_2), \quad \xi \in U_q
\mathfrak{g},\quad v_1,v_2 \in V(0).
\end{equation}
Due to (\ref{uni9pproach.texv1}), (\ref{uni9pproach.texv2}) the maps
\begin{equation}\label{uni9pproach.texcoalstr}
v(0)\mapsto v(0)\otimes v(0),\qquad v(0)\mapsto1
\end{equation}
are extendable up to morphisms of $U_q \mathfrak{g}$-modules
$$
\Delta_-:V(0)\rightarrow V(0)\bigotimes V(0),\qquad
\varepsilon_-:V(0)\rightarrow \mathbb{C}(q^{1/4}).
$$
It can be shown that $\Delta_-$ and $\varepsilon_-$ make $V(0)$ a
coassociative coalgebra with a counit. Thus, the dual maps
\begin{equation}\label{uni9pproach.texalstr}
m=(\Delta_-)^*:\mathbb{C}[\mathfrak{g}_{-1}]_q \bigotimes
\mathbb{C}[\mathfrak{g}_{-1}]_q \rightarrow
\mathbb{C}[\mathfrak{g}_{-1}]_q,\qquad
{1}=(\varepsilon_-)^*:\mathbb{C}(q^{1/4})\rightarrow
\mathbb{C}[\mathfrak{g}_{-1}]_q
\end{equation}
make $\mathbb{C}[\mathfrak{g}_{-1}]_q$ an associative unital algebra.
Moreover, the product map $m$ is a morphism of $U_q \mathfrak{g}$-modules
and the unit $1$ is $U_q \mathfrak{g}$-invariant, i.e.
$\mathbb{C}[\mathfrak{g}_{-1}]_q$ is a $U_q \mathfrak{g}$-module algebra. In
particular, it is a $U_q \mathfrak{g}_0$-module algebra. This $U_q
\mathfrak{g}_0$-module structure is just the one we mentioned in the
Introduction.

\subsection{Second approach}

Now we are going to describe briefly an approach of H. P. Jakobsen \cite{uni9pproach.texJ}
to quantization of $\mathbb{C}[\mathfrak{g}_{-1}]$.

It follows from the definition of $U_q\mathfrak{g}$ that
$$
\mathrm{ad}K_2(E_2)=q^2E_2,\qquad
\mathrm{ad}K_1(E_2)=\mathrm{ad}K_3(E_2)=q^{-1}E_2,
$$
$$\mathrm{ad}F_1(E_2)=\mathrm{ad}F_3(E_2)=0,$$
$$(\mathrm{ad}E_1)^2(E_2)=(\mathrm{ad}E_3)^2(E_2)=0.$$
Thus $\mathrm{ad}U_q \mathfrak{g}_0(E_2)$ is a finite dimensional $U_q
\mathfrak{g}_0$-submodule in $U_q \mathfrak{g}$. Let us denote by
$\mathbb{C}[\mathfrak{g}_{-1}]^{\mathrm{I}}_q$ the minimal subalgebra in
$U_q \mathfrak{g}$ which contains the subspace $\mathrm{ad}U_q
\mathfrak{g}_0(E_2)$. Evidently, it is a $U_q \mathfrak{g}_0$-module
subalgebra in $U_q \mathfrak{g}$. The algebra
$\mathbb{C}[\mathfrak{g}_{-1}]^{\mathrm{I}}_q$ can be treated as a
$q$-analog of $\mathbb{C}[\mathfrak{g}_{-1}]$ (see Remark 2 above).

\subsection{Third approach} Let us turn to description of an approach of A.
Kamita, Y. Morita, and T. Tanisaki \cite{uni9pproach.texKMT}. Note that notation in
\cite{uni9pproach.texKMT} differs from ours.

Let $w'_0 \in \mathrm{W}_{\mathrm{c}}$ be the maximal length element.
Evidently, $w'_0=(1,2)\cdot(3,4)$. Consider the subspace
$\mathbb{C}[\mathfrak{g}_{-1}]^{\mathrm{II}}_q$ in $U_q \mathfrak{g}$
defined by
\begin{equation}\label{uni9pproach.texii}
\mathbb{C}[\mathfrak{g}_{-1}]^{\mathrm{II}}_q=U_q^{+}\cap
T^{-1}_{w'_0}\left(U_q^{+}\right).
\end{equation}

Obviously the subspace $\mathbb{C}[\mathfrak{g}_{-1}]^{\mathrm{II}}_q$ is a
subalgebra in $U_q \mathfrak{g}$. It is shown in \cite{uni9pproach.texKMT} that
$\mathbb{C}[\mathfrak{g}_{-1}]^{\mathrm{II}}_q$ is a $U_q
\mathfrak{g}_0$-module subalgebra in $U_q \mathfrak{g}$ with respect to the
adjoint action. It is one more $q$-analog of the algebra
$\mathbb{C}[\mathfrak{g}_{-1}]$.

\bigskip

\section{Comparison of the approaches}

Note that both the algebras $\mathbb{C}[\mathfrak{g}_{-1}]^{\mathrm{I}}_q$
and $\mathbb{C}[\mathfrak{g}_{-1}]^{\mathrm{II}}_q$ lie within the quantum
universal enveloping algebra $U_q \mathfrak{g}$. Our aim is:

i) to construct an embedding $T$ of the algebra
$\mathbb{C}[\mathfrak{g}_{-1}]_q$ into $U_q \mathfrak{g}$ which intertwines
the $U_q \mathfrak{g}_0$-action in $\mathbb{C}[\mathfrak{g}_{-1}]_q$
mentioned in subsection 3.1 and the adjoint $U_q \mathfrak{g}_0$-action;

ii) to show that the subalgebras $T(\mathbb{C}[\mathfrak{g}_{-1}]_q)$,
$\mathbb{C}[\mathfrak{g}_{-1}]^{\mathrm{I}}_q$, and
$\mathbb{C}[\mathfrak{g}_{-1}]^{\mathrm{II}}_q$ in $U_q \mathfrak{g}$
coincide.

\bigskip

\subsection{$\mathbb{C}[\mathfrak{g}_{-1}]_q \simeq
\mathbb{C}[\mathfrak{g}_{-1}]^{\mathrm{I}}_q$}

In this subsection we construct an $U_q \mathfrak{g}_0$-invariant embedding
$T$ of the algebra $\mathbb{C}[\mathfrak{g}_{-1}]_q$ into $U_q
\mathfrak{g}$, and then we show that $T(\mathbb{C}[\mathfrak{g}_{-1}]_q)=
\mathbb{C}[\mathfrak{g}_{-1}]^{\mathrm{I}}_q$. The embedding $T$ is
constructed via a standard technique due to \cite{uni9pproach.texD}.

Let $R=\sum_ia_i \otimes b_i$ be the universal $R$-matrix for $U_q
\mathfrak{g}$. Consider the linear map $T:\mathbb{C}[\mathfrak{g}_{-1}]_q
\rightarrow U_q \mathfrak{g}$ given by
\begin{equation}\label{uni9pproach.texT}
T(f)\stackrel{\mathrm{def}}=\sum_ia_i \langle b_iv(0),f \rangle
\end{equation}
with $f \in \mathbb{C}[\mathfrak{g}_{-1}]_q$, $v(0)$ being the generator of
the generalized Verma module $V(0)$, $\langle \cdot,\cdot \rangle$ being the
pairing $V(0)\times \mathbb{C}[\mathfrak{g}_{-1}]_q \rightarrow
\mathbb{C}(q^{1/4})$ arising from the equality
$\mathbb{C}[\mathfrak{g}_{-1}]_q=\left(V(0)\right)^*$ (see subsection 3.1).

Let us comment the definition of $T$. Using the multiplicative formula
(\ref{uni9pproach.texmult}) and the definition of the $U_q \mathfrak{g}$-module $V(0)$ one
shows that $\sum_ia_i \otimes b_iv(0)$ is a formal series of elements from
$U^{+}_q \bigotimes V(0)$. We will prove later (see proof of Proposition
\ref{uni9pproach.tex3}) that the right hand side of (\ref{uni9pproach.texT}) is a finite sum for any $f
\in \mathbb{C}[\mathfrak{g}_{-1}]_q$.

\medskip

\begin{proposition}\label{uni9pproach.texhom}
 $T$ is a homomorphism of algebras.
\end{proposition}

\smallskip

{\bf Proof.} Let $f,\varphi \in \mathbb{C}[\mathfrak{g}_{-1}]_q$. Then due
to (\ref{uni9pproach.texalstr}), (\ref{uni9pproach.texcoalstr}), (\ref{uni9pproach.texmodstr})
\begin{multline*}
T(f \cdot \varphi)=\sum_ia_i \langle b_iv(0),f \cdot \varphi
\rangle=\sum_ia_i \langle b_iv(0),m(f \otimes \varphi)\rangle=
\\ =\sum_ia_i \langle \Delta_-(b_iv(0)),f \otimes \varphi \rangle=\sum_ia_i
\langle \Delta^{\mathrm{op}}(b_i)(v(0)\otimes v(0)),f \otimes \varphi
\rangle.
\end{multline*}
By (\ref{uni9pproach.texprop1})
\begin{multline*}
\sum_ia_i \langle \Delta^{\mathrm{op}}(b_i)(v(0)\otimes v(0)),f \otimes
\varphi \rangle=\sum_{i,j}a_i \cdot a_j \langle b_i \otimes b_j (v(0)\otimes
v(0)),f \otimes \varphi \rangle=
\\ =\left(\sum_ia_i \langle b_iv(0),f \rangle \right)\cdot \left(\sum_ja_j
\langle b_iv(0),\varphi \rangle \right)=T(f)\cdot T(\varphi). \hfill \square
\end{multline*}

\pagebreak

\begin{proposition}\label{uni9pproach.texinter}
$T$ is a morphism of $U_q \mathfrak{g}_0$-modules.
\end{proposition}
\smallskip

{\bf Proof.} Denote by $U_q \mathfrak{g}\widehat{\otimes} V(0)$ the vector
space of formal series of the form
$$
\sum_i \xi_i \otimes v_i,\qquad \xi_i \otimes v_i \in U_q
\mathfrak{g}\otimes V(0).
$$
Equip this vector space with a $U_q\mathfrak{g}$-module structure as follows
$$\xi(\eta \otimes v)\stackrel{\mathrm{def}}=\xi_{(3)}\cdot \eta \cdot
S^{-1}(\xi_{(1)})\otimes \xi_{(2)}v,\quad \xi,\eta \in U_q
\mathfrak{g},\quad v \in V(0).
$$

\medskip

\begin{lemma}\label{uni9pproach.texl1} The linear map
$$
V(0)\rightarrow U_q \mathfrak{g}\widehat{\otimes} V(0), \quad v \mapsto
\sum_ia_i \otimes b_iv,
$$
is a morphism of $U_q \mathfrak{g}$-modules (here $\sum_ia_i \otimes b_i$ is
the universal $R$-matrix).
\end{lemma}

\smallskip

{\bf Proof of the lemma.} The map is well defined (it can be explained using
the definition (\ref{uni9pproach.texv1}), (\ref{uni9pproach.texv2}) of the $U_q \mathfrak{g}$-module
$V(0)$ and the multiplicative formula (\ref{uni9pproach.texmult}) for the universal
$R$-matrix). Let $\xi \in U_q \mathfrak{g}$, $v \in V(0)$. Then
\begin{multline*}
\xi \left(\sum_ia_i \otimes b_iv \right)=\sum_i \xi_{(3)}a_i
S^{-1}(\xi_{(1)})\otimes \xi_{(2)}b_iv=\sum_i
\varepsilon(\xi_{(3)(2)})\xi_{(3)(1)}a_i S^{-1}(\xi_{(1)})\otimes
\xi_{(2)}b_iv=
\\ =\sum_i \xi_{(3)(1)}a_i S^{-1}(\xi_{(1)})\otimes \xi_{(2)}b_i
\varepsilon(\xi_{(3)(2)})v= \sum_i \xi_{(2)(2)}a_i S^{-1}(\xi_{(1)})\otimes
\xi_{(2)(1)}b_i \varepsilon(\xi_{(3)})v.
\end{multline*}
Let us make use of the property (\ref{uni9pproach.texprop2}). We get
\begin{multline*}
\xi \left(\sum_ia_i \otimes b_iv \right)=\sum_i \xi_{(2)(2)}a_i
S^{-1}(\xi_{(1)})\otimes \xi_{(2)(1)}b_i \varepsilon(\xi_{(3)})v=
\\ =\sum_ia_i \xi_{(2)(1)}S^{-1}(\xi_{(1)})\otimes b_i
\xi_{(2)(2)}\varepsilon(\xi_{(3)})v= \sum_ia_i
\xi_{(1)(2)}S^{-1}(\xi_{(1)(1)})\otimes b_i
\xi_{(2)(1)}\varepsilon(\xi_{(2)(2)})v=
\\ =\sum_ia_i \varepsilon(\xi_{(1)})\otimes b_i \xi_{(2)}v=\sum_ia_i \otimes
b_i \varepsilon(\xi_{(1)})\xi_{(2)}v=\sum_ia_i \otimes b_i \xi v.
\end{multline*}
\hfill $\square$

\medskip

\begin{lemma}\label{uni9pproach.texl2}
An element $\eta \otimes v \in U_q \mathfrak{g}\widehat{\otimes} V(0)$ is
$U_q \mathfrak{g}_0$-invariant iff for any $\xi \in U_q \mathfrak{g}_0$
\begin{equation}\label{uni9pproach.texinvar}
\xi_{(1)}\eta S(\xi_{(2)})\otimes v=\eta \otimes S(\xi)v.
\end{equation}
\end{lemma}

\smallskip

{\bf Proof of the lemma.} Let $\eta \otimes v \in U_q
\mathfrak{g}\widehat{\otimes} V(0)$ satisfies (\ref{uni9pproach.texinvar}) for any $\xi \in
U_q \mathfrak{g}_0$. Rewrite (\ref{uni9pproach.texinvar}) for $\xi:=S^{-1}(\zeta)$:
\begin{equation}\label{uni9pproach.texinvari}
S^{-1}(\zeta_{(2)})\eta \zeta_{(1)}\otimes v=\eta \otimes \zeta v.
\end{equation}
Using (\ref{uni9pproach.texinvari}) one gets
\begin{multline*}
\xi(\eta \otimes v)=\xi_{(3)}\eta S^{-1}(\xi_{(1)})\otimes
\xi_{(2)}v=\xi_{(3)}S^{-1}(\xi_{(2)(2)})\eta
\xi_{(2)(1)}S^{-1}(\xi_{(1)})\otimes v=
\\ =\xi_{(2)(2)}S^{-1}(\xi_{(2)(1)})\eta
\xi_{(1)(2)}S^{-1}(\xi_{(1)(1)})\otimes v=\varepsilon(\xi_{(2)})\eta
\varepsilon(\xi_{(1)})\otimes v=\varepsilon(\xi)\eta \otimes v.
\end{multline*}
Thus $\eta \otimes v$ is $U_q \mathfrak{g}_0$-invariant.

Conversely, suppose that $\eta \otimes v$ is $U_q \mathfrak{g}_0$-invariant.
Let us prove (\ref{uni9pproach.texinvari}) (obviously, it is equivalent to (\ref{uni9pproach.texinvar})).
$$
S^{-1}(\xi_{(2)})\eta \xi_{(1)}\otimes v=S^{-1}(\xi_{(2)})\eta
\varepsilon(\xi_{(1)(2)})\xi_{(1)(1)}\otimes v=S^{-1}(\xi_{(3)})\eta
\varepsilon(\xi_{(2)})\xi_{(1)}\otimes v.
$$
$U_q \mathfrak{g}_0$-invariance of $\eta \otimes v$ implies
$$
\varepsilon(\xi_{(2)})\eta \otimes v=\xi_{(2)(3)}\eta
S^{-1}(\xi_{(2)(1)})\otimes \xi_{(2)(2)}v.
$$
Thus
\begin{multline*}
S^{-1}(\xi_{(2)})\eta \xi_{(1)}\otimes v=S^{-1}(\xi_{(3)})\xi_{(2)(3)}\eta
S^{-1}(\xi_{(2)(1)})\xi_{(1)}\otimes \xi_{(2)(2)}v=
\\ S^{-1}(\xi_{(3)(2)})\xi_{(3)(1)}\eta
S^{-1}(\xi_{(1)(2)})\xi_{(1)(1)}\otimes
\xi_{(2)}v=\varepsilon(\xi_{(3)})\eta \varepsilon(\xi_{(1)})\otimes
\xi_{(2)}v=\eta \otimes \xi v.
\end{multline*}
\hfill $\square$

\bigskip

Let us complete the proof of Proposition \ref{uni9pproach.texinter}. By (\ref{uni9pproach.texv1}),
(\ref{uni9pproach.texv2}) $v(0)$ is $U_q \mathfrak{g}_0$-invariant. Due to Lemma \ref{uni9pproach.texl1}
the element $\sum_ia_i \otimes b_iv(0)\in U_q \mathfrak{g}\widehat{\otimes}
V(0)$ is $U_q \mathfrak{g}_0$-invariant. By (\ref{uni9pproach.texinvar}) we have: for $f
\in \mathbb{C}[\mathfrak{g}_{-1}]_q$, $\xi \in U_q \mathfrak{g}_0$
\begin{multline*}
T(\xi f)=\sum_ia_i \langle b_iv(0),\xi f \rangle=\sum_ia_i \langle
S(\xi)b_iv(0), f \rangle=
\\ =\sum_i \xi_{(1)}a_iS(\xi_{(2)})\langle b_iv(0),f
\rangle=\mathrm{ad}\xi(T(f)).
\end{multline*}
\hfill $\square$

\bigskip

We have constructed the mapping $T: \mathbb{C}[\mathfrak{g}_{-1}]_q
\rightarrow U_q \mathfrak{g}$ which is a morphism of $U_q
\mathfrak{g}_0$-module algebras. It turns out to be an embedding.

\medskip
\begin{proposition}\label{uni9pproach.tex3}
$T$ is injective.
\end{proposition}

\smallskip

{\bf Proof.} Let $w_0 \in \mathrm{W}$ be the maximal length element. We fix
a reduced expression for $w_0$:
\begin{equation}\label{uni9pproach.texred}
w_0=(1,2)(3,4)(2,3)(1,2)(3,4)(2,3)=s_1s_3s_2s_1s_3s_2.
\end{equation}
Obviously
$
w_0=w'_0s_2s_1s_3s_2
$
with $w'_0$ being the maximal length element in $\mathrm{W}_\mathrm{c}$.
Describe explicitly the order in $\Delta_+$ attached to the expression
(\ref{uni9pproach.texred}) (see section 2):
$$
\beta_1=\alpha_1, \qquad \beta_2=\alpha_3,\qquad
\beta_3=\alpha_1+\alpha_2+\alpha_3,
$$
$$
\beta_4=\alpha_2+\alpha_3,\qquad \beta_5=\alpha_1+\alpha_2, \qquad
\beta_6=\alpha_2.
$$
Thus $\beta_3$, $\beta_4$, $\beta_5$, $\beta_6$ are 'noncompact' roots (they
belong to $\Delta_ \mathrm{n}$), and $\beta_1$, $\beta_2$ are 'compact'
roots (they belong to $\Delta_ \mathrm{c}$). This implies $E_{\beta_1},
F_{\beta_1}, E_{\beta_2}, F_{\beta_2}\in U_q\mathfrak{g}_0$. Indeed, let us
show this, for example, for $E_{\beta_2}$. Consider the
$\mathbb{Z}_+$-grading in $U_q^+$ given by
\begin{equation}\label{uni9pproach.texgrad1}
\deg E_i=\left \{\begin{array}{ccl}1 &,& i=2 \\
                                 0 &,& \mathrm{otherwise} \end{array}\right..
\end{equation}
This grading can be described in another way:
\begin{equation}\label{uni9pproach.texgrad2}
\left(U_q^{+}\right)_j=\{\xi \in U_q^{+}|\quad
\mathrm{ad}\,K_0(\xi)=q^{2j}\xi \}.
\end{equation}
Using (\ref{uni9pproach.texgrad1}) one shows that
\begin{equation}\label{uni9pproach.texu0}
\left(U_q^{+}\right)_0=U_q^{+}\cap U_q \mathfrak{g}_0.
\end{equation}
It follows from the definition of $E_{\beta_2}$ and $K_0$ that
$$\mathrm{ad}K_0(E_{\beta_2})=q^{\beta_2(H_0)}E_{\beta_2}=E_{\beta_2}.$$
Thus, according to (\ref{uni9pproach.texgrad2}) and (\ref{uni9pproach.texu0}) $E_{\beta_2}\in U_q
\mathfrak{g}_0$.

It follows from the relations (\ref{uni9pproach.texv1}) and (\ref{uni9pproach.texv2}) that
$$\sum_ia_i \otimes b_iv(0)=\sum_jA_j \otimes B_jv(0),$$
where $\sum_ia_i \otimes b_i$ is the universal $R$-matrix, \\ $\sum_jA_j
\otimes B_j
\stackrel{\mathrm{def}}=\mathrm{exp}_{q^2}((q^{-1}-q)E_{\beta_6}\otimes
F_{\beta_6})\cdot \mathrm{exp}_{q^2}((q^{-1}-q)E_{\beta_{5}}\otimes
F_{\beta_{5}})\times$
\begin{equation}\label{uni9pproach.texab}
\times \mathrm{exp}_{q^2}((q^{-1}-q)E_{\beta_4}\otimes F_{\beta_4})\cdot
\mathrm{exp}_{q^2}((q^{-1}-q)E_{\beta_3}\otimes F_{\beta_3}).
\end{equation}
It is clear that
\begin{equation}\label{uni9pproach.texab1}
\sum_jA_j \otimes B_j=\sum_{(k_1,\ldots, k_4)\in
\mathbb{Z}_+^4}a_{k_1,\ldots, k_4}E^{k_1}_{\beta_6}\ldots
E^{k_4}_{\beta_3}\otimes F^{k_1}_{\beta_6}\ldots F^{k_4}_{\beta_3},
\end{equation}
where all $a_{k_1,\ldots, k_4}$ are nonzero elements of
$\mathbb{C}(q^{1/4})$. Thus we get the formula
\begin{equation}\label{uni9pproach.texform}
 T(f)=\sum_{(k_1,\ldots, k_4)\in \mathbb{Z}_+^4}a_{k_1,\ldots,
k_4}E^{k_1}_{\beta_6}\ldots E^{k_4}_{\beta_3}\langle F^{k_1}_{\beta_6}\ldots
F^{k_4}_{\beta_3}v(0),f \rangle.
\end{equation}
To complete the proof of Proposition \ref{uni9pproach.tex3} it is sufficient to prove that
the vectors
$$
\{E^{k_1}_{\beta_6}\ldots E^{k_4}_{\beta_3}\}_{(k_1,\ldots, k_4)\in
\mathbb{Z}_+^4}
$$
are linearly independent, and the vectors
$$
\{F^{k_1}_{\beta_6}\ldots F^{k_4}_{\beta_3}v(0)\}_{(k_1,\ldots, k_4)\in
\mathbb{Z}_+^4}
$$
constitute a basis in $V(0)$ (the latter statement implies also the well
definiteness of $T$!). For this purpose we need the following Theorem
\cite[page 14]{uni9pproach.texDC}:

\medskip

\begin{theorem}\label{uni9pproach.texdc} i) The monomials
$$
\{F^{k_1}_{\beta_6}F^{k_2}_{\beta_5}\ldots
F^{k_6}_{\beta_1}K^{m_1}_{1}K^{m_2}_{2}
K^{m_3}_{3}E^{l_1}_{\beta_1}E^{l_2}_{\beta_2}\ldots
E^{l_6}_{\beta_6}\}_{(k_1,\ldots, k_6)\in \mathbb{Z}_+^6, (m_1,m_2, m_3)\in
\mathbb{Z}^3,(l_1,\ldots, l_6)\in \mathbb{Z}_+^6}
$$
constitute a basis in $U_q \mathfrak{g}$;

\medskip

ii) for $i<j$ one has:
$$
E_{\beta_i}E_{\beta_j}-q^{(\beta_i|\beta_j)} E_{\beta_j}E_{\beta_i}=
 \sum_{(k_1,\ldots, k_6)\in \mathbb{Z}_+^6}a_{k_1,\ldots,
k_6}E^{k_1}_{\beta_1}E^{k_2}_{\beta2}\ldots E^{k_6}_{\beta_6},
$$
where $a_{k_1,\ldots, k_6}\ne0$ only when $k_s=0$ for $s \leq i$ or $s \geq
j$.
\end{theorem}

\smallskip

It is not hard to prove that this Theorem implies linear independence of the
vectors $\{E^{k_1}_{\beta_6}\ldots E^{k_4}_{\beta_3}\}_{(k_1,\ldots, k_4)\in
\mathbb{Z}_+^4}$. Let us prove that the vectors $\{F^{k_1}_{\beta_6}\ldots
F^{k_4}_{\beta_3}v(0)\}_{(k_1,\ldots, k_4)\in \mathbb{Z}_+^4}$ constitute a
basis in $V(0)$.

Obviously, $V(0)$ is the linear span of $\{F^{k_1}_{\beta_6}\ldots
F^{k_4}_{\beta_3}v(0)\}_{(k_1,\ldots, k_4)\in \mathbb{Z}_+^4}$. Indeed, this
follows from statement i) of the Theorem, from the observations that the map
$U_q\mathfrak{g}\rightarrow V(0)$, $\xi \mapsto \xi v(0)$ is surjective, and
from the relations
$$E_{\beta_i}v(0)=0,\quad i=\overline{1,6},$$
$$K_iv(0)=v(0),\quad i=\overline{1,3},$$
$$F_{\beta_1}v(0)=F_{\beta_2}v(0)=0.$$
What remains is to prove that $\{F^{k_1}_{\beta_6}\ldots
F^{k_4}_{\beta_3}v(0)\}_{(k_1,\ldots, k_4)\in \mathbb{Z}_+^4}$ are linearly
independent. We prove this statement using its correctness for $q=1$.

Let $A$ be the ring $\mathbb{C}[q,q^{-1}]$. Consider the $A$-algebra
$U^-_{q,A}$ generated by $\{F^A_i \}_{i=\overline{1,3}}$ which satisfy the
same relations as $\{F_i \}_{i=\overline{1,3}}$. Evidently, as
$\mathbb{C}(q^{1/4})$-algebras
\begin{equation}\label{uni9pproach.texisom}
U_q^{-}\simeq \mathbb{C}(q^{1/4})\bigotimes_A U^{-}_{q,A}.
\end{equation}
Let $\{F^A_{\beta_i}\}_{i=\overline{1,6}}$ be the elements of $U^-_{q,A}$
derived from $\{F_{\beta_i}\}_{i=\overline{1,6}}$ via substitution $F_i
\leftrightarrow F^A_i$. Consider the $U^{-}_{q,A}$-module $V(0)^A$ given by
the generator $v(0)^A$ and the relations
$$F^A_iv(0)^A=0, \qquad i \ne2.$$
Similarly to (\ref{uni9pproach.texisom})
\begin{equation}\label{uni9pproach.texisomo}
V(0)\simeq \mathbb{C}(q^{1/4})\bigotimes_A V(0)^A_q
\end{equation}
as $\mathbb{C}(q^{1/4})$-modules.

There is an evident homomorphism of $\mathbb{C}$-algebras
$$J:U^-_{q,A}\rightarrow U^-_{q,A}/(q-1)\cdot U^-_{q,A}\simeq U^-,$$
where $U^-$ is the subalgebra in the classical universal enveloping algebra
$U \mathfrak{g}$ generated by $\{J(F^A_i)\}_{i=\overline{1,3}}$. Denote by
$\widetilde{V(0)}$ a $U^-$-module given by the generator $\widetilde{v(0)}$
and the relations
$$J(F^A_i)\widetilde{v(0)}=0, \qquad i \ne2.$$
It is clear that the map $v(0)^A \mapsto \widetilde{v(0)}$ can be extended
up to a $\mathbb{C}$-linear map
$$J_0: V(0)^A \rightarrow \widetilde{V(0)}$$
such that for any $\xi \in U^-_{q,A}$
\begin{equation}\label{uni9pproach.texlimit}
J(\xi)\widetilde{v(0)}=J_0(\xi v(0)^A).
\end{equation}
It is well known that the vectors $\{J(F^A_{\beta_6})^{k_1}\ldots
J(F^A_{\beta_3})^{k_4}\widetilde{v(0)}\}_{(k_1,\ldots, k_4)\in
\mathbb{Z}_+^4}$ constitute a basis in the $\mathbb{C}$-module
$\widetilde{V(0)}$. Thus by (\ref{uni9pproach.texlimit}) the vectors
$\{(F^A_{\beta_6})^{k_1}\ldots (F^A_{\beta_3})^{k_4}v(0)^A \}_{(k_1,\ldots,
k_4)\in \mathbb{Z}_+^4}$ are linearly independent. Due to (\ref{uni9pproach.texisomo})
$\{F_{\beta_6}^{k_1}\ldots F_{\beta_3}^{k_4}v(0)\}_{(k_1,\ldots, k_4)\in
\mathbb{Z}_+^4}$ are linearly independent.

We have completed the proof of Proposition \ref{uni9pproach.tex3}. $\hfill$ $\square$

\bigskip

The following statement is a direct consequence of Propositions \ref{uni9pproach.texhom},
\ref{uni9pproach.tex3} and formula (\ref{uni9pproach.texform}).

\medskip

\begin{corollary}\label{uni9pproach.texcor}
Linear span of $\{E^{k_1}_{\beta_6}\ldots E^{k_4}_{\beta_3}\}_{(k_1,\ldots,
k_4)\in \mathbb{Z}_+^4}$ is a $\mathrm{ad}\,U_q \mathfrak{g}_0$-invariant
subalgebra in $U_q^+$. It coincides with
$T(\mathbb{C}[\mathfrak{g}_{-1}]_q)$.
\end{corollary}

\medskip

{\sc Remark 3.} The fact that linear span of $\{E^{k_1}_{\beta_6}\ldots
E^{k_4}_{\beta_3}\}_{(k_1,\ldots, k_4)\in \mathbb{Z}_+^4}$ is a subalgebra
in $U_q^+$ easily follows from the statement ii) of the Theorem \ref{uni9pproach.texdc}.

\smallskip

The main statement of this section is
\begin{proposition}\label{uni9pproach.texpropos}
$$
T(\mathbb{C}[\mathfrak{g}_{-1}]_q)=\mathbb{C}[\mathfrak{g}_{-1}]^{\mathrm
I}_q.
$$
\end{proposition}

\smallskip

{\bf Proof.} Let us start with

\smallskip

\begin{lemma}\label{uni9pproach.texlem}
For some $c \in \mathbb{C}(q^{1/4})$
$$E_{\beta_6}=c \cdot E_2.$$
\end{lemma}

\smallskip

{\bf Proof of the Lemma.} The elements $E_{\beta_6}$, $E_2$ of the
$U^0_q$-module $U_q^+$ are weight vectors of the weight $\alpha_2 \in
\mathfrak{h}^*$. But the subspace in $U_q^+$ of weight vectors of that
weight is $1$-dimensional (this follows from linear independence of the
weights $\alpha_1, \alpha_2, \alpha_3 \in \mathfrak{h}^*$ of the generators
$E_1$, $E_2$, $E_3$). \hfill $\square$

\medskip

Remind that we have equipped $U_q^+$ with a $\mathbb{Z}_+$-grading (see
(\ref{uni9pproach.texgrad2})). Obviously,
$$
\deg(E_{\beta_3})=\deg(E_{\beta_4})=\deg(E_{\beta_5})=\deg(E_{\beta_6})=1.
$$
It follows from (\ref{uni9pproach.texho}), (\ref{uni9pproach.texgrad2}) that endomorphisms from
$\mathrm{ad}U_q \mathfrak{g}_0$ preserve this grading: for $\xi \in U_q
\mathfrak{g}_0$ and $\eta \in U_q^+$
$$\deg(\mathrm{ad}\xi(\eta))=\deg(\eta)$$
provided $\mathrm{ad}\xi(\eta)\in U_q^+$. Using this observation and
Corollary \ref{uni9pproach.texcor} we get
\begin{equation}\label{uni9pproach.texincl}
\mathrm{ad}U_q \mathfrak{g}_0(E_{\beta_6})\subseteq \text{linear
span}\{E_{\beta_3}, E_{\beta_4}, E_{\beta_5}, E_{\beta_6}\}.
\end{equation}
Actually the spaces in the both sides of (\ref{uni9pproach.texincl}) coincide: dimension of
$\mathrm{ad}U_q\mathfrak{g}_0(E_{\beta_6})$ should be equal to $4$ just as
in the classical case $q=1$. Thus, by Lemma \ref{uni9pproach.texlem} and by the definition
of the algebra $\mathbb{C}[\mathfrak{g}_{-1}]^{\mathrm I}_q$
$$
\mathbb{C}[\mathfrak{g}_{-1}]^{\mathrm I}_q=\langle
E_{\beta_3},E_{\beta_4},E_{\beta_5},E_{\beta_6}\rangle.
$$
What remains is to use Corollary \ref{uni9pproach.texcor}. \hfill $\square$

\bigskip

\subsection{$\mathbb{C}[\mathfrak{g}_{-1}]^{\mathrm{I}}_q=
\mathbb{C}[\mathfrak{g}_{-1}]^{\mathrm{II}}_q$}

In this subsection we use notation of the previous one.

\begin{proposition}\label{uni9pproach.texpr}
$$
\mathbb{C}[\mathfrak{g}_{-1}]^{\mathrm{I}}_q \subseteq
\mathbb{C}[\mathfrak{g}_{-1}]^{\mathrm{II}}_q.
$$
\end{proposition}

\smallskip

{\bf Proof.} By Corollary \ref{uni9pproach.texcor} and Proposition \ref{uni9pproach.texpropos} the linear
span of $\{E^{k_1}_{\beta_6}\ldots E^{k_4}_{\beta_3}\}_{(k_1,\ldots, k_4)\in
\mathbb{Z}_+^4}$ coincides with
$\mathbb{C}[\mathfrak{g}_{-1}]^{\mathrm{I}}_q$. Thus, due to the definition
(\ref{uni9pproach.texii}) of $\mathbb{C}[\mathfrak{g}_{-1}]^{\mathrm{II}}_q$ it is
sufficient to prove that
$$T^{-1}_{w'_0}(E_{\beta_k})\in U_q^+,\qquad k=\overline{3,6}.$$
Let us prove this, for example, for $E_{\beta_6}$. By definition
$$E_{\beta_6}=T_1T_3T_2T_1T_3(E_2)=T_{w'_0}T_2T_1T_3(E_2).$$
One gets
$$T^{-1}_{w'_0}(E_{\beta_6})=T_2T_1T_3(E_2).$$
It remains to make use of the following well known fact.

\smallskip

\begin{lemma}\label{uni9pproach.texw}
If $w(\alpha_i)\in \Delta_+$ for some $i=\overline{1,3}$ and $w \in
\mathrm{W}$ then $T_w(E_i)\in U_q^+.$
\end{lemma}
\hfill $\square$

\medskip

Now we are ready to prove

\begin{proposition}\label{uni9pproach.texpro}
$$
\mathbb{C}[\mathfrak{g}_{-1}]^{\mathrm{I}}_q=
\mathbb{C}[\mathfrak{g}_{-1}]^{\mathrm{II}}_q.
$$
\end{proposition}

\medskip

{\bf Proof.} Let $\sum a_{k_1,\ldots,k_6}E^{k_1}_{\beta_6}\ldots
E^{k_6}_{\beta_1}\in \mathbb{C}[\mathfrak{g}_{-1}]^{\mathrm{II}}_q,$ i.e.
\begin{equation}\label{uni9pproach.texu+}
T^{-1}_{w'_0}\left(\sum a_{k_1,\ldots,k_6}E^{k_1}_{\beta_6}\ldots
E^{k_6}_{\beta_1}\right)=\sum a_{k_1,\ldots,
k_6}T^{-1}_{w'_0}\left(E^{k_1}_{\beta_6}\ldots E^{k_4}_{\beta_3}\right)
T^{-1}_{w'_0}\left(E^{k_5}_{\beta_2}\cdot E^{k_6}_{\beta_1}\right)\in U_q^+.
\end{equation}
By Proposition \ref{uni9pproach.texpr}
\begin{equation}\label{uni9pproach.texeq}
T^{-1}_{w'_0}\left(E^{k_1}_{\beta_6}\ldots E^{k_4}_{\beta_3}\right)\in
U_q^+.
\end{equation}
\begin{lemma}\label{uni9pproach.texinclu}
$T^{-1}_{w'_0}\left(E_{\beta_2}\right)\in U_q^{\leq0}$,
$T^{-1}_{w'_0}\left(E_{\beta_1}\right)\in U_q^{\leq0}$.
\end{lemma}

\smallskip

{\bf Proof of the Lemma.} Suppose that $\beta_k$ is a 'compact' root
($\beta_k=\beta_1$ or $\beta_k=\beta_2$). Let $s_{i_1}s_{i_2}\ldots s_{i_M}$
be a reduced expression of $w'_0$ (of course, in the case we consider $M=2$
and there are only two different reduced expression for $w'_0$). One has
$$
\beta_k=s_{i_1}s_{i_2}\ldots s_{i_{k-1}}(\alpha_k),\qquad
E_{\beta_k}=T_{i_1}T_{i_2}\ldots T_{i_{k-1}}(E_k).
$$
Since $T_{w'_0}=T_{i_1}T_{i_2}\ldots T_{i_{M}}$ we get
$$
T^{-1}_{w'_0}(E_{\beta_k})=T^{-1}_{i_M}T^{-1}_{i_{M-1}}\ldots
T^{-1}_{i_{k}}(E_k).
$$
So we have to prove that
\begin{equation}\label{uni9pproach.texdoc}
T^{-1}_{i_M}T^{-1}_{i_{M-1}}\ldots T^{-1}_{i_{k}}(E_k)\in U_q^{\leq0}.
\end{equation}

Consider the antiautomorphism $\tau$ of the algebra $U_q\mathfrak{g}$ given
by
$$\tau(K_i)=K^{-1}_i,\quad \tau(E_i)=E_i,\quad \tau(F_i)=F_i.$$
It is not hard to prove that
$$\tau \circ T_i=T^{-1}_i \circ \tau.$$
Thus the inclusion
\begin{equation}\label{uni9pproach.texdoc1}
T_{i_M}T_{i_{M-1}}\ldots T_{i_{k}}(E_k)\in U_q^{\leq0}.
\end{equation}
is equivalent to (\ref{uni9pproach.texdoc}). One has $T_{i_k}(E_k)=-F_kK_k$. Therefore
(\ref{uni9pproach.texdoc1}) is equivalent to
$$
T_{i_M}T_{i_{M-1}}\ldots T_{i_{k+1}}(F_k)\cdot T_{i_M}T_{i_{M-1}}\ldots
T_{i_{k+1}}(K_k)\in U_q^{\leq0}
$$
and thus to
\begin{equation}\label{uni9pproach.texlast}
T_{i_M}T_{i_{M-1}}\ldots T_{i_{k+1}}(F_k)\in U_q^-.
\end{equation}
Applying the antiautomorphism $\it{k}$ (see section 2) to both hand sides of
(\ref{uni9pproach.texlast}) we get the equivalent inclusion
\begin{equation}\label{uni9pproach.texlast1}
T_{i_M}T_{i_{M-1}}\ldots T_{i_{k+1}}(E_k)\in U_q^+.
\end{equation}
But (\ref{uni9pproach.texlast1}) is a direct consequence of Lemma \ref{uni9pproach.texw}. \hfill $\square$

\medskip

The following result is well known.

\smallskip

\begin{lemma}\label{uni9pproach.texisomor}
The multiplication in $U_q \mathfrak{g}$ induces the isomorphism of vector
spaces
$$U_q^+\bigotimes U_q^{\leq0}\rightarrow U_q \mathfrak{g}.$$
\end{lemma}

\smallskip

It follows from (\ref{uni9pproach.texeq}) and Lemmas \ref{uni9pproach.texinclu}, \ref{uni9pproach.texisomor} that
(\ref{uni9pproach.texu+}) holds iff $a_{k_1,\ldots, k_6}=0$ for $k_5\ne0$ or $k_6\ne0$. We
 have completed the proof of Proposition \ref{uni9pproach.texpro}. \hfill $\square$

\bigskip

Comparing results of the two previous sections we get:\\ {\sl the algebras
$\mathbb{C}[\mathfrak{g}_{-1}]_q$, $\mathbb{C}[\mathfrak{g}_{-1}]^{\mathrm
{I}}_q$, and $\mathbb{C}[\mathfrak{g}_{-1}]^{\mathrm{II}}_q$ are isomorphic
to each other as $U_q \mathfrak{g}_0$-module algebras.}





\makeatletter \let\@thefnmark\relax \makeatother
\renewcommand{\theequation}{\thesection.\arabic{equation}}

\title{\bf HIDDEN SYMMETRY OF SOME ALGEBRAS OF q-DIFFERENTIAL OPERATORS}

\author{D. Shklyarov, \and S. Sinel'shchikov\thanks{email:
sinelshchikov@ilt.kharkov.ua}\and L. Vaksman \thanks{email:
vaksman@ilt.kharkov.ua}}

\date{\tt Institute for Low Temperature Physics \& Engineering\\
47 Lenin Avenue, 61103 Kharkov, Ukraine}

\newpage
\setcounter{section}{0}
\large

\makeatletter
\renewcommand{\@oddhead}{HIDDEN SYMMETRY OF SOME ALGEBRAS OF q-DIFFERENTIAL
OPERATORS \hfill \thepage}
\renewcommand{\@evenhead}{\thepage \hfill D. Shklyarov, S. Sinel'shchikov,
and L. Vaksman}
\let\@thefnmark\relax
\@footnotetext{This lecture has been delivered at the NATO Advanced Research
Workshop 'Non-commutative Structures in Mathematics and Physics', Kiev,
September 2000; published in Noncommutative Structures in Mathematical
Physics, S.~Duplij and J. Wess (eds), Kluwer AP, Netherlands, 2001, 309 --
320.}
\addcontentsline{toc}{chapter}{\@title \\ {\sl D. Shklyarov, S.
Sinel'shchikov, and L. Vaksman}\dotfill} \makeatother

\maketitle

\section{Introduction}

Let us explain the meaning of the words "q-differential operators" and
"hidden symmetry". Let $\mathbb{C}[z]_q$ be the algebra of polynomials in
$z$ over the field of rational functions $\mathbb{C}(q^{1/2})$ (we assume
this field to be the ground field throughout the paper). We denote by
$\Lambda^1(\mathbb{C})_q$ the $\mathbb{C}[z]_q$-bimodule with the generator
$dz$ such that
$$z \cdot dz=q^{-2}dz \cdot z.$$
Let $d$ be the linear map $\mathbb{C}[z]_q \rightarrow
\Lambda^1(\mathbb{C})_q$ given by the two conditions:
$$d:z \mapsto dz,$$
$$d(f_1(z)f_2(z))=d(f_1(z))f_2(z)+f_1(z)d(f_2(z)).$$
(The later condition is just the Leibniz rule). The bimodule
$\Lambda^1(\mathbb{C})_q$ (together with the map $d$) is a well known first
order differential calculus over the algebra $\mathbb{C}[z]_q$. The
differential $d$ allows one to introduce an operator of "partial derivative"
$\dfrac{d}{dz}$ in $\mathbb{C}[z]_q$:
$$df=dz \cdot \frac{df}{dz}.$$
Let us introduce also the notation $\widehat{z}$ for the operator in
$\mathbb{C}[z]_q$ of multiplication by $z$:
$$\widehat{z}:f(z)\mapsto zf(z).$$
Let $D(\mathbb{C})_q$ be the subalgebra in the algebra
$\mathrm{End}_{\mathbb{C}(q^{1/2})}(\mathbb{C}[z]_q)$ (of all endomorphisms
of the linear space $\mathbb{C}[z]_q$) containing $1$ and generated by
$\dfrac{d}{dz}$, $\widehat{z}$. It is easy to check that
$$\frac{d}{dz}\cdot \widehat{z}=q^{-2}\widehat{z}\cdot \frac{d}{dz}+1.$$
Thus the algebra $D(\mathbb{C})_q$ is an analogue of the Weyl algebra
$A_1(\mathbb{C})$.

Let $\lambda \in \mathbb{C}(q^{1/2})$. One checks that the map
$$\widehat{z}\mapsto \lambda \cdot \widehat{z},\qquad
\frac{d}{dz}\mapsto{\lambda}^{-1}\cdot \frac{d}{dz}$$ is extendable up to an
automorphism of the algebra $D(\mathbb{C})_q$. Such automorphisms are
"evident" symmetries of $D(\mathbb{C})_q$. It turn out that they belong to a
wider set of symmetries of $D(\mathbb{C})_q$. This set does not consists of
automorphisms only. Let us turn to precise formulations.

To start with, we recall the definition of the quantum universal enveloping
algebra $U_q \mathfrak{sl}_2$ \cite{uni9z_vaksm.texR}. It is

i) the algebra given by the generators $E$, $F$, $K$, $K^{-1}$, and the
relations
$$KK^{-1}=K^{-1}K=1,\quad KE=q^{2}EK,\quad KF=q^{-2}FK,$$
$$EF-FE=\frac{K-K^{-1}}{q-q^{-1}};$$

ii) the Hopf algebra: the comultiplication $\Delta$, the antipode $S$, and
the counit $\varepsilon$ are determined by
$$
\Delta(E)=E \otimes 1+K \otimes E,\quad \Delta(F)=F \otimes K^{-1}+1 \otimes
F,\quad \Delta(K)=K \otimes K,
$$
$$S(E)=-K^{-1}E,\qquad S(F)=-FK,\qquad S(K)=K^{-1}, $$
$$\varepsilon(E)=\varepsilon(F)=0,\qquad \varepsilon(K)=1.$$

There is a well known structure of $U_q \mathfrak{sl}_2$-module in the space
$\mathbb{C}[z]_q$. Let us describe it explicitly:
$$E:f(z)\mapsto -q^{1/2}z^2 \frac{f(z)-f(q^2z)}{z-q^2z},$$
$$F:f(z)\mapsto q^{1/2}\frac{f(z)-f(q^{-2}z)}{z-q^{-2}z},$$
$$K^{\pm 1}:f(z)\mapsto f(q^{\pm 2}z).$$
It can be checked that $\mathbb{C}[z]_q$ is a $U_q \mathfrak{sl}_2$-module
algebra, i.e. for any $\xi \in U_q \mathfrak{sl}_2$, $f_1,f_2 \in
\mathbb{C}[z]_q$
\begin{equation}\label{uni9z_vaksm.texUma1}
\xi(1)=\varepsilon(\xi)\cdot1,
\end{equation}
\begin{equation}\label{uni9z_vaksm.texUma2}
\xi(f_1f_2)=\sum_j \xi'_j(f_1)\xi''_j(f_2),
\end{equation}
with $\Delta(\xi)=\sum \limits_j \xi'_j \otimes \xi''_j$.

\medskip

{\sc Remark.} This observation is an analogue of the following one. The
group $SL_2(\mathbb{C})$ acts on $\mathbb{CP}^1$ via the fractional-linear
transformations. Thus the universal enveloping algebra $U \mathfrak{sl}_2$
acts via differential operators in the space of holomorphic functions on the
open cell $\mathbb{C}\subset \mathbb{CP}^1$.

Let $V$ be a $U_q \mathfrak{sl}_2$-module. Then the algebra
$\mathrm{End}(V)$ admits a "canonical" structure of $U_q
\mathfrak{sl}_2$-module: for $\xi \in U_q \mathfrak{sl}_2$, $T \in
\mathrm{End}(V)$
\begin{equation}\label{uni9z_vaksm.texend}
\xi(T)=\sum_j \xi'_j \cdot T \cdot S(\xi''_j),
\end{equation}
where $\Delta(\xi)=\sum \limits_j \xi'_j \otimes \xi''_j$, $S$ is the
antipode, and the elements in the right-hand side are multiplied within the
algebra $\mathrm{End}(V)$. It is well known that this action of $U_q
\mathfrak{sl}_2$ in $\mathrm{End}(V)$ makes $\mathrm{End}(V)$ into a $U_q
\mathfrak{sl}_2$-module algebra (i.e. for $\xi \in U_q \mathfrak{sl}_2$,
$T_1,T_2 \in \mathrm{End}(V)$ (\ref{uni9z_vaksm.texUma1}), (\ref{uni9z_vaksm.texUma2}) hold with $f_1,f_2$
being replaced by $T_1,T_2$, respectively).

The objects considered above are the simplest among ones we deal with in the
present paper. In this simplest case our main result can be formulated as
follows: the algebra $D(\mathbb{C})_q$ is a $U_q \mathfrak{sl}_2$-module
subalgebra in the $U_q \mathfrak{sl}_2$-module algebra
$\mathrm{End}_{\mathbb{C}(q^{1/2})}(\mathbb{C}[z]_q)$ (where the $U_q
\mathfrak{sl}_2$-action is given by (\ref{uni9z_vaksm.texend})). This $U_q
\mathfrak{sl}_2$-module structure in the algebra $D(\mathbb{C})_q$ is what
we call "hidden symmetry" of $D(\mathbb{C})_q$.

\medskip

{\sc Remark.} In the setting of the previous Remark the analogous fact is
evident: for $\xi \in \mathfrak{sl}_2$ the action (\ref{uni9z_vaksm.texend}) is just the
commutator of the differential operators $\xi$ and $T$ in the space of
holomorphic functions on $\mathbb{C}$. The commutator is again a
differential operator.

We can describe the $U_q \mathfrak{sl}_2$-action in $D(\mathbb{C})_q$
explicitly:
$$
E(\widehat{z})=-q^{1/2}\widehat{z}^2,\quad F(\widehat{z})=q^{1/2}, \quad
K^{\pm 1}(\widehat{z})=q^{\pm 2}\widehat{z},
$$
$$
E(\frac{d}{dz})=q^{-3/2}(q^{-2}+1)\widehat{z}\frac{d}{dz},\quad
F(\frac{d}{dz})=0, \quad K^{\pm 1}(\frac{d}{dz})=q^{\mp 2}\frac{d}{dz}.
$$
(The action of $U_q \mathfrak{sl}_2$ on an arbitrary element of
$D(\mathbb{C})_q$ can be produced via the rule (\ref{uni9z_vaksm.texUma2}).)

The paper is organized as follows.

In Section 2 we recall one definitions of the quantum universal enveloping
algebra $U_q \mathfrak{sl}_N$, a $U_q \mathfrak{sl}_N$-module algebra of
holomorphic polynomials on a quantum space of $m \times n$ matrices
($N=m+n$), and a well known first order differential calculus over this
algebra (in this Introduction the case $m=n=1$ was considered). Then we
introduce an algebra of q-differential operators and formulate a main
theorem concerning a hidden symmetry of this algebra.

Section 3 contains a sketch of the proof of the main theorem.

In Section 4 we discuss briefly q-analogues of the notions of a holomorphic
vector bundle and a differential operator in sections of holomorphic vector
bundles. We formulate an analogue of our main theorem for such differential
operators.

Section 'Concluding notes' discusses one of possible generalizations of our
results, specifically, the case when the matrix space is replaced by an
arbitrary prehomogeneous vector space of commutative parabolic type
\cite{uni9z_vaksm.texRu}.

Appendix deals with q-analogues of constant coefficient differential
operators in functions and in sections of holomorphic vector bundles.

We take this opportunity to thank Prof. H. P. Jakobsen and Prof. T.~
Tanisaki who attracted our attention to other approaches to the notion of
quantum differential operators.

This research was partially supported by Award No.UM1-2091 of the U.S.
Civilian Research and Development Foundation.

\bigskip

\section{The main theorem}

{\sl Throughout the paper (except for the section 'Concluding notes') we
suppose that $\mathfrak{g}=\mathfrak{sl}_N$,
$\mathfrak{k}=\mathfrak{s}(\mathfrak{gl}_n \times \mathfrak{gl}_m)$,
$N=n+m$, and $\mathfrak{p}^-$ is the space of complex $m \times n$
matrices.}

The algebra $\mathbb{C}[\mathfrak{p}^-]_q$ is the unital algebra given by
its generators $z_a^\alpha$, $a=1,\ldots n$, $\alpha=1,\ldots m$, and the
following relations
\begin{equation}\label{uni9z_vaksm.texz}
z_a^\alpha z_b^\beta=\left \{\begin{array}{ccl}qz_b^\beta z_a^\alpha &,& a=b
\;\&\;\alpha<\beta \quad{\rm or}\quad a<b \;\&\;\alpha=\beta \\ z_b^\beta
z_a^\alpha &,& a<b \;\&\;\alpha>\beta \\ z_b^\beta
z_a^\alpha+(q-q^{-1})z_a^\beta z_b^\alpha &,& a<b \;\&\;\alpha<\beta
\end{array} \right.,
\end{equation}

Let $(a_{ij})$ be the Cartan matrix for $\mathfrak{g}$. The Hopf algebra
$U_q \mathfrak{g}$ is determined by the generators $E_i$, $F_i$, $K_i$,
$K_i^{-1}$, $i=1,\ldots,N-1$, and the relations

$$
K_iK_j=K_jK_i,\qquad K_iK_i^{-1}=K_i^{-1}K_i=1,\qquad
K_iE_j=q^{a_{ij}}E_jK_i,
$$
$$
K_iF_j=q^{-a_{ij}}F_jK_i, \qquad
E_iF_j-F_jE_i=\delta_{ij}(K_i-K_i^{-1})/(q-q^{-1})
$$
\begin{equation}
E_i^2E_j-(q+q^{-1})E_iE_jE_i+E_jE_i^2=0,\qquad |i-j|=1
\end{equation}
$$F_i^2F_j-(q+q^{-1})F_iF_jF_i+F_jF_i^2=0,\qquad |i-j|=1$$
$$[E_i,E_j]=[F_i,F_j]=0,\qquad |i-j|\ne 1.$$

The comultiplication $\Delta$, the antipode $S$, and the counit
$\varepsilon$ are determined by
\begin{equation}
\Delta(E_i)=E_i \otimes 1+K_i \otimes E_i,\;\Delta(F_i)=F_i \otimes
K_i^{-1}+1 \otimes F_i,\;\Delta(K_i)=K_i \otimes K_i,
\end{equation}
\begin{equation}
S(E_i)=-K_i^{-1}E_i,\qquad S(F_i)=-F_iK_i,\qquad S(K_i)=K_i^{-1},
\end{equation}
$$\varepsilon(E_i)=\varepsilon(F_i)=0,\qquad \varepsilon(K_i)=1.$$
Denote by $U_q \mathfrak{k}$ the Hopf subalgebra in $U_q \mathfrak{g}$
generated by $E_j$, $F_j$, $K_i$, $K_i^{-1}$, $i,j=1,\ldots,N-1$, $j \ne n$.

The algebra $\mathbb{C}[\mathfrak{p}^-]_q$ possesses the well known
structure of $U_q \mathfrak{k}$-module algebra:
\begin{equation}\label{uni9z_vaksm.texh}
K_nz_a^\alpha=\left \{\begin{array}{ccl}q^2z_a^\alpha &,&a=n \;\&\;\alpha=m
\\ qz_a^\alpha &,&a=n \;\&\;\alpha \ne m \quad{\rm or}\quad a \ne n \;\&\;
\alpha=m \\ z_a^\alpha &,&\mathrm{otherwise}\end{array}\right.,
\end{equation}
and with $k \ne n$
\begin{equation}\label{uni9z_vaksm.texh's}
K_kz_a^\alpha=\left \{\begin{array}{ccl}qz_a^\alpha &,& k<n \;\&\;a=k \quad
\mathrm{or}\quad k>n \;\&\;\alpha=N-k \\ q^{-1}z_a^\alpha &,& k<n
\;\&\;a=k+1 \quad{\rm or}\quad k>n \;\&\;\alpha=N-k+1 \\ z_a^\alpha
&,&\mathrm{otherwise}\end{array}\right.,
\end{equation}
\begin{equation}\label{uni9z_vaksm.texf's}
F_kz_a^\alpha=q^{1/2}\cdot \left \{\begin{array}{ccl}z_{a+1}^\alpha &,& k<n
\;\&\;a=k \\ z_a^{\alpha+1} &,& k>n \;\&\;\alpha=N-k \\ 0 &,&{\rm
otherwise}\end{array}\right.,
\end{equation}
\begin{equation}\label{uni9z_vaksm.texe's}
E_kz_a^\alpha=q^{-1/2}\cdot \left \{\begin{array}{ccl}z_{a-1}^\alpha &,& k<n
\;\&\; a=k+1 \\ z_a^{\alpha-1}&,& k>n \;\&\;\alpha=N-k+1 \\ 0
&,&\mathrm{otherwise}\end{array}\right..
\end{equation}
This $U_q \mathfrak{k}$-module algebra structure in
$\mathbb{C}[\mathfrak{p}^-]_q$ can be extended up to a $U_q
\mathfrak{g}$-module algebra structure as follows (see \cite{uni9z_vaksm.texSSV1}):
\begin{equation}\label{uni9z_vaksm.texf}
F_nz_a^\alpha=q^{1/2}\cdot \left \{\begin{array}{ccl}1 &,& a=n \;\&
\;\alpha=m \\ 0 &,&\mathrm{otherwise}\end{array}\right.,
\end{equation}
\begin{equation}\label{uni9z_vaksm.texe}
E_nz_a^\alpha=-q^{1/2}\cdot \left \{\begin{array}{ccl}q^{-1}z_a^mz_n^\alpha
&,&a \ne n \;\&\;\alpha \ne m \\ (z_n^m)^2 &,& a=n \;\&\;\alpha=m \\
z_n^mz_a^{\alpha} &,&\mathrm{otherwise}\end{array}\right..
\end{equation}

\medskip

{\sc Remarks.} i) In the classical case the corresponding action of $U
\mathfrak{g}$ in the space of holomorphic functions on $\mathfrak{p}^-$ can
be produced via an embedding $\mathfrak{p}^-$ into the Grassmanian
$Gr_{m,N}$ as an open cell (we describe a q-analogue of the embedding in
\cite{uni9z_vaksm.texSSV1}).

Now let us recall a definition of a well known first order differential
calculus over $\mathbb{C}[\mathfrak{p}^-]_q$. Let
$\Lambda^1(\mathfrak{p}^-)_q$ be the $\mathbb{C}[\mathfrak{p}^-]_q$-bimodule
given by its generators $dz_a^\alpha$, $a=1,\ldots n$, $\alpha=1,\ldots m$,
and the relations

\begin{equation}\label{uni9z_vaksm.texzdz}
z_b^\beta dz_a^\alpha=\sum_{\alpha',\beta'=1}^m \sum_{a',b'=1}^n R_{\beta
\alpha}^{\beta'\alpha'}R^{b'a'}_{ba}dz_{a'}^{\alpha'}\cdot z_{b'}^{\beta'},
\end{equation}
with
\begin{equation}\label{uni9z_vaksm.texR}
R^{b'a'}_{ba}=\left \{\begin{array}{ccl}q^{-1} &,& a=b=a'=b' \\ 1 &,& a \ne
b \quad \&\quad a=a'\quad \&\quad b=b'\\ q^{-1}-q &,& a<b \quad \&\quad a=b'
\quad \&\quad b=a'\\ 0 &,& \mathrm{otherwise}\end{array}\right..
\end{equation}

\medskip

The map $d:z_a^\alpha \mapsto dz_a^\alpha$ can be extended up to a linear
operator $d:\mathbb{C}[\mathfrak{p}^-]_q \rightarrow
\Lambda^1(\mathfrak{p}^-)_q$ satisfying the Leibniz rule. It was noted for
the first time in \cite{uni9z_vaksm.texSV1}, that there exists a unique structure of a $U_q
\mathfrak{g}$-module $\mathbb{C}[\mathfrak{p}^-]_q$-bimodule in
$\Lambda^1(\mathfrak{p}^-)_q$ such that the map $d$ is a morphism of $U_q
\mathfrak{g}$-modules. The pair $\left(\Lambda^1(\mathfrak{p}^-)_q, d
\right)$ is the first order differential calculus over
$\mathbb{C}[\mathfrak{p}^-]_q$.

Let us introduce an algebra $D(\mathfrak{p}^-)_q$ of q-differential
operators on $\mathfrak{p}^-$. Define the linear operators
$\dfrac{\partial}{\partial z_a^{\alpha}}$ in $\mathbb{C}[\mathfrak{p}^-]_q$
via the differential $d$:
$$
df=\sum_{a=1}^{n}\sum_{\alpha=1}^{m}dz_a^{\alpha}\cdot \frac{\partial
f}{\partial z_a^{\alpha}},\qquad f\in \mathbb{C}[\mathfrak{p}^-]_q,
$$
and the operators $\widehat{z_a^{\alpha}}$ by
$$
\widehat{z_a^{\alpha}}f=z_a^{\alpha}\cdot f, \qquad f \in \mathbb{C}
[\mathfrak{p}^-]_q.
$$
Then $D(\mathfrak{p}^-)_q$ is the unital subalgebra in
$\mathrm{End}_{\mathbb{C}(q^{1/2})}(\mathbb{C}[\mathfrak{p}^-]_q)$ generated
by the operators $\dfrac{\partial }{\partial z_a^{\alpha}}$,
$\widehat{z_a^{\alpha}}$, $a=1,\ldots n,$ $\alpha=1,\ldots m.$

To start with, we describe $D(\mathfrak{p}^-)_q$ in terms of generators and
relations.

\medskip

\begin{proposition}
The complete list of relations between the generators
$\widehat{z_a^\alpha}$, $\dfrac{\partial}{\partial z_a^\alpha}$, $a=1,\ldots
n$, $\alpha=1,\ldots m$, of $D(\mathfrak{p}^-)_q$ is as follows

\begin{equation}\label{uni9z_vaksm.texz^}
\widehat{z}_a^\alpha \widehat{z}_b^\beta=\left \{\begin{array}{ccl}q
\widehat{z}_b^\beta \widehat{z}_a^\alpha &,& a=b \;\&\;\alpha<\beta \quad
\mathrm{or}\quad a<b \;\&\;\alpha=\beta \\ \widehat{z}_b^\beta
\widehat{z}_a^\alpha &,& a<b \;\&\;\alpha>\beta \\ \widehat{z}_b^\beta
\widehat{z}_a^\alpha+(q-q^{-1})\widehat{z}_a^\beta \widehat{z}_b^\alpha &,&
a<b \;\&\;\alpha<\beta \end{array} \right.,
\end{equation}
$\dfrac{\partial}{\partial z_b^\beta}\dfrac{\partial}{\partial z_a^\alpha}=$
\begin{equation}\label{uni9z_vaksm.texd/dz}
=\left \{\begin{array}{ccl}q \dfrac{\partial}{\partial
z_a^\alpha}\dfrac{\partial}{\partial z_b^\beta} &,& a=b \;\&\;\alpha<\beta
\quad \mathrm{or}\quad a<b \;\&\;\alpha=\beta \\ \dfrac{\partial}{\partial
z_a^\alpha}\dfrac{\partial}{\partial z_b^\beta} &,& a<b \;\&\;\alpha>\beta
\\ \dfrac{\partial}{\partial z_a^\alpha}\dfrac{\partial}{\partial z_b^\beta}
+(q-q^{-1})\dfrac{\partial}{\partial z_a^\beta}\dfrac{\partial}{\partial
z_b^\alpha} &,& a<b \;\&\;\alpha<\beta \end{array} \right.,
\end{equation}
\begin{equation}\label{uni9z_vaksm.texd/dz,z}
\frac{\partial}{\partial z_a^\alpha}\widehat{z}_b^\beta=\sum_{a',b'=1}^n
\sum_{\alpha',\beta'=1}^m R_{ba'}^{b'a} R_{\beta
\alpha'}^{\beta'\alpha}\widehat{z}_{b'}^{\beta'} \frac{\partial}{\partial
z_{a'}^{\alpha'}}+\delta_{ab} \delta^{\alpha \beta},
\end{equation}
with $\delta_{ab},\delta^{\alpha \beta}$ being the Kronecker symbols, and
$R_{ba'}^{b'a}$ given by (\ref{uni9z_vaksm.texR}).
\end{proposition}

\medskip

Using the $U_q \mathfrak{g}$-module structure in
$\mathbb{C}[\mathfrak{p}^-]_q$, we can define the structure of $U_q
\mathfrak{g}$-module algebra in
$\mathrm{End}_{\mathbb{C}(q^{1/2})}(\mathbb{C}[\mathfrak{p}^-]_q)$ via
(\ref{uni9z_vaksm.texend}) with $\xi \in U_q \mathfrak{g}$, $T \in
\mathrm{End}_{\mathbb{C}(q^{1/2})}(\mathbb{C}[\mathfrak{p}^-]_q)$.

Our main result is

\begin{theorem}
i) The algebra $D(\mathfrak{p}^-)_q$ is a $U_q \mathfrak{g}$-module
subalgebra in the $U_q \mathfrak{g}$-module algebra
$\mathrm{End}_{\mathbb{C}(q^{1/2})}(\mathbb{C}[\mathfrak{p}^-]_q)$.

ii) The $U_q \mathfrak{g}$-module structure in $D(\mathfrak{p}^-)_q$ is
described explicitly as follows:

$U_q \mathfrak{g}$ acts on the generators $\widehat{z}_a^{\alpha}$ via
formulae (\ref{uni9z_vaksm.texh})-(\ref{uni9z_vaksm.texe's}) (where $z_a^{\alpha}$ should be replaced by
$\widehat{z}_a^{\alpha}$); for the generators $\dfrac{\partial}{\partial
z_a^\alpha}$ the formulae are
\begin{equation}\label{uni9z_vaksm.texh'}
K_n \frac{\partial}{\partial z_a^\alpha}=\left
\{\begin{array}{ccl}q^{-2}\dfrac{\partial}{\partial z_a^\alpha} &,&a=n
\;\&\;\alpha=m \\ q^{-1}\dfrac{\partial}{\partial z_a^\alpha} &,&a=n
\;\&\;\alpha \ne m \quad \mathrm{or}\quad a \ne n \;\&\; \alpha=m \\
\dfrac{\partial}{\partial z_a^\alpha}
&,&\mathrm{otherwise}\end{array}\right.,
\end{equation}
\begin{equation}\label{uni9z_vaksm.texf'} F_n \dfrac{\partial}{\partial z_a^\alpha}=0
\qquad a=1,\ldots n, \qquad \alpha=1,\ldots m,
\end{equation}
$E_n \dfrac{\partial}{\partial z_a^\alpha}=q^{-3/2}\cdot$
\begin{equation}\label{uni9z_vaksm.texe'}
\cdot \left \{\begin{array}{ccl}\sum \limits_{b=1}^n \widehat{z}_b^m
\dfrac{\partial}{\partial z_b^m}+\sum \limits_{\beta=1}^m
\widehat{z}_n^\beta \dfrac{\partial}{\partial z_n^\beta}+(q^{-2}-1)\sum
\limits_{b=1}^n \sum \limits_{\beta=1}^m \widehat{z}_b^\beta
\dfrac{\partial}{\partial z_b^\beta} &,&a=n \;\&\;\alpha=m \\ \sum
\limits_{\beta=1}^m \widehat{z}_n^\beta \dfrac{\partial}{\partial z_a^\beta}
&,&a \ne n \;\&\;\alpha=m \\ \sum \limits_{b=1}^n \widehat{z}_b^m
\dfrac{\partial}{\partial z_b^\alpha} &,&a=n \;\&\;\alpha \ne m \\ 0
&,&\mathrm{otherwise}\end{array}\right.,
\end{equation}
and with $k \ne n$
\begin{equation}\label{uni9z_vaksm.texh''s}
K_k \dfrac{\partial}{\partial z_a^\alpha}=\left
\{\begin{array}{ccl}q^{-1}\dfrac{\partial}{\partial z_a^\alpha} &,& k<n
\;\&\;a=k \quad \mathrm{or}\quad k>n \;\&\;\alpha=N-k \\ q
\dfrac{\partial}{\partial z_a^\alpha} &,& k<n \;\&\;a=k+1 \quad
\mathrm{or}\quad k>n \;\&\;\alpha=N-k+1 \\ \dfrac{\partial}{\partial
z_a^\alpha} &,&\mathrm{otherwise}\end{array}\right.,
\end{equation}
\begin{equation}\label{uni9z_vaksm.texf''s}
F_k \frac{\partial}{\partial z_a^\alpha}=-q^{3/2}\cdot \left
\{\begin{array}{ccl}\dfrac{\partial}{\partial z_{a-1}^\alpha} &,& k<n
\;\&\;a=k+1
\\ \dfrac{\partial}{\partial z_a^{\alpha-1}} &,& k>n \;\&\;\alpha=N-k+1 \\ 0
&,& \mathrm{otherwise}\end{array}\right.,
\end{equation}
\begin{equation}\label{uni9z_vaksm.texe''s}
E_k \frac{\partial}{\partial z_a^\alpha}=-q^{-3/2}\cdot \left
\{\begin{array}{ccl}\dfrac{\partial}{\partial z_{a+1}^\alpha} &,& k<n \;\&\;
a=k \\ \dfrac{\partial}{\partial z_a^{\alpha+1}}&,& k>n \;\&\;\alpha=N-k \\
0 &,& \mathrm{otherwise}\end{array}\right..
\end{equation}
\end{theorem}

\bigskip

\section{Sketch of the proof}

Let us outline an idea of the proof of the main theorem. To prove the
statement i) of the theorem we have to explain why for arbitrary $\xi \in
U_q \mathfrak{g}$, $T \in D(\mathfrak{p}^-)_q$
\begin{equation}\label{uni9z_vaksm.texin}
\xi(T)\in D(\mathfrak{p}^-)_q.
\end{equation}

The map $z_a^\alpha \mapsto \widehat{z}_a^\alpha,$ $a=1,\ldots n$,
$\alpha=1,\ldots m$, is extendable up to an embedding of algebras
$\mathbb{C}[\mathfrak{p}^-]_q \hookrightarrow
\mathrm{End}_{\mathbb{C}(q^{1/2})}(\mathbb{C}[\mathfrak{p}^-]_q)$. This
embedding intertwines the actions of $U_q \mathfrak{g}$ in
$\mathbb{C}[\mathfrak{p}^-]_q$ and
$\mathrm{End}_{\mathbb{C}(q^{1/2})}(\mathbb{C}[\mathfrak{p}^-]_q)$ (this is
a corollary of the fact that $\mathbb{C}[\mathfrak{p}^-]_q$ is a $U_q
\mathfrak{g}$-module algebra). This observation proves (\ref{uni9z_vaksm.texin}) for $T$ of
the form $\widehat{f}$, $f \in \mathbb{C}[\mathfrak{p}^-]_q,$ as well as the
first part of the statement ii) of the theorem. What remains is to prove
(\ref{uni9z_vaksm.texin}) for $T=\dfrac{\partial}{\partial z_a^\alpha},$ $a=1,\ldots n$,
$\alpha=1,\ldots m$.

The space $\mathrm{End}_{\mathbb{C}(q^{1/2})}(\mathbb{C}[\mathfrak{p}^-]_q)$
can be made into a left $\mathbb{C}[\mathfrak{p}^-]_q$-module as follows:
$$z_a^\alpha(T)=\widehat{z}_a^\alpha \cdot T,$$
with $a=1,\ldots n$, $\alpha=1,\ldots m$, $T \in
\mathrm{End}_{\mathbb{C}(q^{1/2})}(\mathbb{C}[\mathfrak{p}^-]_q)$. This
structure is compatible with the action of $U_q \mathfrak{g}$. Consider the
$U_q \mathfrak{g}$-module $\Lambda^1(\mathfrak{p}^-)_q
\bigotimes_{\mathbb{C}[\mathfrak{p}^-]_q}
\mathrm{End}_{\mathbb{C}(q^{1/2})}(\mathbb{C}[\mathfrak{p}^-]_q).$ The
differential $d:\mathbb{C}[\mathfrak{p}^-]_q \rightarrow
\Lambda^1(\mathfrak{p}^-)_q$ is a morphism of the $U_q
\mathfrak{g}$-modules. This implies $U_q \mathfrak{g}$-invariance of the
element $\sum \limits_{a=1}^n \sum \limits_{\alpha=1}^mdz_a^\alpha \otimes
\dfrac{\partial}{\partial z_a^\alpha} \in \Lambda^1(\mathfrak{p}^-)_q
\bigotimes_{\mathbb{C}[\mathfrak{p}^-]_q}
\mathrm{End}_{\mathbb{C}(q^{1/2})}(\mathbb{C}[\mathfrak{p}^-]_q),$ i.e. for
all $\xi \in U_q \mathfrak{g}$
\begin{equation}\label{uni9z_vaksm.texeq}
\sum_{a=1}^n \sum_{\alpha=1}^m \sum_j \xi'_jdz_a^\alpha \otimes \xi''_j
\frac{\partial}{\partial z_a^\alpha}=\varepsilon(\xi) \sum_{a=1}^n
\sum_{\alpha=1}^mdz_a^\alpha \otimes \frac{\partial}{\partial z_a^\alpha}
\end{equation}
with $\varepsilon$ being the counit of $U_q \mathfrak{g}$, $\Delta(\xi)=\sum
\limits_j \xi'_j \otimes \xi''_j$ ($\Delta$ is the coproduct in $U_q
\mathfrak{g}$). As it was proved in \cite{uni9z_vaksm.texSSV1},
$\Lambda^1(\mathfrak{p}^-)_q$ is a free right
$\mathbb{C}[\mathfrak{p}^-]_q$-module with the generators $dz_a^\alpha$,
$a=1,\ldots n$, $\alpha=1,\ldots m$. Thus, for $\xi \in U_q \mathfrak{g}$
there exists a unique set $f_{\beta,a}^{b,\alpha}(\xi)\in
\mathbb{C}[\mathfrak{p}^-]_q$, $a=1,\ldots n$, $\alpha=1,\ldots m$,
$b=1,\ldots n$, $\beta=1,\ldots m,$ such that
$$
\xi dz_a^\alpha=\sum_{b=1}^n \sum_{\beta=1}^mdz_b^\beta
f_{\beta,a}^{b,\alpha}(\xi).
$$
Using the later equality, we can rewrite (\ref{uni9z_vaksm.texeq}) as follows:
\begin{equation}\label{uni9z_vaksm.texeq1}
\sum_{a,b=1}^n \sum_{\alpha,\beta=1}^m \sum_jdz_b^\beta \otimes
f_{\beta,a}^{b,\alpha}(\xi'_j)\xi''_j \frac{\partial}{\partial
z_a^\alpha}=\varepsilon(\xi) \sum_{a=1}^n \sum_{\alpha=1}^mdz_a^\alpha
\otimes \frac{\partial}{\partial z_a^\alpha}.
\end{equation}
Now one can obtain formulae (\ref{uni9z_vaksm.texh'}) - (\ref{uni9z_vaksm.texe''s}) (and thus prove
(\ref{uni9z_vaksm.texin}) for $T=\dfrac{\partial}{\partial z_a^\alpha},$ $a=1,\ldots n$,
$\alpha=1,\ldots m$) via applying (\ref{uni9z_vaksm.texeq1}) to the generators $E_i, F_i,
K_i, K^{-1}_i$ of $U_q \mathfrak{g}$.

\bigskip

\section{A generalization: q-differential operators in holomorphic
q-bundles.}

Let $\Gamma$ be a finitely generated free right
$\mathbb{C}[\mathfrak{p}^-]_q$-module, i.e. there exists an isomorphism of
the right $\mathbb{C}[\mathfrak{p}^-]_q$-modules
$$\pi:\Gamma \rightarrow V \bigotimes \mathbb{C}[\mathfrak{p}^-]_q,$$
with $V$ being a finite dimensional vector space. Elements of $\Gamma$ are
q-analogues of sections of a holomorphic vector bundle over
$\mathfrak{p}^-$. {\sl Finitely generated free right
$\mathbb{C}[\mathfrak{p}^-]_q$-modules will be called vector q-bundles.} The
isomorphism $\pi$ will be called a trivialization of $\Gamma$.

Let $\Gamma_1$, $\Gamma_2$ be two vector q-bundles, $\pi_1:\Gamma_1
\rightarrow V_1 \bigotimes \mathbb{C}[\mathfrak{p}^-]_q$, $\pi_2:\Gamma_2
\rightarrow V_2 \bigotimes \mathbb{C}[\mathfrak{p}^-]_q$ their
trivializations. Set
$$
D(\Gamma_1,\Gamma_2)_q=\left \{D \in
\mathrm{Hom}(\Gamma_1,\Gamma_2)\left|\pi_2 \cdot D \cdot{\pi_1}^{-1}\in
\mathrm{Hom}(V_1,V_2)\bigotimes D(\mathfrak{p}^-)_q \right.\right \}.
$$
Elements of $D(\Gamma_1,\Gamma_2)_q$ can be treated as q-analogues of
differential operators in sections of holomorphic vector bundles.

We need to verify that $D(\Gamma_1,\Gamma_2)_q$ is independent of the choice
of trivializations. This follows from the observation that for two
trivializations $\pi_1:\Gamma \rightarrow V_1
\bigotimes\mathbb{C}[\mathfrak{p}^-]_q$, $\pi_2:\Gamma \rightarrow V_2
\bigotimes\mathbb{C}[\mathfrak{p}^-]_q$ of a q-bundle $\Gamma$
$$
\pi_2 \cdot(\pi_1)^{-1}\in \mathrm{Hom}(V_1,V_2)\bigotimes
\widehat{\mathbb{C}[\mathfrak{p}^-]_q}
$$
with $\widehat{\mathbb{C}[\mathfrak{p}^-]_q}$ being the unital subalgebra in
$D(\mathfrak{p}^-)_q$ generated by $\widehat{z}_a^\alpha,$ $a=1,\ldots n$,
$\alpha=1,\ldots m$.

Suppose that $\Gamma$ is a {\sl $U_q \mathfrak{g}$-module} q-bundle. It
means that $\Gamma$ is a vector q-bundle and a $U_q \mathfrak{g}$-module,
and the multiplication map $\Gamma \bigotimes \mathbb{C}[\mathfrak{p}^-]_q
\rightarrow \Gamma$ is a morphism of the $U_q \mathfrak{g}$-modules.

For $U_q \mathfrak{g}$-module vector q-bundles a result analogous to the
main theorem (Section 2) can be obtained. Let us formulate it.

If $V_1$, $V_2$ are modules over a Hopf algebra $A$ then the space
$\mathrm{Hom}(V_1,V_2)$ admits the following "canonical" structure of an
$A$-module: for $\xi \in A$, $T \in \mathrm{Hom}(V_1,V_2)$
\begin{equation}\label{uni9z_vaksm.texhom}
\xi(T)=\sum_j \xi'_j \cdot T \cdot S(\xi''_j),
\end{equation}
where $\Delta(\xi)=\sum \limits_j \xi'_j \otimes \xi''_j$ ($\Delta$ is the
coproduct), $S$ is the antipode, and the product in the right-hand side
means the composition of the maps $S(\xi''_j)\in \mathrm{End}(V_1)$, $T \in
\mathrm{Hom}(V_1,V_2)$, $\xi'_j \in \mathrm{End}(V_2)$. It is well known
that this action makes $\mathrm{Hom}(V_1,V_2)$ into an $A$-module left
$\mathrm{End}(V_2)$-module and an $A$-module right
$\mathrm{End}(V_1)$-module, i.e. the composition map
$$
\mathrm{End}(V_2)\bigotimes \mathrm{Hom}(V_1,V_2)\bigotimes
\mathrm{End}(V_1) \rightarrow \mathrm{Hom}(V_1,V_2)
$$
is a morphism of the $A$-modules.

Let $\Gamma_1$, $\Gamma_2$ be $U_q \mathfrak{g}$-module vector q-bundles.
Using our main theorem, one can prove the following

\medskip

\begin{proposition}
The subspace $D(\Gamma_1,\Gamma_2)_q \subset
\mathrm{Hom}(\Gamma_1,\Gamma_2)$ is $U_q \mathfrak{g}$-invariant; thus, the
composition map
$$
D(\Gamma_2)_q \bigotimes D(\Gamma_1,\Gamma_2)_q \bigotimes D(\Gamma_1)_q
\rightarrow D(\Gamma_1,\Gamma_2)_q
$$
(here $D(\Gamma)_q$ denotes $D(\Gamma,\Gamma)_q$) makes
$D(\Gamma_1,\Gamma_2)_q$ into a $U_q \mathfrak{g}$-module left
$D(\Gamma_2)_q$-module and a $U_q \mathfrak{g}$-module right
$D(\Gamma_1)_q$-module.
\end{proposition}

\bigskip

\section{Concluding notes}

The space of $m \times n$ matrices considered in the present paper is the
simplest example of a prehomogeneous vector space of commutative parabolic
type \cite{uni9z_vaksm.texRu}. Such vector spaces are closely related to non-compact
Hermitian symmetric spaces. Specifically, any non-compact Hermitian
symmetric space can be realized (via the Harish-Chandra embedding) as a
bounded symmetric domain in some prehomogeneous vector space of commutative
parabolic type.

In \cite{uni9z_vaksm.texSV2} a q-analogue of an arbitrary prehomogeneous vector space of
commutative parabolic type was constructed. More precisely, let $U$ be a
bounded symmetric domain, $\mathfrak{p}^-$ the corresponding prehomogeneous
vector space, and $\mathfrak{g}$ the complexification of the Lie algebra of
the automorphism group of $U$. In the paper \cite{uni9z_vaksm.texSV2} a $U_q
\mathfrak{g}$-module algebra $\mathbb{C}[\mathfrak{p}^-]_q$ and a covariant
first order differential calculus $\left(\Lambda^1(\mathfrak{p}^-),d
\right)$ over $\mathbb{C}[\mathfrak{p}^-]_q$ were introduced (the notation
$\mathfrak{g}_{-1}$ was used in \cite{uni9z_vaksm.texSV2} instead of $\mathfrak{p}^-$).
Using the first order differential calculus, one can define an algebra
$D(\mathfrak{p}^-)_q$ of q-differential operators in
$\mathbb{C}[\mathfrak{p}^-]_q$ just as it was done in Section 2 in the case
of the matrix space.

In this general setting it can also be proved that $D(\mathfrak{p}^-)_q$ is
a $U_q \mathfrak{g}$-module subalgebra in the $U_q\mathfrak{g}$-module
algebra ${\rm End}(\mathbb{C}[\mathfrak{p}^-]_q)$. Indeed, it easy to see
that the proof of our main theorem (Section 3) does not use a specific
nature of the case when $\mathfrak{p}^-$ is the matrix space.

\bigskip

\section{Appendix: Constant coefficient q-differential operators.}

Let $D(\mathfrak{p}^-)^\mathrm{const}_q$ be the unital subalgebra in
$D(\mathfrak{p}^-)_q$ generated by $\dfrac{\partial}{\partial z_a^\alpha},$
$a=1,\ldots n$, $\alpha=1,\ldots m$. By (\ref{uni9z_vaksm.texh'}), (\ref{uni9z_vaksm.texh''s}),
(\ref{uni9z_vaksm.texf''s}), (\ref{uni9z_vaksm.texe''s}), it is a $U_q \mathfrak{k}$-module subalgebra in
$D(\mathfrak{p}^-)_q$. Its elements are q-analogues of the constant
coefficient differential operators in the space of holomorphic polynomials.
It turn out that the natural action of $D(\mathfrak{p}^-)^\mathrm{const}_q$
in $\mathbb{C}[\mathfrak{p}^-]_q$ is related to the action of the quantized
universal enveloping algebra. Turn to precise formulations.

Denote by $\check{U}_q \mathfrak{g}$ the Hopf algebra derived from $U_q
\mathfrak{g}$ by adjoining the pairwise commuting generators $L^{\pm 1}_i$,
$i=1,\ldots N-1$, such that $K_i=\prod_{j=1}^{N-1}L_j^{a_{ij}}$ (with
$(a_{ij})$ being the Cartan matrix for $\mathfrak{g}$), and
$$
L_iL_i^{-1}=L_i^{-1}L_i=1,\qquad L_iE_j=q^{\delta_{ij}}E_jL_i,\qquad
L_iF_j=q^{-\delta_{ij}}F_jL_i,
$$
$$
\Delta(L_i)=L_i \otimes L_i,\qquad S(L_i)=L_i^{-1},\qquad
\varepsilon(L_i)=1.
$$
Endow $\check{U}_q \mathfrak{g}$ with a structure of $U_q
\mathfrak{g}$-module algebra via the usual quantum adjoint action:
$\mathrm{ad}\xi(\eta)=\sum \limits_j \xi'_j \cdot \eta \cdot S(\xi''_j)$,
with $\Delta(\xi)=\sum \limits_j \xi'_j \otimes \xi''_j.$

As it was proved in \cite{uni9z_vaksm.texJL}, the $U_q \mathfrak{g}$-submodule
$\mathrm{ad}U_q \mathfrak{g}(L_1^{-2})\subset \check{U}_q \mathfrak{g}$ is
finite-dimensional. Let $F$ be its primitive vector
$\mathrm{ad}F_{N-1}\mathrm{ad}F_{N-2}\ldots \mathrm{ad}F_1(L_1^{-2})$.
Evidently, $\mathrm{ad}U_q \mathfrak{k}(F)$ is a finite dimensional $U_q
\mathfrak{k}$-submodule in $\check{U}_q \mathfrak{g}$. Let $U_q
\mathfrak{p}^-$ be the minimal subalgebra in $\check{U}_q \mathfrak{g}$
containing ${\rm ad}U_q \mathfrak{k}(F)$. The following statement can be
proved.

\medskip

\begin{proposition}
For any $\xi \in U_q \mathfrak{p}^-$ there exists a unique $\partial_ \xi
\in D(\mathfrak{p}^-)^\mathrm{const}_q$ such that
$$\xi(f)=\partial_ \xi(f),$$
for any $f \in \mathbb{C}[\mathfrak{p}^-]_q$. The map $\phi:U_q
\mathfrak{p}^-\rightarrow D(\mathfrak{p}^-)^\mathrm{const}_q$, $\phi:\xi
\mapsto \partial_ \xi$, is an isomorphism of the $U_q \mathfrak{k}$-module
algebras.
\end{proposition}

\medskip

Let us produce a notion of constant coefficient differential operators for
vector q-bundles. First of all, we need to distinguish a class of vector
q-bundles for which this notion is well defined.

Let $U_q(\mathfrak{p}^-+\mathfrak{k})$ be the Hopf subalgebra in $U_q
\mathfrak{g}$ generated by $F_n$ and $U_q \mathfrak{k}$. Suppose that
$\Gamma$ is a $U_q(\mathfrak{p}^-+\mathfrak{k})$-module vector q-bundle (it
means that $\Gamma$ is a vector q-bundle and a
$U_q(\mathfrak{p}^-+\mathfrak{k} )$-module, and the multiplication map
$\Gamma \bigotimes \mathbb{C}[\mathfrak{p}^-]_q \rightarrow \Gamma$ is a
morphism of the $U_q(\mathfrak{p}^-+\mathfrak{k} )$-modules). A
trivialization $\pi:\Gamma \rightarrow V \bigotimes
\mathbb{C}[\mathfrak{p}^-]_q$ is called good trivialization if it satisfies
the following conditions:

i) $V$ is a finite dimensional $U_q(\mathfrak{p}^-+\mathfrak{k})$-module
with the property $F_nv=0$ for any $v \in V$;

ii) $\pi$ is a morphism of the $U_q(\mathfrak{p}^-+\mathfrak{k} )$-modules
(here $V \bigotimes \mathbb{C}[\mathfrak{p}^-]_q$ is endowed with
$U_q(\mathfrak{p}^-+\mathfrak{k} )$-module structure via the coproduct
$\Delta:U_q(\mathfrak{p}^-+\mathfrak{k} )\rightarrow
U_q(\mathfrak{p}^-+\mathfrak{k} )\bigotimes
U_q(\mathfrak{p}^-+\mathfrak{k})$).

It can be proved that for any two good trivializations $\pi_1:\Gamma
\rightarrow V_1 \bigotimes \mathbb{C}[\mathfrak{p}^-]_q$, $\pi_2:\Gamma
\rightarrow V_2 \bigotimes \mathbb{C}[\mathfrak{p}^-]_q$ of a vector
q-bundle $\Gamma$
\begin{equation}\label{uni9z_vaksm.textriv}
\pi_2 \cdot \pi_1^{-1}=T \otimes 1
\end{equation}
with $T \in \mathrm{Hom}_{U_q(\mathfrak{p}^-+\mathfrak{k})}(V_1,V_2).$

The set of $U_q(\mathfrak{p}^-+\mathfrak{k})$-module vector q-bundles
admitting good trivializations is the class of vector q-bundles for which
the notion of a q-differential operator with constant coefficients is
well-defined: if $\Gamma_1$, $\Gamma_2$ admit good trivializations
$\pi_1:\Gamma_1 \rightarrow V_1 \bigotimes \mathbb{C}[\mathfrak{p}^-]_q$,
$\pi_2:\Gamma_2 \rightarrow V_2 \bigotimes \mathbb{C}[\mathfrak{p}^-]_q$,
then the elements of the space
$$
\left \{D \in D(\Gamma_1,\Gamma_2)_q \left|\pi_2 \cdot D
\cdot{\pi_1}^{-1}\in \mathrm{Hom}(V_1,V_2)\bigotimes
D(\mathfrak{p}^-)^\mathrm{const}_q \right.\right \}
$$
can be treated as q-analogues of the constant coefficient differential
operators in sections of holomorphic vector bundles. By (\ref{uni9z_vaksm.textriv}), this
space is independent of good trivializations.

\end{document}